\documentclass[english]{article}
\usepackage[T1]{fontenc}
\usepackage[utf8]{inputenc}
\usepackage[a4paper]{geometry}
\usepackage{color}
\usepackage{babel}
\usepackage{bbding}
\usepackage{url}
\usepackage{amsmath}
\usepackage{amsthm}
\usepackage{amssymb}
\usepackage{setspace}
\usepackage[all]{xy}
\usepackage[unicode=true,pdfusetitle,
 bookmarks=true,bookmarksnumbered=false,bookmarksopen=false,
 breaklinks=false,pdfborder={0 0 0},pdfborderstyle={},backref=false,colorlinks=true]
 {hyperref}
\hypersetup{
 pdfencoding=auto, psdextra}

\theoremstyle{plain}
\newtheorem{thm}{\protect\theoremname}
\theoremstyle{plain}
\newtheorem{lem}[thm]{\protect\lemmaname}
\theoremstyle{plain}
\newtheorem{cor}[thm]{\protect\corollaryname}
\theoremstyle{plain}
\newtheorem{prop}[thm]{\protect\propositionname}
\theoremstyle{plain}
\newtheorem{defn}[thm]{\protect\definitionname}

\usepackage[leqno]{nccmath}
\usepackage{mathtools}
\usepackage{amssymb}
\usepackage{enumitem}
\usepackage{booktabs}
\usepackage{mathdots}
\usepackage[labelformat=empty,belowskip=2pt]{caption}
\usepackage{float}

\hypersetup{linkcolor=blue}

\usepackage{color}
\usepackage{tikz}
\usetikzlibrary{arrows, fit, shapes.geometric}
\tikzset{>=latex}

\def \ARDots {{\footnotesize{}$\bullet\bullet\bullet$}}

\def \ARPshiftAlAA {(0.4,-0.2)}
\def \ARRshiftAlAA {(0.0,-0.0)}
\def \ARIshiftAlAA {(0.4,-0.2)}

\def \ARPshiftAlAAA {(0.4,-0.2)}
\def \ARRshiftAlAAA {(0.0,-0.0)}
\def \ARIshiftAlAAA {(0.4,-0.2)}

\def \ARPshiftAlAAAA {(0.4,-0.2)}
\def \ARRshiftAlAAAA {(0.0,-0.0)}
\def \ARIshiftAlAAAA {(0.4,-0.2)}

\def \ARPshiftAAlAA {(0.4,-0.2)}
\def \ARRshiftAAlAA {(0.0,-0.0)}
\def \ARIshiftAAlAA {(0.4,-0.2)}

\def \ARPshiftAAlAAA {(0.4,-0.2)}
\def \ARRshiftAAlAAA {(0.0,-0.0)}
\def \ARIshiftAAlAAA {(0.4,-0.2)}

\def \ARPshiftAAlAAAA {(0.4,-0.2)}
\def \ARRshiftAAlAAAA {(0.0,-0.0)}
\def \ARIshiftAAlAAAA {(0.4,-0.2)}

\def \ARPshiftAAAlAAA {(0.4,-0.2)}
\def \ARRshiftAAAlAAA {(0.0,-0.0)}
\def \ARIshiftAAAlAAA {(0.4,-0.2)}

\def \ARPshiftAAAlAAAA {(0.4,-0.2)}
\def \ARRshiftAAAlAAAA {(0.0,-0.0)}
\def \ARIshiftAAAlAAAA {(0.4,-0.2)}

\def \ARPshiftAAAlAAAAA {(0.4,-0.2)}
\def \ARRshiftAAAlAAAAA {(0.0,-0.0)}
\def \ARIshiftAAAlAAAAA {(0.4,-0.2)}

\def \ARPshiftDDDDD {(0.0,-0.10)}
\def \ARRshiftDDDDD {(0.0,-0.10)}
\def \ARIshiftDDDDD {(0.0,-0.10)}

\def \ARPshiftDDDDDD {(0.0,-0.10)}
\def \ARRshiftDDDDDD {(0.0,-0.10)}
\def \ARIshiftDDDDDD {(0.0,-0.10)}

\def \ARPshiftDDDDDDD {(0.0,-0.25)}
\def \ARRshiftDDDDDDD {(0.0,-0.25)}
\def \ARIshiftDDDDDDD {(0.0,-0.25)}

\def \ARPshiftEEEEEEE {(0.5,-0.4)}
\def \ARRshiftEEEEEEE {(0.5,-0.4)}
\def \ARIshiftEEEEEEE {(0.5,-0.4)}

\def \ARPshiftEEEEEEEE {(0.45,-0.35)}
\def \ARRshiftEEEEEEEE {(0.45,-0.35)}
\def \ARIshiftEEEEEEEE {(0.45,-0.35)}

\def \ARPshiftEEEEEEEEE {(0.25,-0.35)}
\def \ARRshiftEEEEEEEEE {(0.25,-0.35)}
\def \ARIshiftEEEEEEEEE {(0.25,-0.35)}

\tikzset{
	seloneAlAA/.style = {scale=.8,fill=magenta!20, inner sep=0pt,outer sep=4pt},
	seltwoAlAA/.style = {scale=.8,fill=green!20, inner sep=0pt,outer sep=4pt},
	noselAlAA/.style = {scale=.8,inner sep=0pt,outer sep=4pt},
	seloneAlAAA/.style = {scale=.8,fill=magenta!20, inner sep=0pt,outer sep=4pt},
	seltwoAlAAA/.style = {scale=.8,fill=green!20, inner sep=0pt,outer sep=4pt},
	noselAlAAA/.style = {scale=.8,inner sep=0pt,outer sep=4pt},
	seloneAlAAAA/.style = {scale=.8,fill=magenta!20, inner sep=0pt,outer sep=4pt},
	seltwoAlAAAA/.style = {scale=.8,fill=green!20, inner sep=0pt,outer sep=4pt},
	noselAlAAAA/.style = {scale=.8,inner sep=0pt,outer sep=4pt},
	seloneAAlAA/.style = {scale=.8,fill=magenta!20, inner sep=0pt,outer sep=4pt},
	seltwoAAlAA/.style = {scale=.8,fill=green!20, inner sep=0pt,outer sep=4pt},
	noselAAlAA/.style = {scale=.8,inner sep=0pt,outer sep=4pt},
	seloneAAlAAA/.style = {scale=.8,fill=magenta!20, inner sep=0pt,outer sep=4pt},
	seltwoAAlAAA/.style = {scale=.8,fill=green!20, inner sep=0pt,outer sep=4pt},
	noselAAlAAA/.style = {scale=.8,inner sep=0pt,outer sep=4pt},
	seloneAAlAAAA/.style = {scale=.8,fill=magenta!20, inner sep=0pt,outer sep=4pt},
	seltwoAAlAAAA/.style = {scale=.8,fill=green!20, inner sep=0pt,outer sep=4pt},
	noselAAlAAAA/.style = {scale=.8,inner sep=0pt,outer sep=4pt},
	seloneAAAlAAA/.style = {scale=.8,fill=magenta!20, inner sep=0pt,outer sep=4pt},
	seltwoAAAlAAA/.style = {scale=.8,fill=green!20, inner sep=0pt,outer sep=4pt},
	noselAAAlAAA/.style = {scale=.8,inner sep=0pt,outer sep=4pt},
	seloneAAAlAAAA/.style = {scale=.8,fill=magenta!20, inner sep=0pt,outer sep=4pt},
	seltwoAAAlAAAA/.style = {scale=.8,fill=green!20, inner sep=0pt,outer sep=4pt},
	noselAAAlAAAA/.style = {scale=.8,inner sep=0pt,outer sep=4pt},
	seloneAAAlAAAAA/.style = {scale=.8,fill=magenta!20, inner sep=0pt,outer sep=4pt},
	seltwoAAAlAAAAA/.style = {scale=.8,fill=green!20, inner sep=0pt,outer sep=4pt},
	noselAAAlAAAAA/.style = {scale=.8,inner sep=0pt,outer sep=4pt},
	seloneDDDDD/.style = {scale=.8,fill=magenta!20, inner sep=0pt,outer sep=4pt},
	seltwoDDDDD/.style = {scale=.8,fill=green!20, inner sep=0pt,outer sep=4pt},
	noselDDDDD/.style = {scale=.8,inner sep=0pt,outer sep=4pt},
	seloneDDDDDD/.style = {scale=.8,fill=magenta!20, inner sep=0pt,outer sep=4pt},
	seltwoDDDDDD/.style = {scale=.8,fill=green!20, inner sep=0pt,outer sep=4pt},
	noselDDDDDD/.style = {scale=.8,inner sep=0pt,outer sep=4pt},
	seloneDDDDDDD/.style = {scale=.8,fill=magenta!20, inner sep=0pt,outer sep=6pt},
	seltwoDDDDDDD/.style = {scale=.8,fill=green!20, inner sep=0pt,outer sep=6pt},
	noselDDDDDDD/.style = {scale=.8,inner sep=0pt,outer sep=6pt},
	seloneDDDDDDDD/.style = {scale=.7,fill=magenta!20, inner sep=0pt,outer sep=6pt},
	seltwoDDDDDDDD/.style = {scale=.7,fill=green!20, inner sep=0pt,outer sep=6pt},
	noselDDDDDDDD/.style = {scale=.7,inner sep=0pt,outer sep=6pt},
	seloneDDDDDDDDD/.style = {scale=.7,fill=magenta!20, inner sep=0pt,outer sep=6pt},
	seltwoDDDDDDDDD/.style = {scale=.7,fill=green!20, inner sep=0pt,outer sep=6pt},
	noselDDDDDDDDD/.style = {scale=.7,inner sep=0pt,outer sep=6pt},
	seloneDDDDDDDDDD/.style = {scale=.7,fill=magenta!20, inner sep=0pt,outer sep=6pt},
	seltwoDDDDDDDDDD/.style = {scale=.7,fill=green!20, inner sep=0pt,outer sep=6pt},
	noselDDDDDDDDDD/.style = {scale=.7,inner sep=0pt,outer sep=6pt},
	seloneDDDDDDDDDDD/.style = {scale=.7,fill=magenta!20, inner sep=0pt,outer sep=6pt},
	seltwoDDDDDDDDDDD/.style = {scale=.7,fill=green!20, inner sep=0pt,outer sep=6pt},
	noselDDDDDDDDDDD/.style = {scale=.7,inner sep=0pt,outer sep=6pt},
	seloneEEEEEEE/.style = {scale=.8,fill=magenta!20, inner sep=0pt,outer sep=3pt},
	seltwoEEEEEEE/.style = {scale=.8,fill=green!20, inner sep=0pt,outer sep=3pt},
	noselEEEEEEE/.style = {scale=.8,inner sep=0pt,outer sep=3pt},
	seloneEEEEEEEE/.style = {scale=.8,fill=magenta!20, inner sep=0pt,outer sep=4pt},
	seltwoEEEEEEEE/.style = {scale=.8,fill=green!20, inner sep=0pt,outer sep=4pt},
	noselEEEEEEEE/.style = {scale=.8,inner sep=0pt,outer sep=4pt},
	seloneEEEEEEEEE/.style = {scale=.8,fill=magenta!20, inner sep=0pt,outer sep=4pt},
	seltwoEEEEEEEEE/.style = {scale=.8,fill=green!20, inner sep=0pt,outer sep=4pt},
	noselEEEEEEEEE/.style = {scale=.8,inner sep=0pt,outer sep=4pt}
}

\usepackage{environ}
\makeatletter
\newsavebox{\measure@tikzpicture}
\NewEnviron{scaletikzpicturetowidth}[1]{%
\def\tikz@width{#1}%
\def\tikzscale{1}\begin{lrbox}{\measure@tikzpicture}%
\BODY
\end{lrbox}%
\pgfmathparse{#1/\wd\measure@tikzpicture}%
\edef\tikzscale{\pgfmathresult}%
\BODY
}
\makeatother

\allowdisplaybreaks

\newsavebox{\boxAlAA}
\newsavebox{\boxAlAAA}
\newsavebox{\boxAlAAAA}
\newsavebox{\boxAAlAA}
\newsavebox{\boxAAlAAA}
\newsavebox{\boxAAlAAAA}
\newsavebox{\boxAAAlAAA}
\newsavebox{\boxAAAlAAAA}
\newsavebox{\boxAAAlAAAAA}
\newsavebox{\boxDDDDD}
\newsavebox{\boxDDDDDD}
\newsavebox{\boxDDDDDDD}
\newsavebox{\boxDDDDDDDD}
\newsavebox{\boxDDDDDDDDD}
\newsavebox{\boxDDDDDDDDDD}
\newsavebox{\boxDDDDDDDDDDD}
\newsavebox{\boxEEEEEEE}
\newsavebox{\boxEEEEEEEE}
\newsavebox{\boxEEEEEEEEE}
\newsavebox{\boxDm}

\newcommand{\Mod}{\text{\rm mod-}}

\newcommand{\D}{\widetilde{\mathbb{D}}}
\newcommand{\A}{\widetilde{\mathbb{A}}}
\newcommand{\E}{\widetilde{\mathbb{E}}}

\newcommand{\dimz}{\underline\dim}
\newcommand{\End}{\operatorname{End}}
\newcommand{\Hom}{\operatorname{Hom}}
\newcommand{\Homk}{\operatorname{Hom}_k}
\newcommand{\Ext}{\operatorname{Ext}}

\newcommand{\Ima}{\operatorname{Im}}
\newcommand{\Ker}{\operatorname{Ker}}

\providecommand{\corollaryname}{Corollary}
\providecommand{\lemmaname}{Lemma}
\providecommand{\theoremname}{Theorem}
\providecommand{\propositionname}{Proposition}
\providecommand{\definitionname}{Definition}

\begin{document}
\title{Schofield sequences in the Euclidean case}
\author{Csaba Sz\'ant\'o\thanks{\Envelope{} \protect\url{szanto.cs@gmail.com}, \textbf{Faculty of
Mathematics and Computer Science, Babe\textcommabelow{s}-Bolyai University}
(400084 Cluj-Napoca, str. M. Kog\u{a}lniceanu, nr. 1, Romania)}, Istv\'an Sz\"oll\H{o}si\thanks{\Envelope{} \protect\url{szollosi@gmail.com} (corresponding author),
\textbf{Faculty of Mathematics and Computer Science, Babe\textcommabelow{s}-Bolyai
University} (400084 Cluj-Napoca, str. M. Kog\u{a}lniceanu, nr. 1,
Romania), \textbf{E\"otv\"os Lor\'and University, Faculty of Informatics}
(H-1117 Budapest, P\'azm\'any P. sny 1/C, Hungary)}}
\maketitle
\begin{abstract}
Let $k$ be a field and consider the path algebra $kQ$ of the quiver
$Q$. A pair of indecomposable $kQ$-modules $(Y,X)$ is called an orthogonal
	exceptional pair if the modules are exceptional and $\Hom(X,Y)=\Hom(Y,X)=\Ext^{1}(X,Y)=0$.
Denote by $\mathcal{F}(X,Y)$ the full subcategory of objects having a 
filtration with factors $X$ and $Y$. By a theorem of Schofield if $Z$ is exceptional but not simple, then $Z\in\mathcal{F}(X,Y)$
for some orthogonal exceptional pair $(Y,X)$, and $Z$ is not a simple
object in $\mathcal{F}(X,Y)$. In fact, there are precisely $s(Z)-1$
such pairs, where $s(Z)$ is the support of $Z$ (i.e. the number of nonzero components
in ${\dimz}Z$).  Whereas it is easy to construct
$Z$ given $X$ and $Y$, there is no convenient procedure yet to
determine the possible modules $X$ (called Schofield submodules of $Z$) and then $Y$ (called Schofield factors of $Z$), when $Z$ is given. We present such an explicit  procedure in the tame case, i.e. when $Q$ is Euclidean.
\end{abstract}

\medskip


\medskip


\section{Introduction}
Let $k$ be a field and consider the path algebra $kQ$ of the quiver
$Q$. A pair of indecomposable $kQ$-modules $(Y,X)$ is called an \emph{orthogonal
	exceptional pair} if the modules are exceptional (i.e. without self extension) and $\Hom(X,Y)=\Hom(Y,X)=\Ext^{1}(X,Y)=0$.
Denote by $\mathcal{F}(X,Y)$ the full subcategory of objects having a 
filtration with factors $X$ and $Y$. Observe that $\mathcal{F}(X,Y)$
is an exact, hereditary, abelian subcategory equivalent to the category
of finite dimensional $k$-representations of the quiver having two
vertices and $d=\dim_k\Ext^{1}(Y,X)$ identically oriented arrows.

We know due to Hubery (see \cite{Hubery}) that in the tame case we
have $\dim_k\Ext^{1}(Y,X)\leq2$; moreover, if equality holds then
${\dimz}(X\oplus Y)=\delta$ and $\partial Y=1$. Thus in
the case when $\dim_k\Ext^{1}(Y,X)=2$ we have $X=P$ indecomposable
preprojective of defect $-1$ and $Y=I$ indecomposable preinjective
having dimension vector $\delta-{\dimz}P$ and defect $1$. This pair
$(I,P)$ is then called a \emph{Kronecker pair} since the category $\mathcal{F}(P,I)$
is equivalent to the category of finite dimensional $k$-representations
of the Kronecker quiver $K$.

The following theorem by Schofield (see \cite{Ringel2,Schofield})
makes it possible to construct exceptional modules as extensions of smaller
exceptional ones. This procedure is generally called \emph{Schofield
	induction} and it is a very useful tool in representation theory of algebras.
\begin{thm}
	\label{thm:Schofield}(Schofield, \cite{Hubery,Ringel2,Schofield})
	If $Z$ is exceptional but not simple, then $Z\in\mathcal{F}(X,Y)$
	for some orthogonal exceptional pair $(Y,X)$, and $Z$ is not a simple
	object in $\mathcal{F}(X,Y)$. In fact, there are precisely $s(Z)-1$
	such pairs, where $s(Z)$ is the support of $Z$ (i.e.  $s(Z)=\sharp \{i\in Q_0\ | \ Z_i\neq 0\}$). 
\end{thm}

Note that in the theorem above the condition requiring $Z$ not to
be a simple object in $\mathcal{F}(Y,X)$ is equivalent to the existence
of an exact sequence of the form ~$\xymatrix{0\ar[r] & uX\ar[r] & Z\ar[r] & vY\ar[r] & 0}
$ with $u,v$ positive and uniquely determined by $X,Y,Z$.

\begin{defn}
Consider an orthogonal exceptional pair $(Y,X)$. A short exact sequence of the form $\xymatrix{0\ar[r] & uX\ar[r] & Z\ar[r] & vY\ar[r] & 0}$ with $u,v$ positive integers and $Z$ exceptional is called a \emph{Schofield sequence} and the pair $(Y,X)$ a \emph{Schofield pair} associated to $Z$. 
\end{defn}

As Ringel states in \cite{Ringel3} whereas it is easy to construct
$Z$ given $X$ and $Y$, there is no convenient procedure yet to
determine the possible modules $X$ (called \emph{Schofield submodules}
of $Z$) and then $Y$ (called \emph{Schofield factors of $Z$}), when $Z$ is given. In the Dynkin case Bo
Chen's Theorem provides a method to find at least some of these modules
$X$, namely the Gabriel-Roiter submodules of $Z$ (see \cite{Chen, Ringel3}).

In the present paper we describe all the tame Schofield sequences
using some numerical criteria for the characterization of these sequences.
As a result we give in the tame cases an explicit procedure to obtain
all the possible Schofield submodules and factor modules for a given exceptional module.
We will notice that all the tame Schofield sequences are field independent,
thus of purely combinatorial nature, depending only on the oriented
quiver.

\section{Tame quivers and roots}

We begin with a short survey of the problem context. For further details
we refer to \cite{Skow1,Aus,Kac,Skow2}.

Let $Q=(Q_{0},Q_{1})$ be a simply-laced tame quiver without oriented
cycles. This means $Q$ is a directed acyclic graph with vertex set
$Q_{0}$ and arrow set $Q_{1}$ and with its underlying undirected
graph a so-called Euclidean (affine) diagram of type $\A_{m}$, $\D_{m}$, $\E_{6}$, $\E_{7}$, $\E_{8}$.
Suppose that the vertex set $Q_{0}$ has $n$ elements and for an
arrow $\alpha\in Q_{1}$ we denote by $t(\alpha),h(\alpha)\in Q_{0}$
the tail, respectively the head of $\alpha$.

The \emph{Euler form} of $Q$ is a bilinear form on $\mathbb{Z}Q_{0}\cong\mathbb{Z}^{n}$
given by $\langle x,y\rangle=\sum_{i\in Q_{0}}x_{i}y_{i}-\sum_{\alpha\in Q_{1}}x_{t(\alpha)}y_{h(\alpha)}$.
Its quadratic form $q_{Q}$ (called \emph{Tits form}) is independent
of the orientation of $Q$ and in the tame case it is positive semi-definite
with radical $\{a\in\mathbb{Z}Q_{0}|q_{Q}(a)=0\}=\mathbb{Z}\delta$,
where $\delta$ is called \emph{minimal positive imaginary root} of
the corresponding Kac-Moody root system (see \cite{Kac}). The \emph{defect}
of $x\in\mathbb{Z}Q_{0}$ is then $\partial x=\langle\delta,x\rangle$,
the \emph{absolute defect} being the absolute value $|\partial x|$.
The vectors $a\in\mathbb{N}Q_{0}$ with $q_{Q}(a)=1$ are called \emph{positive
real roots}.

Below is the list of all the Euclidean diagrams marking each vertex
$i\in Q_{0}$ with the component $\delta_{i}$ (the number of vertices
being $n$).

Type $\A_m$, with $m\geq1$ (the case $m=0$ having only
cyclic orientation) and $n=m+1$:

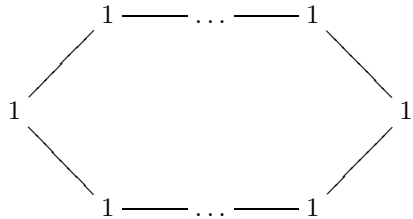
\[
\xymatrix{ & 1\ar@{-}[r] & \dots\ar@{-}[r] & 1\ar@{-}[dr]\\
1\ar@{-}[ur]\ar@{-}[dr] &  &  &  & 1\\
 & 1\ar@{-}[r] & \dots\ar@{-}[r] & 1\ar@{-}[ur]
}
\]

The \emph{Kronecker quiver} is the particular quiver of type $\A_{1}$
with non cyclic orientation 
\[
\xymatrix{{\bullet} & {\bullet}\ar@<1ex>[l]\ar@<-1ex>[l]}
.
\]

Type $\D_{m}$, with $m\geq4$ and $n=m+1$:

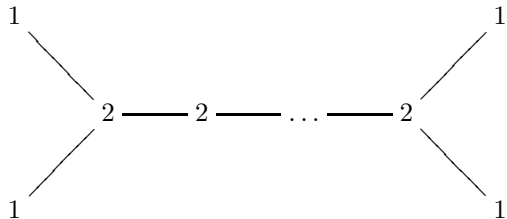
\[
\xymatrix{1\ar@{-}[dr] &  &  &  &  & 1\\
 & 2\ar@{-}[r] & 2\ar@{-}[r] & \dots\ar@{-}[r] & 2\ar@{-}[ur]\\
1\ar@{-}[ur] &  &  &  &  & 1\ar@{-}[ul]
}
\]

Type $\E_{6}$, with $n=7$: 
\[
\xymatrix{ &  & 1\ar@{-}[d]\\
 &  & 2\ar@{-}[d]\\
1\ar@{-}[r] & 2\ar@{-}[r] & 3\ar@{-}[r] & 2\ar@{-}[r] & 1
}
\]

Type $\E_{7}$, with $n=8$: 
\[
\xymatrix{ &  &  & 2\ar@{-}[d]\\
1\ar@{-}[r] & 2\ar@{-}[r] & 3\ar@{-}[r] & 4\ar@{-}[r] & 3\ar@{-}[r] & 2\ar@{-}[r] & 1
}
\]

Type $\E_{8}$, with $n=9$: 
\[
\xymatrix{ &  & 3\ar@{-}[d]\\
2\ar@{-}[r] & 4\ar@{-}[r] & 6\ar@{-}[r] & 5\ar@{-}[r] & 4\ar@{-}[r] & 3\ar@{-}[r] & 2\ar@{-}[r] & 1
}
\]

\section{Tame hereditary algebras. Representations of tame quivers}

Let $k$ be a field and consider the path algebra $kQ$ of the quiver
$Q$. The algebra $kQ$ is a finite dimensional
\emph{tame hereditary algebra} and in fact all connected elementary
tame hereditary algebras have this form.

The category of finite dimensional right modules (which is abelian
and Krull--Schmidt) is denoted by $\Mod kQ$. For a module $M\in\Mod kQ$
, $[M]$ denotes its isomorphism class, $|M|$ its length and $uM=M\oplus\dots\oplus M$ ($u$-times).

The category $\Mod kQ$ can and will be identified with the category
$\text{{\rm rep-}}kQ$ of  finite dimensional $k$-representations
of the quiver $Q$. Recall that a $k$-representation of $Q$ is defined
as a set of finite dimensional $k$-spaces $\{V_{i}|i\in Q_{0}\}$
corresponding to the vertices, together with $k$-linear maps $V_{\alpha}:V_{t(\alpha)}\to V_{h(\alpha)}$
corresponding to the arrows $\alpha\in Q_{1}$. The dimension vector of a
module $M=(V_{i},V_{\alpha})\in\Mod kQ=\text{{\rm rep-}}kQ$ is then
${\dimz}M=(\dim_kV_{i})_{i\in Q_{0}}\in\mathbb{Z}Q_{0}$. For two (dimension) vectors $d,d'\in \mathbb{Z}Q_0$ we say that $d\leq d'$ if $d_i \leq d'_i$ for all $i\in Q_0$.

Let $S(i)$, $P(i)$ and $I(i)$ be the indecomposable simple, projective
and injective modules corresponding to the vertex $i$ and consider
the \emph{Cartan matrix of the algebra $kQ$} $C_{Q}$ with the $j$-th column being ${\dimz}P(j)$  (see Chapter III of \cite{Skow1}).
The \emph{Coxeter matrix} is defined as $\Phi_{Q}=-C_{Q}^{t}C_{Q}^{-1}$.
Then $\Phi_{Q}\delta=\delta$ and the Euler form satisfies $\langle a,b\rangle=aC_{Q}^{-t}b^{t}=-\langle b,\Phi_{Q}a\rangle$,
where $a,b\in\mathbb{Z}Q_{0}$. Moreover, (because our algebra is hereditary)
for two modules $M,N\in\Mod kQ$ we get 
\[
\langle{\dimz}M,{\dimz}N\rangle=\dim_k\Hom(M,N)-\dim_k\Ext^{1}(M,N).
\]
Let $\partial M=\partial(\dimz M)=\langle\delta,{\dimz}M\rangle=-\langle{\dimz}M,\delta\rangle$
be the \emph{defect} of the module $M$.

Consider the Auslander--Reiten translates $\tau=D\Ext^{1}(-,kQ)$ and
$\tau^{-1}=\Ext^{1}(D(kQ),-)$, where $D=\Homk(-,k)$.
An indecomposable module $M$ is \emph{preprojective} (\emph{preinjective})
if there exists a positive integer $m$ such that $\tau^{m}(M)=0$
($\tau^{-m}(M)=0$). Otherwise $M$ is said to be \emph{regular}.
Note that an indecomposable module $M$ is preprojective (preinjective,
regular) if and only if $\partial M<0$ ($\partial M>0$, $\partial M=0$).
A module is called preprojective (preinjective, regular) if all its
indecomposable components are preprojective (preinjective, regular).
We will use the notation $P$, $I$ and $R$ in these cases.

The category $\Mod kQ$ is well known, its Auslander--Reiten quiver
being completely described. The indecomposable modules are either
preprojective, preinjective or regular. Thus up to isomorphism the
indecomposable preprojective modules are $P(m,i)=\tau^{-m}P(i)$ and the indecomposable preinjectives are  $I(m,i)=\tau^{m}I(i)$ where $m\in\mathbb{N}$ and $i\in Q_{0}$. Note
that ${\dimz}P(m,i)=\Phi_{Q}^{-m}{\dimz}P(i)$
and $\partial P(m,i)=\partial P(i)=-\delta_{i}$, respectively  ${\dimz}I(m,i)=\Phi_{Q}^{m}{\dimz}I(i)$
and $\partial I(m,i)=\partial I(i)=\delta_{i}$. Thus the possible
defects of indecomposable preprojectives are $-1$ in the $\A_{m}$
case, $-1$ or $-2$ in the $\D_{m}$ case, $-1$, $-2$ or
$-3$ in the $\E_{6}$ case, ranging from $-4$ to $-1$ in
the $\E_{7}$ case and ranging from $-6$ to $-1$ in the $\E_{8}$
case. Indecomposable preprojectives and preinjectives are \emph{exceptional}
modules (i.e. they are indecomposable modules without self extension
and with one dimensional endomorphism space), are \emph{directing}
(see \cite{Skow1} pages 357, 358 for details) and also uniquely determined
up to isomorphism by their dimension vectors which are positive real
roots, so they can be considered independently from the field $k$.

The category of regular modules is an abelian, exact subcategory which
decomposes into a direct sum of serial categories with Auslander--Reiten
quiver of the form $\mathbb ZA_{\infty}/m$, called \emph{tubes of rank} $m$. These tubes are indexed by the scheme theoretic, closed points of the projective line $\mathbb{P}_{k}^{1}$,
the degree of such a point being $d_{a}=[\kappa(a): k]$ (where $\kappa(a)$ is the residue field at the point $a$).
A tube of rank 1 is called \emph{homogeneous}, otherwise it is called
\emph{non-homogeneous} or \emph{exceptional}. We have at most 3 non-homogeneous
tubes indexed by points $a$ of degree $d_{a}=1$. All the other tubes
are homogeneous. We assume that the non-homogeneous tubes are labeled
by the elements of some subset $E$ of $\{0,1,\infty\}$, whereas the homogeneous
tubes are labeled by the points of the scheme $\mathbb{H}_{k}=\mathbb{H}_{\mathbb{Z}}\otimes k$
for some open subscheme $\mathbb{H}_{\mathbb{Z}}\subset\mathbb{P}_{\mathbb{Z}}^{1}$. For $x\in\mathbb{H}_{k}$ the indecomposables on a homogeneous
tube $\mathcal{T}_{x}$ are denoted by $R_{x}(1)\subset R_{x}(2)\subset\dots$.
Note that these modules depend on the field $k$ (sometimes this being
emphasized by the notation $R_{x}^{k}(t)$) and ${\dimz}R_{x}(t)=td_{x}\delta$.
The module $R_{x}(1)$ is called \emph{homogeneous regular simple}
and $t$ is called the \emph{regular length} of $R_{x}(t)$. The homogeneous
regular simple $R_{x}(1)$ is the \emph{regular socle} and also the \emph{regular top}
of $R_{x}(t)$. For a partition $\lambda=(\lambda_{1},\dots,\lambda_{n})$
let $R_{x}(\lambda)=R_{x}(\lambda_{1})\oplus\dots\oplus R_{x}(\lambda_{n})$.
The indecomposable modules on the mouth of the non-homogeneous tube
$\mathcal{T}_{e}$ of rank $m$ (where $e\in E$) are labeled $R_{e}^{l}(1)$
($l\in\{1,\dots,m\}$) and are called \emph{non-homogeneous regular
simples}. We have $\tau R_{e}^{l}(1)=R_{e}^{l-1}(1)$ for $l\geq2$,
$\tau R_{e}^{1}(1)=R_{e}^{m}(1)$ and $\sum_{l=1}^{m}{\dimz}R_{e}^{l}(1)=\delta$.
We denote by $R_{e}^{l}(t)$ the non-homogeneous indecomposable regular
with \emph{regular socle} $R_{e}^{l}(1)$ and \emph{regular length} $t$
(and \emph{regular top} $R_{e}^{l}(t)/R_{e}^{l}(t-1)=R_{e}^{l'}(1)$).
Note that $R_{e}^{l}(t)$ is uniquely determined by the values $e$,
$l$ and $t$ thus can be treated field independently; moreover, ${\dimz}R_{e}^{l}(t)$
is a positive real root (uniquely corresponding to the module) unless
$t$ is a multiple of $m$. Finally note that ${\dimz}R_{e}^{l}(sm)=s\delta$.

We should remark here that for any positive real root $a$ we have
a unique indecomposable module of dimension vector $a$ (see \cite{Kac,Dlab1}).

We list below some known facts on morphisms and
extensions (see \cite{Szanto} for all the details). 
\begin{lem}
\label{lem:basic}
\begin{enumerate}
\item[i)] For $P$ preprojective, $I$ preinjective, $R$ regular module we
have $\Hom({R},{P})=\Hom({I},{P})=\Hom({I},{R})=\Ext^{1}({P},{R})=\Ext^{1}({P},{I})=\Ext^{1}({R},{I})=0.$
\item[ii)] If $x\neq x'$ and $R_{x}$ (respectively $R_{x'}$) is a regular
with components from the tube $\mathcal{T}_{x}$ (respectively $\mathcal{T}_{x'}$),
then $\Hom(R_{x},R_{x'})=\Ext^{1}(R_{x},R_{x'})=0.$
\item[iii)]  For $\mathcal{T}_{x}$ homogeneous and $R_{x}(t)$, $R_{x}(t')$
indecomposables from $\mathcal{T}_{x}$ we have $\dim_k\Hom(R_{x}(t),R_{x}(t'))=\dim_k\Ext^{1}(R_{x}(t),R_{x}(t'))=\min(t,t')d_{x}$.
\item[iv)]  For $\mathcal{T}_{e}$ non-homogeneous of rank $m$ and $R_{e}^{l}(t)$
an indecomposable from $\mathcal{T}_{e}$ we have $\dim_k\End(R_{e}^{l}(t))=s+1$
for $sm<t\leq(s+1)m$ and $\dim_k\Ext^{1}(R_{e}^{l}(t),R_{e}^{l}(t))=s$
for $sm\leq t<(s+1)m$, thus $R_{e}^{l}(t)$ is exceptional if and
only if $0<t<m$. Moreover, in this case $R_{e}^{l}(t)$ is uniquely
determined by its unique regular composition series, with distinct
regular composition factors.
\item[v)]  For $P$ indecomposable preprojective and $I$ indecomposable preinjective
modules we have $\dim_k\End(P)=\dim_k\End(I)=1$ and $\Ext^{1}(P,P)=\Ext^{1}(I,I)=0$,
so they are exceptional modules. Moreover, they are also directing.
\end{enumerate}
\end{lem}
\begin{lem} \label{lem:preproj} (\cite{Szanto}) Let $P$ be a preprojective indecomposable with defect
	$\partial P=-1$, $P'$ a preprojective module and $R$ a regular
	indecomposable. Then we have:
\begin{enumerate}
    \item[i)] Every nonzero morphism $f:P\to P'$ is a monomorphism.
    \item[ii)] For every nonzero morphism $f:P\to R$, $f$ is either a
	monomorphism or $\Ima f$ is regular. In particular, if $R$ is
	regular simple and $\Ima f$ is regular then $f$ is an epimorphism.
	\item[iii)]Suppose that $\underline\dim P>\delta$. Then $P$ projects to a unique	regular simple $R^P_e(1)$ from the mouth of each non-homogeneous tube
	$\mathcal{T}_e$; moreover, $\dim_k\Hom(P,R^P_e(1))=\langle \dimz P, \dimz R^P_e(1)\rangle=1$.
	\item[iv)] Suppose that $\underline\dim P<\delta$. Then, depending on its dimension vector, $P$ embeds in or projects onto 
	a unique regular simple $R^P_e(1)$ from each non-homogeneous
	tube $\mathcal{T}_e$; moreover, $\dim_k\Hom(P,R^P_e(1))=\langle \dimz  P, \dimz R^P_e(1)\rangle=1$.
\end{enumerate}
\end{lem}
\begin{lem}\label{lem:preinjextension}(\cite{Szanto}) Let $I$ be an indecomposable preinjective of
		defect 1, $P$ an indecomposable preprojective of defect $-1$ and
		$X\ncong P\oplus I$. There exists a short exact sequence $0\to P\to X\to I\to 0$ iff $X$ satisfies the
		following conditions:
\begin{enumerate}
    \item[i)] $X$ is a regular module with $\underline\dim
		X=\underline\dim I+\underline\dim P$;
	\item[ii)] if $X$ has an indecomposable component from a non-homogeneous tube
		$\mathcal{T}_e$, then the regular top of this component is $R^{P}_e(1)$;
	\item[iii)] the indecomposable components of $X$ are taken from
		pairwise different tubes.
\end{enumerate}
\end{lem}

We end this section with some facts on reflections.

Let $i$ be a sink in the quiver $Q$. Denote by $s_i$ the reflection induced by the vertex $i$ and by $\sigma_i Q$ the quiver obtained by reversing all arrows involving $i$. Let $\Mod kQ\langle i\rangle$ be the full subcategory of modules not containing the simple module corresponding to the vertex $i$ as a direct summand.

We consider the reflection functors $S^+_i:\Mod kQ\to\Mod k\sigma_iQ$ and  $S^-_i:\Mod k\sigma_i Q\to\Mod kQ$. For all details concerning reflection functors we refer the reader to \cite{bernstein,Skow1}. It is well known that for an indecomposable $M$ we have $S^+_iM\neq 0$ iff $M\ncong S(i)$, where $S(i)$ is the projective simple corresponding to the sink $i$. Moreover, in this case $S^+_iM$ is indecomposable and $\dimz S^+_iM=s_i(\dimz M)$. Also the functors $S^+_i$, $S^-_i$ induce quasi-inverse equivalences between $\Mod kQ\langle i\rangle$ and $\Mod k\sigma Q\langle i\rangle$. It is easy to see that for $M\in\Mod kQ\langle i\rangle$ indecomposable we have that $\partial S^+_iM=\partial M$ and if $R$ is a simple homogeneous (respectively non-homogeneous) regular, then so is $S^+_iR$.

The following lemma is taken from \cite{Skow1}.
\begin{lem}\label{lem:reflections} Let $Q$, $Q'$ be trees having the same underlying graph. There exists a sequence $i_1,...,i_t$ of vertices of $Q$ such that for each $s\in\{1,\dots,t\}$ the vertex $i_s$ is a sink in $\sigma_{i_{s-1}}\dots\sigma_{i_1}Q$ and $\sigma_{i_{t}}\dots\sigma_{i_1}Q=Q'$.
\end{lem}

All the information above can be dualized for sources instead of sinks.

\section{Tame Schofield sequences}

We begin with a numerical criterion characterizing
tame Schofield sequences.
\begin{prop}
\label{prop:SchofieldCharaterization}Suppose $X$, $Y$, $Z$ are
exceptional indecomposables such that $u\dimz X+v\dimz Y=\dimz Z$.
Then we have a Schofield sequence 
\[
\xymatrix{0\ar[r] & uX\ar[r] & Z\ar[r] & vY\ar[r] & 0}
\]
if and only if $\left\langle \dimz X,\dimz Y\right\rangle =0$. Moreover,
in this case either $u=v=1$ or $\left|u-v\right|=1$ with $\partial X=-1$,
$\partial Y=1$, $\dimz X+\dimz Y=\delta$ and $\partial Z=\pm1$.
In the latter case we will speak about a \emph{special Schofield sequence}. Finally $u$, $v$ and thus $Y$ are uniquely determined by $X$, $Z$ and also $u$, $v$ and thus $X$ are uniquely determined by $Z$, $Y$ .
\end{prop}

\begin{proof}
The necessity is trivial. For the converse, note that we cannot
have simultaneously $\Hom(X,Y)\neq0$ and $\Ext^{1}(X,Y)\neq0$. By
Lemma \ref{lem:basic} this is trivial if $X$ and $Y$ are not of
the same type (preprojective, regular or preinjective). If both $X$
and $Y$ are preprojective (respectively preinjective), they are directing,
hence the assertion is true. If both $X$ and $Y$ are regular (non-homogeneous
exceptional) then we must have $u=v=1$ and also $\Hom(X,Y)=0$. This
holds because $Z$ being a regular exceptional, by Lemma \ref{lem:basic}
iii) it has a unique regular composition series with distinct regular
composition factors. Hence this is true also for $X$ and $Y$; moreover,
there is no common regular composition factor for $X$ and $Y$. In
case of a nonzero morphism $f:X\to Y$, $X$ would project onto the
regular module $f(X)\subseteq Y$, so the regular composition factors of $f(X)$ would appear in both $X$ and $Y$, which is impossible.

Since $0=\left\langle \dimz X,\dimz Y\right\rangle =\dim_{k}\Hom(X,Y)-\dim_{k}\Ext^{1}(X,Y)$,
it follows that $\Hom(X,Y)=\Ext^{1}(X,Y)=0$ using the previous observation.

Since $\dimz Z$, $\dimz X$ and $\dimz Y$ are all positive roots
and $u\dimz X+v\dimz Y=\dimz Z$ it follows that 
\begin{align*}
1 & =\left\langle \dimz Z,\dimz Z\right\rangle \\
 & =u^{2}\left\langle \dimz X,\dimz X\right\rangle +uv\left\langle \dimz X,\dimz Y\right\rangle +uv\left\langle \dimz Y,\dimz X\right\rangle +v^{2}\left\langle \dimz Y,\dimz Y\right\rangle \\
 & =u^{2}+v^{2}+uv\left\langle \dimz X,\dimz Y\right\rangle +uv\left\langle \dimz Y,\dimz X\right\rangle \\
 & =u^{2}+v^{2}+uv\left\langle \dimz Y,\dimz X\right\rangle .
\end{align*}
So $\left\langle \dimz Y,\dimz X\right\rangle =\frac{1-u^{2}-v^{2}}{uv}<0$
(here $u,v\ge1$), which implies that $\Hom(Y,X)=0$ (using the fact proven above that $\Hom(Y,X)$ and $\Ext^1(Y,X)$ can't both be nonzero) and  $\dim_k\Ext^{1}(Y,X)=\frac{-1+u^{2}+v^{2}}{uv}$.
So we already have that $X$ and $Y$ form an orthogonal exceptional
pair; therefore, by using Lemma 2 from \cite{Hubery1}
we get that the value of $\dim_{k}\Ext^{1}(Y,X)$ is either $1$ or
$2$. Moreover, if the value is 2 we must have $\partial X=-1$ and $\dimz X+\dimz Y=\delta$.

In the case when $\dim_k\Ext^{1}(Y,X)=1=\frac{-1+u^{2}+v^{2}}{uv}$,
we must have $u=v=1$, so $\dimz Z=\dimz X+\dimz Y$ and we know (by \cite{Hubery1}) that $\mathcal{F}(Y,X)$ is an exact, hereditary, abelian subcategory equivalent to the category of modules over the finite dimensional path algebra of the Dynkin quiver of type $A_1$ ($\xymatrix{{\bullet} & {\bullet}\ar[l]}$) over $k$ ($X,Y$ corresponding to the simple Dynkin modules). Knowing the Auslander-Reiten quiver of this Dynkin algebra and the fact that $\dimz Z=\dimz X+\dimz Y$ it follows that there is an Auslander-Reiten sequence $0\to X\to Z\to Y\to 0$, which will be our Schofield sequence.

In case $\dim_k\Ext^{1}(Y,X)=2$, we have $(u-v)^{2}=1$, so $\left|u-v\right|=1$,
$\partial X=-1$, $\dimz X+\dimz Y=\delta$. We also know (by \cite{Hubery1}) that in this case $\mathcal{F}(Y,X)$ is an exact, hereditary, abelian subcategory equivalent to the category of modules over the Kronecker algebra $kK$, where $K$ is the Kronecker quiver ($X,Y$ corresponding to the simple Kronecker modules). Knowing the Auslander-Reiten quiver of the Kronecker algebra and the fact that $\dimz Z=u\dimz X+v\dimz Y$ it follows that there is a short exact sequence $0\to uX\to Z\to vY\to 0$, which will be our Schofield sequence.

Note that $u$, $v$ and thus $Y$ are uniquely determined by $X$, $Z$. Indeed using the results above it is enough to consider the case $\partial X=-1$ and $\partial Z=\pm 1$. If we have two Schofield sequences, a special one $0\to (u\mp 1)X\to Z\to uY'\to 0$ with $\dimz Y'=\delta-\dimz X$ and a non-special one $0\to X\to Z\to Y\to 0$, then $\partial Y=2$ or 0 and
\begin{align*}
1 & =\left\langle \dimz Z,\dimz Z\right\rangle \\
&=\left\langle\dimz Z,(u\mp 1)\dimz X+u\dimz Y'\right\rangle\\
&=\left\langle\dimz X+\dimz Y,u\delta\mp\dimz X\right\rangle\\
& =u\left\langle \dimz X,\delta\right\rangle\mp\left\langle \dimz X,\dimz X\right\rangle +u\left\langle \dimz Y,\delta\right\rangle \mp\left\langle \dimz Y,\dimz X\right\rangle \\
& =u(1-\partial Y).
\end{align*}
It follows that either $u=-1$, a contradiction or $u=1$ which is again impossible (since then $\dimz Y=\dimz Z-\dimz X=\delta$).

Similarly follows that $u,v$ and thus $X$ are uniquely determined by $Z$ and $Y$.
\end{proof}

We can deduce from the proof above that:

\begin{cor} \label{cor:Schofield conditions}Suppose $X$, $Y$, $Z$ are
	exceptional indecomposables such that $\dimz X+\dimz Y=\dimz Z$. Then the following are equivalent:
	\begin{enumerate}
    
    \item[i)] we have a non-special Schofield sequence 
	\[
	\xymatrix{0\ar[r] & X\ar[r] & Z\ar[r] & Y\ar[r] & 0;}
	\]
     \item[ii)] $\dim_k\Ext^1(Y,X)=1$;
	 
	 \item[iii)] in case $X$ and $Z$ are not both regular $\dim_k\Hom(X,Z)=1$ and in case $X$ and $Z$ are both regular they share the same regular socle;
	 
	 \item[iv)] in case $Y$ and $Z$ are not both regular $\dim_k\Hom(Z,Y)=1$ and in case $Y$ and $Z$ are both regular they share the same regular top.
\end{enumerate}
\end{cor}
\begin{proof} $i)\Rightarrow ii)$  Use the proof of the previous proposition. 

$ii)\Rightarrow i)$ Note that $\dim_k\Ext^1(Y,X)=1$ implies $\left\langle \dimz Y,\dimz X\right\rangle =-1$, so because $1=\left\langle \dimz Z,\dimz Z\right\rangle =\left\langle \dimz X,\dimz X\right\rangle +\left\langle \dimz X,\dimz Y\right\rangle +\left\langle \dimz Y,\dimz X\right\rangle +\left\langle \dimz Y,\dimz Y\right\rangle=1+\left\langle \dimz X,\dimz Y\right\rangle$, we obtain $\left\langle \dimz X,\dimz Y\right\rangle =0$. 

$i)\Rightarrow iii)$ In case $X$,$Z$ are not both regular, from $0=\left\langle \dimz X,\dimz Y\right\rangle=\left\langle \dimz X,\dimz Z-\dimz X\right\rangle$ it follows that $\left\langle \dimz X,\dimz Z\right\rangle=\left\langle \dimz X,\dimz X\right\rangle=1$. Also by Lemma \ref{lem:basic} $\Hom(X,Z)$ and $\Ext^1(X,Z)$ can't be both nonzero, so the assertion follows.

In case $X$, $Z$ are both regular the assertion follows, since we have a monomorphism $X\to Z$.

$iii)\Rightarrow i)$ If $\dim_k\Hom(X,Z)=1$ and $X$, $Z$ are not both regular exceptional modules,  then by Lemma \ref{lem:basic} $\Ext^1(X,Z)=0$, which means that $1=\langle\dimz X,\dimz Z\rangle=\langle \dimz X,\dimz X+\dimz Y\rangle=1+\langle \dimz X,\dimz Y\rangle$, that is $\langle \dimz X,\dimz Y\rangle=0$. 

If both $X$ and $Z$ are regular sharing the same regular socle then we will have a unique regular short exact sequence (due to uniseriality) 
$0\to X\to Z\to Y\to 0$, thus we have the long Hom sequences
$$0\to\Hom(Z,X)\to\Hom(Z,Z)\to\Hom(Z,Y)\to\Ext^1(Z,X)\to\Ext^1(Z,Z)\to\Ext^1(Z,Y)\to 0,$$
$$0\to\Hom(Y,Y)\to\Hom(Z,Y)\to\Hom(X,Y)\to\Ext^1(Y,Y)\to\Ext^1(Z,Y)\to\Ext^1(X,Y)\to 0.$$
Knowing that the modules are exceptional we obtain that $\Ext^1(Z,Y)=\Ext^1(X,Y)=0$. Using the previous proof we have $\Hom(X,Y)=0$, so we get $\langle\dimz X,\dimz Y\rangle=0$.

Dually we obtain that $i)\Leftrightarrow iv)$.
\end{proof}

Concerning the special Schofield sequences we can prove the following:
\begin{prop}
\label{prop:SchofieldSpecial}Let $Z$ be an exceptional indecomposable
of defect $-1$ (i.e. preprojective) with $\dimz Z>\delta$. Then
it has a unique special Schofield sequence which is of the form 
\[
\xymatrix{0\ar[r] & (u+1)X\ar[r] & Z\ar[r] & uY\ar[r] & 0}
\]
 with $\partial X=-1$, $\dimz X+\dimz Y=\delta$ and thus $\dimz Z=u\delta +\dimz X$.

Dually, let $Z$ be an exceptional indecomposable of defect $1$ (i.e.
preinjective) with $\dimz Z>\delta$. Then it has a unique special
Schofield sequence which is of the form 
\[
\xymatrix{0\ar[r] & vX\ar[r] & Z\ar[r] & (v+1)Y\ar[r] & 0}
\]
 with $\partial X=-1$, $\dimz X+\dimz Y=\delta$ and thus $\dimz Z=v\delta +\dimz Y$.
\end{prop}

\begin{proof}
For the existence, write $\dimz Z$ as $\dimz Z=u\delta+x$, where
$x<\delta$. Since $\left\langle x,x\right\rangle =\left\langle \dimz Z-u\delta,\dimz Z-u\delta\right\rangle =1$,
$x$ is a positive real root, so there exists a unique indecomposable $X$
such that $\dimz X=x$; moreover, since $\partial x=\left\langle \delta,x\right\rangle =\left\langle \delta,\dimz Z-u\delta\right\rangle =\partial Z=-1$
it follows that $X$ is an exceptional preprojective (of defect $-1$).
Since $\delta-x$ is also a positive real root, let $Y$ be the unique
indecomposable with $\dimz Y=\delta-x$. We have that $\partial Y=1$,
hence $Y$ is a preinjective indecomposable, thus exceptional. Of
course we have $\dimz X+\dimz Y=\delta$; moreover, $\left\langle x,\delta-x\right\rangle =-\partial x-1=1-1=0$.
Therefore, $\left\langle \dimz X,\dimz Y\right\rangle =0$ and $\dimz Z=(u+1)\cdot\dimz X+u\cdot\dimz Y$
and by the previous theorem we have the required short exact sequence.

For the uniqueness, if there is another sequence $0\to(u'+1)X'\to Z\to u'Y'\to0$
fulfilling the requirements, then $\dimz Z=u'\delta+\dimz X'$ with
$\dimz X'<\delta$, $\partial X'=-1$, so $u'=u$ and $\dimz X'=x=\dimz X$.
\end{proof}

We will describe now some non-special Schofield sequences associated to $Z$ of defect $\pm 1$ or 0.
\begin{prop} If $Z$ is a preprojective indecomposable of defect $-1$, $R$ is an exceptional regular module having regular top $R^Z_e(1)$ and $\dimz Z>\dimz R$, then $R$ is a non-special Schofield factor of $Z$. A dual statement can be formulated for the case when $Z$ is preinjective of defect 1.
\end{prop}
\begin{proof}
By Lemma \ref{lem:preproj} we obtain that there is a projection $Z\to R^Z_e(1)$ and $\dim_k\Hom (Z,R^Z_e(1))=\dim_k\Hom (Z,R)=1$. If $f:Z\to R$ is a nonzero morphism then by Lemma \ref{lem:preproj} $\Ima f$ is regular exceptional. In case $\Ima f\neq R$, then the regular composition series of $\Ima f$ does not contain  $R^Z_e(1)$, so $\dim_k\Hom(Z,\Ima f)=0$, a contradiction. Therefore, $f$ must be an epimorphism, but then $\partial\Ker f=-1$, so $\Ker f$ is a preprojective indecomposable of defect $-1$. Applying Corollary \ref{cor:Schofield conditions} iv) we are done.
\end{proof}
\begin{cor}\label{cor: defect 1}
\begin{enumerate}
    \item[i)] In case $\partial Z=-1$ and $\dimz Z>\delta$, all its non-special Schofield factors are the regular exceptionals having regular top $R^Z_e(1)$ for $e\in E$. A dual statement is true for $\partial Z=1$.
    \item[ii)]If Q is of $\A_m$ type then the statement formulated in a) is valid for any indecomposable $Z$ of defect $\partial Z=\pm 1$.
\end{enumerate}
\end{cor}
\begin{proof} i) We know that we have $s(Z)-2=|Q_0|-2$ different non-special Schofield factors, since $\dimz Z>\delta$. But the number of the specified regular exceptionals matches this number (just look at the rank of the non-homogeneous tubes in all the cases, the exceptional regulars being all the non-homogeneous regular indecomposables of dimension vector less than $\delta$).
	
ii) Just note that in the $\A_m$ case the absolute defect is at most 1, so a Schofield factor in case of $\partial Z=-1$ must be regular exceptional with top $R^Z_e(1)$.
\end{proof}

\begin{prop} \label{prop:regular} Suppose $Z$ is an exceptional regular from the non-homogeneous tube $\mathcal{T}_e$ and $P$ is an indecomposable preprojective of defect $-1$. Then $P$ is a non-special Schofield submodule in $Z$ iff the regular top of $Z$ is $R^P_e(1)$ and $\dimz P<\dimz Z$.
\end{prop}

\begin{proof} Using Lemma \ref{lem:preinjextension} the necessity follows. For the sufficiency suppose the regular top of $Z$ is $R^P_e(1)$ and $\dimz P<\dimz Z$. Then $\dim_k\Hom (P,Z)=1$ and if $Z'\subset Z$ is the regular radical of $Z$ (i.e. $Z/Z'=R^P_e(1)$) we have that $\Hom(P,Z')=0$. By Lemma \ref{lem:preproj} this means that nonzero morphisms from $\Hom (P,Z)$ are monomorphisms. So consider $f:P\to Z$ such a monomorphism. Then $Z/P$ has defect 1 and we show that we can't have regular components in it. Indeed the projection from $Z$ to such a component would have a proper regular exceptional kernel containing $\Ima f=P$, which is impossible due to $\Hom(P,Z')=0$. So $Z/P$ is preinjective indecomposable and we are done.
\end{proof}
\begin{cor} \label{cor:regular} In the $\A_m$ case consider the regular exceptional $R_e(t)$ with regular top $R_0$. Then its	 non-special Schofield submodules  are the regulars $R_e(t')$ with $t'<t$ and the indecomposable preprojectives $P$ of defect $-1$, having dimension vector less than $\dimz R_e(t)$ and satisfying $R^P_e(1)=R_0$ (i.e. $\dim_k\Hom(P,R_0)=1$).
\end{cor}

By Corollaries \ref{cor: defect 1} and \ref{cor:regular} we already know all the Schofield pairs associated to exceptional modules in the $\A_m$ case. 

In the other tame cases (when the quiver is a tree) we will use reflections and AR-translations in order to obtain all the Schofield pairs over any orientation.

The following proposition is the result of the general compatibility of exact sequences with reflections.
\begin{prop}
\label{prop:Reflection} Consider the sink $i$, its reflection functor $S^+_i$ and suppose $X$ is an indecomposable module which is not the projective simple corresponding to the sink $i$. If $0\to uX\to Z\to vY\to 0$ is a Schofield sequence, then so is $0\to uS^+_iX\to S^+_iZ\to vS^+_iY\to 0$. The assertion remains valid also for the functor $S^-_i$, in the case $i$ is a source and $Y$ is not the injective simple corresponding to $i$. Moreover, we obtain the same result for the AR-translations $\tau$ (in the case $X$ is not an indecomposable projective) and $\tau^{-1}$ (in the case $Y$ is not an indecomposable injective)
\end{prop}

As we can see from the proposition above a Schofield sequence $0\to uX\to Z\to vY\to 0$ vanishes under the reflection $S^+_i$ in the case $X$ is the projective simple corresponding to the sink $i$. This might suggest that the set of all the Schofield pairs associated to $Z$ is not compatible with reflections. However, for non-special Schofield pairs the contrary is true.
\begin{prop}
\label{prop:ReflectionSupport} Consider the sink $i$ and its reflection functor $S^+_i$. If $\mathcal{M}=\{(Y_j,X_j)|j\in J\}$ is the set of all non-special Schofield pairs associated to $Z$, then $S^+_i\mathcal{M}=\{(S^+_iY_j,S^+_iX_j)|j\in J, S^+_iX_j\neq 0\}$ is the set of all non-special Schofield pairs associated to $S^+_iZ$. Dually if $i$ is a source and $\mathcal{M}=\{(Y_j,X_j)|j\in J\}$ is the set of all non-special Schofield pairs associated to $Z$, then $S^-_i\mathcal{M}=\{(S^-_iY_j,S^-_iX_j)|j\in J, S^-_iY_j\neq 0\}$ is the set of all non-special Schofield pairs associated to $S^-_iZ$.
\end{prop}
\begin{proof} Suppose $S^+_iX_j=0$ for some $j$. It follows that $X_j=S(i)$ is the projective simple corresponding to the sink $i$ and, moreover, due to Proposition \ref{prop:SchofieldCharaterization} we have a single Schofield pair with this second term. 
Suppose this pair $(Y_j,X_j)$ with $X_j=S(i)$ is not special. Since $\langle\dimz Y_j,\dimz S(i)\rangle=-1$, we have $(\dimz Y_j)_i-\sum_{t}(\dimz Y_j)_t=-1$ (where $t$ are the neighbours of $i$ in the quiver), so $(\dimz S^+_iZ)_i=-(\dimz Z)_i+\sum_{t}(\dimz Z)_t=0<(\dimz Z)_i$ and $(\dimz S^+_iZ)_l=(\dimz Z)_l$ (for $l\neq i$). This means that $s(S^+_iZ)=s(Z)-1$, so the number of Schofield pairs associated to $S^+_iZ$ (which equals $s(S^+_iZ)-1$ by Theorem \ref{thm:Schofield}) is one less than the number of Schofield pairs associated to $Z$ (which equals $s(Z)-1$). 
\end{proof}

\begin{thm}
\label{thm:SchofieldGeneralCanonical}In case of a canonically oriented
tame tree quiver, for a given exceptional module $Z$ all its Schofield pairs are listed in the Appendix (Section \ref{sec:Appendix}).
\end{thm}

\begin{proof}
Let $Q$ be a canonically oriented Euclidean quiver (the canonical
orientations correspond to those given in \cite{Skow2} and are also
shown in this document). We determine the dimension vectors of all
exceptional modules, falling below (at least) $3\delta$. Say $M$
is such an exceptional module, which is not simple. Then, we perform
an ``exhaustive search'' for non-special Schofield pairs $(Y,X)$
using the criterion $\left\langle \dimz X,\dimz Y\right\rangle =0$
from Proposition \ref{prop:SchofieldCharaterization}, finding in this way
exactly $s(Z)-1$ pairs (if $|\partial Z|>1$
or $|\partial Z|=1$ and $\dimz Z<\delta$), respectively $s(Z)-2$
pairs (if $|\partial Z|=1$ and $\dim Z>\delta$).
But this is exactly the number of different Schofield pairs possible
under Theorem \ref{thm:Schofield}, with the additional information
on special Schofield pairs provided by Proposition \ref{prop:SchofieldSpecial}. The cases $\dimz Z>3\delta$ will easily follow using AR-translations and Proposition \ref{prop:Reflection}.
\end{proof}
The computations required in the previous proof were carried out using
the computer algebra system GAP (see \cite{GAP}) but any other computer
algebra system, standalone implementation in any programming language,
or considerable amounts of paper, spare time and patience may be used
to verify the results, since only vector and matrix operations are
needed, over small integers.

\medskip

We are ready now to describe the procedure of obtaining all the Schofield pairs of any exceptional module $Z$ over any tame quiver not of type $\A_m$ (this type being already discussed). 

Consider the tame tree quiver $Q$ and denote by $Q'$ the canonically oriented quiver having the same underlying graph. 

{\bf Case 1} If $Z$ is preinjective or regular.

We know by Lemma \ref{lem:reflections} that there exists a sequence $i_1,\dots,i_t$ of vertices in $Q$ such that for each $s\in\{1,\dots,t\}$ the vertex $i_s$ is a sink in $\sigma_{i_{s-1}}\dots\sigma_{i_1}Q$ and $\sigma_{i_{t}}\dots\sigma_{i_1}Q=Q'$. We  perform the following steps:
\begin{description}
\item[{Step 1.}]  We fix the special Schofield pair associated to $Z$ (if $\dimz Z>\delta$) using Proposition \ref{prop:SchofieldSpecial}.
\item[{Step 2.}] To determine the non-special Schofield pairs associated to $Z$ consider the reflected exceptional module 
$S^+_{i_t}\dots S^+_{i_1} Z$ in the canonically oriented quiver $Q'=\sigma_{i_{t}}\dots\sigma_{i_1}Q$. Note that since $Z$ is not preprojective $S^+_{i_t}\dots S^+_{i_1} Z\neq 0$.
\item[{Step 3.}] Using Theorem \ref{thm:SchofieldGeneralCanonical} and the Appendix we determine the non-special Schofield pairs of $S^+_{i_t}\dots S^+_{i_1} Z$ in $Q'$.
\item[{Step 4.}] Reflect back these Schofield pairs using the functor $S^-_{i_1}\dots S^-_{i_t}$. By Proposition \ref{prop:ReflectionSupport} these pairs will be all the non-special Schofield pairs of $Z$.
\end{description}

{\bf Case 2} If $Z$ is preprojective.

We know by Lemma \ref{lem:reflections} that there exists a sequence $i_1,\dots,i_t$ of vertices in $Q$ such that for each $s\in\{1,\dots,t\}$ the vertex $i_s$ is a source in $\sigma_{i_{s-1}}\dots\sigma_{i_1}Q$ and $\sigma_{i_{t}}\dots\sigma_{i_1}Q=Q'$. We  perform the following steps:
\begin{description}
\item[{Step 1.}]  We fix the special Schofield pair associated to $Z$ (if $\dimz Z>\delta$) using Proposition \ref{prop:SchofieldSpecial}.
\item[{Step 2.}] To determine the non-special Schofield pairs associated to $Z$ consider the reflected exceptional module 
$S^-_{i_t}\dots S^-_{i_1} Z$ in the canonically oriented quiver $Q'=\sigma_{i_{t}}\dots\sigma_{i_1}Q$. Note that since $Z$ is not preinjective $S^-_{i_t}\dots S^-_{i_1} Z\neq 0$.
\item[{Step 3.}] Using Theorem \ref{thm:SchofieldGeneralCanonical} and the Appendix we determine the non-special Schofield pairs of $S^-_{i_t}\dots S^-_{i_1} Z$ in $Q'$.
\item[{Step 4.}] Reflect back these Schofield pairs using the functor $S^+_{i_1}\dots S^+_{i_t}$. By Proposition \ref{prop:ReflectionSupport} these pairs will be all the non-special Schofield pairs of $Z$.
\end{description}

\subsection*{Acknowledgment}

We are extremely grateful for the thorough and professional report
of the anonymous referee and for the advices which made the presentation
more precise.


\newpage
\section{Appendix}\label{sec:Appendix}
In this part we present the list of Schofield pairs of canonically oriented tame quivers (see Theorem \ref{thm:SchofieldGeneralCanonical}). Following \cite{Skow2} we denote by $\Delta(\A_{p,q})$, $\Delta(\D_{m})$, $\Delta(\E_{6})$, $\Delta(\E_{7})$, $\Delta(\E_{8})$ the canonically orientated quivers.

For every considered canonically oriented quiver $\Delta(Q)$, the corresponding data about the Schofield pairs is to be found in the subsection with the title  ``Schofield pairs for the quiver $\Delta(Q)$ -- $\delta = \dots$'' where $\delta$ is the minimal radical vector in the case of $Q$. Each subsection begins with a drawing of $\Delta(Q)$ to specify the node labeling and we also give the Cartan and Coxeter matrices denoted by $C_{\Delta(Q)}$, respectively $\Phi_{\Delta(Q)}$. Then there will be three subdivisions (``Schofield pairs associated to preprojective exceptional modules'', ``Schofield pairs associated to preinjective exceptional modules'' and ``Schofield pairs associated to regular exceptional modules'') containing the drawing of the corresponding part of the Auslander--Reiten quiver and the computed list of Schofield pairs associated to preprojective, preinjective, respectively regular exceptionals. On the graphical representation of the AR quiver, blue arrows show the existence
of a so-called irreducible monomorphism, while red arrows represent
the existence of irreducible epimorphisms between suitable indecomposable
modules (for details see \cite{Skow1}). In case of each indecomposable on the Auslander--Reiten quiver,
we show its dimension vector, since this is the only information needed
to determine the Schofield sequences.

For any given $M\in\Mod kQ$ exceptional with $s(M)=s$ we enlist
all the $s-1$ Schofield pairs associated to $M$ as 
\[
M:\ \bigl(Y_{1},X_{1}\bigr),\ \bigl(Y_{2},X_{2}\bigr),\dots \bigl(Y_{s-1},X_{s-1}\bigr),
\]
when there are no special pairs associated to $M$, respectively
\[
M:\ \bigl(Y_{1},X_{1}\bigr),\ \dots,\ \bigl(Y_{s-2},X_{s-2}\bigr),\  \bigl(uI,(u+1)P\bigr),
\]
when $M$ is a preprojective indecomposable with $\partial M=-1$ and
\[
M:\ \bigl(Y_{1},X_{1}\bigr),\ \dots,\ \bigl(Y_{s-2},X_{s-2}\bigr),\  \bigl((v+1)I,vP\bigr),
\]
when $M$ is a preinjective indecomposable with $\partial M=1$. Note that for every $M$ with $\partial M=\pm1$, its associated special Schofield pair may be obtained using Proposition \ref{prop:SchofieldSpecial} (therefore, in the general formulas we just mark their presence).

We group the modules $M$ according to the vertices of the Auslander--Reiten quiver. Thus for every vertex $i\in Q_0$ we will have groups entitled ``Modules of the form $P(n,i)$'' enlisting Schofield pairs associated to families of preprojective exceptionals and ``Modules of the form $I(n,i)$'' with Schofield pairs associated to families of preinjective exceptionals. The regular exceptional modules are grouped by their non-homogeneous tubes. On the drawings depicting non-homogeneous tubes, the exceptional regulars are marked with green background color. We have marked with pink background the preprojectives and the preinjectives beyond which all Schofield sequences may be obtained by Auslander--Reiten translation (inverse translation in the case of preprojectives).

Although using Corollaries \ref{cor: defect 1} and \ref{cor:regular} the Schofield pairs are completely described in the case of quiver of type $\mathbb{A}_m$, we compute and enlist concretely the Schofield pairs for the quivers $\Delta(\mathbb{A}_{1,2})$, $\Delta(\mathbb{A}_{1,3})$, $\Delta(\mathbb{A}_{1,4})$, $\Delta(\mathbb{A}_{2,2})$, $\Delta(\mathbb{A}_{2,3})$, $\Delta(\mathbb{A}_{2,4})$, $\Delta(\mathbb{A}_{3,3})$, $\Delta(\mathbb{A}_{3,4})$ and $\Delta(\mathbb{A}_{3,5})$. We also list the pairs for $\Delta(\mathbb{D}_{4})$, $\Delta(\mathbb{D}_{5})$ and $\Delta(\mathbb{D}_{6})$ as examples before giving the general formulas for $\Delta(\mathbb{D}_{m})$, $m\ge4$. The Appendix ends with the complete lists for $\Delta(\mathbb{E}_{6})$, $\Delta(\mathbb{E}_{7})$ and $\Delta(\mathbb{E}_{8})$.

The computations were performed by a software  implemented in the computer algebra system GAP (see \cite{GAP}), which also generated the \LaTeX{} source for this part (including the drawings of the various Auslander--Reiten quivers).

\begin{lrbox}{\boxAlAA}$\delta =\begin{smallmatrix}1&&1\\&1\end{smallmatrix} $\end{lrbox}
\subsection[Schofield pairs for the quiver $\Delta(\A_{1,2})$]{Schofield pairs for the quiver $\Delta(\A_{1,2})$ -- {\usebox{\boxAlAA}}}
\[
\vcenter{\hbox{\xymatrix{1 &   & 3\ar[ll]\ar[dl]\\
  & 2\ar[ul] &  }}}\qquad C_{\Delta(\widetilde{\mathbb{A}}_{1,2})} = \begin{bmatrix}1 & 1 & 2\\
0 & 1 & 1\\
0 & 0 & 1\end{bmatrix}\quad\Phi_{\Delta(\widetilde{\mathbb{A}}_{1,2})} = \begin{bmatrix}-1 & 1 & 1\\
-1 & 0 & 2\\
-2 & 1 & 2\end{bmatrix}
\]
\subsubsection{Schofield pairs associated to preprojective exceptional modules}

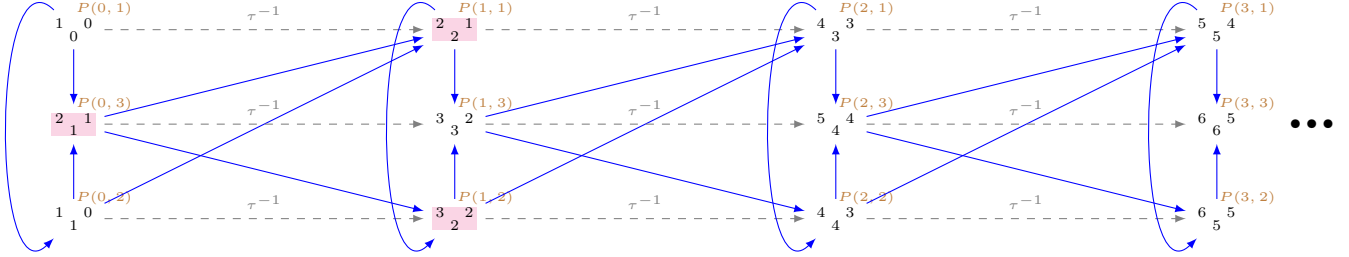
\begin{figure}[ht]


\begin{center}
\begin{scaletikzpicturetowidth}{\textwidth}

\end{scaletikzpicturetowidth}
\end{center}


\end{figure}
\medskip{}
\subsubsection*{Modules of the form $P(n,1)$}

Defect: $\partial P(n,1) = -1$, for $n\ge 0$.

\begin{fleqn}
\begin{align*}
P(0,1):\ & - \\
P(1,1):\ & \bigl(R_{0}^{1}(1),P(0,3)\bigr),\ \bigl(I(0,3),2P(0,2)\bigr)\\
P(n,1):\ & \bigl(R_{0}^{(n-1)\bmod 2+1}(1),P(n-1,3)\bigr),\ \bigl(uI,(u+1)P\bigr),\ n>1\\
\end{align*}
\end{fleqn}
\subsubsection*{Modules of the form $P(n,2)$}

Defect: $\partial P(n,2) = -1$, for $n\ge 0$.

\begin{fleqn}
\begin{align*}
P(0,2):\ & \bigl(R_{0}^{1}(1),P(0,1)\bigr)\\
P(1,2):\ & \bigl(R_{0}^{2}(1),P(1,1)\bigr),\ \bigl(2I(0,2),3P(0,1)\bigr)\\
P(n,2):\ & \bigl(R_{0}^{n\bmod 2+1}(1),P(n,1)\bigr),\ \bigl(uI,(u+1)P\bigr),\ n>1\\
\end{align*}
\end{fleqn}
\subsubsection*{Modules of the form $P(n,3)$}

Defect: $\partial P(n,3) = -1$, for $n\ge 0$.

\begin{fleqn}
\begin{align*}
P(0,3):\ & \bigl(R_{0}^{2}(1),P(0,2)\bigr),\ \bigl(I(0,2),2P(0,1)\bigr)\\
P(n,3):\ & \bigl(R_{0}^{(n+1)\bmod 2+1}(1),P(n,2)\bigr),\ \bigl(uI,(u+1)P\bigr),\ n>0\\
\end{align*}
\end{fleqn}
\subsubsection{Schofield pairs associated to preinjective exceptional modules}

\begin{figure}[ht]


\begin{center}
\begin{scaletikzpicturetowidth}{\textwidth}

\end{scaletikzpicturetowidth}
\end{center}

\end{figure}
\medskip{}
\subsubsection*{Modules of the form $I(n,1)$}

Defect: $\partial I(n,1) = 1$, for $n\ge 0$.

\begin{fleqn}
\begin{align*}
I(0,1):\ & \bigl(I(0,2),R_{0}^{2}(1)\bigr),\ \bigl(2I(0,3),P(0,2)\bigr)\\
I(n,1):\ & \bigl(I(n,2),R_{0}^{(-n+1)\bmod 2+1}(1)\bigr),\ \bigl((v+1)I,vP\bigr),\ n>0\\
\end{align*}
\end{fleqn}
\subsubsection*{Modules of the form $I(n,2)$}

Defect: $\partial I(n,2) = 1$, for $n\ge 0$.

\begin{fleqn}
\begin{align*}
I(0,2):\ & \bigl(I(0,3),R_{0}^{1}(1)\bigr)\\
I(1,2):\ & \bigl(I(1,3),R_{0}^{2}(1)\bigr),\ \bigl(3I(0,3),2P(0,2)\bigr)\\
I(n,2):\ & \bigl(I(n,3),R_{0}^{(-n+2)\bmod 2+1}(1)\bigr),\ \bigl((v+1)I,vP\bigr),\ n>1\\
\end{align*}
\end{fleqn}
\subsubsection*{Modules of the form $I(n,3)$}

Defect: $\partial I(n,3) = 1$, for $n\ge 0$.

\begin{fleqn}
\begin{align*}
I(0,3):\ & - \\
I(1,3):\ & \bigl(I(0,1),R_{0}^{1}(1)\bigr),\ \bigl(2I(0,2),P(0,1)\bigr)\\
I(n,3):\ & \bigl(I(n-1,1),R_{0}^{(-n+1)\bmod 2+1}(1)\bigr),\ \bigl((v+1)I,vP\bigr),\ n>1\\
\end{align*}
\end{fleqn}
\subsubsection{Schofield pairs associated to regular exceptional modules}

\subsubsection*{The non-homogeneous tube $\mathcal{T}_{0}^{\Delta(\widetilde{\mathbb{A}}_{1,2})}$}

\begin{figure}[ht]



\begin{center}
\captionof{figure}{\vspace*{-10pt}$\mathcal{T}_{0}^{\Delta(\widetilde{\mathbb{A}}_{1,2})}$}
\begin{scaletikzpicturetowidth}{0.40000000000000002\textwidth}
 $\end{lrbox}
\subsection[Schofield pairs for the quiver $\Delta(\A_{1,3})$]{Schofield pairs for the quiver $\Delta(\A_{1,3})$ -- {\usebox{\boxAlAAA}}}
\[
\vcenter{\hbox{\xymatrix{1 &   &   & 4\ar[lll]\ar[dl]\\
  & 2\ar[ul] & 3\ar[l] &  }}}\qquad C_{\Delta(\widetilde{\mathbb{A}}_{1,3})} = %
\]
\subsubsection{Schofield pairs associated to preprojective exceptional modules}

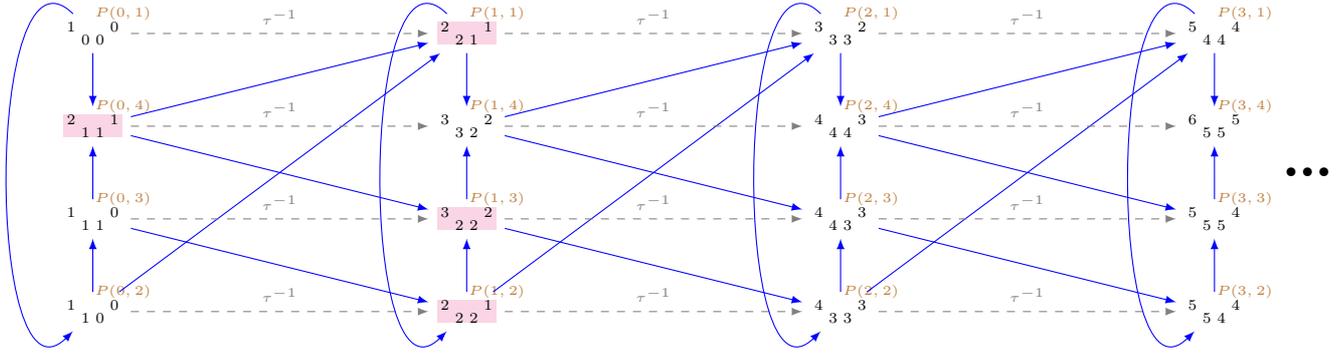
\begin{figure}[ht]


\begin{center}
\begin{scaletikzpicturetowidth}{\textwidth}

\end{scaletikzpicturetowidth}
\end{center}


\end{figure}
\medskip{}
\subsubsection*{Modules of the form $P(n,1)$}

Defect: $\partial P(n,1) = -1$, for $n\ge 0$.

\begin{fleqn}
\begin{align*}
P(0,1):\ & - \\
P(1,1):\ & \bigl(R_{0}^{2}(2),P(0,3)\bigr),\ \bigl(R_{0}^{3}(1),P(0,4)\bigr),\ \bigl(I(0,3),2P(0,2)\bigr)\\
P(n,1):\ & \bigl(R_{0}^{n\bmod 3+1}(2),P(n-1,3)\bigr),\ \bigl(R_{0}^{(n+1)\bmod 3+1}(1),P(n-1,4)\bigr),\ \bigl(uI,(u+1)P\bigr),\ n>1\\
\end{align*}
\end{fleqn}
\subsubsection*{Modules of the form $P(n,2)$}

Defect: $\partial P(n,2) = -1$, for $n\ge 0$.

\begin{fleqn}
\begin{align*}
P(0,2):\ & \bigl(R_{0}^{3}(1),P(0,1)\bigr)\\
P(1,2):\ & \bigl(R_{0}^{3}(2),P(0,4)\bigr),\ \bigl(R_{0}^{1}(1),P(1,1)\bigr),\ \bigl(I(0,4),2P(0,3)\bigr)\\
P(n,2):\ & \bigl(R_{0}^{(n+1)\bmod 3+1}(2),P(n-1,4)\bigr),\ \bigl(R_{0}^{(n-1)\bmod 3+1}(1),P(n,1)\bigr),\ \bigl(uI,(u+1)P\bigr),\ n>1\\
\end{align*}
\end{fleqn}
\subsubsection*{Modules of the form $P(n,3)$}

Defect: $\partial P(n,3) = -1$, for $n\ge 0$.

\begin{fleqn}
\begin{align*}
P(0,3):\ & \bigl(R_{0}^{3}(2),P(0,1)\bigr),\ \bigl(R_{0}^{1}(1),P(0,2)\bigr)\\
P(1,3):\ & \bigl(R_{0}^{1}(2),P(1,1)\bigr),\ \bigl(R_{0}^{2}(1),P(1,2)\bigr),\ \bigl(2I(0,2),3P(0,1)\bigr)\\
P(n,3):\ & \bigl(R_{0}^{(n-1)\bmod 3+1}(2),P(n,1)\bigr),\ \bigl(R_{0}^{n\bmod 3+1}(1),P(n,2)\bigr),\ \bigl(uI,(u+1)P\bigr),\ n>1\\
\end{align*}
\end{fleqn}
\subsubsection*{Modules of the form $P(n,4)$}

Defect: $\partial P(n,4) = -1$, for $n\ge 0$.

\begin{fleqn}
\begin{align*}
P(0,4):\ & \bigl(R_{0}^{1}(2),P(0,2)\bigr),\ \bigl(R_{0}^{2}(1),P(0,3)\bigr),\ \bigl(I(0,2),2P(0,1)\bigr)\\
P(n,4):\ & \bigl(R_{0}^{n\bmod 3+1}(2),P(n,2)\bigr),\ \bigl(R_{0}^{(n+1)\bmod 3+1}(1),P(n,3)\bigr),\ \bigl(uI,(u+1)P\bigr),\ n>0\\
\end{align*}
\end{fleqn}
\subsubsection{Schofield pairs associated to preinjective exceptional modules}

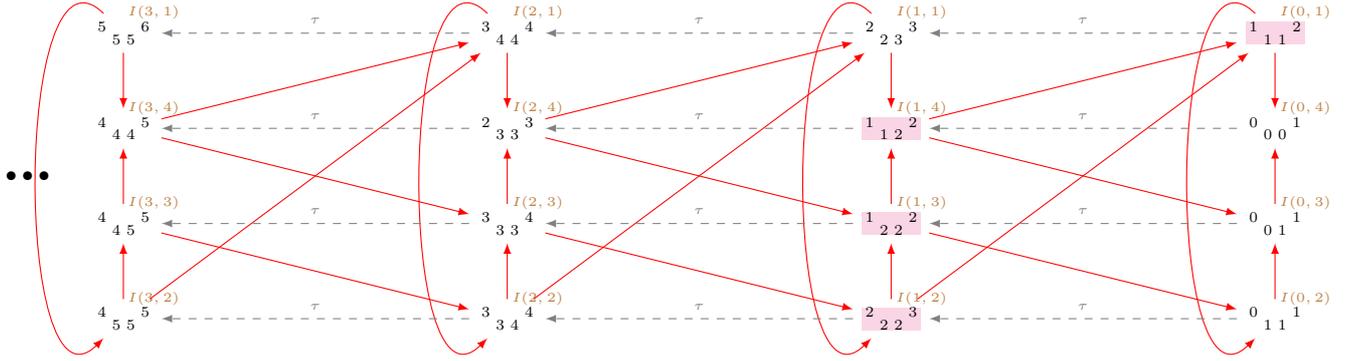
\begin{figure}[ht]


\begin{center}
\begin{scaletikzpicturetowidth}{\textwidth}

\end{scaletikzpicturetowidth}
\end{center}

\end{figure}
\medskip{}
\subsubsection*{Modules of the form $I(n,1)$}

Defect: $\partial I(n,1) = 1$, for $n\ge 0$.

\begin{fleqn}
\begin{align*}
I(0,1):\ & \bigl(I(0,2),R_{0}^{2}(1)\bigr),\ \bigl(I(0,3),R_{0}^{2}(2)\bigr),\ \bigl(2I(0,4),P(0,3)\bigr)\\
I(n,1):\ & \bigl(I(n,2),R_{0}^{(-n+1)\bmod 3+1}(1)\bigr),\ \bigl(I(n,3),R_{0}^{(-n+1)\bmod 3+1}(2)\bigr),\ \bigl((v+1)I,vP\bigr),\ n>0\\
\end{align*}
\end{fleqn}
\subsubsection*{Modules of the form $I(n,2)$}

Defect: $\partial I(n,2) = 1$, for $n\ge 0$.

\begin{fleqn}
\begin{align*}
I(0,2):\ & \bigl(I(0,3),R_{0}^{3}(1)\bigr),\ \bigl(I(0,4),R_{0}^{3}(2)\bigr)\\
I(1,2):\ & \bigl(I(1,3),R_{0}^{2}(1)\bigr),\ \bigl(I(1,4),R_{0}^{2}(2)\bigr),\ \bigl(3I(0,4),2P(0,3)\bigr)\\
I(n,2):\ & \bigl(I(n,3),R_{0}^{(-n+2)\bmod 3+1}(1)\bigr),\ \bigl(I(n,4),R_{0}^{(-n+2)\bmod 3+1}(2)\bigr),\ \bigl((v+1)I,vP\bigr),\ n>1\\
\end{align*}
\end{fleqn}
\subsubsection*{Modules of the form $I(n,3)$}

Defect: $\partial I(n,3) = 1$, for $n\ge 0$.

\begin{fleqn}
\begin{align*}
I(0,3):\ & \bigl(I(0,4),R_{0}^{1}(1)\bigr)\\
I(1,3):\ & \bigl(I(0,1),R_{0}^{3}(2)\bigr),\ \bigl(I(1,4),R_{0}^{3}(1)\bigr),\ \bigl(2I(0,2),P(0,1)\bigr)\\
I(n,3):\ & \bigl(I(n-1,1),R_{0}^{(-n+3)\bmod 3+1}(2)\bigr),\ \bigl(I(n,4),R_{0}^{(-n+3)\bmod 3+1}(1)\bigr),\ \bigl((v+1)I,vP\bigr),\ n>1\\
\end{align*}
\end{fleqn}
\subsubsection*{Modules of the form $I(n,4)$}

Defect: $\partial I(n,4) = 1$, for $n\ge 0$.

\begin{fleqn}
\begin{align*}
I(0,4):\ & - \\
I(1,4):\ & \bigl(I(0,1),R_{0}^{1}(1)\bigr),\ \bigl(I(0,2),R_{0}^{1}(2)\bigr),\ \bigl(2I(0,3),P(0,2)\bigr)\\
I(n,4):\ & \bigl(I(n-1,1),R_{0}^{(-n+1)\bmod 3+1}(1)\bigr),\ \bigl(I(n-1,2),R_{0}^{(-n+1)\bmod 3+1}(2)\bigr),\ \bigl((v+1)I,vP\bigr),\ n>1\\
\end{align*}
\end{fleqn}
\subsubsection{Schofield pairs associated to regular exceptional modules}

\subsubsection*{The non-homogeneous tube $\mathcal{T}_{0}^{\Delta(\widetilde{\mathbb{A}}_{1,3})}$}

\begin{figure}[ht]



\begin{center}
\captionof{figure}{\vspace*{-10pt}$\mathcal{T}_{0}^{\Delta(\widetilde{\mathbb{A}}_{1,3})}$}
\begin{scaletikzpicturetowidth}{0.69999999999999996\textwidth}
 $\end{lrbox}
\subsection[Schofield pairs for the quiver $\Delta(\A_{1,4})$]{Schofield pairs for the quiver $\Delta(\A_{1,4})$ -- {\usebox{\boxAlAAAA}}}
\[
\vcenter{\hbox{\xymatrix{1 &   &   &   & 5\ar[llll]\ar[dl]\\
  & 2\ar[ul] & 3\ar[l] & 4\ar[l] &  }}}\qquad C_{\Delta(\widetilde{\mathbb{A}}_{1,4})} = %
\]
\subsubsection{Schofield pairs associated to preprojective exceptional modules}

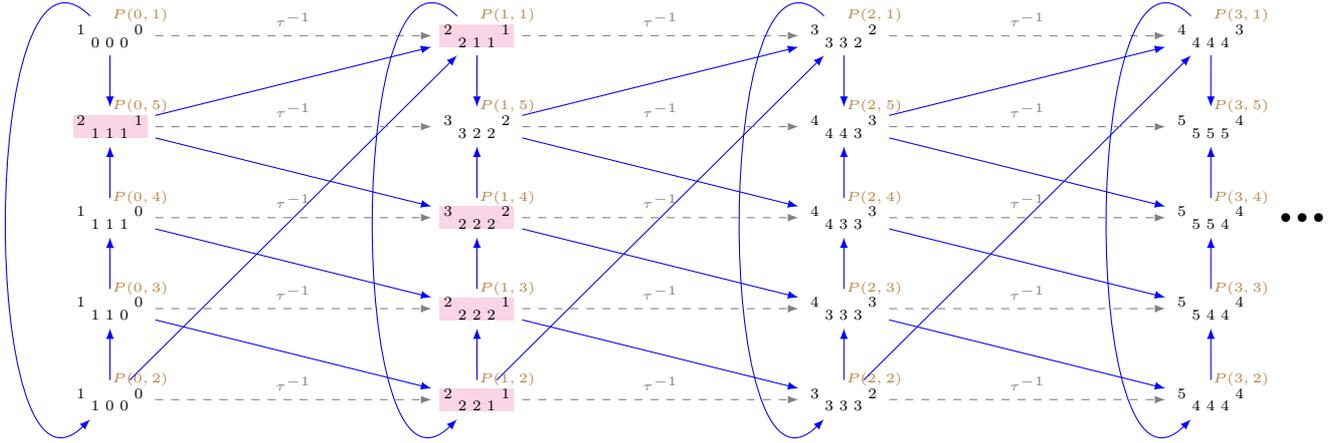
\begin{figure}[ht]


\begin{center}
\begin{scaletikzpicturetowidth}{\textwidth}

\end{scaletikzpicturetowidth}
\end{center}


\end{figure}
\medskip{}
\subsubsection*{Modules of the form $P(n,1)$}

Defect: $\partial P(n,1) = -1$, for $n\ge 0$.

\begin{fleqn}
\begin{align*}
P(0,1):\ & - \\
P(1,1):\ & \bigl(R_{0}^{1}(3),P(0,3)\bigr),\ \bigl(R_{0}^{2}(2),P(0,4)\bigr),\ \bigl(R_{0}^{3}(1),P(0,5)\bigr),\ \bigl(I(0,3),2P(0,2)\bigr)\\
P(n,1):\ & \bigl(R_{0}^{(n-1)\bmod 4+1}(3),P(n-1,3)\bigr),\ \bigl(R_{0}^{n\bmod 4+1}(2),P(n-1,4)\bigr),\ \bigl(R_{0}^{(n+1)\bmod 4+1}(1),P(n-1,5)\bigr)\\
 & \bigl(uI,(u+1)P\bigr),\ n>1\\
\end{align*}
\end{fleqn}
\subsubsection*{Modules of the form $P(n,2)$}

Defect: $\partial P(n,2) = -1$, for $n\ge 0$.

\begin{fleqn}
\begin{align*}
P(0,2):\ & \bigl(R_{0}^{3}(1),P(0,1)\bigr)\\
P(1,2):\ & \bigl(R_{0}^{2}(3),P(0,4)\bigr),\ \bigl(R_{0}^{3}(2),P(0,5)\bigr),\ \bigl(R_{0}^{4}(1),P(1,1)\bigr),\ \bigl(I(0,4),2P(0,3)\bigr)\\
P(n,2):\ & \bigl(R_{0}^{n\bmod 4+1}(3),P(n-1,4)\bigr),\ \bigl(R_{0}^{(n+1)\bmod 4+1}(2),P(n-1,5)\bigr),\ \bigl(R_{0}^{(n+2)\bmod 4+1}(1),P(n,1)\bigr)\\
 & \bigl(uI,(u+1)P\bigr),\ n>1\\
\end{align*}
\end{fleqn}
\subsubsection*{Modules of the form $P(n,3)$}

Defect: $\partial P(n,3) = -1$, for $n\ge 0$.

\begin{fleqn}
\begin{align*}
P(0,3):\ & \bigl(R_{0}^{3}(2),P(0,1)\bigr),\ \bigl(R_{0}^{4}(1),P(0,2)\bigr)\\
P(1,3):\ & \bigl(R_{0}^{3}(3),P(0,5)\bigr),\ \bigl(R_{0}^{4}(2),P(1,1)\bigr),\ \bigl(R_{0}^{1}(1),P(1,2)\bigr),\ \bigl(I(0,5),2P(0,4)\bigr)\\
P(n,3):\ & \bigl(R_{0}^{(n+1)\bmod 4+1}(3),P(n-1,5)\bigr),\ \bigl(R_{0}^{(n+2)\bmod 4+1}(2),P(n,1)\bigr),\ \bigl(R_{0}^{(n-1)\bmod 4+1}(1),P(n,2)\bigr)\\
 & \bigl(uI,(u+1)P\bigr),\ n>1\\
\end{align*}
\end{fleqn}
\subsubsection*{Modules of the form $P(n,4)$}

Defect: $\partial P(n,4) = -1$, for $n\ge 0$.

\begin{fleqn}
\begin{align*}
P(0,4):\ & \bigl(R_{0}^{3}(3),P(0,1)\bigr),\ \bigl(R_{0}^{4}(2),P(0,2)\bigr),\ \bigl(R_{0}^{1}(1),P(0,3)\bigr)\\
P(1,4):\ & \bigl(R_{0}^{4}(3),P(1,1)\bigr),\ \bigl(R_{0}^{1}(2),P(1,2)\bigr),\ \bigl(R_{0}^{2}(1),P(1,3)\bigr),\ \bigl(2I(0,2),3P(0,1)\bigr)\\
P(n,4):\ & \bigl(R_{0}^{(n+2)\bmod 4+1}(3),P(n,1)\bigr),\ \bigl(R_{0}^{(n-1)\bmod 4+1}(2),P(n,2)\bigr),\ \bigl(R_{0}^{n\bmod 4+1}(1),P(n,3)\bigr)\\
 & \bigl(uI,(u+1)P\bigr),\ n>1\\
\end{align*}
\end{fleqn}
\subsubsection*{Modules of the form $P(n,5)$}

Defect: $\partial P(n,5) = -1$, for $n\ge 0$.

\begin{fleqn}
\begin{align*}
P(0,5):\ & \bigl(R_{0}^{4}(3),P(0,2)\bigr),\ \bigl(R_{0}^{1}(2),P(0,3)\bigr),\ \bigl(R_{0}^{2}(1),P(0,4)\bigr),\ \bigl(I(0,2),2P(0,1)\bigr)\\
P(n,5):\ & \bigl(R_{0}^{(n+3)\bmod 4+1}(3),P(n,2)\bigr),\ \bigl(R_{0}^{n\bmod 4+1}(2),P(n,3)\bigr),\ \bigl(R_{0}^{(n+1)\bmod 4+1}(1),P(n,4)\bigr)\\
 & \bigl(uI,(u+1)P\bigr),\ n>0\\
\end{align*}
\end{fleqn}
\subsubsection{Schofield pairs associated to preinjective exceptional modules}

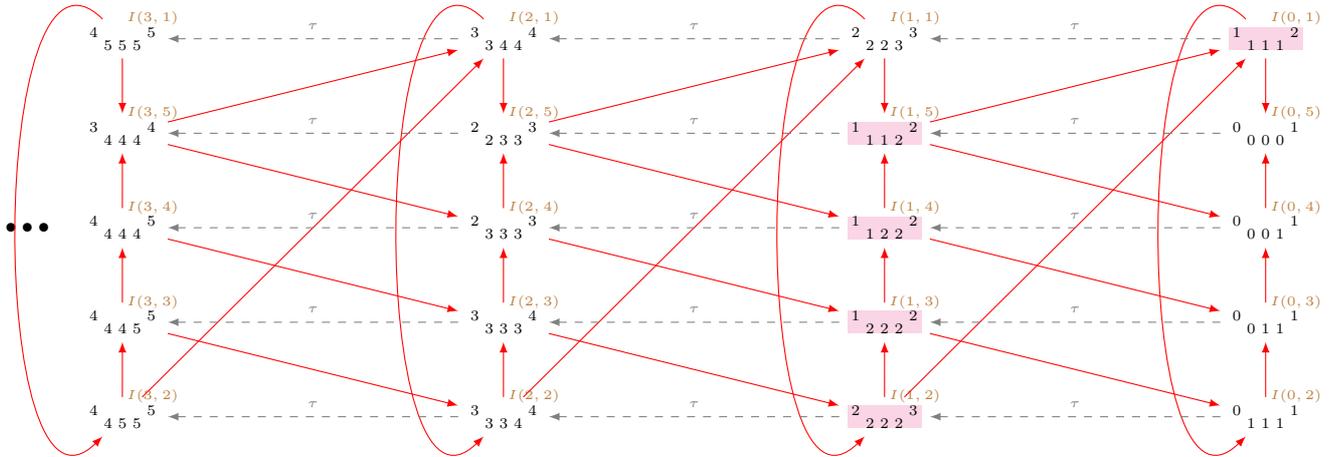
\begin{figure}[ht]


\begin{center}
\begin{scaletikzpicturetowidth}{\textwidth}

\end{scaletikzpicturetowidth}
\end{center}

\end{figure}
\medskip{}
\subsubsection*{Modules of the form $I(n,1)$}

Defect: $\partial I(n,1) = 1$, for $n\ge 0$.

\begin{fleqn}
\begin{align*}
I(0,1):\ & \bigl(I(0,2),R_{0}^{2}(1)\bigr),\ \bigl(I(0,3),R_{0}^{2}(2)\bigr),\ \bigl(I(0,4),R_{0}^{2}(3)\bigr),\ \bigl(2I(0,5),P(0,4)\bigr)\\
I(n,1):\ & \bigl(I(n,2),R_{0}^{(-n+1)\bmod 4+1}(1)\bigr),\ \bigl(I(n,3),R_{0}^{(-n+1)\bmod 4+1}(2)\bigr),\ \bigl(I(n,4),R_{0}^{(-n+1)\bmod 4+1}(3)\bigr)\\
 & \bigl((v+1)I,vP\bigr),\ n>0\\
\end{align*}
\end{fleqn}
\subsubsection*{Modules of the form $I(n,2)$}

Defect: $\partial I(n,2) = 1$, for $n\ge 0$.

\begin{fleqn}
\begin{align*}
I(0,2):\ & \bigl(I(0,3),R_{0}^{3}(1)\bigr),\ \bigl(I(0,4),R_{0}^{3}(2)\bigr),\ \bigl(I(0,5),R_{0}^{3}(3)\bigr)\\
I(1,2):\ & \bigl(I(1,3),R_{0}^{2}(1)\bigr),\ \bigl(I(1,4),R_{0}^{2}(2)\bigr),\ \bigl(I(1,5),R_{0}^{2}(3)\bigr),\ \bigl(3I(0,5),2P(0,4)\bigr)\\
I(n,2):\ & \bigl(I(n,3),R_{0}^{(-n+2)\bmod 4+1}(1)\bigr),\ \bigl(I(n,4),R_{0}^{(-n+2)\bmod 4+1}(2)\bigr),\ \bigl(I(n,5),R_{0}^{(-n+2)\bmod 4+1}(3)\bigr)\\
 & \bigl((v+1)I,vP\bigr),\ n>1\\
\end{align*}
\end{fleqn}
\subsubsection*{Modules of the form $I(n,3)$}

Defect: $\partial I(n,3) = 1$, for $n\ge 0$.

\begin{fleqn}
\begin{align*}
I(0,3):\ & \bigl(I(0,4),R_{0}^{4}(1)\bigr),\ \bigl(I(0,5),R_{0}^{4}(2)\bigr)\\
I(1,3):\ & \bigl(I(0,1),R_{0}^{3}(3)\bigr),\ \bigl(I(1,4),R_{0}^{3}(1)\bigr),\ \bigl(I(1,5),R_{0}^{3}(2)\bigr),\ \bigl(2I(0,2),P(0,1)\bigr)\\
I(n,3):\ & \bigl(I(n-1,1),R_{0}^{(-n+3)\bmod 4+1}(3)\bigr),\ \bigl(I(n,4),R_{0}^{(-n+3)\bmod 4+1}(1)\bigr),\ \bigl(I(n,5),R_{0}^{(-n+3)\bmod 4+1}(2)\bigr)\\
 & \bigl((v+1)I,vP\bigr),\ n>1\\
\end{align*}
\end{fleqn}
\subsubsection*{Modules of the form $I(n,4)$}

Defect: $\partial I(n,4) = 1$, for $n\ge 0$.

\begin{fleqn}
\begin{align*}
I(0,4):\ & \bigl(I(0,5),R_{0}^{1}(1)\bigr)\\
I(1,4):\ & \bigl(I(0,1),R_{0}^{4}(2)\bigr),\ \bigl(I(0,2),R_{0}^{4}(3)\bigr),\ \bigl(I(1,5),R_{0}^{4}(1)\bigr),\ \bigl(2I(0,3),P(0,2)\bigr)\\
I(n,4):\ & \bigl(I(n-1,1),R_{0}^{(-n+4)\bmod 4+1}(2)\bigr),\ \bigl(I(n-1,2),R_{0}^{(-n+4)\bmod 4+1}(3)\bigr),\ \bigl(I(n,5),R_{0}^{(-n+4)\bmod 4+1}(1)\bigr)\\
 & \bigl((v+1)I,vP\bigr),\ n>1\\
\end{align*}
\end{fleqn}
\subsubsection*{Modules of the form $I(n,5)$}

Defect: $\partial I(n,5) = 1$, for $n\ge 0$.

\begin{fleqn}
\begin{align*}
I(0,5):\ & - \\
I(1,5):\ & \bigl(I(0,1),R_{0}^{1}(1)\bigr),\ \bigl(I(0,2),R_{0}^{1}(2)\bigr),\ \bigl(I(0,3),R_{0}^{1}(3)\bigr),\ \bigl(2I(0,4),P(0,3)\bigr)\\
I(n,5):\ & \bigl(I(n-1,1),R_{0}^{(-n+1)\bmod 4+1}(1)\bigr),\ \bigl(I(n-1,2),R_{0}^{(-n+1)\bmod 4+1}(2)\bigr),\ \bigl(I(n-1,3),R_{0}^{(-n+1)\bmod 4+1}(3)\bigr)\\
 & \bigl((v+1)I,vP\bigr),\ n>1\\
\end{align*}
\end{fleqn}
\subsubsection{Schofield pairs associated to regular exceptional modules}

\subsubsection*{The non-homogeneous tube $\mathcal{T}_{0}^{\Delta(\widetilde{\mathbb{A}}_{1,4})}$}

\begin{figure}[ht]



\begin{center}
\captionof{figure}{\vspace*{-10pt}$\mathcal{T}_{0}^{\Delta(\widetilde{\mathbb{A}}_{1,4})}$}
\begin{scaletikzpicturetowidth}{1.\textwidth}

\end{scaletikzpicturetowidth}
\end{center}

\end{figure}
\begin{fleqn}
\begin{align*}
R_{0}^{1}(1):\ & - \\
R_{0}^{1}(2):\ & \bigl(R_{0}^{2}(1),R_{0}^{1}(1)\bigr),\ \bigl(I(0,4),P(0,1)\bigr)\\
R_{0}^{2}(1):\ & \bigl(I(0,5),P(0,1)\bigr)\\
R_{0}^{2}(2):\ & \bigl(R_{0}^{3}(1),R_{0}^{2}(1)\bigr),\ \bigl(I(0,5),P(0,2)\bigr)\\
R_{0}^{3}(1):\ & - \\
R_{0}^{3}(2):\ & \bigl(R_{0}^{4}(1),R_{0}^{3}(1)\bigr)\\
R_{0}^{4}(1):\ & - \\
R_{0}^{4}(2):\ & \bigl(R_{0}^{1}(1),R_{0}^{4}(1)\bigr)\\
R_{0}^{4}(3):\ & \bigl(R_{0}^{1}(2),R_{0}^{4}(1)\bigr),\ \bigl(R_{0}^{2}(1),R_{0}^{4}(2)\bigr),\ \bigl(I(0,3),P(0,1)\bigr)\\
R_{0}^{1}(3):\ & \bigl(R_{0}^{2}(2),R_{0}^{1}(1)\bigr),\ \bigl(R_{0}^{3}(1),R_{0}^{1}(2)\bigr),\ \bigl(I(0,4),P(0,2)\bigr)\\
R_{0}^{2}(3):\ & \bigl(R_{0}^{3}(2),R_{0}^{2}(1)\bigr),\ \bigl(R_{0}^{4}(1),R_{0}^{2}(2)\bigr),\ \bigl(I(0,5),P(0,3)\bigr)\\
R_{0}^{3}(3):\ & \bigl(R_{0}^{4}(2),R_{0}^{3}(1)\bigr),\ \bigl(R_{0}^{1}(1),R_{0}^{3}(2)\bigr)\\
\end{align*}
\end{fleqn}
\clearpage

\begin{lrbox}{\boxAAlAA}$\delta =\begin{smallmatrix}&1\\1&&1\\&1\end{smallmatrix} $\end{lrbox}
\subsection[Schofield pairs for the quiver $\Delta(\A_{2,2})$]{Schofield pairs for the quiver $\Delta(\A_{2,2})$ -- {\usebox{\boxAAlAA}}}
\[
\vcenter{\hbox{\xymatrix{  & 2\ar[dl] &   &  \\
1 &   & 4\ar[ul]\ar[dl] &  \\
  & 3\ar[ul] &   &  }}}\qquad C_{\Delta(\widetilde{\mathbb{A}}_{2,2})} = \begin{bmatrix}1 & 1 & 1 & 2\\
0 & 1 & 0 & 1\\
0 & 0 & 1 & 1\\
0 & 0 & 0 & 1\end{bmatrix}\quad\Phi_{\Delta(\widetilde{\mathbb{A}}_{2,2})} = \begin{bmatrix}-1 & 1 & 1 & 0\\
-1 & 0 & 1 & 1\\
-1 & 1 & 0 & 1\\
-2 & 1 & 1 & 1\end{bmatrix}
\]
\subsubsection{Schofield pairs associated to preprojective exceptional modules}

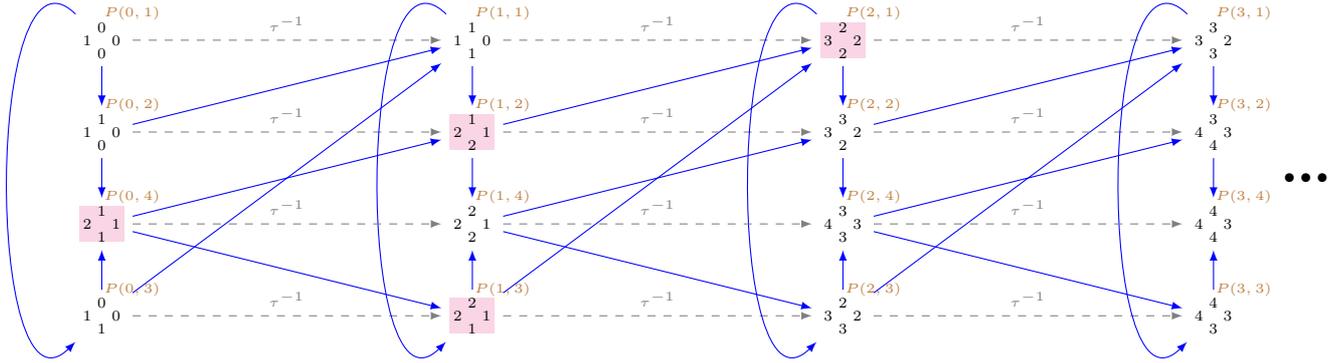
\begin{figure}[ht]


\begin{center}
\begin{scaletikzpicturetowidth}{\textwidth}

\end{scaletikzpicturetowidth}
\end{center}


\end{figure}
\medskip{}
\subsubsection*{Modules of the form $P(n,1)$}

Defect: $\partial P(n,1) = -1$, for $n\ge 0$.

\begin{fleqn}
\begin{align*}
P(0,1):\ & - \\
P(1,1):\ & \bigl(R_{0}^{1}(1),P(0,2)\bigr),\ \bigl(R_{\infty}^{1}(1),P(0,3)\bigr)\\
P(2,1):\ & \bigl(R_{0}^{2}(1),P(1,2)\bigr),\ \bigl(R_{\infty}^{2}(1),P(1,3)\bigr),\ \bigl(2I(1,4),3P(0,1)\bigr)\\
P(n,1):\ & \bigl(R_{0}^{(n-1)\bmod 2+1}(1),P(n-1,2)\bigr),\ \bigl(R_{\infty}^{(n-1)\bmod 2+1}(1),P(n-1,3)\bigr),\ \bigl(uI,(u+1)P\bigr),\ n>2\\
\end{align*}
\end{fleqn}
\subsubsection*{Modules of the form $P(n,2)$}

Defect: $\partial P(n,2) = -1$, for $n\ge 0$.

\begin{fleqn}
\begin{align*}
P(0,2):\ & \bigl(R_{\infty}^{1}(1),P(0,1)\bigr)\\
P(1,2):\ & \bigl(R_{0}^{1}(1),P(0,4)\bigr),\ \bigl(R_{\infty}^{2}(1),P(1,1)\bigr),\ \bigl(I(0,2),2P(0,3)\bigr)\\
P(n,2):\ & \bigl(R_{0}^{(n-1)\bmod 2+1}(1),P(n-1,4)\bigr),\ \bigl(R_{\infty}^{n\bmod 2+1}(1),P(n,1)\bigr),\ \bigl(uI,(u+1)P\bigr),\ n>1\\
\end{align*}
\end{fleqn}
\subsubsection*{Modules of the form $P(n,3)$}

Defect: $\partial P(n,3) = -1$, for $n\ge 0$.

\begin{fleqn}
\begin{align*}
P(0,3):\ & \bigl(R_{0}^{1}(1),P(0,1)\bigr)\\
P(1,3):\ & \bigl(R_{\infty}^{1}(1),P(0,4)\bigr),\ \bigl(R_{0}^{2}(1),P(1,1)\bigr),\ \bigl(I(0,3),2P(0,2)\bigr)\\
P(n,3):\ & \bigl(R_{\infty}^{(n-1)\bmod 2+1}(1),P(n-1,4)\bigr),\ \bigl(R_{0}^{n\bmod 2+1}(1),P(n,1)\bigr),\ \bigl(uI,(u+1)P\bigr),\ n>1\\
\end{align*}
\end{fleqn}
\subsubsection*{Modules of the form $P(n,4)$}

Defect: $\partial P(n,4) = -1$, for $n\ge 0$.

\begin{fleqn}
\begin{align*}
P(0,4):\ & \bigl(R_{\infty}^{2}(1),P(0,2)\bigr),\ \bigl(R_{0}^{2}(1),P(0,3)\bigr),\ \bigl(I(1,4),2P(0,1)\bigr)\\
P(n,4):\ & \bigl(R_{\infty}^{(n+1)\bmod 2+1}(1),P(n,2)\bigr),\ \bigl(R_{0}^{(n+1)\bmod 2+1}(1),P(n,3)\bigr),\ \bigl(uI,(u+1)P\bigr),\ n>0\\
\end{align*}
\end{fleqn}
\subsubsection{Schofield pairs associated to preinjective exceptional modules}

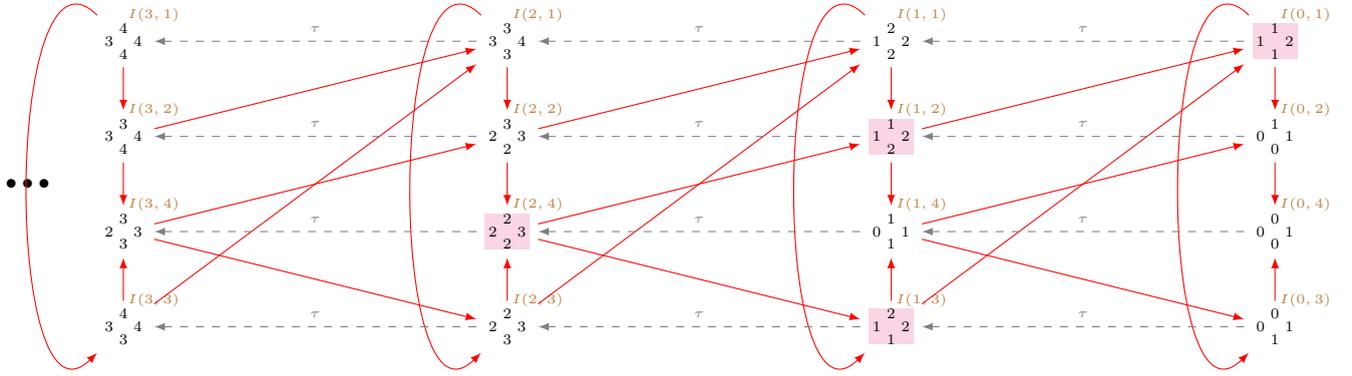
\begin{figure}[ht]


\begin{center}
\begin{scaletikzpicturetowidth}{\textwidth}

\end{scaletikzpicturetowidth}
\end{center}

\end{figure}
\medskip{}
\subsubsection*{Modules of the form $I(n,1)$}

Defect: $\partial I(n,1) = 1$, for $n\ge 0$.

\begin{fleqn}
\begin{align*}
I(0,1):\ & \bigl(I(0,2),R_{\infty}^{2}(1)\bigr),\ \bigl(I(0,3),R_{0}^{2}(1)\bigr),\ \bigl(2I(0,4),P(1,1)\bigr)\\
I(n,1):\ & \bigl(I(n,2),R_{\infty}^{(-n+1)\bmod 2+1}(1)\bigr),\ \bigl(I(n,3),R_{0}^{(-n+1)\bmod 2+1}(1)\bigr),\ \bigl((v+1)I,vP\bigr),\ n>0\\
\end{align*}
\end{fleqn}
\subsubsection*{Modules of the form $I(n,2)$}

Defect: $\partial I(n,2) = 1$, for $n\ge 0$.

\begin{fleqn}
\begin{align*}
I(0,2):\ & \bigl(I(0,4),R_{\infty}^{1}(1)\bigr)\\
I(1,2):\ & \bigl(I(0,1),R_{0}^{1}(1)\bigr),\ \bigl(I(1,4),R_{\infty}^{2}(1)\bigr),\ \bigl(2I(0,3),P(0,2)\bigr)\\
I(n,2):\ & \bigl(I(n-1,1),R_{0}^{(-n+1)\bmod 2+1}(1)\bigr),\ \bigl(I(n,4),R_{\infty}^{(-n+2)\bmod 2+1}(1)\bigr),\ \bigl((v+1)I,vP\bigr),\ n>1\\
\end{align*}
\end{fleqn}
\subsubsection*{Modules of the form $I(n,3)$}

Defect: $\partial I(n,3) = 1$, for $n\ge 0$.

\begin{fleqn}
\begin{align*}
I(0,3):\ & \bigl(I(0,4),R_{0}^{1}(1)\bigr)\\
I(1,3):\ & \bigl(I(0,1),R_{\infty}^{1}(1)\bigr),\ \bigl(I(1,4),R_{0}^{2}(1)\bigr),\ \bigl(2I(0,2),P(0,3)\bigr)\\
I(n,3):\ & \bigl(I(n-1,1),R_{\infty}^{(-n+1)\bmod 2+1}(1)\bigr),\ \bigl(I(n,4),R_{0}^{(-n+2)\bmod 2+1}(1)\bigr),\ \bigl((v+1)I,vP\bigr),\ n>1\\
\end{align*}
\end{fleqn}
\subsubsection*{Modules of the form $I(n,4)$}

Defect: $\partial I(n,4) = 1$, for $n\ge 0$.

\begin{fleqn}
\begin{align*}
I(0,4):\ & - \\
I(1,4):\ & \bigl(I(0,2),R_{0}^{1}(1)\bigr),\ \bigl(I(0,3),R_{\infty}^{1}(1)\bigr)\\
I(2,4):\ & \bigl(I(1,2),R_{0}^{2}(1)\bigr),\ \bigl(I(1,3),R_{\infty}^{2}(1)\bigr),\ \bigl(3I(0,4),2P(1,1)\bigr)\\
I(n,4):\ & \bigl(I(n-1,2),R_{0}^{(-n+3)\bmod 2+1}(1)\bigr),\ \bigl(I(n-1,3),R_{\infty}^{(-n+3)\bmod 2+1}(1)\bigr),\ \bigl((v+1)I,vP\bigr),\ n>2\\
\end{align*}
\end{fleqn}
\subsubsection{Schofield pairs associated to regular exceptional modules}

\subsubsection*{The non-homogeneous tube $\mathcal{T}_{0}^{\Delta(\widetilde{\mathbb{A}}_{2,2})}$}

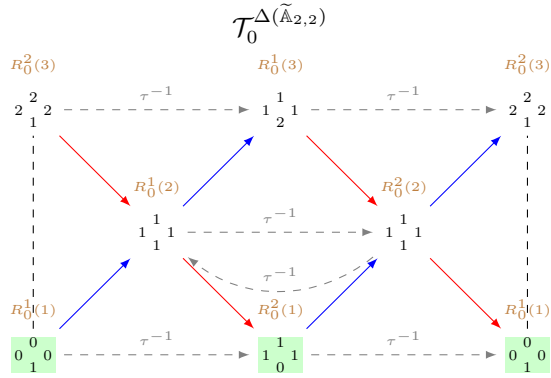
\begin{figure}[ht]



\begin{center}
\captionof{figure}{\vspace*{-10pt}$\mathcal{T}_{0}^{\Delta(\widetilde{\mathbb{A}}_{2,2})}$}
\begin{scaletikzpicturetowidth}{0.40000000000000002\textwidth}

\end{scaletikzpicturetowidth}
\end{center}

\end{figure}
\begin{fleqn}
\begin{align*}
R_{0}^{1}(1):\ & - \\
R_{0}^{2}(1):\ & \bigl(I(0,2),P(0,1)\bigr),\ \bigl(I(0,4),P(0,2)\bigr)\\
\end{align*}
\end{fleqn}
\subsubsection*{The non-homogeneous tube $\mathcal{T}_{\infty}^{\Delta(\widetilde{\mathbb{A}}_{2,2})}$}

\begin{figure}[ht]



\begin{center}
\captionof{figure}{\vspace*{-10pt}$\mathcal{T}_{\infty}^{\Delta(\widetilde{\mathbb{A}}_{2,2})}$}
\begin{scaletikzpicturetowidth}{0.40000000000000002\textwidth}
 $\end{lrbox}
\subsection[Schofield pairs for the quiver $\Delta(\A_{2,3})$]{Schofield pairs for the quiver $\Delta(\A_{2,3})$ -- {\usebox{\boxAAlAAA}}}
\[
\vcenter{\hbox{\xymatrix{  & 2\ar[dl] &   &   &  \\
1 &   &   & 5\ar[ull]\ar[dl] &  \\
  & 3\ar[ul] & 4\ar[l] &   &  }}}\qquad C_{\Delta(\widetilde{\mathbb{A}}_{2,3})} = %
\]
\subsubsection{Schofield pairs associated to preprojective exceptional modules}

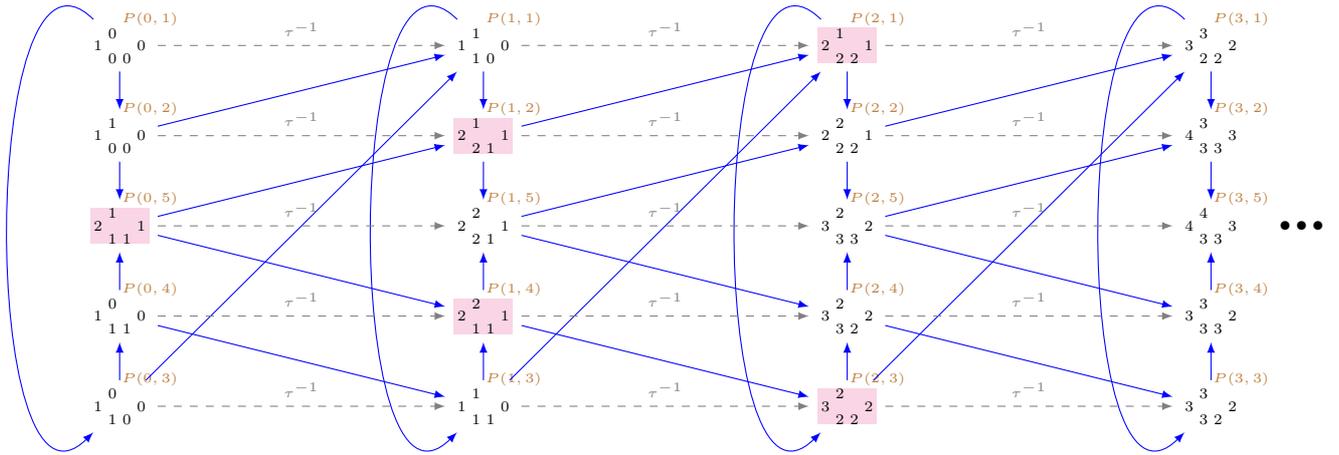
\begin{figure}[ht]


\begin{center}
\begin{scaletikzpicturetowidth}{\textwidth}

\end{scaletikzpicturetowidth}
\end{center}


\end{figure}
\medskip{}
\subsubsection*{Modules of the form $P(n,1)$}

Defect: $\partial P(n,1) = -1$, for $n\ge 0$.

\begin{fleqn}
\begin{align*}
P(0,1):\ & - \\
P(1,1):\ & \bigl(R_{0}^{3}(1),P(0,2)\bigr),\ \bigl(R_{\infty}^{1}(1),P(0,3)\bigr)\\
P(2,1):\ & \bigl(R_{0}^{3}(2),P(0,5)\bigr),\ \bigl(R_{0}^{1}(1),P(1,2)\bigr),\ \bigl(R_{\infty}^{2}(1),P(1,3)\bigr),\ \bigl(I(0,2),2P(0,4)\bigr)\\
P(n,1):\ & \bigl(R_{0}^{n\bmod 3+1}(2),P(n-2,5)\bigr),\ \bigl(R_{0}^{(n-2)\bmod 3+1}(1),P(n-1,2)\bigr),\ \bigl(R_{\infty}^{(n-1)\bmod 2+1}(1),P(n-1,3)\bigr)\\
 & \bigl(uI,(u+1)P\bigr),\ n>2\\
\end{align*}
\end{fleqn}
\subsubsection*{Modules of the form $P(n,2)$}

Defect: $\partial P(n,2) = -1$, for $n\ge 0$.

\begin{fleqn}
\begin{align*}
P(0,2):\ & \bigl(R_{\infty}^{1}(1),P(0,1)\bigr)\\
P(1,2):\ & \bigl(R_{0}^{2}(2),P(0,4)\bigr),\ \bigl(R_{0}^{3}(1),P(0,5)\bigr),\ \bigl(R_{\infty}^{2}(1),P(1,1)\bigr),\ \bigl(I(1,5),2P(0,3)\bigr)\\
P(n,2):\ & \bigl(R_{0}^{n\bmod 3+1}(2),P(n-1,4)\bigr),\ \bigl(R_{0}^{(n+1)\bmod 3+1}(1),P(n-1,5)\bigr),\ \bigl(R_{\infty}^{n\bmod 2+1}(1),P(n,1)\bigr)\\
 & \bigl(uI,(u+1)P\bigr),\ n>1\\
\end{align*}
\end{fleqn}
\subsubsection*{Modules of the form $P(n,3)$}

Defect: $\partial P(n,3) = -1$, for $n\ge 0$.

\begin{fleqn}
\begin{align*}
P(0,3):\ & \bigl(R_{0}^{3}(1),P(0,1)\bigr)\\
P(1,3):\ & \bigl(R_{0}^{3}(2),P(0,2)\bigr),\ \bigl(R_{\infty}^{1}(1),P(0,4)\bigr),\ \bigl(R_{0}^{1}(1),P(1,1)\bigr)\\
P(2,3):\ & \bigl(R_{0}^{1}(2),P(1,2)\bigr),\ \bigl(R_{\infty}^{2}(1),P(1,4)\bigr),\ \bigl(R_{0}^{2}(1),P(2,1)\bigr),\ \bigl(2I(1,4),3P(0,1)\bigr)\\
P(n,3):\ & \bigl(R_{0}^{(n-2)\bmod 3+1}(2),P(n-1,2)\bigr),\ \bigl(R_{\infty}^{(n-1)\bmod 2+1}(1),P(n-1,4)\bigr),\ \bigl(R_{0}^{(n-1)\bmod 3+1}(1),P(n,1)\bigr)\\
 & \bigl(uI,(u+1)P\bigr),\ n>2\\
\end{align*}
\end{fleqn}
\subsubsection*{Modules of the form $P(n,4)$}

Defect: $\partial P(n,4) = -1$, for $n\ge 0$.

\begin{fleqn}
\begin{align*}
P(0,4):\ & \bigl(R_{0}^{3}(2),P(0,1)\bigr),\ \bigl(R_{0}^{1}(1),P(0,3)\bigr)\\
P(1,4):\ & \bigl(R_{\infty}^{1}(1),P(0,5)\bigr),\ \bigl(R_{0}^{1}(2),P(1,1)\bigr),\ \bigl(R_{0}^{2}(1),P(1,3)\bigr),\ \bigl(I(0,3),2P(0,2)\bigr)\\
P(n,4):\ & \bigl(R_{\infty}^{(n-1)\bmod 2+1}(1),P(n-1,5)\bigr),\ \bigl(R_{0}^{(n-1)\bmod 3+1}(2),P(n,1)\bigr),\ \bigl(R_{0}^{n\bmod 3+1}(1),P(n,3)\bigr)\\
 & \bigl(uI,(u+1)P\bigr),\ n>1\\
\end{align*}
\end{fleqn}
\subsubsection*{Modules of the form $P(n,5)$}

Defect: $\partial P(n,5) = -1$, for $n\ge 0$.

\begin{fleqn}
\begin{align*}
P(0,5):\ & \bigl(R_{\infty}^{2}(1),P(0,2)\bigr),\ \bigl(R_{0}^{1}(2),P(0,3)\bigr),\ \bigl(R_{0}^{2}(1),P(0,4)\bigr),\ \bigl(I(1,4),2P(0,1)\bigr)\\
P(n,5):\ & \bigl(R_{\infty}^{(n+1)\bmod 2+1}(1),P(n,2)\bigr),\ \bigl(R_{0}^{n\bmod 3+1}(2),P(n,3)\bigr),\ \bigl(R_{0}^{(n+1)\bmod 3+1}(1),P(n,4)\bigr)\\
 & \bigl(uI,(u+1)P\bigr),\ n>0\\
\end{align*}
\end{fleqn}
\subsubsection{Schofield pairs associated to preinjective exceptional modules}

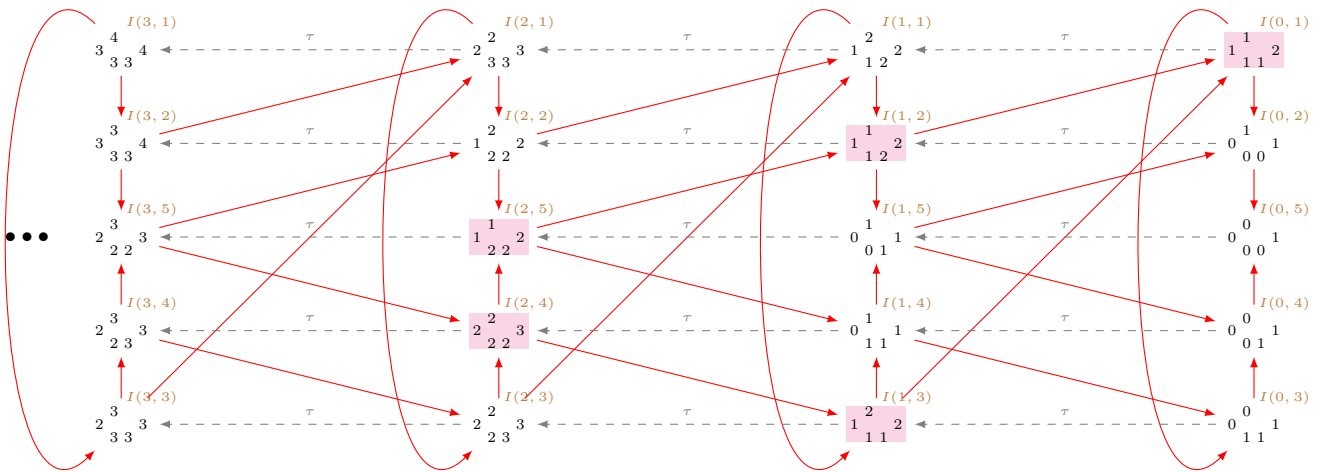
\begin{figure}[ht]


\begin{center}
\begin{scaletikzpicturetowidth}{\textwidth}

\end{scaletikzpicturetowidth}
\end{center}

\end{figure}
\medskip{}
\subsubsection*{Modules of the form $I(n,1)$}

Defect: $\partial I(n,1) = 1$, for $n\ge 0$.

\begin{fleqn}
\begin{align*}
I(0,1):\ & \bigl(I(0,2),R_{\infty}^{2}(1)\bigr),\ \bigl(I(0,3),R_{0}^{2}(1)\bigr),\ \bigl(I(0,4),R_{0}^{2}(2)\bigr),\ \bigl(2I(0,5),P(1,3)\bigr)\\
I(n,1):\ & \bigl(I(n,2),R_{\infty}^{(-n+1)\bmod 2+1}(1)\bigr),\ \bigl(I(n,3),R_{0}^{(-n+1)\bmod 3+1}(1)\bigr),\ \bigl(I(n,4),R_{0}^{(-n+1)\bmod 3+1}(2)\bigr)\\
 & \bigl((v+1)I,vP\bigr),\ n>0\\
\end{align*}
\end{fleqn}
\subsubsection*{Modules of the form $I(n,2)$}

Defect: $\partial I(n,2) = 1$, for $n\ge 0$.

\begin{fleqn}
\begin{align*}
I(0,2):\ & \bigl(I(0,5),R_{\infty}^{1}(1)\bigr)\\
I(1,2):\ & \bigl(I(0,1),R_{0}^{1}(1)\bigr),\ \bigl(I(0,3),R_{0}^{1}(2)\bigr),\ \bigl(I(1,5),R_{\infty}^{2}(1)\bigr),\ \bigl(2I(0,4),P(1,1)\bigr)\\
I(n,2):\ & \bigl(I(n-1,1),R_{0}^{(-n+1)\bmod 3+1}(1)\bigr),\ \bigl(I(n-1,3),R_{0}^{(-n+1)\bmod 3+1}(2)\bigr),\ \bigl(I(n,5),R_{\infty}^{(-n+2)\bmod 2+1}(1)\bigr)\\
 & \bigl((v+1)I,vP\bigr),\ n>1\\
\end{align*}
\end{fleqn}
\subsubsection*{Modules of the form $I(n,3)$}

Defect: $\partial I(n,3) = 1$, for $n\ge 0$.

\begin{fleqn}
\begin{align*}
I(0,3):\ & \bigl(I(0,4),R_{0}^{3}(1)\bigr),\ \bigl(I(0,5),R_{0}^{3}(2)\bigr)\\
I(1,3):\ & \bigl(I(0,1),R_{\infty}^{1}(1)\bigr),\ \bigl(I(1,4),R_{0}^{2}(1)\bigr),\ \bigl(I(1,5),R_{0}^{2}(2)\bigr),\ \bigl(2I(0,2),P(0,4)\bigr)\\
I(n,3):\ & \bigl(I(n-1,1),R_{\infty}^{(-n+1)\bmod 2+1}(1)\bigr),\ \bigl(I(n,4),R_{0}^{(-n+2)\bmod 3+1}(1)\bigr),\ \bigl(I(n,5),R_{0}^{(-n+2)\bmod 3+1}(2)\bigr)\\
 & \bigl((v+1)I,vP\bigr),\ n>1\\
\end{align*}
\end{fleqn}
\subsubsection*{Modules of the form $I(n,4)$}

Defect: $\partial I(n,4) = 1$, for $n\ge 0$.

\begin{fleqn}
\begin{align*}
I(0,4):\ & \bigl(I(0,5),R_{0}^{1}(1)\bigr)\\
I(1,4):\ & \bigl(I(0,2),R_{0}^{3}(2)\bigr),\ \bigl(I(0,3),R_{\infty}^{1}(1)\bigr),\ \bigl(I(1,5),R_{0}^{3}(1)\bigr)\\
I(2,4):\ & \bigl(I(1,2),R_{0}^{2}(2)\bigr),\ \bigl(I(1,3),R_{\infty}^{2}(1)\bigr),\ \bigl(I(2,5),R_{0}^{2}(1)\bigr),\ \bigl(3I(0,5),2P(1,3)\bigr)\\
I(n,4):\ & \bigl(I(n-1,2),R_{0}^{(-n+3)\bmod 3+1}(2)\bigr),\ \bigl(I(n-1,3),R_{\infty}^{(-n+3)\bmod 2+1}(1)\bigr),\ \bigl(I(n,5),R_{0}^{(-n+3)\bmod 3+1}(1)\bigr)\\
 & \bigl((v+1)I,vP\bigr),\ n>2\\
\end{align*}
\end{fleqn}
\subsubsection*{Modules of the form $I(n,5)$}

Defect: $\partial I(n,5) = 1$, for $n\ge 0$.

\begin{fleqn}
\begin{align*}
I(0,5):\ & - \\
I(1,5):\ & \bigl(I(0,2),R_{0}^{1}(1)\bigr),\ \bigl(I(0,4),R_{\infty}^{1}(1)\bigr)\\
I(2,5):\ & \bigl(I(0,1),R_{0}^{3}(2)\bigr),\ \bigl(I(1,2),R_{0}^{3}(1)\bigr),\ \bigl(I(1,4),R_{\infty}^{2}(1)\bigr),\ \bigl(2I(0,3),P(0,2)\bigr)\\
I(n,5):\ & \bigl(I(n-2,1),R_{0}^{(-n+4)\bmod 3+1}(2)\bigr),\ \bigl(I(n-1,2),R_{0}^{(-n+4)\bmod 3+1}(1)\bigr),\ \bigl(I(n-1,4),R_{\infty}^{(-n+3)\bmod 2+1}(1)\bigr)\\
 & \bigl((v+1)I,vP\bigr),\ n>2\\
\end{align*}
\end{fleqn}
\subsubsection{Schofield pairs associated to regular exceptional modules}

\subsubsection*{The non-homogeneous tube $\mathcal{T}_{0}^{\Delta(\widetilde{\mathbb{A}}_{2,3})}$}

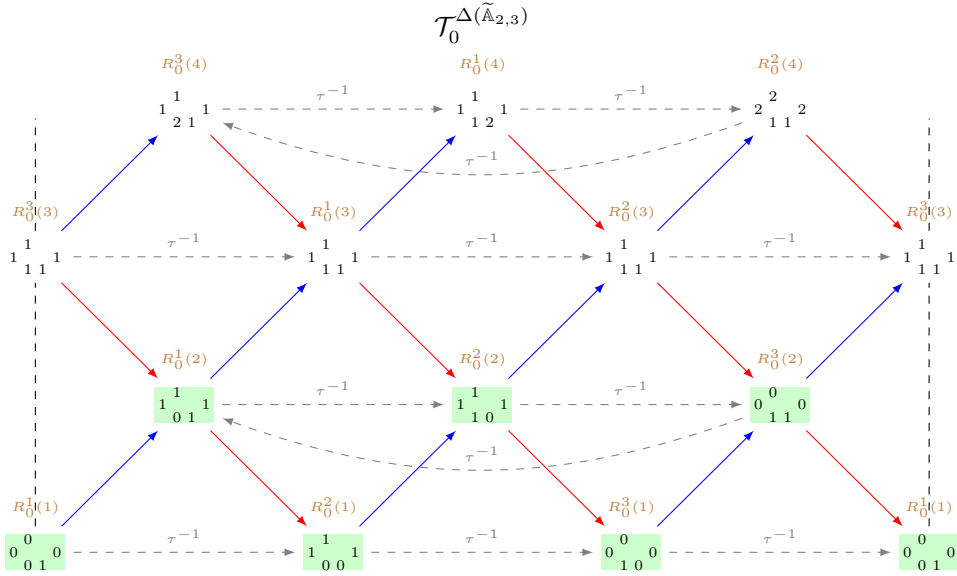
\begin{figure}[ht]



\begin{center}
\captionof{figure}{\vspace*{-10pt}$\mathcal{T}_{0}^{\Delta(\widetilde{\mathbb{A}}_{2,3})}$}
\begin{scaletikzpicturetowidth}{0.69999999999999996\textwidth}

\end{scaletikzpicturetowidth}
\end{center}

\end{figure}
\begin{fleqn}
\begin{align*}
R_{0}^{1}(1):\ & - \\
R_{0}^{1}(2):\ & \bigl(R_{0}^{2}(1),R_{0}^{1}(1)\bigr),\ \bigl(I(1,5),P(0,1)\bigr),\ \bigl(I(0,4),P(0,2)\bigr)\\
R_{0}^{2}(1):\ & \bigl(I(0,2),P(0,1)\bigr),\ \bigl(I(0,5),P(0,2)\bigr)\\
R_{0}^{2}(2):\ & \bigl(R_{0}^{3}(1),R_{0}^{2}(1)\bigr),\ \bigl(I(0,2),P(0,3)\bigr),\ \bigl(I(0,5),P(1,1)\bigr)\\
R_{0}^{3}(1):\ & - \\
R_{0}^{3}(2):\ & \bigl(R_{0}^{1}(1),R_{0}^{3}(1)\bigr)\\
\end{align*}
\end{fleqn}
\subsubsection*{The non-homogeneous tube $\mathcal{T}_{\infty}^{\Delta(\widetilde{\mathbb{A}}_{2,3})}$}

\begin{figure}[ht]



\begin{center}
\captionof{figure}{\vspace*{-10pt}$\mathcal{T}_{\infty}^{\Delta(\widetilde{\mathbb{A}}_{2,3})}$}
\begin{scaletikzpicturetowidth}{0.5\textwidth}
 $\end{lrbox}
\subsection[Schofield pairs for the quiver $\Delta(\A_{2,4})$]{Schofield pairs for the quiver $\Delta(\A_{2,4})$ -- {\usebox{\boxAAlAAAA}}}
\[
\vcenter{\hbox{\xymatrix{  & 2\ar[dl] &   &   &   &  \\
1 &   &   &   & 6\ar[ulll]\ar[dl] &  \\
  & 3\ar[ul] & 4\ar[l] & 5\ar[l] &   &  }}}\qquad C_{\Delta(\widetilde{\mathbb{A}}_{2,4})} = %
\]
\subsubsection{Schofield pairs associated to preprojective exceptional modules}

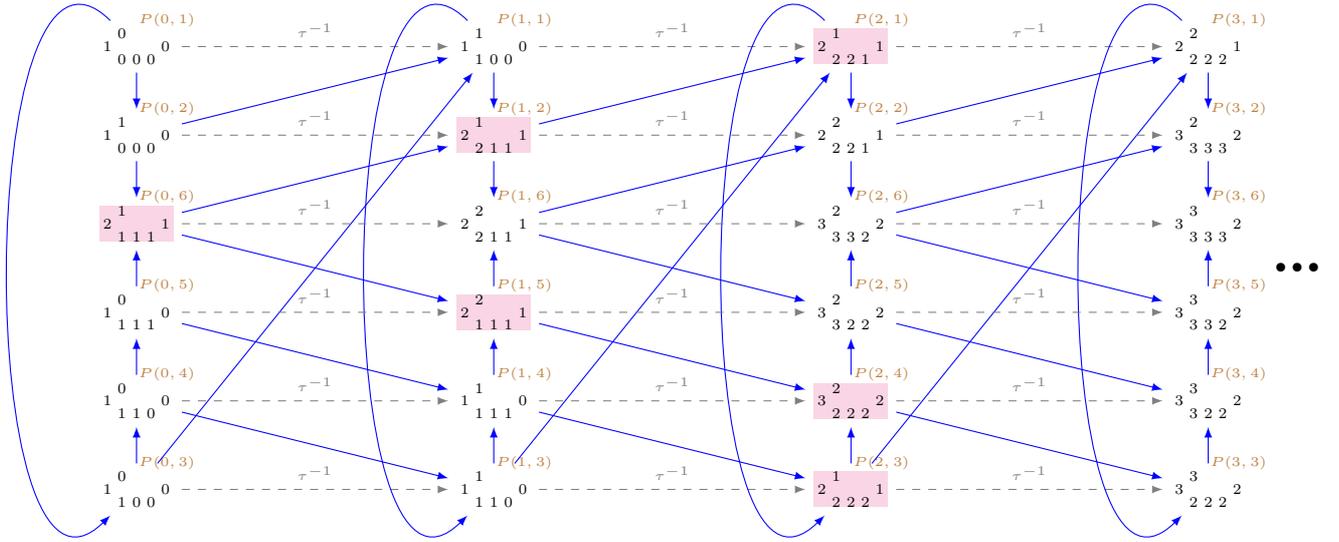
\begin{figure}[ht]


\begin{center}
\begin{scaletikzpicturetowidth}{\textwidth}

\end{scaletikzpicturetowidth}
\end{center}


\end{figure}
\medskip{}
\subsubsection*{Modules of the form $P(n,1)$}

Defect: $\partial P(n,1) = -1$, for $n\ge 0$.

\begin{fleqn}
\begin{align*}
P(0,1):\ & - \\
P(1,1):\ & \bigl(R_{0}^{3}(1),P(0,2)\bigr),\ \bigl(R_{\infty}^{1}(1),P(0,3)\bigr)\\
P(2,1):\ & \bigl(R_{0}^{2}(3),P(0,5)\bigr),\ \bigl(R_{0}^{3}(2),P(0,6)\bigr),\ \bigl(R_{0}^{4}(1),P(1,2)\bigr),\ \bigl(R_{\infty}^{2}(1),P(1,3)\bigr),\ \bigl(I(1,6),2P(0,4)\bigr)\\
P(n,1):\ & \bigl(R_{0}^{(n-1)\bmod 4+1}(3),P(n-2,5)\bigr),\ \bigl(R_{0}^{n\bmod 4+1}(2),P(n-2,6)\bigr),\ \bigl(R_{0}^{(n+1)\bmod 4+1}(1),P(n-1,2)\bigr)\\
 & \bigl(R_{\infty}^{(n-1)\bmod 2+1}(1),P(n-1,3)\bigr),\ \bigl(uI,(u+1)P\bigr),\ n>2\\
\end{align*}
\end{fleqn}
\subsubsection*{Modules of the form $P(n,2)$}

Defect: $\partial P(n,2) = -1$, for $n\ge 0$.

\begin{fleqn}
\begin{align*}
P(0,2):\ & \bigl(R_{\infty}^{1}(1),P(0,1)\bigr)\\
P(1,2):\ & \bigl(R_{0}^{1}(3),P(0,4)\bigr),\ \bigl(R_{0}^{2}(2),P(0,5)\bigr),\ \bigl(R_{0}^{3}(1),P(0,6)\bigr),\ \bigl(R_{\infty}^{2}(1),P(1,1)\bigr),\ \bigl(I(1,5),2P(0,3)\bigr)\\
P(n,2):\ & \bigl(R_{0}^{(n-1)\bmod 4+1}(3),P(n-1,4)\bigr),\ \bigl(R_{0}^{n\bmod 4+1}(2),P(n-1,5)\bigr),\ \bigl(R_{0}^{(n+1)\bmod 4+1}(1),P(n-1,6)\bigr)\\
 & \bigl(R_{\infty}^{n\bmod 2+1}(1),P(n,1)\bigr),\ \bigl(uI,(u+1)P\bigr),\ n>1\\
\end{align*}
\end{fleqn}
\subsubsection*{Modules of the form $P(n,3)$}

Defect: $\partial P(n,3) = -1$, for $n\ge 0$.

\begin{fleqn}
\begin{align*}
P(0,3):\ & \bigl(R_{0}^{3}(1),P(0,1)\bigr)\\
P(1,3):\ & \bigl(R_{0}^{3}(2),P(0,2)\bigr),\ \bigl(R_{\infty}^{1}(1),P(0,4)\bigr),\ \bigl(R_{0}^{4}(1),P(1,1)\bigr)\\
P(2,3):\ & \bigl(R_{0}^{3}(3),P(0,6)\bigr),\ \bigl(R_{0}^{4}(2),P(1,2)\bigr),\ \bigl(R_{\infty}^{2}(1),P(1,4)\bigr),\ \bigl(R_{0}^{1}(1),P(2,1)\bigr),\ \bigl(I(0,2),2P(0,5)\bigr)\\
P(n,3):\ & \bigl(R_{0}^{n\bmod 4+1}(3),P(n-2,6)\bigr),\ \bigl(R_{0}^{(n+1)\bmod 4+1}(2),P(n-1,2)\bigr),\ \bigl(R_{\infty}^{(n-1)\bmod 2+1}(1),P(n-1,4)\bigr)\\
 & \bigl(R_{0}^{(n-2)\bmod 4+1}(1),P(n,1)\bigr),\ \bigl(uI,(u+1)P\bigr),\ n>2\\
\end{align*}
\end{fleqn}
\subsubsection*{Modules of the form $P(n,4)$}

Defect: $\partial P(n,4) = -1$, for $n\ge 0$.

\begin{fleqn}
\begin{align*}
P(0,4):\ & \bigl(R_{0}^{3}(2),P(0,1)\bigr),\ \bigl(R_{0}^{4}(1),P(0,3)\bigr)\\
P(1,4):\ & \bigl(R_{0}^{3}(3),P(0,2)\bigr),\ \bigl(R_{\infty}^{1}(1),P(0,5)\bigr),\ \bigl(R_{0}^{4}(2),P(1,1)\bigr),\ \bigl(R_{0}^{1}(1),P(1,3)\bigr)\\
P(2,4):\ & \bigl(R_{0}^{4}(3),P(1,2)\bigr),\ \bigl(R_{\infty}^{2}(1),P(1,5)\bigr),\ \bigl(R_{0}^{1}(2),P(2,1)\bigr),\ \bigl(R_{0}^{2}(1),P(2,3)\bigr),\ \bigl(2I(1,4),3P(0,1)\bigr)\\
P(n,4):\ & \bigl(R_{0}^{(n+1)\bmod 4+1}(3),P(n-1,2)\bigr),\ \bigl(R_{\infty}^{(n-1)\bmod 2+1}(1),P(n-1,5)\bigr),\ \bigl(R_{0}^{(n-2)\bmod 4+1}(2),P(n,1)\bigr)\\
 & \bigl(R_{0}^{(n-1)\bmod 4+1}(1),P(n,3)\bigr),\ \bigl(uI,(u+1)P\bigr),\ n>2\\
\end{align*}
\end{fleqn}
\subsubsection*{Modules of the form $P(n,5)$}

Defect: $\partial P(n,5) = -1$, for $n\ge 0$.

\begin{fleqn}
\begin{align*}
P(0,5):\ & \bigl(R_{0}^{3}(3),P(0,1)\bigr),\ \bigl(R_{0}^{4}(2),P(0,3)\bigr),\ \bigl(R_{0}^{1}(1),P(0,4)\bigr)\\
P(1,5):\ & \bigl(R_{\infty}^{1}(1),P(0,6)\bigr),\ \bigl(R_{0}^{4}(3),P(1,1)\bigr),\ \bigl(R_{0}^{1}(2),P(1,3)\bigr),\ \bigl(R_{0}^{2}(1),P(1,4)\bigr),\ \bigl(I(0,3),2P(0,2)\bigr)\\
P(n,5):\ & \bigl(R_{\infty}^{(n-1)\bmod 2+1}(1),P(n-1,6)\bigr),\ \bigl(R_{0}^{(n+2)\bmod 4+1}(3),P(n,1)\bigr),\ \bigl(R_{0}^{(n-1)\bmod 4+1}(2),P(n,3)\bigr)\\
 & \bigl(R_{0}^{n\bmod 4+1}(1),P(n,4)\bigr),\ \bigl(uI,(u+1)P\bigr),\ n>1\\
\end{align*}
\end{fleqn}
\subsubsection*{Modules of the form $P(n,6)$}

Defect: $\partial P(n,6) = -1$, for $n\ge 0$.

\begin{fleqn}
\begin{align*}
P(0,6):\ & \bigl(R_{\infty}^{2}(1),P(0,2)\bigr),\ \bigl(R_{0}^{4}(3),P(0,3)\bigr),\ \bigl(R_{0}^{1}(2),P(0,4)\bigr),\ \bigl(R_{0}^{2}(1),P(0,5)\bigr),\ \bigl(I(1,4),2P(0,1)\bigr)\\
P(n,6):\ & \bigl(R_{\infty}^{(n+1)\bmod 2+1}(1),P(n,2)\bigr),\ \bigl(R_{0}^{(n+3)\bmod 4+1}(3),P(n,3)\bigr),\ \bigl(R_{0}^{n\bmod 4+1}(2),P(n,4)\bigr)\\
 & \bigl(R_{0}^{(n+1)\bmod 4+1}(1),P(n,5)\bigr),\ \bigl(uI,(u+1)P\bigr),\ n>0\\
\end{align*}
\end{fleqn}
\subsubsection{Schofield pairs associated to preinjective exceptional modules}

\begin{figure}[ht]


\begin{center}
\begin{scaletikzpicturetowidth}{\textwidth}

\end{scaletikzpicturetowidth}
\end{center}

\end{figure}
\medskip{}
\subsubsection*{Modules of the form $I(n,1)$}

Defect: $\partial I(n,1) = 1$, for $n\ge 0$.

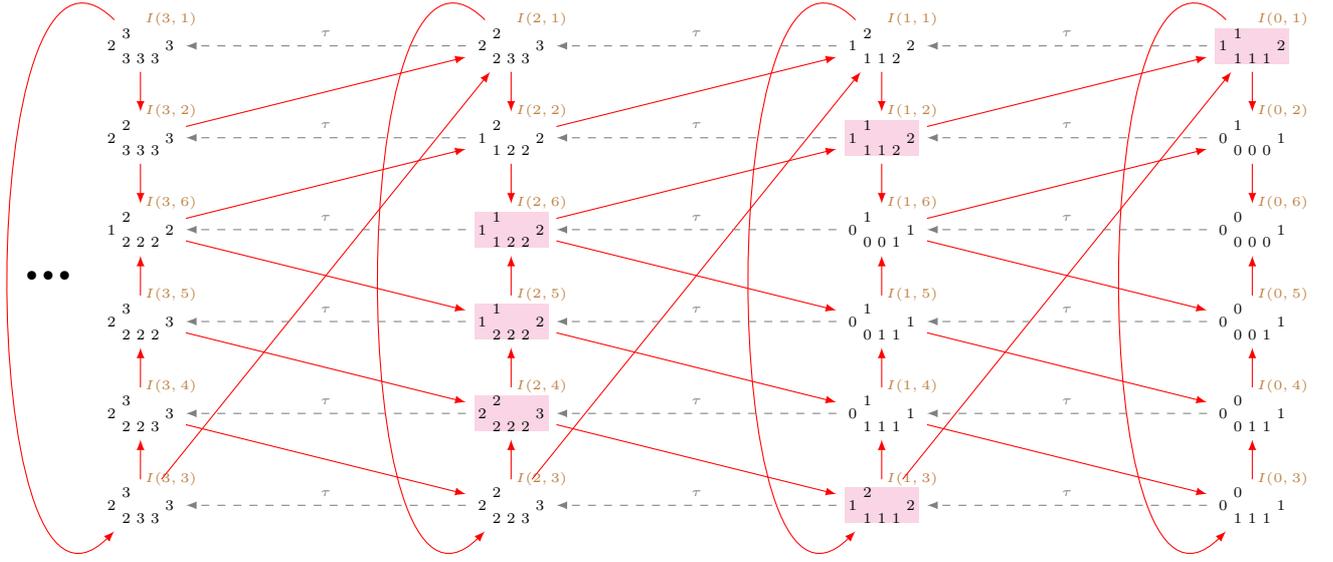
\begin{fleqn}
\begin{align*}
I(0,1):\ & \bigl(I(0,2),R_{\infty}^{2}(1)\bigr),\ \bigl(I(0,3),R_{0}^{2}(1)\bigr),\ \bigl(I(0,4),R_{0}^{2}(2)\bigr),\ \bigl(I(0,5),R_{0}^{2}(3)\bigr),\ \bigl(2I(0,6),P(1,4)\bigr)\\
I(n,1):\ & \bigl(I(n,2),R_{\infty}^{(-n+1)\bmod 2+1}(1)\bigr),\ \bigl(I(n,3),R_{0}^{(-n+1)\bmod 4+1}(1)\bigr),\ \bigl(I(n,4),R_{0}^{(-n+1)\bmod 4+1}(2)\bigr)\\
 & \bigl(I(n,5),R_{0}^{(-n+1)\bmod 4+1}(3)\bigr),\ \bigl((v+1)I,vP\bigr),\ n>0\\
\end{align*}
\end{fleqn}
\subsubsection*{Modules of the form $I(n,2)$}

Defect: $\partial I(n,2) = 1$, for $n\ge 0$.

\begin{fleqn}
\begin{align*}
I(0,2):\ & \bigl(I(0,6),R_{\infty}^{1}(1)\bigr)\\
I(1,2):\ & \bigl(I(0,1),R_{0}^{1}(1)\bigr),\ \bigl(I(0,3),R_{0}^{1}(2)\bigr),\ \bigl(I(0,4),R_{0}^{1}(3)\bigr),\ \bigl(I(1,6),R_{\infty}^{2}(1)\bigr),\ \bigl(2I(0,5),P(1,3)\bigr)\\
I(n,2):\ & \bigl(I(n-1,1),R_{0}^{(-n+1)\bmod 4+1}(1)\bigr),\ \bigl(I(n-1,3),R_{0}^{(-n+1)\bmod 4+1}(2)\bigr),\ \bigl(I(n-1,4),R_{0}^{(-n+1)\bmod 4+1}(3)\bigr)\\
 & \bigl(I(n,6),R_{\infty}^{(-n+2)\bmod 2+1}(1)\bigr),\ \bigl((v+1)I,vP\bigr),\ n>1\\
\end{align*}
\end{fleqn}
\subsubsection*{Modules of the form $I(n,3)$}

Defect: $\partial I(n,3) = 1$, for $n\ge 0$.

\begin{fleqn}
\begin{align*}
I(0,3):\ & \bigl(I(0,4),R_{0}^{3}(1)\bigr),\ \bigl(I(0,5),R_{0}^{3}(2)\bigr),\ \bigl(I(0,6),R_{0}^{3}(3)\bigr)\\
I(1,3):\ & \bigl(I(0,1),R_{\infty}^{1}(1)\bigr),\ \bigl(I(1,4),R_{0}^{2}(1)\bigr),\ \bigl(I(1,5),R_{0}^{2}(2)\bigr),\ \bigl(I(1,6),R_{0}^{2}(3)\bigr),\ \bigl(2I(0,2),P(0,5)\bigr)\\
I(n,3):\ & \bigl(I(n-1,1),R_{\infty}^{(-n+1)\bmod 2+1}(1)\bigr),\ \bigl(I(n,4),R_{0}^{(-n+2)\bmod 4+1}(1)\bigr),\ \bigl(I(n,5),R_{0}^{(-n+2)\bmod 4+1}(2)\bigr)\\
 & \bigl(I(n,6),R_{0}^{(-n+2)\bmod 4+1}(3)\bigr),\ \bigl((v+1)I,vP\bigr),\ n>1\\
\end{align*}
\end{fleqn}
\subsubsection*{Modules of the form $I(n,4)$}

Defect: $\partial I(n,4) = 1$, for $n\ge 0$.

\begin{fleqn}
\begin{align*}
I(0,4):\ & \bigl(I(0,5),R_{0}^{4}(1)\bigr),\ \bigl(I(0,6),R_{0}^{4}(2)\bigr)\\
I(1,4):\ & \bigl(I(0,2),R_{0}^{3}(3)\bigr),\ \bigl(I(0,3),R_{\infty}^{1}(1)\bigr),\ \bigl(I(1,5),R_{0}^{3}(1)\bigr),\ \bigl(I(1,6),R_{0}^{3}(2)\bigr)\\
I(2,4):\ & \bigl(I(1,2),R_{0}^{2}(3)\bigr),\ \bigl(I(1,3),R_{\infty}^{2}(1)\bigr),\ \bigl(I(2,5),R_{0}^{2}(1)\bigr),\ \bigl(I(2,6),R_{0}^{2}(2)\bigr),\ \bigl(3I(0,6),2P(1,4)\bigr)\\
I(n,4):\ & \bigl(I(n-1,2),R_{0}^{(-n+3)\bmod 4+1}(3)\bigr),\ \bigl(I(n-1,3),R_{\infty}^{(-n+3)\bmod 2+1}(1)\bigr),\ \bigl(I(n,5),R_{0}^{(-n+3)\bmod 4+1}(1)\bigr)\\
 & \bigl(I(n,6),R_{0}^{(-n+3)\bmod 4+1}(2)\bigr),\ \bigl((v+1)I,vP\bigr),\ n>2\\
\end{align*}
\end{fleqn}
\subsubsection*{Modules of the form $I(n,5)$}

Defect: $\partial I(n,5) = 1$, for $n\ge 0$.

\begin{fleqn}
\begin{align*}
I(0,5):\ & \bigl(I(0,6),R_{0}^{1}(1)\bigr)\\
I(1,5):\ & \bigl(I(0,2),R_{0}^{4}(2)\bigr),\ \bigl(I(0,4),R_{\infty}^{1}(1)\bigr),\ \bigl(I(1,6),R_{0}^{4}(1)\bigr)\\
I(2,5):\ & \bigl(I(0,1),R_{0}^{3}(3)\bigr),\ \bigl(I(1,2),R_{0}^{3}(2)\bigr),\ \bigl(I(1,4),R_{\infty}^{2}(1)\bigr),\ \bigl(I(2,6),R_{0}^{3}(1)\bigr),\ \bigl(2I(0,3),P(0,2)\bigr)\\
I(n,5):\ & \bigl(I(n-2,1),R_{0}^{(-n+4)\bmod 4+1}(3)\bigr),\ \bigl(I(n-1,2),R_{0}^{(-n+4)\bmod 4+1}(2)\bigr),\ \bigl(I(n-1,4),R_{\infty}^{(-n+3)\bmod 2+1}(1)\bigr)\\
 & \bigl(I(n,6),R_{0}^{(-n+4)\bmod 4+1}(1)\bigr),\ \bigl((v+1)I,vP\bigr),\ n>2\\
\end{align*}
\end{fleqn}
\subsubsection*{Modules of the form $I(n,6)$}

Defect: $\partial I(n,6) = 1$, for $n\ge 0$.

\begin{fleqn}
\begin{align*}
I(0,6):\ & - \\
I(1,6):\ & \bigl(I(0,2),R_{0}^{1}(1)\bigr),\ \bigl(I(0,5),R_{\infty}^{1}(1)\bigr)\\
I(2,6):\ & \bigl(I(0,1),R_{0}^{4}(2)\bigr),\ \bigl(I(0,3),R_{0}^{4}(3)\bigr),\ \bigl(I(1,2),R_{0}^{4}(1)\bigr),\ \bigl(I(1,5),R_{\infty}^{2}(1)\bigr),\ \bigl(2I(0,4),P(1,1)\bigr)\\
I(n,6):\ & \bigl(I(n-2,1),R_{0}^{(-n+5)\bmod 4+1}(2)\bigr),\ \bigl(I(n-2,3),R_{0}^{(-n+5)\bmod 4+1}(3)\bigr),\ \bigl(I(n-1,2),R_{0}^{(-n+5)\bmod 4+1}(1)\bigr)\\
 & \bigl(I(n-1,5),R_{\infty}^{(-n+3)\bmod 2+1}(1)\bigr),\ \bigl((v+1)I,vP\bigr),\ n>2\\
\end{align*}
\end{fleqn}
\subsubsection{Schofield pairs associated to regular exceptional modules}

\subsubsection*{The non-homogeneous tube $\mathcal{T}_{0}^{\Delta(\widetilde{\mathbb{A}}_{2,4})}$}

\begin{figure}[ht]



\begin{center}
\captionof{figure}{\vspace*{-10pt}$\mathcal{T}_{0}^{\Delta(\widetilde{\mathbb{A}}_{2,4})}$}
\begin{scaletikzpicturetowidth}{1.\textwidth}

\end{scaletikzpicturetowidth}
\end{center}

\end{figure}
\begin{fleqn}
\begin{align*}
R_{0}^{1}(1):\ & - \\
R_{0}^{1}(2):\ & \bigl(R_{0}^{2}(1),R_{0}^{1}(1)\bigr),\ \bigl(I(1,6),P(0,1)\bigr),\ \bigl(I(0,5),P(0,2)\bigr)\\
R_{0}^{2}(1):\ & \bigl(I(0,2),P(0,1)\bigr),\ \bigl(I(0,6),P(0,2)\bigr)\\
R_{0}^{2}(2):\ & \bigl(R_{0}^{3}(1),R_{0}^{2}(1)\bigr),\ \bigl(I(0,2),P(0,3)\bigr),\ \bigl(I(0,6),P(1,1)\bigr)\\
R_{0}^{3}(1):\ & - \\
R_{0}^{3}(2):\ & \bigl(R_{0}^{4}(1),R_{0}^{3}(1)\bigr)\\
R_{0}^{4}(1):\ & - \\
R_{0}^{4}(2):\ & \bigl(R_{0}^{1}(1),R_{0}^{4}(1)\bigr)\\
R_{0}^{4}(3):\ & \bigl(R_{0}^{1}(2),R_{0}^{4}(1)\bigr),\ \bigl(R_{0}^{2}(1),R_{0}^{4}(2)\bigr),\ \bigl(I(1,5),P(0,1)\bigr),\ \bigl(I(0,4),P(0,2)\bigr)\\
R_{0}^{1}(3):\ & \bigl(R_{0}^{2}(2),R_{0}^{1}(1)\bigr),\ \bigl(R_{0}^{3}(1),R_{0}^{1}(2)\bigr),\ \bigl(I(1,6),P(0,3)\bigr),\ \bigl(I(0,5),P(1,1)\bigr)\\
R_{0}^{2}(3):\ & \bigl(R_{0}^{3}(2),R_{0}^{2}(1)\bigr),\ \bigl(R_{0}^{4}(1),R_{0}^{2}(2)\bigr),\ \bigl(I(0,2),P(0,4)\bigr),\ \bigl(I(0,6),P(1,3)\bigr)\\
R_{0}^{3}(3):\ & \bigl(R_{0}^{4}(2),R_{0}^{3}(1)\bigr),\ \bigl(R_{0}^{1}(1),R_{0}^{3}(2)\bigr)\\
\end{align*}
\end{fleqn}
\subsubsection*{The non-homogeneous tube $\mathcal{T}_{\infty}^{\Delta(\widetilde{\mathbb{A}}_{2,4})}$}

\begin{figure}[ht]



\begin{center}
\captionof{figure}{\vspace*{-10pt}$\mathcal{T}_{\infty}^{\Delta(\widetilde{\mathbb{A}}_{2,4})}$}
\begin{scaletikzpicturetowidth}{0.5\textwidth}
 $\end{lrbox}
\subsection[Schofield pairs for the quiver $\Delta(\A_{3,3})$]{Schofield pairs for the quiver $\Delta(\A_{3,3})$ -- {\usebox{\boxAAAlAAA}}}
\[
\vcenter{\hbox{\xymatrix{  & 2\ar[dl] & 3\ar[l] &   &   &  \\
1 &   &   & 6\ar[ul]\ar[dl] &   &  \\
  & 4\ar[ul] & 5\ar[l] &   &   &  }}}\qquad C_{\Delta(\widetilde{\mathbb{A}}_{3,3})} = %
\]
\subsubsection{Schofield pairs associated to preprojective exceptional modules}

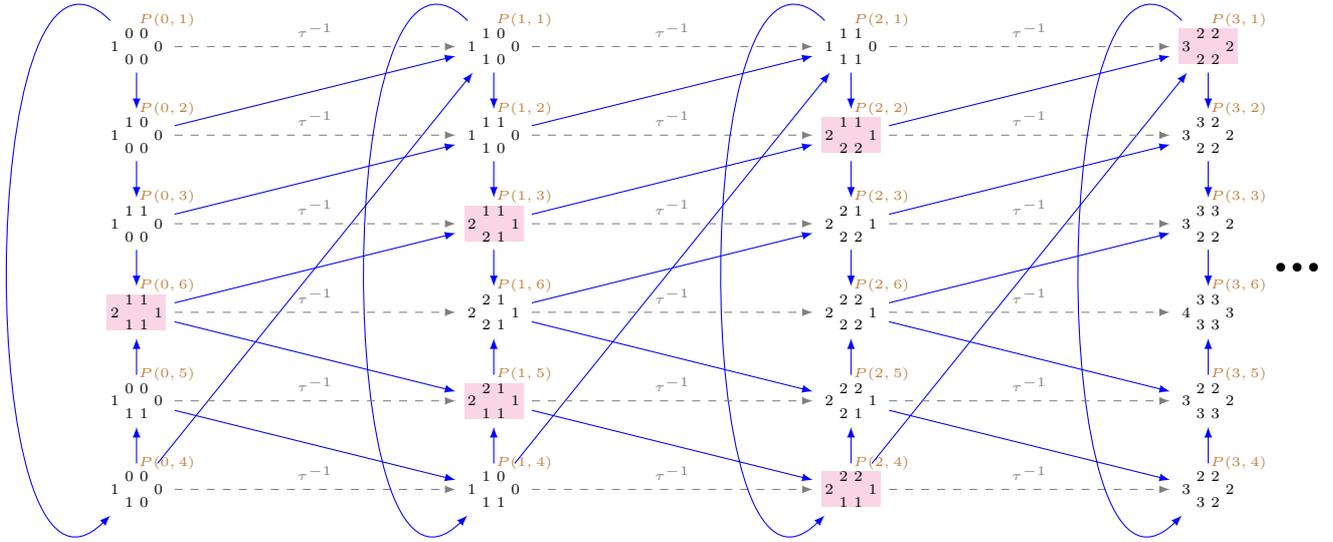
\begin{figure}[ht]


\begin{center}
\begin{scaletikzpicturetowidth}{\textwidth}

\end{scaletikzpicturetowidth}
\end{center}


\end{figure}
\medskip{}
\subsubsection*{Modules of the form $P(n,1)$}

Defect: $\partial P(n,1) = -1$, for $n\ge 0$.

\begin{fleqn}
\begin{align*}
P(0,1):\ & - \\
P(1,1):\ & \bigl(R_{0}^{3}(1),P(0,2)\bigr),\ \bigl(R_{\infty}^{3}(1),P(0,4)\bigr)\\
P(2,1):\ & \bigl(R_{0}^{3}(2),P(0,3)\bigr),\ \bigl(R_{\infty}^{3}(2),P(0,5)\bigr),\ \bigl(R_{0}^{1}(1),P(1,2)\bigr),\ \bigl(R_{\infty}^{1}(1),P(1,4)\bigr)\\
P(3,1):\ & \bigl(R_{0}^{1}(2),P(1,3)\bigr),\ \bigl(R_{\infty}^{1}(2),P(1,5)\bigr),\ \bigl(R_{0}^{2}(1),P(2,2)\bigr),\ \bigl(R_{\infty}^{2}(1),P(2,4)\bigr),\ \bigl(2I(2,6),3P(0,1)\bigr)\\
P(n,1):\ & \bigl(R_{0}^{(n-3)\bmod 3+1}(2),P(n-2,3)\bigr),\ \bigl(R_{\infty}^{(n-3)\bmod 3+1}(2),P(n-2,5)\bigr),\ \bigl(R_{0}^{(n-2)\bmod 3+1}(1),P(n-1,2)\bigr)\\
 & \bigl(R_{\infty}^{(n-2)\bmod 3+1}(1),P(n-1,4)\bigr),\ \bigl(uI,(u+1)P\bigr),\ n>3\\
\end{align*}
\end{fleqn}
\subsubsection*{Modules of the form $P(n,2)$}

Defect: $\partial P(n,2) = -1$, for $n\ge 0$.

\begin{fleqn}
\begin{align*}
P(0,2):\ & \bigl(R_{\infty}^{3}(1),P(0,1)\bigr)\\
P(1,2):\ & \bigl(R_{0}^{3}(1),P(0,3)\bigr),\ \bigl(R_{\infty}^{3}(2),P(0,4)\bigr),\ \bigl(R_{\infty}^{1}(1),P(1,1)\bigr)\\
P(2,2):\ & \bigl(R_{0}^{3}(2),P(0,6)\bigr),\ \bigl(R_{0}^{1}(1),P(1,3)\bigr),\ \bigl(R_{\infty}^{1}(2),P(1,4)\bigr),\ \bigl(R_{\infty}^{2}(1),P(2,1)\bigr),\ \bigl(I(0,2),2P(0,5)\bigr)\\
P(n,2):\ & \bigl(R_{0}^{n\bmod 3+1}(2),P(n-2,6)\bigr),\ \bigl(R_{0}^{(n-2)\bmod 3+1}(1),P(n-1,3)\bigr),\ \bigl(R_{\infty}^{(n-2)\bmod 3+1}(2),P(n-1,4)\bigr)\\
 & \bigl(R_{\infty}^{(n-1)\bmod 3+1}(1),P(n,1)\bigr),\ \bigl(uI,(u+1)P\bigr),\ n>2\\
\end{align*}
\end{fleqn}
\subsubsection*{Modules of the form $P(n,3)$}

Defect: $\partial P(n,3) = -1$, for $n\ge 0$.

\begin{fleqn}
\begin{align*}
P(0,3):\ & \bigl(R_{\infty}^{3}(2),P(0,1)\bigr),\ \bigl(R_{\infty}^{1}(1),P(0,2)\bigr)\\
P(1,3):\ & \bigl(R_{0}^{2}(2),P(0,5)\bigr),\ \bigl(R_{0}^{3}(1),P(0,6)\bigr),\ \bigl(R_{\infty}^{1}(2),P(1,1)\bigr),\ \bigl(R_{\infty}^{2}(1),P(1,2)\bigr),\ \bigl(I(1,3),2P(0,4)\bigr)\\
P(n,3):\ & \bigl(R_{0}^{n\bmod 3+1}(2),P(n-1,5)\bigr),\ \bigl(R_{0}^{(n+1)\bmod 3+1}(1),P(n-1,6)\bigr),\ \bigl(R_{\infty}^{(n-1)\bmod 3+1}(2),P(n,1)\bigr)\\
 & \bigl(R_{\infty}^{n\bmod 3+1}(1),P(n,2)\bigr),\ \bigl(uI,(u+1)P\bigr),\ n>1\\
\end{align*}
\end{fleqn}
\subsubsection*{Modules of the form $P(n,4)$}

Defect: $\partial P(n,4) = -1$, for $n\ge 0$.

\begin{fleqn}
\begin{align*}
P(0,4):\ & \bigl(R_{0}^{3}(1),P(0,1)\bigr)\\
P(1,4):\ & \bigl(R_{0}^{3}(2),P(0,2)\bigr),\ \bigl(R_{\infty}^{3}(1),P(0,5)\bigr),\ \bigl(R_{0}^{1}(1),P(1,1)\bigr)\\
P(2,4):\ & \bigl(R_{\infty}^{3}(2),P(0,6)\bigr),\ \bigl(R_{0}^{1}(2),P(1,2)\bigr),\ \bigl(R_{\infty}^{1}(1),P(1,5)\bigr),\ \bigl(R_{0}^{2}(1),P(2,1)\bigr),\ \bigl(I(0,4),2P(0,3)\bigr)\\
P(n,4):\ & \bigl(R_{\infty}^{n\bmod 3+1}(2),P(n-2,6)\bigr),\ \bigl(R_{0}^{(n-2)\bmod 3+1}(2),P(n-1,2)\bigr),\ \bigl(R_{\infty}^{(n-2)\bmod 3+1}(1),P(n-1,5)\bigr)\\
 & \bigl(R_{0}^{(n-1)\bmod 3+1}(1),P(n,1)\bigr),\ \bigl(uI,(u+1)P\bigr),\ n>2\\
\end{align*}
\end{fleqn}
\subsubsection*{Modules of the form $P(n,5)$}

Defect: $\partial P(n,5) = -1$, for $n\ge 0$.

\begin{fleqn}
\begin{align*}
P(0,5):\ & \bigl(R_{0}^{3}(2),P(0,1)\bigr),\ \bigl(R_{0}^{1}(1),P(0,4)\bigr)\\
P(1,5):\ & \bigl(R_{\infty}^{2}(2),P(0,3)\bigr),\ \bigl(R_{\infty}^{3}(1),P(0,6)\bigr),\ \bigl(R_{0}^{1}(2),P(1,1)\bigr),\ \bigl(R_{0}^{2}(1),P(1,4)\bigr),\ \bigl(I(1,5),2P(0,2)\bigr)\\
P(n,5):\ & \bigl(R_{\infty}^{n\bmod 3+1}(2),P(n-1,3)\bigr),\ \bigl(R_{\infty}^{(n+1)\bmod 3+1}(1),P(n-1,6)\bigr),\ \bigl(R_{0}^{(n-1)\bmod 3+1}(2),P(n,1)\bigr)\\
 & \bigl(R_{0}^{n\bmod 3+1}(1),P(n,4)\bigr),\ \bigl(uI,(u+1)P\bigr),\ n>1\\
\end{align*}
\end{fleqn}
\subsubsection*{Modules of the form $P(n,6)$}

Defect: $\partial P(n,6) = -1$, for $n\ge 0$.

\begin{fleqn}
\begin{align*}
P(0,6):\ & \bigl(R_{\infty}^{1}(2),P(0,2)\bigr),\ \bigl(R_{\infty}^{2}(1),P(0,3)\bigr),\ \bigl(R_{0}^{1}(2),P(0,4)\bigr),\ \bigl(R_{0}^{2}(1),P(0,5)\bigr),\ \bigl(I(2,6),2P(0,1)\bigr)\\
P(n,6):\ & \bigl(R_{\infty}^{n\bmod 3+1}(2),P(n,2)\bigr),\ \bigl(R_{\infty}^{(n+1)\bmod 3+1}(1),P(n,3)\bigr),\ \bigl(R_{0}^{n\bmod 3+1}(2),P(n,4)\bigr)\\
 & \bigl(R_{0}^{(n+1)\bmod 3+1}(1),P(n,5)\bigr),\ \bigl(uI,(u+1)P\bigr),\ n>0\\
\end{align*}
\end{fleqn}
\subsubsection{Schofield pairs associated to preinjective exceptional modules}

\begin{figure}[ht]


\begin{center}
\begin{scaletikzpicturetowidth}{\textwidth}

\end{scaletikzpicturetowidth}
\end{center}

\end{figure}
\medskip{}
\subsubsection*{Modules of the form $I(n,1)$}

Defect: $\partial I(n,1) = 1$, for $n\ge 0$.

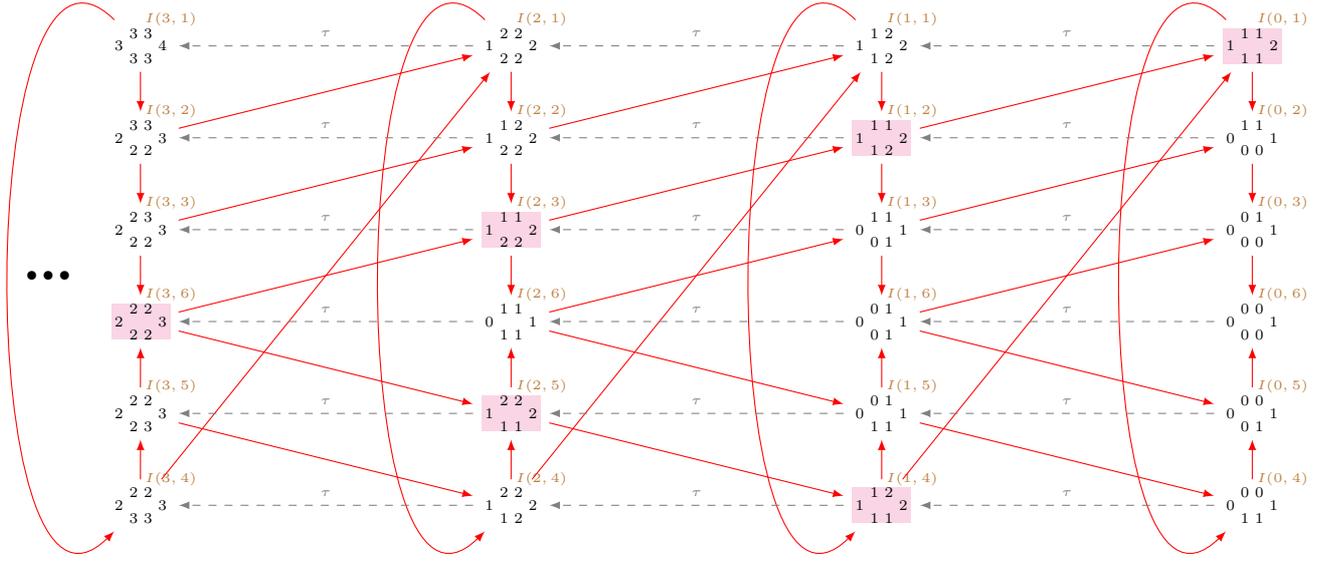
\begin{fleqn}
\begin{align*}
I(0,1):\ & \bigl(I(0,2),R_{\infty}^{2}(1)\bigr),\ \bigl(I(0,3),R_{\infty}^{2}(2)\bigr),\ \bigl(I(0,4),R_{0}^{2}(1)\bigr),\ \bigl(I(0,5),R_{0}^{2}(2)\bigr),\ \bigl(2I(0,6),P(2,1)\bigr)\\
I(n,1):\ & \bigl(I(n,2),R_{\infty}^{(-n+1)\bmod 3+1}(1)\bigr),\ \bigl(I(n,3),R_{\infty}^{(-n+1)\bmod 3+1}(2)\bigr),\ \bigl(I(n,4),R_{0}^{(-n+1)\bmod 3+1}(1)\bigr)\\
 & \bigl(I(n,5),R_{0}^{(-n+1)\bmod 3+1}(2)\bigr),\ \bigl((v+1)I,vP\bigr),\ n>0\\
\end{align*}
\end{fleqn}
\subsubsection*{Modules of the form $I(n,2)$}

Defect: $\partial I(n,2) = 1$, for $n\ge 0$.

\begin{fleqn}
\begin{align*}
I(0,2):\ & \bigl(I(0,3),R_{\infty}^{3}(1)\bigr),\ \bigl(I(0,6),R_{\infty}^{3}(2)\bigr)\\
I(1,2):\ & \bigl(I(0,1),R_{0}^{1}(1)\bigr),\ \bigl(I(0,4),R_{0}^{1}(2)\bigr),\ \bigl(I(1,3),R_{\infty}^{2}(1)\bigr),\ \bigl(I(1,6),R_{\infty}^{2}(2)\bigr),\ \bigl(2I(0,5),P(1,2)\bigr)\\
I(n,2):\ & \bigl(I(n-1,1),R_{0}^{(-n+1)\bmod 3+1}(1)\bigr),\ \bigl(I(n-1,4),R_{0}^{(-n+1)\bmod 3+1}(2)\bigr),\ \bigl(I(n,3),R_{\infty}^{(-n+2)\bmod 3+1}(1)\bigr)\\
 & \bigl(I(n,6),R_{\infty}^{(-n+2)\bmod 3+1}(2)\bigr),\ \bigl((v+1)I,vP\bigr),\ n>1\\
\end{align*}
\end{fleqn}
\subsubsection*{Modules of the form $I(n,3)$}

Defect: $\partial I(n,3) = 1$, for $n\ge 0$.

\begin{fleqn}
\begin{align*}
I(0,3):\ & \bigl(I(0,6),R_{\infty}^{1}(1)\bigr)\\
I(1,3):\ & \bigl(I(0,2),R_{0}^{1}(1)\bigr),\ \bigl(I(0,5),R_{\infty}^{3}(2)\bigr),\ \bigl(I(1,6),R_{\infty}^{3}(1)\bigr)\\
I(2,3):\ & \bigl(I(0,1),R_{0}^{3}(2)\bigr),\ \bigl(I(1,2),R_{0}^{3}(1)\bigr),\ \bigl(I(1,5),R_{\infty}^{2}(2)\bigr),\ \bigl(I(2,6),R_{\infty}^{2}(1)\bigr),\ \bigl(2I(0,4),P(0,3)\bigr)\\
I(n,3):\ & \bigl(I(n-2,1),R_{0}^{(-n+4)\bmod 3+1}(2)\bigr),\ \bigl(I(n-1,2),R_{0}^{(-n+4)\bmod 3+1}(1)\bigr),\ \bigl(I(n-1,5),R_{\infty}^{(-n+3)\bmod 3+1}(2)\bigr)\\
 & \bigl(I(n,6),R_{\infty}^{(-n+3)\bmod 3+1}(1)\bigr),\ \bigl((v+1)I,vP\bigr),\ n>2\\
\end{align*}
\end{fleqn}
\subsubsection*{Modules of the form $I(n,4)$}

Defect: $\partial I(n,4) = 1$, for $n\ge 0$.

\begin{fleqn}
\begin{align*}
I(0,4):\ & \bigl(I(0,5),R_{0}^{3}(1)\bigr),\ \bigl(I(0,6),R_{0}^{3}(2)\bigr)\\
I(1,4):\ & \bigl(I(0,1),R_{\infty}^{1}(1)\bigr),\ \bigl(I(0,2),R_{\infty}^{1}(2)\bigr),\ \bigl(I(1,5),R_{0}^{2}(1)\bigr),\ \bigl(I(1,6),R_{0}^{2}(2)\bigr),\ \bigl(2I(0,3),P(1,4)\bigr)\\
I(n,4):\ & \bigl(I(n-1,1),R_{\infty}^{(-n+1)\bmod 3+1}(1)\bigr),\ \bigl(I(n-1,2),R_{\infty}^{(-n+1)\bmod 3+1}(2)\bigr),\ \bigl(I(n,5),R_{0}^{(-n+2)\bmod 3+1}(1)\bigr)\\
 & \bigl(I(n,6),R_{0}^{(-n+2)\bmod 3+1}(2)\bigr),\ \bigl((v+1)I,vP\bigr),\ n>1\\
\end{align*}
\end{fleqn}
\subsubsection*{Modules of the form $I(n,5)$}

Defect: $\partial I(n,5) = 1$, for $n\ge 0$.

\begin{fleqn}
\begin{align*}
I(0,5):\ & \bigl(I(0,6),R_{0}^{1}(1)\bigr)\\
I(1,5):\ & \bigl(I(0,3),R_{0}^{3}(2)\bigr),\ \bigl(I(0,4),R_{\infty}^{1}(1)\bigr),\ \bigl(I(1,6),R_{0}^{3}(1)\bigr)\\
I(2,5):\ & \bigl(I(0,1),R_{\infty}^{3}(2)\bigr),\ \bigl(I(1,3),R_{0}^{2}(2)\bigr),\ \bigl(I(1,4),R_{\infty}^{3}(1)\bigr),\ \bigl(I(2,6),R_{0}^{2}(1)\bigr),\ \bigl(2I(0,2),P(0,5)\bigr)\\
I(n,5):\ & \bigl(I(n-2,1),R_{\infty}^{(-n+4)\bmod 3+1}(2)\bigr),\ \bigl(I(n-1,3),R_{0}^{(-n+3)\bmod 3+1}(2)\bigr),\ \bigl(I(n-1,4),R_{\infty}^{(-n+4)\bmod 3+1}(1)\bigr)\\
 & \bigl(I(n,6),R_{0}^{(-n+3)\bmod 3+1}(1)\bigr),\ \bigl((v+1)I,vP\bigr),\ n>2\\
\end{align*}
\end{fleqn}
\subsubsection*{Modules of the form $I(n,6)$}

Defect: $\partial I(n,6) = 1$, for $n\ge 0$.

\begin{fleqn}
\begin{align*}
I(0,6):\ & - \\
I(1,6):\ & \bigl(I(0,3),R_{0}^{1}(1)\bigr),\ \bigl(I(0,5),R_{\infty}^{1}(1)\bigr)\\
I(2,6):\ & \bigl(I(0,2),R_{0}^{3}(2)\bigr),\ \bigl(I(0,4),R_{\infty}^{3}(2)\bigr),\ \bigl(I(1,3),R_{0}^{3}(1)\bigr),\ \bigl(I(1,5),R_{\infty}^{3}(1)\bigr)\\
I(3,6):\ & \bigl(I(1,2),R_{0}^{2}(2)\bigr),\ \bigl(I(1,4),R_{\infty}^{2}(2)\bigr),\ \bigl(I(2,3),R_{0}^{2}(1)\bigr),\ \bigl(I(2,5),R_{\infty}^{2}(1)\bigr),\ \bigl(3I(0,6),2P(2,1)\bigr)\\
I(n,6):\ & \bigl(I(n-2,2),R_{0}^{(-n+4)\bmod 3+1}(2)\bigr),\ \bigl(I(n-2,4),R_{\infty}^{(-n+4)\bmod 3+1}(2)\bigr),\ \bigl(I(n-1,3),R_{0}^{(-n+4)\bmod 3+1}(1)\bigr)\\
 & \bigl(I(n-1,5),R_{\infty}^{(-n+4)\bmod 3+1}(1)\bigr),\ \bigl((v+1)I,vP\bigr),\ n>3\\
\end{align*}
\end{fleqn}
\subsubsection{Schofield pairs associated to regular exceptional modules}

\subsubsection*{The non-homogeneous tube $\mathcal{T}_{0}^{\Delta(\widetilde{\mathbb{A}}_{3,3})}$}

\begin{figure}[ht]



\begin{center}
\captionof{figure}{\vspace*{-10pt}$\mathcal{T}_{0}^{\Delta(\widetilde{\mathbb{A}}_{3,3})}$}
\begin{scaletikzpicturetowidth}{0.80000000000000004\textwidth}

\end{scaletikzpicturetowidth}
\end{center}

\end{figure}
\begin{fleqn}
\begin{align*}
R_{0}^{1}(1):\ & - \\
R_{0}^{1}(2):\ & \bigl(R_{0}^{2}(1),R_{0}^{1}(1)\bigr),\ \bigl(I(1,3),P(0,1)\bigr),\ \bigl(I(1,6),P(0,2)\bigr),\ \bigl(I(0,5),P(0,3)\bigr)\\
R_{0}^{2}(1):\ & \bigl(I(0,2),P(0,1)\bigr),\ \bigl(I(0,3),P(0,2)\bigr),\ \bigl(I(0,6),P(0,3)\bigr)\\
R_{0}^{2}(2):\ & \bigl(R_{0}^{3}(1),R_{0}^{2}(1)\bigr),\ \bigl(I(0,2),P(0,4)\bigr),\ \bigl(I(0,3),P(1,1)\bigr),\ \bigl(I(0,6),P(1,2)\bigr)\\
R_{0}^{3}(1):\ & - \\
R_{0}^{3}(2):\ & \bigl(R_{0}^{1}(1),R_{0}^{3}(1)\bigr)\\
\end{align*}
\end{fleqn}
\subsubsection*{The non-homogeneous tube $\mathcal{T}_{\infty}^{\Delta(\widetilde{\mathbb{A}}_{3,3})}$}

\begin{figure}[ht]



\begin{center}
\captionof{figure}{\vspace*{-10pt}$\mathcal{T}_{\infty}^{\Delta(\widetilde{\mathbb{A}}_{3,3})}$}
\begin{scaletikzpicturetowidth}{0.80000000000000004\textwidth}

\end{scaletikzpicturetowidth}
\end{center}

\end{figure}
\begin{fleqn}
\begin{align*}
R_{\infty}^{1}(1):\ & - \\
R_{\infty}^{1}(2):\ & \bigl(R_{\infty}^{2}(1),R_{\infty}^{1}(1)\bigr),\ \bigl(I(1,5),P(0,1)\bigr),\ \bigl(I(1,6),P(0,4)\bigr),\ \bigl(I(0,3),P(0,5)\bigr)\\
R_{\infty}^{2}(1):\ & \bigl(I(0,4),P(0,1)\bigr),\ \bigl(I(0,5),P(0,4)\bigr),\ \bigl(I(0,6),P(0,5)\bigr)\\
R_{\infty}^{2}(2):\ & \bigl(R_{\infty}^{3}(1),R_{\infty}^{2}(1)\bigr),\ \bigl(I(0,4),P(0,2)\bigr),\ \bigl(I(0,5),P(1,1)\bigr),\ \bigl(I(0,6),P(1,4)\bigr)\\
R_{\infty}^{3}(1):\ & - \\
R_{\infty}^{3}(2):\ & \bigl(R_{\infty}^{1}(1),R_{\infty}^{3}(1)\bigr)\\
\end{align*}
\end{fleqn}
\clearpage

\begin{lrbox}{\boxAAAlAAAA}$\delta =\begin{smallmatrix}&1&1\\1&&&&1\\&1&1&1\end{smallmatrix} $\end{lrbox}
\subsection[Schofield pairs for the quiver $\Delta(\A_{3,4})$]{Schofield pairs for the quiver $\Delta(\A_{3,4})$ -- {\usebox{\boxAAAlAAAA}}}
\[
\xymatrix{  & 2\ar[dl] & 3\ar[l] &   &   &   &  \\
1 &   &   &   & 7\ar[ull]\ar[dl] &   &  \\
  & 4\ar[ul] & 5\ar[l] & 6\ar[l] &   &   &  }
\]
\medskip{}
\[
 C_{\Delta(\widetilde{\mathbb{A}}_{3,4})} = \begin{bmatrix}1 & 1 & 1 & 1 & 1 & 1 & 2\\
0 & 1 & 1 & 0 & 0 & 0 & 1\\
0 & 0 & 1 & 0 & 0 & 0 & 1\\
0 & 0 & 0 & 1 & 1 & 1 & 1\\
0 & 0 & 0 & 0 & 1 & 1 & 1\\
0 & 0 & 0 & 0 & 0 & 1 & 1\\
0 & 0 & 0 & 0 & 0 & 0 & 1\end{bmatrix}\qquad\Phi_{\Delta(\widetilde{\mathbb{A}}_{3,4})} = \begin{bmatrix}-1 & 1 & 0 & 1 & 0 & 0 & 0\\
-1 & 0 & 1 & 1 & 0 & 0 & 0\\
-1 & 0 & 0 & 1 & 0 & 0 & 1\\
-1 & 1 & 0 & 0 & 1 & 0 & 0\\
-1 & 1 & 0 & 0 & 0 & 1 & 0\\
-1 & 1 & 0 & 0 & 0 & 0 & 1\\
-2 & 1 & 0 & 1 & 0 & 0 & 1\end{bmatrix}
\]

\medskip{}\subsubsection{Schofield pairs associated to preprojective exceptional modules}

\begin{figure}[ht]


\begin{center}
\begin{scaletikzpicturetowidth}{\textwidth}

\end{scaletikzpicturetowidth}
\end{center}


\end{figure}
\medskip{}
\subsubsection*{Modules of the form $P(n,1)$}

Defect: $\partial P(n,1) = -1$, for $n\ge 0$.

\begin{fleqn}
\begin{align*}
P(0,1):\ & - \\
P(1,1):\ & \bigl(R_{0}^{3}(1),P(0,2)\bigr),\ \bigl(R_{\infty}^{3}(1),P(0,4)\bigr)\\
P(2,1):\ & \bigl(R_{0}^{3}(2),P(0,3)\bigr),\ \bigl(R_{\infty}^{3}(2),P(0,5)\bigr),\ \bigl(R_{0}^{4}(1),P(1,2)\bigr),\ \bigl(R_{\infty}^{1}(1),P(1,4)\bigr)\\
P(3,1):\ & \bigl(R_{0}^{3}(3),P(0,7)\bigr),\ \bigl(R_{0}^{4}(2),P(1,3)\bigr),\ \bigl(R_{\infty}^{1}(2),P(1,5)\bigr),\ \bigl(R_{0}^{1}(1),P(2,2)\bigr),\ \bigl(R_{\infty}^{2}(1),P(2,4)\bigr)\\
 & \bigl(I(0,2),2P(0,6)\bigr)\\
P(n,1):\ & \bigl(R_{0}^{(n-1)\bmod 4+1}(3),P(n-3,7)\bigr),\ \bigl(R_{0}^{n\bmod 4+1}(2),P(n-2,3)\bigr),\ \bigl(R_{\infty}^{(n-3)\bmod 3+1}(2),P(n-2,5)\bigr)\\
 & \bigl(R_{0}^{(n-3)\bmod 4+1}(1),P(n-1,2)\bigr),\ \bigl(R_{\infty}^{(n-2)\bmod 3+1}(1),P(n-1,4)\bigr),\ \bigl(uI,(u+1)P\bigr),\ n>3\\
\end{align*}
\end{fleqn}
\subsubsection*{Modules of the form $P(n,2)$}

Defect: $\partial P(n,2) = -1$, for $n\ge 0$.

\begin{fleqn}
\begin{align*}
P(0,2):\ & \bigl(R_{\infty}^{3}(1),P(0,1)\bigr)\\
P(1,2):\ & \bigl(R_{0}^{3}(1),P(0,3)\bigr),\ \bigl(R_{\infty}^{3}(2),P(0,4)\bigr),\ \bigl(R_{\infty}^{1}(1),P(1,1)\bigr)\\
P(2,2):\ & \bigl(R_{0}^{2}(3),P(0,6)\bigr),\ \bigl(R_{0}^{3}(2),P(0,7)\bigr),\ \bigl(R_{0}^{4}(1),P(1,3)\bigr),\ \bigl(R_{\infty}^{1}(2),P(1,4)\bigr),\ \bigl(R_{\infty}^{2}(1),P(2,1)\bigr)\\
 & \bigl(I(1,3),2P(0,5)\bigr)\\
P(n,2):\ & \bigl(R_{0}^{(n-1)\bmod 4+1}(3),P(n-2,6)\bigr),\ \bigl(R_{0}^{n\bmod 4+1}(2),P(n-2,7)\bigr),\ \bigl(R_{0}^{(n+1)\bmod 4+1}(1),P(n-1,3)\bigr)\\
 & \bigl(R_{\infty}^{(n-2)\bmod 3+1}(2),P(n-1,4)\bigr),\ \bigl(R_{\infty}^{(n-1)\bmod 3+1}(1),P(n,1)\bigr),\ \bigl(uI,(u+1)P\bigr),\ n>2\\
\end{align*}
\end{fleqn}
\subsubsection*{Modules of the form $P(n,3)$}

Defect: $\partial P(n,3) = -1$, for $n\ge 0$.

\begin{fleqn}
\begin{align*}
P(0,3):\ & \bigl(R_{\infty}^{3}(2),P(0,1)\bigr),\ \bigl(R_{\infty}^{1}(1),P(0,2)\bigr)\\
P(1,3):\ & \bigl(R_{0}^{1}(3),P(0,5)\bigr),\ \bigl(R_{0}^{2}(2),P(0,6)\bigr),\ \bigl(R_{0}^{3}(1),P(0,7)\bigr),\ \bigl(R_{\infty}^{1}(2),P(1,1)\bigr),\ \bigl(R_{\infty}^{2}(1),P(1,2)\bigr)\\
 & \bigl(I(2,7),2P(0,4)\bigr)\\
P(n,3):\ & \bigl(R_{0}^{(n-1)\bmod 4+1}(3),P(n-1,5)\bigr),\ \bigl(R_{0}^{n\bmod 4+1}(2),P(n-1,6)\bigr),\ \bigl(R_{0}^{(n+1)\bmod 4+1}(1),P(n-1,7)\bigr)\\
 & \bigl(R_{\infty}^{(n-1)\bmod 3+1}(2),P(n,1)\bigr),\ \bigl(R_{\infty}^{n\bmod 3+1}(1),P(n,2)\bigr),\ \bigl(uI,(u+1)P\bigr),\ n>1\\
\end{align*}
\end{fleqn}
\subsubsection*{Modules of the form $P(n,4)$}

Defect: $\partial P(n,4) = -1$, for $n\ge 0$.

\begin{fleqn}
\begin{align*}
P(0,4):\ & \bigl(R_{0}^{3}(1),P(0,1)\bigr)\\
P(1,4):\ & \bigl(R_{0}^{3}(2),P(0,2)\bigr),\ \bigl(R_{\infty}^{3}(1),P(0,5)\bigr),\ \bigl(R_{0}^{4}(1),P(1,1)\bigr)\\
P(2,4):\ & \bigl(R_{0}^{3}(3),P(0,3)\bigr),\ \bigl(R_{\infty}^{3}(2),P(0,6)\bigr),\ \bigl(R_{0}^{4}(2),P(1,2)\bigr),\ \bigl(R_{\infty}^{1}(1),P(1,5)\bigr),\ \bigl(R_{0}^{1}(1),P(2,1)\bigr)\\
P(3,4):\ & \bigl(R_{0}^{4}(3),P(1,3)\bigr),\ \bigl(R_{\infty}^{1}(2),P(1,6)\bigr),\ \bigl(R_{0}^{1}(2),P(2,2)\bigr),\ \bigl(R_{\infty}^{2}(1),P(2,5)\bigr),\ \bigl(R_{0}^{2}(1),P(3,1)\bigr)\\
 & \bigl(2I(2,6),3P(0,1)\bigr)\\
P(n,4):\ & \bigl(R_{0}^{n\bmod 4+1}(3),P(n-2,3)\bigr),\ \bigl(R_{\infty}^{(n-3)\bmod 3+1}(2),P(n-2,6)\bigr),\ \bigl(R_{0}^{(n-3)\bmod 4+1}(2),P(n-1,2)\bigr)\\
 & \bigl(R_{\infty}^{(n-2)\bmod 3+1}(1),P(n-1,5)\bigr),\ \bigl(R_{0}^{(n-2)\bmod 4+1}(1),P(n,1)\bigr),\ \bigl(uI,(u+1)P\bigr),\ n>3\\
\end{align*}
\end{fleqn}
\subsubsection*{Modules of the form $P(n,5)$}

Defect: $\partial P(n,5) = -1$, for $n\ge 0$.

\begin{fleqn}
\begin{align*}
P(0,5):\ & \bigl(R_{0}^{3}(2),P(0,1)\bigr),\ \bigl(R_{0}^{4}(1),P(0,4)\bigr)\\
P(1,5):\ & \bigl(R_{0}^{3}(3),P(0,2)\bigr),\ \bigl(R_{\infty}^{3}(1),P(0,6)\bigr),\ \bigl(R_{0}^{4}(2),P(1,1)\bigr),\ \bigl(R_{0}^{1}(1),P(1,4)\bigr)\\
P(2,5):\ & \bigl(R_{\infty}^{3}(2),P(0,7)\bigr),\ \bigl(R_{0}^{4}(3),P(1,2)\bigr),\ \bigl(R_{\infty}^{1}(1),P(1,6)\bigr),\ \bigl(R_{0}^{1}(2),P(2,1)\bigr),\ \bigl(R_{0}^{2}(1),P(2,4)\bigr)\\
 & \bigl(I(0,4),2P(0,3)\bigr)\\
P(n,5):\ & \bigl(R_{\infty}^{n\bmod 3+1}(2),P(n-2,7)\bigr),\ \bigl(R_{0}^{(n+1)\bmod 4+1}(3),P(n-1,2)\bigr),\ \bigl(R_{\infty}^{(n-2)\bmod 3+1}(1),P(n-1,6)\bigr)\\
 & \bigl(R_{0}^{(n-2)\bmod 4+1}(2),P(n,1)\bigr),\ \bigl(R_{0}^{(n-1)\bmod 4+1}(1),P(n,4)\bigr),\ \bigl(uI,(u+1)P\bigr),\ n>2\\
\end{align*}
\end{fleqn}
\subsubsection*{Modules of the form $P(n,6)$}

Defect: $\partial P(n,6) = -1$, for $n\ge 0$.

\begin{fleqn}
\begin{align*}
P(0,6):\ & \bigl(R_{0}^{3}(3),P(0,1)\bigr),\ \bigl(R_{0}^{4}(2),P(0,4)\bigr),\ \bigl(R_{0}^{1}(1),P(0,5)\bigr)\\
P(1,6):\ & \bigl(R_{\infty}^{2}(2),P(0,3)\bigr),\ \bigl(R_{\infty}^{3}(1),P(0,7)\bigr),\ \bigl(R_{0}^{4}(3),P(1,1)\bigr),\ \bigl(R_{0}^{1}(2),P(1,4)\bigr),\ \bigl(R_{0}^{2}(1),P(1,5)\bigr)\\
 & \bigl(I(1,5),2P(0,2)\bigr)\\
P(n,6):\ & \bigl(R_{\infty}^{n\bmod 3+1}(2),P(n-1,3)\bigr),\ \bigl(R_{\infty}^{(n+1)\bmod 3+1}(1),P(n-1,7)\bigr),\ \bigl(R_{0}^{(n+2)\bmod 4+1}(3),P(n,1)\bigr)\\
 & \bigl(R_{0}^{(n-1)\bmod 4+1}(2),P(n,4)\bigr),\ \bigl(R_{0}^{n\bmod 4+1}(1),P(n,5)\bigr),\ \bigl(uI,(u+1)P\bigr),\ n>1\\
\end{align*}
\end{fleqn}
\subsubsection*{Modules of the form $P(n,7)$}

Defect: $\partial P(n,7) = -1$, for $n\ge 0$.

\begin{fleqn}
\begin{align*}
P(0,7):\ & \bigl(R_{\infty}^{1}(2),P(0,2)\bigr),\ \bigl(R_{\infty}^{2}(1),P(0,3)\bigr),\ \bigl(R_{0}^{4}(3),P(0,4)\bigr),\ \bigl(R_{0}^{1}(2),P(0,5)\bigr),\ \bigl(R_{0}^{2}(1),P(0,6)\bigr)\\
 & \bigl(I(2,6),2P(0,1)\bigr)\\
P(n,7):\ & \bigl(R_{\infty}^{n\bmod 3+1}(2),P(n,2)\bigr),\ \bigl(R_{\infty}^{(n+1)\bmod 3+1}(1),P(n,3)\bigr),\ \bigl(R_{0}^{(n+3)\bmod 4+1}(3),P(n,4)\bigr)\\
 & \bigl(R_{0}^{n\bmod 4+1}(2),P(n,5)\bigr),\ \bigl(R_{0}^{(n+1)\bmod 4+1}(1),P(n,6)\bigr),\ \bigl(uI,(u+1)P\bigr),\ n>0\\
\end{align*}
\end{fleqn}
\subsubsection{Schofield pairs associated to preinjective exceptional modules}

\begin{figure}[ht]


\begin{center}
\begin{scaletikzpicturetowidth}{\textwidth}

\end{scaletikzpicturetowidth}
\end{center}

\end{figure}
\medskip{}
\subsubsection*{Modules of the form $I(n,1)$}

Defect: $\partial I(n,1) = 1$, for $n\ge 0$.

\begin{fleqn}
\begin{align*}
I(0,1):\ & \bigl(I(0,2),R_{\infty}^{2}(1)\bigr),\ \bigl(I(0,3),R_{\infty}^{2}(2)\bigr),\ \bigl(I(0,4),R_{0}^{2}(1)\bigr),\ \bigl(I(0,5),R_{0}^{2}(2)\bigr),\ \bigl(I(0,6),R_{0}^{2}(3)\bigr)\\
 & \bigl(2I(0,7),P(2,4)\bigr)\\
I(n,1):\ & \bigl(I(n,2),R_{\infty}^{(-n+1)\bmod 3+1}(1)\bigr),\ \bigl(I(n,3),R_{\infty}^{(-n+1)\bmod 3+1}(2)\bigr),\ \bigl(I(n,4),R_{0}^{(-n+1)\bmod 4+1}(1)\bigr)\\
 & \bigl(I(n,5),R_{0}^{(-n+1)\bmod 4+1}(2)\bigr),\ \bigl(I(n,6),R_{0}^{(-n+1)\bmod 4+1}(3)\bigr),\ \bigl((v+1)I,vP\bigr),\ n>0\\
\end{align*}
\end{fleqn}
\subsubsection*{Modules of the form $I(n,2)$}

Defect: $\partial I(n,2) = 1$, for $n\ge 0$.

\begin{fleqn}
\begin{align*}
I(0,2):\ & \bigl(I(0,3),R_{\infty}^{3}(1)\bigr),\ \bigl(I(0,7),R_{\infty}^{3}(2)\bigr)\\
I(1,2):\ & \bigl(I(0,1),R_{0}^{1}(1)\bigr),\ \bigl(I(0,4),R_{0}^{1}(2)\bigr),\ \bigl(I(0,5),R_{0}^{1}(3)\bigr),\ \bigl(I(1,3),R_{\infty}^{2}(1)\bigr),\ \bigl(I(1,7),R_{\infty}^{2}(2)\bigr)\\
 & \bigl(2I(0,6),P(2,1)\bigr)\\
I(n,2):\ & \bigl(I(n-1,1),R_{0}^{(-n+1)\bmod 4+1}(1)\bigr),\ \bigl(I(n-1,4),R_{0}^{(-n+1)\bmod 4+1}(2)\bigr),\ \bigl(I(n-1,5),R_{0}^{(-n+1)\bmod 4+1}(3)\bigr)\\
 & \bigl(I(n,3),R_{\infty}^{(-n+2)\bmod 3+1}(1)\bigr),\ \bigl(I(n,7),R_{\infty}^{(-n+2)\bmod 3+1}(2)\bigr),\ \bigl((v+1)I,vP\bigr),\ n>1\\
\end{align*}
\end{fleqn}
\subsubsection*{Modules of the form $I(n,3)$}

Defect: $\partial I(n,3) = 1$, for $n\ge 0$.

\begin{fleqn}
\begin{align*}
I(0,3):\ & \bigl(I(0,7),R_{\infty}^{1}(1)\bigr)\\
I(1,3):\ & \bigl(I(0,2),R_{0}^{1}(1)\bigr),\ \bigl(I(0,6),R_{\infty}^{3}(2)\bigr),\ \bigl(I(1,7),R_{\infty}^{3}(1)\bigr)\\
I(2,3):\ & \bigl(I(0,1),R_{0}^{4}(2)\bigr),\ \bigl(I(0,4),R_{0}^{4}(3)\bigr),\ \bigl(I(1,2),R_{0}^{4}(1)\bigr),\ \bigl(I(1,6),R_{\infty}^{2}(2)\bigr),\ \bigl(I(2,7),R_{\infty}^{2}(1)\bigr)\\
 & \bigl(2I(0,5),P(1,2)\bigr)\\
I(n,3):\ & \bigl(I(n-2,1),R_{0}^{(-n+5)\bmod 4+1}(2)\bigr),\ \bigl(I(n-2,4),R_{0}^{(-n+5)\bmod 4+1}(3)\bigr),\ \bigl(I(n-1,2),R_{0}^{(-n+5)\bmod 4+1}(1)\bigr)\\
 & \bigl(I(n-1,6),R_{\infty}^{(-n+3)\bmod 3+1}(2)\bigr),\ \bigl(I(n,7),R_{\infty}^{(-n+3)\bmod 3+1}(1)\bigr),\ \bigl((v+1)I,vP\bigr),\ n>2\\
\end{align*}
\end{fleqn}
\subsubsection*{Modules of the form $I(n,4)$}

Defect: $\partial I(n,4) = 1$, for $n\ge 0$.

\begin{fleqn}
\begin{align*}
I(0,4):\ & \bigl(I(0,5),R_{0}^{3}(1)\bigr),\ \bigl(I(0,6),R_{0}^{3}(2)\bigr),\ \bigl(I(0,7),R_{0}^{3}(3)\bigr)\\
I(1,4):\ & \bigl(I(0,1),R_{\infty}^{1}(1)\bigr),\ \bigl(I(0,2),R_{\infty}^{1}(2)\bigr),\ \bigl(I(1,5),R_{0}^{2}(1)\bigr),\ \bigl(I(1,6),R_{0}^{2}(2)\bigr),\ \bigl(I(1,7),R_{0}^{2}(3)\bigr)\\
 & \bigl(2I(0,3),P(1,5)\bigr)\\
I(n,4):\ & \bigl(I(n-1,1),R_{\infty}^{(-n+1)\bmod 3+1}(1)\bigr),\ \bigl(I(n-1,2),R_{\infty}^{(-n+1)\bmod 3+1}(2)\bigr),\ \bigl(I(n,5),R_{0}^{(-n+2)\bmod 4+1}(1)\bigr)\\
 & \bigl(I(n,6),R_{0}^{(-n+2)\bmod 4+1}(2)\bigr),\ \bigl(I(n,7),R_{0}^{(-n+2)\bmod 4+1}(3)\bigr),\ \bigl((v+1)I,vP\bigr),\ n>1\\
\end{align*}
\end{fleqn}
\subsubsection*{Modules of the form $I(n,5)$}

Defect: $\partial I(n,5) = 1$, for $n\ge 0$.

\begin{fleqn}
\begin{align*}
I(0,5):\ & \bigl(I(0,6),R_{0}^{4}(1)\bigr),\ \bigl(I(0,7),R_{0}^{4}(2)\bigr)\\
I(1,5):\ & \bigl(I(0,3),R_{0}^{3}(3)\bigr),\ \bigl(I(0,4),R_{\infty}^{1}(1)\bigr),\ \bigl(I(1,6),R_{0}^{3}(1)\bigr),\ \bigl(I(1,7),R_{0}^{3}(2)\bigr)\\
I(2,5):\ & \bigl(I(0,1),R_{\infty}^{3}(2)\bigr),\ \bigl(I(1,3),R_{0}^{2}(3)\bigr),\ \bigl(I(1,4),R_{\infty}^{3}(1)\bigr),\ \bigl(I(2,6),R_{0}^{2}(1)\bigr),\ \bigl(I(2,7),R_{0}^{2}(2)\bigr)\\
 & \bigl(2I(0,2),P(0,6)\bigr)\\
I(n,5):\ & \bigl(I(n-2,1),R_{\infty}^{(-n+4)\bmod 3+1}(2)\bigr),\ \bigl(I(n-1,3),R_{0}^{(-n+3)\bmod 4+1}(3)\bigr),\ \bigl(I(n-1,4),R_{\infty}^{(-n+4)\bmod 3+1}(1)\bigr)\\
 & \bigl(I(n,6),R_{0}^{(-n+3)\bmod 4+1}(1)\bigr),\ \bigl(I(n,7),R_{0}^{(-n+3)\bmod 4+1}(2)\bigr),\ \bigl((v+1)I,vP\bigr),\ n>2\\
\end{align*}
\end{fleqn}
\subsubsection*{Modules of the form $I(n,6)$}

Defect: $\partial I(n,6) = 1$, for $n\ge 0$.

\begin{fleqn}
\begin{align*}
I(0,6):\ & \bigl(I(0,7),R_{0}^{1}(1)\bigr)\\
I(1,6):\ & \bigl(I(0,3),R_{0}^{4}(2)\bigr),\ \bigl(I(0,5),R_{\infty}^{1}(1)\bigr),\ \bigl(I(1,7),R_{0}^{4}(1)\bigr)\\
I(2,6):\ & \bigl(I(0,2),R_{0}^{3}(3)\bigr),\ \bigl(I(0,4),R_{\infty}^{3}(2)\bigr),\ \bigl(I(1,3),R_{0}^{3}(2)\bigr),\ \bigl(I(1,5),R_{\infty}^{3}(1)\bigr),\ \bigl(I(2,7),R_{0}^{3}(1)\bigr)\\
I(3,6):\ & \bigl(I(1,2),R_{0}^{2}(3)\bigr),\ \bigl(I(1,4),R_{\infty}^{2}(2)\bigr),\ \bigl(I(2,3),R_{0}^{2}(2)\bigr),\ \bigl(I(2,5),R_{\infty}^{2}(1)\bigr),\ \bigl(I(3,7),R_{0}^{2}(1)\bigr)\\
 & \bigl(3I(0,7),2P(2,4)\bigr)\\
I(n,6):\ & \bigl(I(n-2,2),R_{0}^{(-n+4)\bmod 4+1}(3)\bigr),\ \bigl(I(n-2,4),R_{\infty}^{(-n+4)\bmod 3+1}(2)\bigr),\ \bigl(I(n-1,3),R_{0}^{(-n+4)\bmod 4+1}(2)\bigr)\\
 & \bigl(I(n-1,5),R_{\infty}^{(-n+4)\bmod 3+1}(1)\bigr),\ \bigl(I(n,7),R_{0}^{(-n+4)\bmod 4+1}(1)\bigr),\ \bigl((v+1)I,vP\bigr),\ n>3\\
\end{align*}
\end{fleqn}
\subsubsection*{Modules of the form $I(n,7)$}

Defect: $\partial I(n,7) = 1$, for $n\ge 0$.

\begin{fleqn}
\begin{align*}
I(0,7):\ & - \\
I(1,7):\ & \bigl(I(0,3),R_{0}^{1}(1)\bigr),\ \bigl(I(0,6),R_{\infty}^{1}(1)\bigr)\\
I(2,7):\ & \bigl(I(0,2),R_{0}^{4}(2)\bigr),\ \bigl(I(0,5),R_{\infty}^{3}(2)\bigr),\ \bigl(I(1,3),R_{0}^{4}(1)\bigr),\ \bigl(I(1,6),R_{\infty}^{3}(1)\bigr)\\
I(3,7):\ & \bigl(I(0,1),R_{0}^{3}(3)\bigr),\ \bigl(I(1,2),R_{0}^{3}(2)\bigr),\ \bigl(I(1,5),R_{\infty}^{2}(2)\bigr),\ \bigl(I(2,3),R_{0}^{3}(1)\bigr),\ \bigl(I(2,6),R_{\infty}^{2}(1)\bigr)\\
 & \bigl(2I(0,4),P(0,3)\bigr)\\
I(n,7):\ & \bigl(I(n-3,1),R_{0}^{(-n+5)\bmod 4+1}(3)\bigr),\ \bigl(I(n-2,2),R_{0}^{(-n+5)\bmod 4+1}(2)\bigr),\ \bigl(I(n-2,5),R_{\infty}^{(-n+4)\bmod 3+1}(2)\bigr)\\
 & \bigl(I(n-1,3),R_{0}^{(-n+5)\bmod 4+1}(1)\bigr),\ \bigl(I(n-1,6),R_{\infty}^{(-n+4)\bmod 3+1}(1)\bigr),\ \bigl((v+1)I,vP\bigr),\ n>3\\
\end{align*}
\end{fleqn}
\subsubsection{Schofield pairs associated to regular exceptional modules}

\subsubsection*{The non-homogeneous tube $\mathcal{T}_{0}^{\Delta(\widetilde{\mathbb{A}}_{3,4})}$}

\begin{figure}[ht]



\begin{center}
\captionof{figure}{\vspace*{-10pt}$\mathcal{T}_{0}^{\Delta(\widetilde{\mathbb{A}}_{3,4})}$}
\begin{scaletikzpicturetowidth}{1.\textwidth}

\end{scaletikzpicturetowidth}
\end{center}

\end{figure}
\begin{fleqn}
\begin{align*}
R_{0}^{1}(1):\ & - \\
R_{0}^{1}(2):\ & \bigl(R_{0}^{2}(1),R_{0}^{1}(1)\bigr),\ \bigl(I(1,3),P(0,1)\bigr),\ \bigl(I(1,7),P(0,2)\bigr),\ \bigl(I(0,6),P(0,3)\bigr)\\
R_{0}^{2}(1):\ & \bigl(I(0,2),P(0,1)\bigr),\ \bigl(I(0,3),P(0,2)\bigr),\ \bigl(I(0,7),P(0,3)\bigr)\\
R_{0}^{2}(2):\ & \bigl(R_{0}^{3}(1),R_{0}^{2}(1)\bigr),\ \bigl(I(0,2),P(0,4)\bigr),\ \bigl(I(0,3),P(1,1)\bigr),\ \bigl(I(0,7),P(1,2)\bigr)\\
R_{0}^{3}(1):\ & - \\
R_{0}^{3}(2):\ & \bigl(R_{0}^{4}(1),R_{0}^{3}(1)\bigr)\\
R_{0}^{4}(1):\ & - \\
R_{0}^{4}(2):\ & \bigl(R_{0}^{1}(1),R_{0}^{4}(1)\bigr)\\
R_{0}^{4}(3):\ & \bigl(R_{0}^{1}(2),R_{0}^{4}(1)\bigr),\ \bigl(R_{0}^{2}(1),R_{0}^{4}(2)\bigr),\ \bigl(I(2,7),P(0,1)\bigr),\ \bigl(I(1,6),P(0,2)\bigr),\ \bigl(I(0,5),P(0,3)\bigr)\\
R_{0}^{1}(3):\ & \bigl(R_{0}^{2}(2),R_{0}^{1}(1)\bigr),\ \bigl(R_{0}^{3}(1),R_{0}^{1}(2)\bigr),\ \bigl(I(1,3),P(0,4)\bigr),\ \bigl(I(1,7),P(1,1)\bigr),\ \bigl(I(0,6),P(1,2)\bigr)\\
R_{0}^{2}(3):\ & \bigl(R_{0}^{3}(2),R_{0}^{2}(1)\bigr),\ \bigl(R_{0}^{4}(1),R_{0}^{2}(2)\bigr),\ \bigl(I(0,2),P(0,5)\bigr),\ \bigl(I(0,3),P(1,4)\bigr),\ \bigl(I(0,7),P(2,1)\bigr)\\
R_{0}^{3}(3):\ & \bigl(R_{0}^{4}(2),R_{0}^{3}(1)\bigr),\ \bigl(R_{0}^{1}(1),R_{0}^{3}(2)\bigr)\\
\end{align*}
\end{fleqn}
\subsubsection*{The non-homogeneous tube $\mathcal{T}_{\infty}^{\Delta(\widetilde{\mathbb{A}}_{3,4})}$}

\begin{figure}[ht]



\begin{center}
\captionof{figure}{\vspace*{-10pt}$\mathcal{T}_{\infty}^{\Delta(\widetilde{\mathbb{A}}_{3,4})}$}
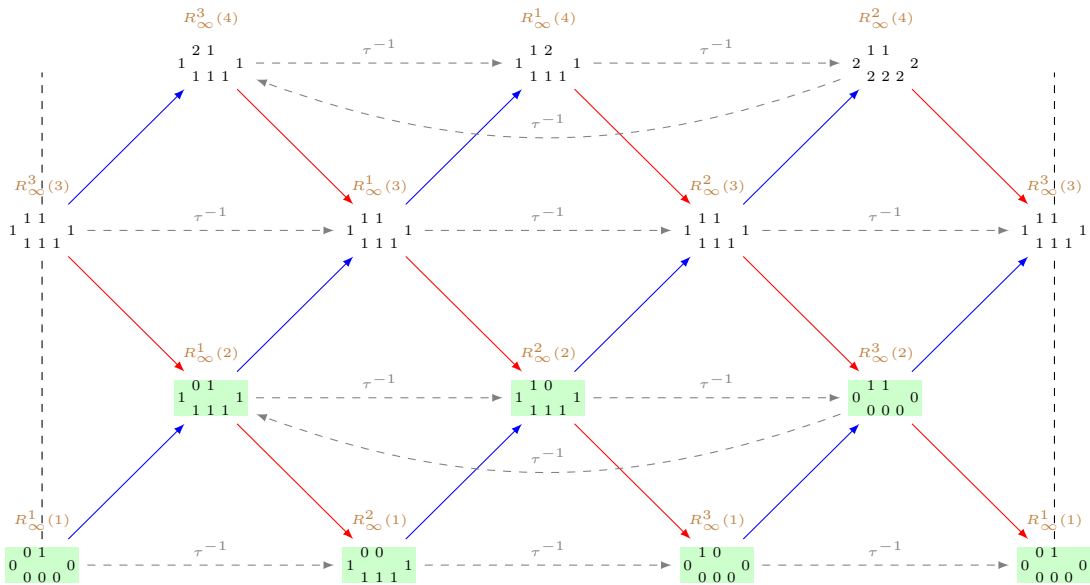
\begin{scaletikzpicturetowidth}{0.80000000000000004\textwidth}

\end{scaletikzpicturetowidth}
\end{center}

\end{figure}
\begin{fleqn}
\begin{align*}
R_{\infty}^{1}(1):\ & - \\
R_{\infty}^{1}(2):\ & \bigl(R_{\infty}^{2}(1),R_{\infty}^{1}(1)\bigr),\ \bigl(I(1,5),P(0,1)\bigr),\ \bigl(I(1,6),P(0,4)\bigr),\ \bigl(I(1,7),P(0,5)\bigr),\ \bigl(I(0,3),P(0,6)\bigr)\\
R_{\infty}^{2}(1):\ & \bigl(I(0,4),P(0,1)\bigr),\ \bigl(I(0,5),P(0,4)\bigr),\ \bigl(I(0,6),P(0,5)\bigr),\ \bigl(I(0,7),P(0,6)\bigr)\\
R_{\infty}^{2}(2):\ & \bigl(R_{\infty}^{3}(1),R_{\infty}^{2}(1)\bigr),\ \bigl(I(0,4),P(0,2)\bigr),\ \bigl(I(0,5),P(1,1)\bigr),\ \bigl(I(0,6),P(1,4)\bigr),\ \bigl(I(0,7),P(1,5)\bigr)\\
R_{\infty}^{3}(1):\ & - \\
R_{\infty}^{3}(2):\ & \bigl(R_{\infty}^{1}(1),R_{\infty}^{3}(1)\bigr)\\
\end{align*}
\end{fleqn}
\clearpage

\begin{lrbox}{\boxAAAlAAAAA}$\delta =\begin{smallmatrix}&1&1\\1&&&&&1\\&1&1&1&1\end{smallmatrix} $\end{lrbox}
\subsection[Schofield pairs for the quiver $\Delta(\A_{3,5})$]{Schofield pairs for the quiver $\Delta(\A_{3,5})$ -- {\usebox{\boxAAAlAAAAA}}}
\[
\xymatrix{  & 2\ar[dl] & 3\ar[l] &   &   &   &   &  \\
1 &   &   &   &   & 8\ar[ulll]\ar[dl] &   &  \\
  & 4\ar[ul] & 5\ar[l] & 6\ar[l] & 7\ar[l] &   &   &  }
\]
\medskip{}
\[
 C_{\Delta(\widetilde{\mathbb{A}}_{3,5})} = \begin{bmatrix}1 & 1 & 1 & 1 & 1 & 1 & 1 & 2\\
0 & 1 & 1 & 0 & 0 & 0 & 0 & 1\\
0 & 0 & 1 & 0 & 0 & 0 & 0 & 1\\
0 & 0 & 0 & 1 & 1 & 1 & 1 & 1\\
0 & 0 & 0 & 0 & 1 & 1 & 1 & 1\\
0 & 0 & 0 & 0 & 0 & 1 & 1 & 1\\
0 & 0 & 0 & 0 & 0 & 0 & 1 & 1\\
0 & 0 & 0 & 0 & 0 & 0 & 0 & 1\end{bmatrix}\qquad\Phi_{\Delta(\widetilde{\mathbb{A}}_{3,5})} = \begin{bmatrix}-1 & 1 & 0 & 1 & 0 & 0 & 0 & 0\\
-1 & 0 & 1 & 1 & 0 & 0 & 0 & 0\\
-1 & 0 & 0 & 1 & 0 & 0 & 0 & 1\\
-1 & 1 & 0 & 0 & 1 & 0 & 0 & 0\\
-1 & 1 & 0 & 0 & 0 & 1 & 0 & 0\\
-1 & 1 & 0 & 0 & 0 & 0 & 1 & 0\\
-1 & 1 & 0 & 0 & 0 & 0 & 0 & 1\\
-2 & 1 & 0 & 1 & 0 & 0 & 0 & 1\end{bmatrix}
\]

\medskip{}\subsubsection{Schofield pairs associated to preprojective exceptional modules}

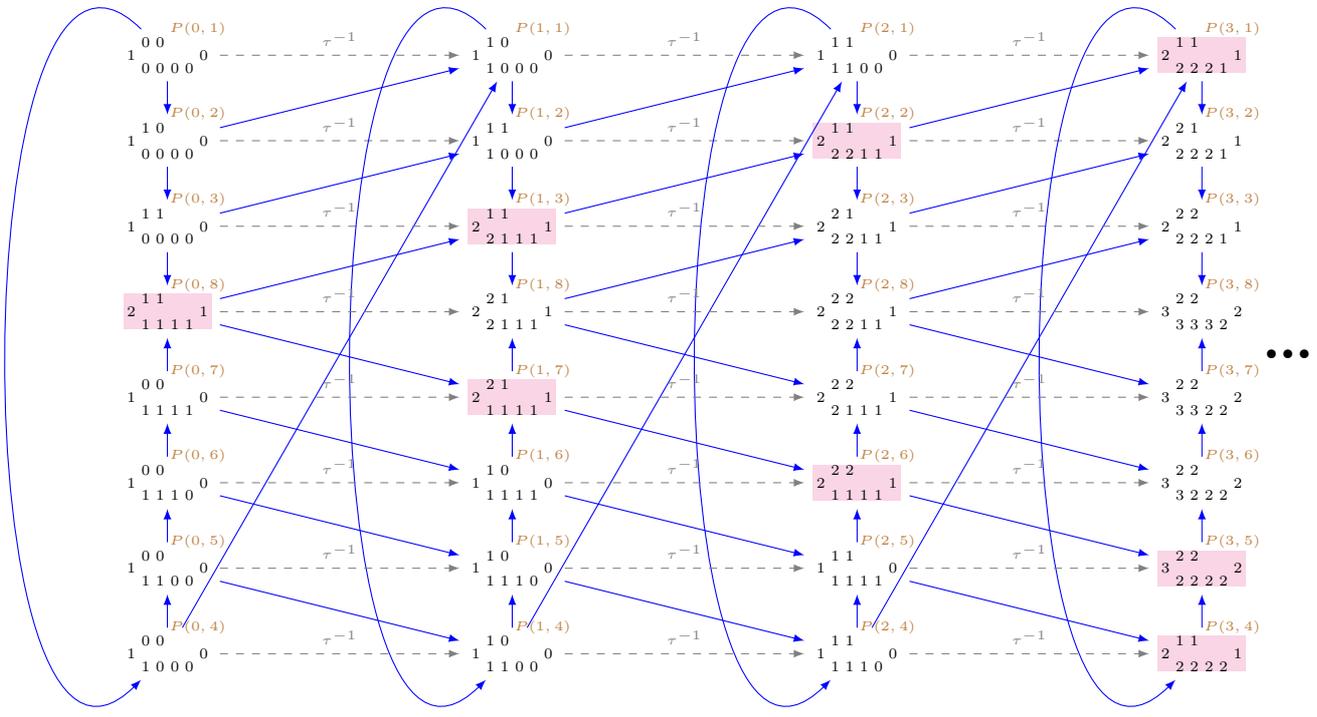
\begin{figure}[ht]


\begin{center}
\begin{scaletikzpicturetowidth}{\textwidth}

\end{scaletikzpicturetowidth}
\end{center}


\end{figure}
\medskip{}
\subsubsection*{Modules of the form $P(n,1)$}

Defect: $\partial P(n,1) = -1$, for $n\ge 0$.

\begin{fleqn}
\begin{align*}
P(0,1):\ & - \\
P(1,1):\ & \bigl(R_{0}^{3}(1),P(0,2)\bigr),\ \bigl(R_{\infty}^{3}(1),P(0,4)\bigr)\\
P(2,1):\ & \bigl(R_{0}^{3}(2),P(0,3)\bigr),\ \bigl(R_{\infty}^{3}(2),P(0,5)\bigr),\ \bigl(R_{0}^{4}(1),P(1,2)\bigr),\ \bigl(R_{\infty}^{1}(1),P(1,4)\bigr)\\
P(3,1):\ & \bigl(R_{0}^{2}(4),P(0,7)\bigr),\ \bigl(R_{0}^{3}(3),P(0,8)\bigr),\ \bigl(R_{0}^{4}(2),P(1,3)\bigr),\ \bigl(R_{\infty}^{1}(2),P(1,5)\bigr),\ \bigl(R_{0}^{5}(1),P(2,2)\bigr)\\
 & \bigl(R_{\infty}^{2}(1),P(2,4)\bigr),\ \bigl(I(1,3),2P(0,6)\bigr)\\
P(n,1):\ & \bigl(R_{0}^{(n-2)\bmod 5+1}(4),P(n-3,7)\bigr),\ \bigl(R_{0}^{(n-1)\bmod 5+1}(3),P(n-3,8)\bigr),\ \bigl(R_{0}^{n\bmod 5+1}(2),P(n-2,3)\bigr)\\
 & \bigl(R_{\infty}^{(n-3)\bmod 3+1}(2),P(n-2,5)\bigr),\ \bigl(R_{0}^{(n+1)\bmod 5+1}(1),P(n-1,2)\bigr),\ \bigl(R_{\infty}^{(n-2)\bmod 3+1}(1),P(n-1,4)\bigr)\\
 & \bigl(uI,(u+1)P\bigr),\ n>3\\
\end{align*}
\end{fleqn}
\subsubsection*{Modules of the form $P(n,2)$}

Defect: $\partial P(n,2) = -1$, for $n\ge 0$.

\begin{fleqn}
\begin{align*}
P(0,2):\ & \bigl(R_{\infty}^{3}(1),P(0,1)\bigr)\\
P(1,2):\ & \bigl(R_{0}^{3}(1),P(0,3)\bigr),\ \bigl(R_{\infty}^{3}(2),P(0,4)\bigr),\ \bigl(R_{\infty}^{1}(1),P(1,1)\bigr)\\
P(2,2):\ & \bigl(R_{0}^{1}(4),P(0,6)\bigr),\ \bigl(R_{0}^{2}(3),P(0,7)\bigr),\ \bigl(R_{0}^{3}(2),P(0,8)\bigr),\ \bigl(R_{0}^{4}(1),P(1,3)\bigr),\ \bigl(R_{\infty}^{1}(2),P(1,4)\bigr)\\
 & \bigl(R_{\infty}^{2}(1),P(2,1)\bigr),\ \bigl(I(2,8),2P(0,5)\bigr)\\
P(n,2):\ & \bigl(R_{0}^{(n-2)\bmod 5+1}(4),P(n-2,6)\bigr),\ \bigl(R_{0}^{(n-1)\bmod 5+1}(3),P(n-2,7)\bigr),\ \bigl(R_{0}^{n\bmod 5+1}(2),P(n-2,8)\bigr)\\
 & \bigl(R_{0}^{(n+1)\bmod 5+1}(1),P(n-1,3)\bigr),\ \bigl(R_{\infty}^{(n-2)\bmod 3+1}(2),P(n-1,4)\bigr),\ \bigl(R_{\infty}^{(n-1)\bmod 3+1}(1),P(n,1)\bigr)\\
 & \bigl(uI,(u+1)P\bigr),\ n>2\\
\end{align*}
\end{fleqn}
\subsubsection*{Modules of the form $P(n,3)$}

Defect: $\partial P(n,3) = -1$, for $n\ge 0$.

\begin{fleqn}
\begin{align*}
P(0,3):\ & \bigl(R_{\infty}^{3}(2),P(0,1)\bigr),\ \bigl(R_{\infty}^{1}(1),P(0,2)\bigr)\\
P(1,3):\ & \bigl(R_{0}^{5}(4),P(0,5)\bigr),\ \bigl(R_{0}^{1}(3),P(0,6)\bigr),\ \bigl(R_{0}^{2}(2),P(0,7)\bigr),\ \bigl(R_{0}^{3}(1),P(0,8)\bigr),\ \bigl(R_{\infty}^{1}(2),P(1,1)\bigr)\\
 & \bigl(R_{\infty}^{2}(1),P(1,2)\bigr),\ \bigl(I(2,7),2P(0,4)\bigr)\\
P(n,3):\ & \bigl(R_{0}^{(n+3)\bmod 5+1}(4),P(n-1,5)\bigr),\ \bigl(R_{0}^{(n-1)\bmod 5+1}(3),P(n-1,6)\bigr),\ \bigl(R_{0}^{n\bmod 5+1}(2),P(n-1,7)\bigr)\\
 & \bigl(R_{0}^{(n+1)\bmod 5+1}(1),P(n-1,8)\bigr),\ \bigl(R_{\infty}^{(n-1)\bmod 3+1}(2),P(n,1)\bigr),\ \bigl(R_{\infty}^{n\bmod 3+1}(1),P(n,2)\bigr)\\
 & \bigl(uI,(u+1)P\bigr),\ n>1\\
\end{align*}
\end{fleqn}
\subsubsection*{Modules of the form $P(n,4)$}

Defect: $\partial P(n,4) = -1$, for $n\ge 0$.

\begin{fleqn}
\begin{align*}
P(0,4):\ & \bigl(R_{0}^{3}(1),P(0,1)\bigr)\\
P(1,4):\ & \bigl(R_{0}^{3}(2),P(0,2)\bigr),\ \bigl(R_{\infty}^{3}(1),P(0,5)\bigr),\ \bigl(R_{0}^{4}(1),P(1,1)\bigr)\\
P(2,4):\ & \bigl(R_{0}^{3}(3),P(0,3)\bigr),\ \bigl(R_{\infty}^{3}(2),P(0,6)\bigr),\ \bigl(R_{0}^{4}(2),P(1,2)\bigr),\ \bigl(R_{\infty}^{1}(1),P(1,5)\bigr),\ \bigl(R_{0}^{5}(1),P(2,1)\bigr)\\
P(3,4):\ & \bigl(R_{0}^{3}(4),P(0,8)\bigr),\ \bigl(R_{0}^{4}(3),P(1,3)\bigr),\ \bigl(R_{\infty}^{1}(2),P(1,6)\bigr),\ \bigl(R_{0}^{5}(2),P(2,2)\bigr),\ \bigl(R_{\infty}^{2}(1),P(2,5)\bigr)\\
 & \bigl(R_{0}^{1}(1),P(3,1)\bigr),\ \bigl(I(0,2),2P(0,7)\bigr)\\
P(n,4):\ & \bigl(R_{0}^{(n-1)\bmod 5+1}(4),P(n-3,8)\bigr),\ \bigl(R_{0}^{n\bmod 5+1}(3),P(n-2,3)\bigr),\ \bigl(R_{\infty}^{(n-3)\bmod 3+1}(2),P(n-2,6)\bigr)\\
 & \bigl(R_{0}^{(n+1)\bmod 5+1}(2),P(n-1,2)\bigr),\ \bigl(R_{\infty}^{(n-2)\bmod 3+1}(1),P(n-1,5)\bigr),\ \bigl(R_{0}^{(n-3)\bmod 5+1}(1),P(n,1)\bigr)\\
 & \bigl(uI,(u+1)P\bigr),\ n>3\\
\end{align*}
\end{fleqn}
\subsubsection*{Modules of the form $P(n,5)$}

Defect: $\partial P(n,5) = -1$, for $n\ge 0$.

\begin{fleqn}
\begin{align*}
P(0,5):\ & \bigl(R_{0}^{3}(2),P(0,1)\bigr),\ \bigl(R_{0}^{4}(1),P(0,4)\bigr)\\
P(1,5):\ & \bigl(R_{0}^{3}(3),P(0,2)\bigr),\ \bigl(R_{\infty}^{3}(1),P(0,6)\bigr),\ \bigl(R_{0}^{4}(2),P(1,1)\bigr),\ \bigl(R_{0}^{5}(1),P(1,4)\bigr)\\
P(2,5):\ & \bigl(R_{0}^{3}(4),P(0,3)\bigr),\ \bigl(R_{\infty}^{3}(2),P(0,7)\bigr),\ \bigl(R_{0}^{4}(3),P(1,2)\bigr),\ \bigl(R_{\infty}^{1}(1),P(1,6)\bigr),\ \bigl(R_{0}^{5}(2),P(2,1)\bigr)\\
 & \bigl(R_{0}^{1}(1),P(2,4)\bigr)\\
P(3,5):\ & \bigl(R_{0}^{4}(4),P(1,3)\bigr),\ \bigl(R_{\infty}^{1}(2),P(1,7)\bigr),\ \bigl(R_{0}^{5}(3),P(2,2)\bigr),\ \bigl(R_{\infty}^{2}(1),P(2,6)\bigr),\ \bigl(R_{0}^{1}(2),P(3,1)\bigr)\\
 & \bigl(R_{0}^{2}(1),P(3,4)\bigr),\ \bigl(2I(2,6),3P(0,1)\bigr)\\
P(n,5):\ & \bigl(R_{0}^{n\bmod 5+1}(4),P(n-2,3)\bigr),\ \bigl(R_{\infty}^{(n-3)\bmod 3+1}(2),P(n-2,7)\bigr),\ \bigl(R_{0}^{(n+1)\bmod 5+1}(3),P(n-1,2)\bigr)\\
 & \bigl(R_{\infty}^{(n-2)\bmod 3+1}(1),P(n-1,6)\bigr),\ \bigl(R_{0}^{(n-3)\bmod 5+1}(2),P(n,1)\bigr),\ \bigl(R_{0}^{(n-2)\bmod 5+1}(1),P(n,4)\bigr)\\
 & \bigl(uI,(u+1)P\bigr),\ n>3\\
\end{align*}
\end{fleqn}
\subsubsection*{Modules of the form $P(n,6)$}

Defect: $\partial P(n,6) = -1$, for $n\ge 0$.

\begin{fleqn}
\begin{align*}
P(0,6):\ & \bigl(R_{0}^{3}(3),P(0,1)\bigr),\ \bigl(R_{0}^{4}(2),P(0,4)\bigr),\ \bigl(R_{0}^{5}(1),P(0,5)\bigr)\\
P(1,6):\ & \bigl(R_{0}^{3}(4),P(0,2)\bigr),\ \bigl(R_{\infty}^{3}(1),P(0,7)\bigr),\ \bigl(R_{0}^{4}(3),P(1,1)\bigr),\ \bigl(R_{0}^{5}(2),P(1,4)\bigr),\ \bigl(R_{0}^{1}(1),P(1,5)\bigr)\\
P(2,6):\ & \bigl(R_{\infty}^{3}(2),P(0,8)\bigr),\ \bigl(R_{0}^{4}(4),P(1,2)\bigr),\ \bigl(R_{\infty}^{1}(1),P(1,7)\bigr),\ \bigl(R_{0}^{5}(3),P(2,1)\bigr),\ \bigl(R_{0}^{1}(2),P(2,4)\bigr)\\
 & \bigl(R_{0}^{2}(1),P(2,5)\bigr),\ \bigl(I(0,4),2P(0,3)\bigr)\\
P(n,6):\ & \bigl(R_{\infty}^{n\bmod 3+1}(2),P(n-2,8)\bigr),\ \bigl(R_{0}^{(n+1)\bmod 5+1}(4),P(n-1,2)\bigr),\ \bigl(R_{\infty}^{(n-2)\bmod 3+1}(1),P(n-1,7)\bigr)\\
 & \bigl(R_{0}^{(n+2)\bmod 5+1}(3),P(n,1)\bigr),\ \bigl(R_{0}^{(n-2)\bmod 5+1}(2),P(n,4)\bigr),\ \bigl(R_{0}^{(n-1)\bmod 5+1}(1),P(n,5)\bigr)\\
 & \bigl(uI,(u+1)P\bigr),\ n>2\\
\end{align*}
\end{fleqn}
\subsubsection*{Modules of the form $P(n,7)$}

Defect: $\partial P(n,7) = -1$, for $n\ge 0$.

\begin{fleqn}
\begin{align*}
P(0,7):\ & \bigl(R_{0}^{3}(4),P(0,1)\bigr),\ \bigl(R_{0}^{4}(3),P(0,4)\bigr),\ \bigl(R_{0}^{5}(2),P(0,5)\bigr),\ \bigl(R_{0}^{1}(1),P(0,6)\bigr)\\
P(1,7):\ & \bigl(R_{\infty}^{2}(2),P(0,3)\bigr),\ \bigl(R_{\infty}^{3}(1),P(0,8)\bigr),\ \bigl(R_{0}^{4}(4),P(1,1)\bigr),\ \bigl(R_{0}^{5}(3),P(1,4)\bigr),\ \bigl(R_{0}^{1}(2),P(1,5)\bigr)\\
 & \bigl(R_{0}^{2}(1),P(1,6)\bigr),\ \bigl(I(1,5),2P(0,2)\bigr)\\
P(n,7):\ & \bigl(R_{\infty}^{n\bmod 3+1}(2),P(n-1,3)\bigr),\ \bigl(R_{\infty}^{(n+1)\bmod 3+1}(1),P(n-1,8)\bigr),\ \bigl(R_{0}^{(n+2)\bmod 5+1}(4),P(n,1)\bigr)\\
 & \bigl(R_{0}^{(n+3)\bmod 5+1}(3),P(n,4)\bigr),\ \bigl(R_{0}^{(n-1)\bmod 5+1}(2),P(n,5)\bigr),\ \bigl(R_{0}^{n\bmod 5+1}(1),P(n,6)\bigr)\\
 & \bigl(uI,(u+1)P\bigr),\ n>1\\
\end{align*}
\end{fleqn}
\subsubsection*{Modules of the form $P(n,8)$}

Defect: $\partial P(n,8) = -1$, for $n\ge 0$.

\begin{fleqn}
\begin{align*}
P(0,8):\ & \bigl(R_{\infty}^{1}(2),P(0,2)\bigr),\ \bigl(R_{\infty}^{2}(1),P(0,3)\bigr),\ \bigl(R_{0}^{4}(4),P(0,4)\bigr),\ \bigl(R_{0}^{5}(3),P(0,5)\bigr),\ \bigl(R_{0}^{1}(2),P(0,6)\bigr)\\
 & \bigl(R_{0}^{2}(1),P(0,7)\bigr),\ \bigl(I(2,6),2P(0,1)\bigr)\\
P(n,8):\ & \bigl(R_{\infty}^{n\bmod 3+1}(2),P(n,2)\bigr),\ \bigl(R_{\infty}^{(n+1)\bmod 3+1}(1),P(n,3)\bigr),\ \bigl(R_{0}^{(n+3)\bmod 5+1}(4),P(n,4)\bigr)\\
 & \bigl(R_{0}^{(n+4)\bmod 5+1}(3),P(n,5)\bigr),\ \bigl(R_{0}^{n\bmod 5+1}(2),P(n,6)\bigr),\ \bigl(R_{0}^{(n+1)\bmod 5+1}(1),P(n,7)\bigr)\\
 & \bigl(uI,(u+1)P\bigr),\ n>0\\
\end{align*}
\end{fleqn}
\subsubsection{Schofield pairs associated to preinjective exceptional modules}

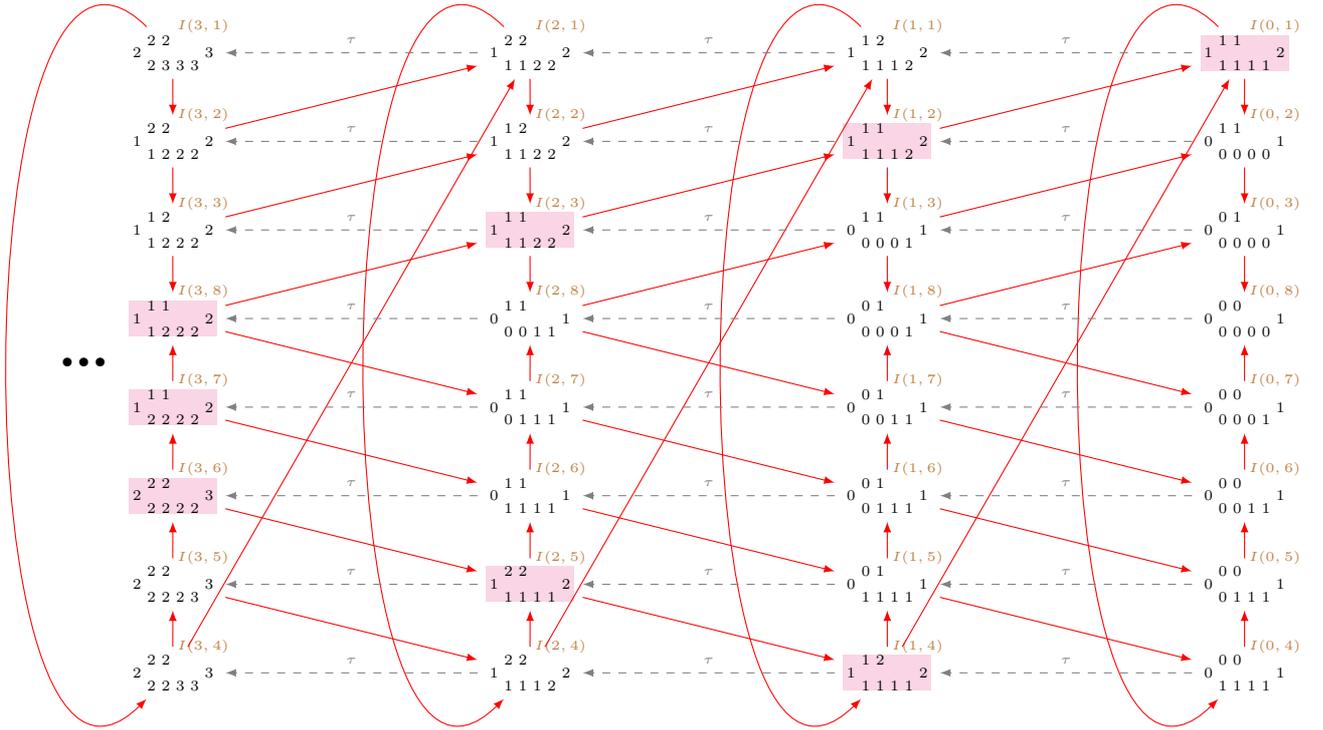
\begin{figure}[ht]


\begin{center}
\begin{scaletikzpicturetowidth}{\textwidth}

\end{scaletikzpicturetowidth}
\end{center}

\end{figure}
\medskip{}
\subsubsection*{Modules of the form $I(n,1)$}

Defect: $\partial I(n,1) = 1$, for $n\ge 0$.

\begin{fleqn}
\begin{align*}
I(0,1):\ & \bigl(I(0,2),R_{\infty}^{2}(1)\bigr),\ \bigl(I(0,3),R_{\infty}^{2}(2)\bigr),\ \bigl(I(0,4),R_{0}^{2}(1)\bigr),\ \bigl(I(0,5),R_{0}^{2}(2)\bigr),\ \bigl(I(0,6),R_{0}^{2}(3)\bigr)\\
 & \bigl(I(0,7),R_{0}^{2}(4)\bigr),\ \bigl(2I(0,8),P(2,5)\bigr)\\
I(n,1):\ & \bigl(I(n,2),R_{\infty}^{(-n+1)\bmod 3+1}(1)\bigr),\ \bigl(I(n,3),R_{\infty}^{(-n+1)\bmod 3+1}(2)\bigr),\ \bigl(I(n,4),R_{0}^{(-n+1)\bmod 5+1}(1)\bigr)\\
 & \bigl(I(n,5),R_{0}^{(-n+1)\bmod 5+1}(2)\bigr),\ \bigl(I(n,6),R_{0}^{(-n+1)\bmod 5+1}(3)\bigr),\ \bigl(I(n,7),R_{0}^{(-n+1)\bmod 5+1}(4)\bigr)\\
 & \bigl((v+1)I,vP\bigr),\ n>0\\
\end{align*}
\end{fleqn}
\subsubsection*{Modules of the form $I(n,2)$}

Defect: $\partial I(n,2) = 1$, for $n\ge 0$.

\begin{fleqn}
\begin{align*}
I(0,2):\ & \bigl(I(0,3),R_{\infty}^{3}(1)\bigr),\ \bigl(I(0,8),R_{\infty}^{3}(2)\bigr)\\
I(1,2):\ & \bigl(I(0,1),R_{0}^{1}(1)\bigr),\ \bigl(I(0,4),R_{0}^{1}(2)\bigr),\ \bigl(I(0,5),R_{0}^{1}(3)\bigr),\ \bigl(I(0,6),R_{0}^{1}(4)\bigr),\ \bigl(I(1,3),R_{\infty}^{2}(1)\bigr)\\
 & \bigl(I(1,8),R_{\infty}^{2}(2)\bigr),\ \bigl(2I(0,7),P(2,4)\bigr)\\
I(n,2):\ & \bigl(I(n-1,1),R_{0}^{(-n+1)\bmod 5+1}(1)\bigr),\ \bigl(I(n-1,4),R_{0}^{(-n+1)\bmod 5+1}(2)\bigr),\ \bigl(I(n-1,5),R_{0}^{(-n+1)\bmod 5+1}(3)\bigr)\\
 & \bigl(I(n-1,6),R_{0}^{(-n+1)\bmod 5+1}(4)\bigr),\ \bigl(I(n,3),R_{\infty}^{(-n+2)\bmod 3+1}(1)\bigr),\ \bigl(I(n,8),R_{\infty}^{(-n+2)\bmod 3+1}(2)\bigr)\\
 & \bigl((v+1)I,vP\bigr),\ n>1\\
\end{align*}
\end{fleqn}
\subsubsection*{Modules of the form $I(n,3)$}

Defect: $\partial I(n,3) = 1$, for $n\ge 0$.

\begin{fleqn}
\begin{align*}
I(0,3):\ & \bigl(I(0,8),R_{\infty}^{1}(1)\bigr)\\
I(1,3):\ & \bigl(I(0,2),R_{0}^{1}(1)\bigr),\ \bigl(I(0,7),R_{\infty}^{3}(2)\bigr),\ \bigl(I(1,8),R_{\infty}^{3}(1)\bigr)\\
I(2,3):\ & \bigl(I(0,1),R_{0}^{5}(2)\bigr),\ \bigl(I(0,4),R_{0}^{5}(3)\bigr),\ \bigl(I(0,5),R_{0}^{5}(4)\bigr),\ \bigl(I(1,2),R_{0}^{5}(1)\bigr),\ \bigl(I(1,7),R_{\infty}^{2}(2)\bigr)\\
 & \bigl(I(2,8),R_{\infty}^{2}(1)\bigr),\ \bigl(2I(0,6),P(2,1)\bigr)\\
I(n,3):\ & \bigl(I(n-2,1),R_{0}^{(-n+6)\bmod 5+1}(2)\bigr),\ \bigl(I(n-2,4),R_{0}^{(-n+6)\bmod 5+1}(3)\bigr),\ \bigl(I(n-2,5),R_{0}^{(-n+6)\bmod 5+1}(4)\bigr)\\
 & \bigl(I(n-1,2),R_{0}^{(-n+6)\bmod 5+1}(1)\bigr),\ \bigl(I(n-1,7),R_{\infty}^{(-n+3)\bmod 3+1}(2)\bigr),\ \bigl(I(n,8),R_{\infty}^{(-n+3)\bmod 3+1}(1)\bigr)\\
 & \bigl((v+1)I,vP\bigr),\ n>2\\
\end{align*}
\end{fleqn}
\subsubsection*{Modules of the form $I(n,4)$}

Defect: $\partial I(n,4) = 1$, for $n\ge 0$.

\begin{fleqn}
\begin{align*}
I(0,4):\ & \bigl(I(0,5),R_{0}^{3}(1)\bigr),\ \bigl(I(0,6),R_{0}^{3}(2)\bigr),\ \bigl(I(0,7),R_{0}^{3}(3)\bigr),\ \bigl(I(0,8),R_{0}^{3}(4)\bigr)\\
I(1,4):\ & \bigl(I(0,1),R_{\infty}^{1}(1)\bigr),\ \bigl(I(0,2),R_{\infty}^{1}(2)\bigr),\ \bigl(I(1,5),R_{0}^{2}(1)\bigr),\ \bigl(I(1,6),R_{0}^{2}(2)\bigr),\ \bigl(I(1,7),R_{0}^{2}(3)\bigr)\\
 & \bigl(I(1,8),R_{0}^{2}(4)\bigr),\ \bigl(2I(0,3),P(1,6)\bigr)\\
I(n,4):\ & \bigl(I(n-1,1),R_{\infty}^{(-n+1)\bmod 3+1}(1)\bigr),\ \bigl(I(n-1,2),R_{\infty}^{(-n+1)\bmod 3+1}(2)\bigr),\ \bigl(I(n,5),R_{0}^{(-n+2)\bmod 5+1}(1)\bigr)\\
 & \bigl(I(n,6),R_{0}^{(-n+2)\bmod 5+1}(2)\bigr),\ \bigl(I(n,7),R_{0}^{(-n+2)\bmod 5+1}(3)\bigr),\ \bigl(I(n,8),R_{0}^{(-n+2)\bmod 5+1}(4)\bigr)\\
 & \bigl((v+1)I,vP\bigr),\ n>1\\
\end{align*}
\end{fleqn}
\subsubsection*{Modules of the form $I(n,5)$}

Defect: $\partial I(n,5) = 1$, for $n\ge 0$.

\begin{fleqn}
\begin{align*}
I(0,5):\ & \bigl(I(0,6),R_{0}^{4}(1)\bigr),\ \bigl(I(0,7),R_{0}^{4}(2)\bigr),\ \bigl(I(0,8),R_{0}^{4}(3)\bigr)\\
I(1,5):\ & \bigl(I(0,3),R_{0}^{3}(4)\bigr),\ \bigl(I(0,4),R_{\infty}^{1}(1)\bigr),\ \bigl(I(1,6),R_{0}^{3}(1)\bigr),\ \bigl(I(1,7),R_{0}^{3}(2)\bigr),\ \bigl(I(1,8),R_{0}^{3}(3)\bigr)\\
I(2,5):\ & \bigl(I(0,1),R_{\infty}^{3}(2)\bigr),\ \bigl(I(1,3),R_{0}^{2}(4)\bigr),\ \bigl(I(1,4),R_{\infty}^{3}(1)\bigr),\ \bigl(I(2,6),R_{0}^{2}(1)\bigr),\ \bigl(I(2,7),R_{0}^{2}(2)\bigr)\\
 & \bigl(I(2,8),R_{0}^{2}(3)\bigr),\ \bigl(2I(0,2),P(0,7)\bigr)\\
I(n,5):\ & \bigl(I(n-2,1),R_{\infty}^{(-n+4)\bmod 3+1}(2)\bigr),\ \bigl(I(n-1,3),R_{0}^{(-n+3)\bmod 5+1}(4)\bigr),\ \bigl(I(n-1,4),R_{\infty}^{(-n+4)\bmod 3+1}(1)\bigr)\\
 & \bigl(I(n,6),R_{0}^{(-n+3)\bmod 5+1}(1)\bigr),\ \bigl(I(n,7),R_{0}^{(-n+3)\bmod 5+1}(2)\bigr),\ \bigl(I(n,8),R_{0}^{(-n+3)\bmod 5+1}(3)\bigr)\\
 & \bigl((v+1)I,vP\bigr),\ n>2\\
\end{align*}
\end{fleqn}
\subsubsection*{Modules of the form $I(n,6)$}

Defect: $\partial I(n,6) = 1$, for $n\ge 0$.

\begin{fleqn}
\begin{align*}
I(0,6):\ & \bigl(I(0,7),R_{0}^{5}(1)\bigr),\ \bigl(I(0,8),R_{0}^{5}(2)\bigr)\\
I(1,6):\ & \bigl(I(0,3),R_{0}^{4}(3)\bigr),\ \bigl(I(0,5),R_{\infty}^{1}(1)\bigr),\ \bigl(I(1,7),R_{0}^{4}(1)\bigr),\ \bigl(I(1,8),R_{0}^{4}(2)\bigr)\\
I(2,6):\ & \bigl(I(0,2),R_{0}^{3}(4)\bigr),\ \bigl(I(0,4),R_{\infty}^{3}(2)\bigr),\ \bigl(I(1,3),R_{0}^{3}(3)\bigr),\ \bigl(I(1,5),R_{\infty}^{3}(1)\bigr),\ \bigl(I(2,7),R_{0}^{3}(1)\bigr)\\
 & \bigl(I(2,8),R_{0}^{3}(2)\bigr)\\
I(3,6):\ & \bigl(I(1,2),R_{0}^{2}(4)\bigr),\ \bigl(I(1,4),R_{\infty}^{2}(2)\bigr),\ \bigl(I(2,3),R_{0}^{2}(3)\bigr),\ \bigl(I(2,5),R_{\infty}^{2}(1)\bigr),\ \bigl(I(3,7),R_{0}^{2}(1)\bigr)\\
 & \bigl(I(3,8),R_{0}^{2}(2)\bigr),\ \bigl(3I(0,8),2P(2,5)\bigr)\\
I(n,6):\ & \bigl(I(n-2,2),R_{0}^{(-n+4)\bmod 5+1}(4)\bigr),\ \bigl(I(n-2,4),R_{\infty}^{(-n+4)\bmod 3+1}(2)\bigr),\ \bigl(I(n-1,3),R_{0}^{(-n+4)\bmod 5+1}(3)\bigr)\\
 & \bigl(I(n-1,5),R_{\infty}^{(-n+4)\bmod 3+1}(1)\bigr),\ \bigl(I(n,7),R_{0}^{(-n+4)\bmod 5+1}(1)\bigr),\ \bigl(I(n,8),R_{0}^{(-n+4)\bmod 5+1}(2)\bigr)\\
 & \bigl((v+1)I,vP\bigr),\ n>3\\
\end{align*}
\end{fleqn}
\subsubsection*{Modules of the form $I(n,7)$}

Defect: $\partial I(n,7) = 1$, for $n\ge 0$.

\begin{fleqn}
\begin{align*}
I(0,7):\ & \bigl(I(0,8),R_{0}^{1}(1)\bigr)\\
I(1,7):\ & \bigl(I(0,3),R_{0}^{5}(2)\bigr),\ \bigl(I(0,6),R_{\infty}^{1}(1)\bigr),\ \bigl(I(1,8),R_{0}^{5}(1)\bigr)\\
I(2,7):\ & \bigl(I(0,2),R_{0}^{4}(3)\bigr),\ \bigl(I(0,5),R_{\infty}^{3}(2)\bigr),\ \bigl(I(1,3),R_{0}^{4}(2)\bigr),\ \bigl(I(1,6),R_{\infty}^{3}(1)\bigr),\ \bigl(I(2,8),R_{0}^{4}(1)\bigr)\\
I(3,7):\ & \bigl(I(0,1),R_{0}^{3}(4)\bigr),\ \bigl(I(1,2),R_{0}^{3}(3)\bigr),\ \bigl(I(1,5),R_{\infty}^{2}(2)\bigr),\ \bigl(I(2,3),R_{0}^{3}(2)\bigr),\ \bigl(I(2,6),R_{\infty}^{2}(1)\bigr)\\
 & \bigl(I(3,8),R_{0}^{3}(1)\bigr),\ \bigl(2I(0,4),P(0,3)\bigr)\\
I(n,7):\ & \bigl(I(n-3,1),R_{0}^{(-n+5)\bmod 5+1}(4)\bigr),\ \bigl(I(n-2,2),R_{0}^{(-n+5)\bmod 5+1}(3)\bigr),\ \bigl(I(n-2,5),R_{\infty}^{(-n+4)\bmod 3+1}(2)\bigr)\\
 & \bigl(I(n-1,3),R_{0}^{(-n+5)\bmod 5+1}(2)\bigr),\ \bigl(I(n-1,6),R_{\infty}^{(-n+4)\bmod 3+1}(1)\bigr),\ \bigl(I(n,8),R_{0}^{(-n+5)\bmod 5+1}(1)\bigr)\\
 & \bigl((v+1)I,vP\bigr),\ n>3\\
\end{align*}
\end{fleqn}
\subsubsection*{Modules of the form $I(n,8)$}

Defect: $\partial I(n,8) = 1$, for $n\ge 0$.

\begin{fleqn}
\begin{align*}
I(0,8):\ & - \\
I(1,8):\ & \bigl(I(0,3),R_{0}^{1}(1)\bigr),\ \bigl(I(0,7),R_{\infty}^{1}(1)\bigr)\\
I(2,8):\ & \bigl(I(0,2),R_{0}^{5}(2)\bigr),\ \bigl(I(0,6),R_{\infty}^{3}(2)\bigr),\ \bigl(I(1,3),R_{0}^{5}(1)\bigr),\ \bigl(I(1,7),R_{\infty}^{3}(1)\bigr)\\
I(3,8):\ & \bigl(I(0,1),R_{0}^{4}(3)\bigr),\ \bigl(I(0,4),R_{0}^{4}(4)\bigr),\ \bigl(I(1,2),R_{0}^{4}(2)\bigr),\ \bigl(I(1,6),R_{\infty}^{2}(2)\bigr),\ \bigl(I(2,3),R_{0}^{4}(1)\bigr)\\
 & \bigl(I(2,7),R_{\infty}^{2}(1)\bigr),\ \bigl(2I(0,5),P(1,2)\bigr)\\
I(n,8):\ & \bigl(I(n-3,1),R_{0}^{(-n+6)\bmod 5+1}(3)\bigr),\ \bigl(I(n-3,4),R_{0}^{(-n+6)\bmod 5+1}(4)\bigr),\ \bigl(I(n-2,2),R_{0}^{(-n+6)\bmod 5+1}(2)\bigr)\\
 & \bigl(I(n-2,6),R_{\infty}^{(-n+4)\bmod 3+1}(2)\bigr),\ \bigl(I(n-1,3),R_{0}^{(-n+6)\bmod 5+1}(1)\bigr),\ \bigl(I(n-1,7),R_{\infty}^{(-n+4)\bmod 3+1}(1)\bigr)\\
 & \bigl((v+1)I,vP\bigr),\ n>3\\
\end{align*}
\end{fleqn}
\subsubsection{Schofield pairs associated to regular exceptional modules}

\subsubsection*{The non-homogeneous tube $\mathcal{T}_{0}^{\Delta(\widetilde{\mathbb{A}}_{3,5})}$}

\begin{figure}[ht]



\begin{center}
\captionof{figure}{\vspace*{-10pt}$\mathcal{T}_{0}^{\Delta(\widetilde{\mathbb{A}}_{3,5})}$}
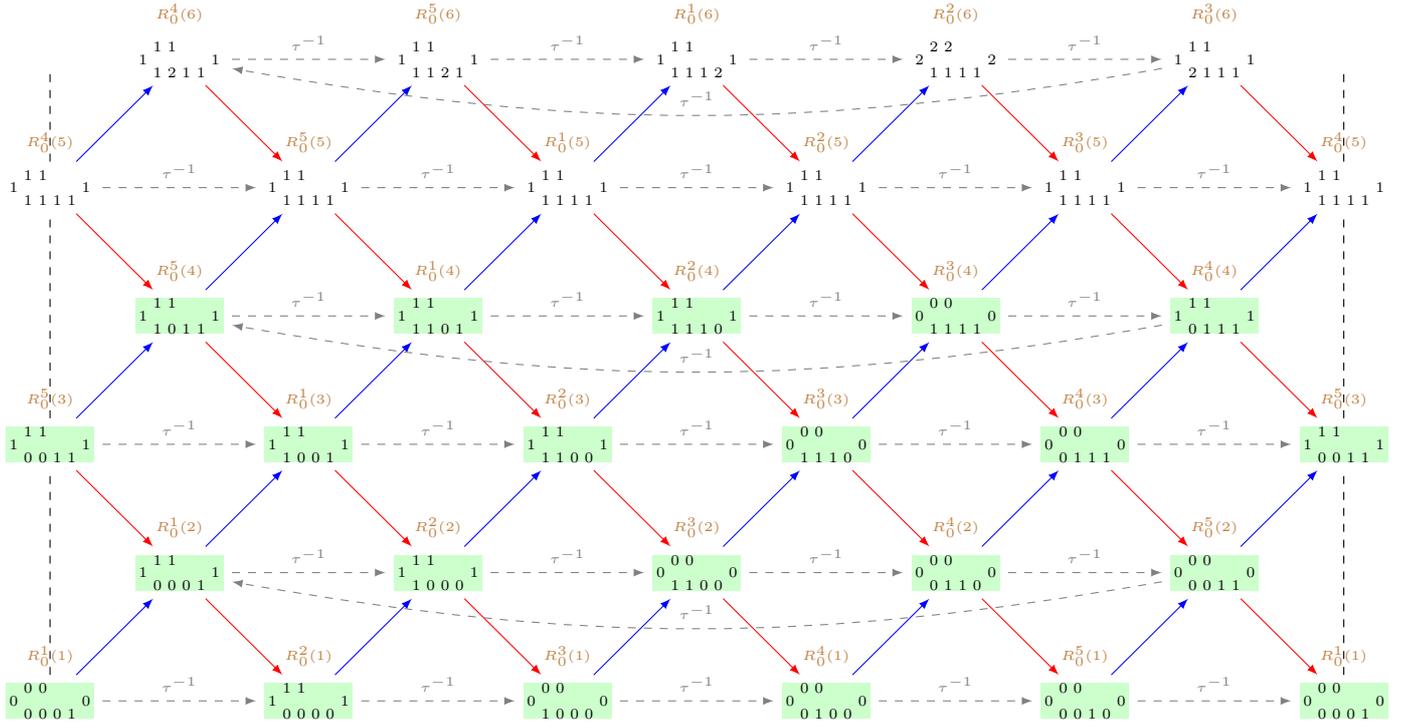
\begin{scaletikzpicturetowidth}{1.\textwidth}

\end{scaletikzpicturetowidth}
\end{center}

\end{figure}
\begin{fleqn}
\begin{align*}
R_{0}^{1}(1):\ & - \\
R_{0}^{1}(2):\ & \bigl(R_{0}^{2}(1),R_{0}^{1}(1)\bigr),\ \bigl(I(1,3),P(0,1)\bigr),\ \bigl(I(1,8),P(0,2)\bigr),\ \bigl(I(0,7),P(0,3)\bigr)\\
R_{0}^{2}(1):\ & \bigl(I(0,2),P(0,1)\bigr),\ \bigl(I(0,3),P(0,2)\bigr),\ \bigl(I(0,8),P(0,3)\bigr)\\
R_{0}^{2}(2):\ & \bigl(R_{0}^{3}(1),R_{0}^{2}(1)\bigr),\ \bigl(I(0,2),P(0,4)\bigr),\ \bigl(I(0,3),P(1,1)\bigr),\ \bigl(I(0,8),P(1,2)\bigr)\\
R_{0}^{3}(1):\ & - \\
R_{0}^{3}(2):\ & \bigl(R_{0}^{4}(1),R_{0}^{3}(1)\bigr)\\
R_{0}^{4}(1):\ & - \\
R_{0}^{4}(2):\ & \bigl(R_{0}^{5}(1),R_{0}^{4}(1)\bigr)\\
R_{0}^{5}(1):\ & - \\
R_{0}^{5}(2):\ & \bigl(R_{0}^{1}(1),R_{0}^{5}(1)\bigr)\\
R_{0}^{5}(3):\ & \bigl(R_{0}^{1}(2),R_{0}^{5}(1)\bigr),\ \bigl(R_{0}^{2}(1),R_{0}^{5}(2)\bigr),\ \bigl(I(2,8),P(0,1)\bigr),\ \bigl(I(1,7),P(0,2)\bigr),\ \bigl(I(0,6),P(0,3)\bigr)\\
R_{0}^{5}(4):\ & \bigl(R_{0}^{1}(3),R_{0}^{5}(1)\bigr),\ \bigl(R_{0}^{2}(2),R_{0}^{5}(2)\bigr),\ \bigl(R_{0}^{3}(1),R_{0}^{5}(3)\bigr),\ \bigl(I(2,8),P(0,4)\bigr),\ \bigl(I(1,7),P(1,1)\bigr)\\
 & \bigl(I(0,6),P(1,2)\bigr)\\
R_{0}^{1}(3):\ & \bigl(R_{0}^{2}(2),R_{0}^{1}(1)\bigr),\ \bigl(R_{0}^{3}(1),R_{0}^{1}(2)\bigr),\ \bigl(I(1,3),P(0,4)\bigr),\ \bigl(I(1,8),P(1,1)\bigr),\ \bigl(I(0,7),P(1,2)\bigr)\\
R_{0}^{1}(4):\ & \bigl(R_{0}^{2}(3),R_{0}^{1}(1)\bigr),\ \bigl(R_{0}^{3}(2),R_{0}^{1}(2)\bigr),\ \bigl(R_{0}^{4}(1),R_{0}^{1}(3)\bigr),\ \bigl(I(1,3),P(0,5)\bigr),\ \bigl(I(1,8),P(1,4)\bigr)\\
 & \bigl(I(0,7),P(2,1)\bigr)\\
R_{0}^{2}(3):\ & \bigl(R_{0}^{3}(2),R_{0}^{2}(1)\bigr),\ \bigl(R_{0}^{4}(1),R_{0}^{2}(2)\bigr),\ \bigl(I(0,2),P(0,5)\bigr),\ \bigl(I(0,3),P(1,4)\bigr),\ \bigl(I(0,8),P(2,1)\bigr)\\
R_{0}^{2}(4):\ & \bigl(R_{0}^{3}(3),R_{0}^{2}(1)\bigr),\ \bigl(R_{0}^{4}(2),R_{0}^{2}(2)\bigr),\ \bigl(R_{0}^{5}(1),R_{0}^{2}(3)\bigr),\ \bigl(I(0,2),P(0,6)\bigr),\ \bigl(I(0,3),P(1,5)\bigr)\\
 & \bigl(I(0,8),P(2,4)\bigr)\\
R_{0}^{3}(3):\ & \bigl(R_{0}^{4}(2),R_{0}^{3}(1)\bigr),\ \bigl(R_{0}^{5}(1),R_{0}^{3}(2)\bigr)\\
R_{0}^{3}(4):\ & \bigl(R_{0}^{4}(3),R_{0}^{3}(1)\bigr),\ \bigl(R_{0}^{5}(2),R_{0}^{3}(2)\bigr),\ \bigl(R_{0}^{1}(1),R_{0}^{3}(3)\bigr)\\
R_{0}^{4}(3):\ & \bigl(R_{0}^{5}(2),R_{0}^{4}(1)\bigr),\ \bigl(R_{0}^{1}(1),R_{0}^{4}(2)\bigr)\\
R_{0}^{4}(4):\ & \bigl(R_{0}^{5}(3),R_{0}^{4}(1)\bigr),\ \bigl(R_{0}^{1}(2),R_{0}^{4}(2)\bigr),\ \bigl(R_{0}^{2}(1),R_{0}^{4}(3)\bigr),\ \bigl(I(2,7),P(0,1)\bigr),\ \bigl(I(1,6),P(0,2)\bigr)\\
 & \bigl(I(0,5),P(0,3)\bigr)\\
\end{align*}
\end{fleqn}
\subsubsection*{The non-homogeneous tube $\mathcal{T}_{\infty}^{\Delta(\widetilde{\mathbb{A}}_{3,5})}$}

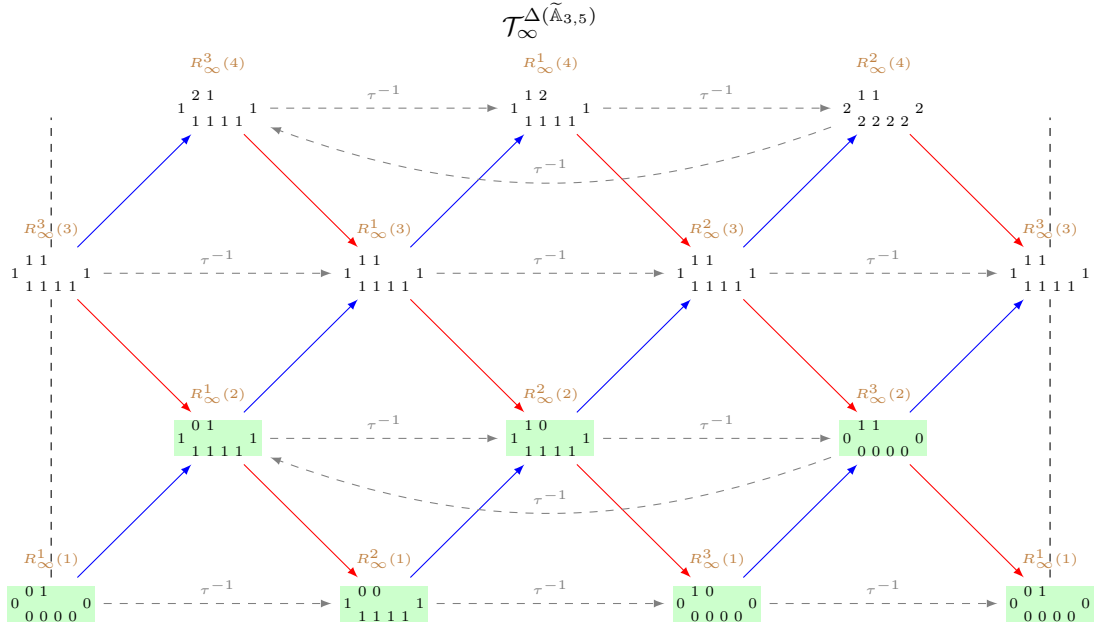
\begin{figure}[ht]



\begin{center}
\captionof{figure}{\vspace*{-10pt}$\mathcal{T}_{\infty}^{\Delta(\widetilde{\mathbb{A}}_{3,5})}$}
\begin{scaletikzpicturetowidth}{0.80000000000000004\textwidth}

\end{scaletikzpicturetowidth}
\end{center}

\end{figure}
\begin{fleqn}
\begin{align*}
R_{\infty}^{1}(1):\ & - \\
R_{\infty}^{1}(2):\ & \bigl(R_{\infty}^{2}(1),R_{\infty}^{1}(1)\bigr),\ \bigl(I(1,5),P(0,1)\bigr),\ \bigl(I(1,6),P(0,4)\bigr),\ \bigl(I(1,7),P(0,5)\bigr),\ \bigl(I(1,8),P(0,6)\bigr)\\
 & \bigl(I(0,3),P(0,7)\bigr)\\
R_{\infty}^{2}(1):\ & \bigl(I(0,4),P(0,1)\bigr),\ \bigl(I(0,5),P(0,4)\bigr),\ \bigl(I(0,6),P(0,5)\bigr),\ \bigl(I(0,7),P(0,6)\bigr),\ \bigl(I(0,8),P(0,7)\bigr)\\
R_{\infty}^{2}(2):\ & \bigl(R_{\infty}^{3}(1),R_{\infty}^{2}(1)\bigr),\ \bigl(I(0,4),P(0,2)\bigr),\ \bigl(I(0,5),P(1,1)\bigr),\ \bigl(I(0,6),P(1,4)\bigr),\ \bigl(I(0,7),P(1,5)\bigr)\\
 & \bigl(I(0,8),P(1,6)\bigr)\\
R_{\infty}^{3}(1):\ & - \\
R_{\infty}^{3}(2):\ & \bigl(R_{\infty}^{1}(1),R_{\infty}^{3}(1)\bigr)\\
\end{align*}
\end{fleqn}
\clearpage

\begin{lrbox}{\boxDDDDD}$\delta =\begin{smallmatrix}1&&1\\&2\\1&&1\end{smallmatrix} $\end{lrbox}
\subsection[Schofield pairs for the quiver $\Delta(\D_4)$]{Schofield pairs for the quiver $\Delta(\D_4)$ -- {\usebox{\boxDDDDD}}}
\[
\vcenter{\hbox{\xymatrix{2 &   & 3\ar[dl]\\
  & 5\ar[dl]\ar[ul] &  \\
1 &   & 4\ar[ul]}}}\qquad C_{\Delta(\widetilde{\mathbb{D}}_{4})} = \begin{bmatrix}1 & 0 & 1 & 1 & 1\\
0 & 1 & 1 & 1 & 1\\
0 & 0 & 1 & 0 & 0\\
0 & 0 & 0 & 1 & 0\\
0 & 0 & 1 & 1 & 1\end{bmatrix}\quad\Phi_{\Delta(\widetilde{\mathbb{D}}_{4})} = \begin{bmatrix}-1 & 0 & 0 & 0 & 1\\
0 & -1 & 0 & 0 & 1\\
-1 & -1 & 0 & 1 & 1\\
-1 & -1 & 1 & 0 & 1\\
-1 & -1 & 1 & 1 & 1\end{bmatrix}
\]
\subsubsection{Schofield pairs associated to preprojective exceptional modules}

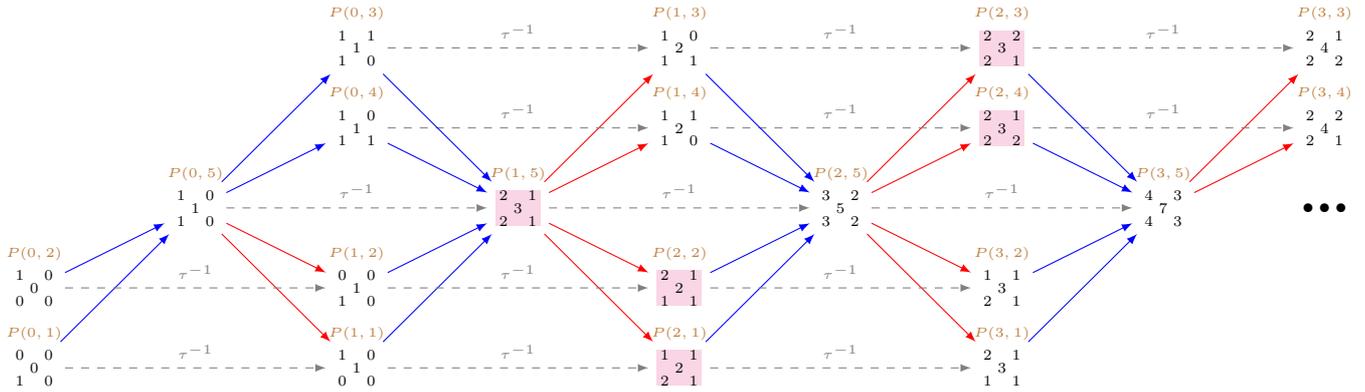
\begin{figure}[ht]


\begin{center}
\begin{scaletikzpicturetowidth}{\textwidth}

\end{scaletikzpicturetowidth}
\end{center}


\end{figure}
\medskip{}
\subsubsection*{Modules of the form $P(n,1)$}

Defect: $\partial P(n,1) = -1$, for $n\ge 0$.

\begin{fleqn}
\begin{align*}
P(0,1):\ & - \\
P(1,1):\ & \bigl(R_{1}^{1}(1),P(0,2)\bigr)\\
P(2,1):\ & \bigl(R_{0}^{2}(1),P(0,3)\bigr),\ \bigl(R_{\infty}^{2}(1),P(0,4)\bigr),\ \bigl(R_{1}^{2}(1),P(1,2)\bigr),\ \bigl(I(1,1),2P(0,1)\bigr)\\
P(n,1):\ & \bigl(R_{0}^{(n-1)\bmod 2+1}(1),P(n-2,3)\bigr),\ \bigl(R_{\infty}^{(n-1)\bmod 2+1}(1),P(n-2,4)\bigr),\ \bigl(R_{1}^{(n-1)\bmod 2+1}(1),P(n-1,2)\bigr)\\
 & \bigl(uI,(u+1)P\bigr),\ n>2\\
\end{align*}
\end{fleqn}
\subsubsection*{Modules of the form $P(n,2)$}

Defect: $\partial P(n,2) = -1$, for $n\ge 0$.

\begin{fleqn}
\begin{align*}
P(0,2):\ & - \\
P(1,2):\ & \bigl(R_{1}^{1}(1),P(0,1)\bigr)\\
P(2,2):\ & \bigl(R_{\infty}^{1}(1),P(0,3)\bigr),\ \bigl(R_{0}^{1}(1),P(0,4)\bigr),\ \bigl(R_{1}^{2}(1),P(1,1)\bigr),\ \bigl(I(1,2),2P(0,2)\bigr)\\
P(n,2):\ & \bigl(R_{\infty}^{(n-2)\bmod 2+1}(1),P(n-2,3)\bigr),\ \bigl(R_{0}^{(n-2)\bmod 2+1}(1),P(n-2,4)\bigr),\ \bigl(R_{1}^{(n-1)\bmod 2+1}(1),P(n-1,1)\bigr)\\
 & \bigl(uI,(u+1)P\bigr),\ n>2\\
\end{align*}
\end{fleqn}
\subsubsection*{Modules of the form $P(n,3)$}

Defect: $\partial P(n,3) = -1$, for $n\ge 0$.

\begin{fleqn}
\begin{align*}
P(0,3):\ & \bigl(R_{0}^{1}(1),P(0,1)\bigr),\ \bigl(R_{\infty}^{2}(1),P(0,2)\bigr),\ \bigl(I(0,3),P(0,5)\bigr)\\
P(1,3):\ & \bigl(R_{1}^{1}(1),P(0,4)\bigr),\ \bigl(R_{0}^{2}(1),P(1,1)\bigr),\ \bigl(R_{\infty}^{1}(1),P(1,2)\bigr)\\
P(2,3):\ & \bigl(R_{1}^{2}(1),P(1,4)\bigr),\ \bigl(R_{0}^{1}(1),P(2,1)\bigr),\ \bigl(R_{\infty}^{2}(1),P(2,2)\bigr),\ \bigl(I(1,3),2P(0,3)\bigr)\\
P(n,3):\ & \bigl(R_{1}^{(n-1)\bmod 2+1}(1),P(n-1,4)\bigr),\ \bigl(R_{0}^{(n-2)\bmod 2+1}(1),P(n,1)\bigr),\ \bigl(R_{\infty}^{(n-1)\bmod 2+1}(1),P(n,2)\bigr)\\
 & \bigl(uI,(u+1)P\bigr),\ n>2\\
\end{align*}
\end{fleqn}
\subsubsection*{Modules of the form $P(n,4)$}

Defect: $\partial P(n,4) = -1$, for $n\ge 0$.

\begin{fleqn}
\begin{align*}
P(0,4):\ & \bigl(R_{\infty}^{1}(1),P(0,1)\bigr),\ \bigl(R_{0}^{2}(1),P(0,2)\bigr),\ \bigl(I(0,4),P(0,5)\bigr)\\
P(1,4):\ & \bigl(R_{1}^{1}(1),P(0,3)\bigr),\ \bigl(R_{\infty}^{2}(1),P(1,1)\bigr),\ \bigl(R_{0}^{1}(1),P(1,2)\bigr)\\
P(2,4):\ & \bigl(R_{1}^{2}(1),P(1,3)\bigr),\ \bigl(R_{\infty}^{1}(1),P(2,1)\bigr),\ \bigl(R_{0}^{2}(1),P(2,2)\bigr),\ \bigl(I(1,4),2P(0,4)\bigr)\\
P(n,4):\ & \bigl(R_{1}^{(n-1)\bmod 2+1}(1),P(n-1,3)\bigr),\ \bigl(R_{\infty}^{(n-2)\bmod 2+1}(1),P(n,1)\bigr),\ \bigl(R_{0}^{(n-1)\bmod 2+1}(1),P(n,2)\bigr)\\
 & \bigl(uI,(u+1)P\bigr),\ n>2\\
\end{align*}
\end{fleqn}
\subsubsection*{Modules of the form $P(n,5)$}

Defect: $\partial P(n,5) = -2$, for $n\ge 0$.

\begin{fleqn}
\begin{align*}
P(0,5):\ & \bigl(P(1,1),P(0,1)\bigr),\ \bigl(P(1,2),P(0,2)\bigr)\\
P(1,5):\ & \bigl(P(1,3),P(0,3)\bigr),\ \bigl(P(1,4),P(0,4)\bigr),\ \bigl(P(2,1),P(1,1)\bigr),\ \bigl(P(2,2),P(1,2)\bigr)\\
P(n,5):\ & \bigl(P(n,3),P(n-1,3)\bigr),\ \bigl(P(n,4),P(n-1,4)\bigr),\ \bigl(P(n+1,1),P(n,1)\bigr)\\
 & \bigl(P(n+1,2),P(n,2)\bigr),\ n>1\\
\end{align*}
\end{fleqn}
\subsubsection{Schofield pairs associated to preinjective exceptional modules}

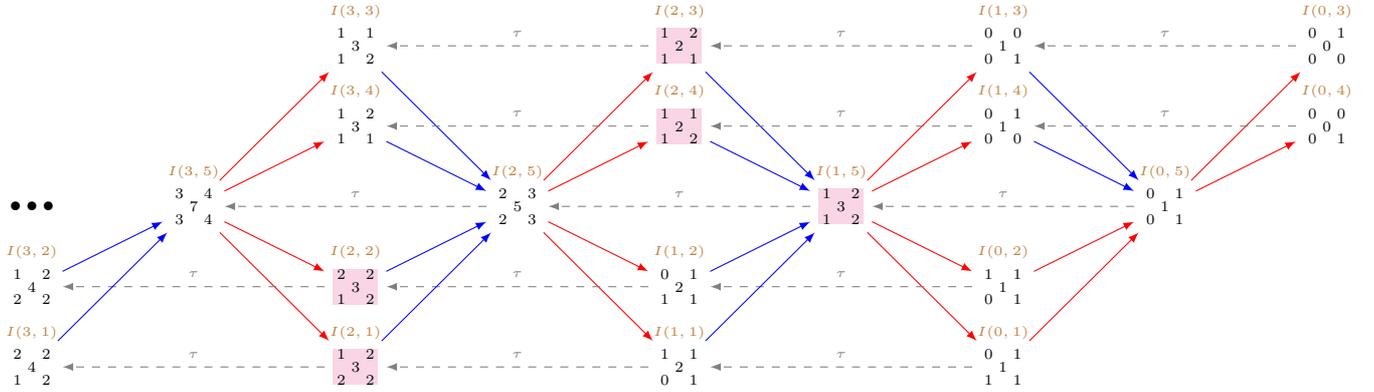
\begin{figure}[ht]


\begin{center}
\begin{scaletikzpicturetowidth}{\textwidth}

\end{scaletikzpicturetowidth}
\end{center}

\end{figure}
\medskip{}
\subsubsection*{Modules of the form $I(n,1)$}

Defect: $\partial I(n,1) = 1$, for $n\ge 0$.

\begin{fleqn}
\begin{align*}
I(0,1):\ & \bigl(I(0,3),R_{0}^{2}(1)\bigr),\ \bigl(I(0,4),R_{\infty}^{2}(1)\bigr),\ \bigl(I(0,5),P(0,1)\bigr)\\
I(1,1):\ & \bigl(I(0,2),R_{1}^{1}(1)\bigr),\ \bigl(I(1,3),R_{0}^{1}(1)\bigr),\ \bigl(I(1,4),R_{\infty}^{1}(1)\bigr)\\
I(2,1):\ & \bigl(I(1,2),R_{1}^{2}(1)\bigr),\ \bigl(I(2,3),R_{0}^{2}(1)\bigr),\ \bigl(I(2,4),R_{\infty}^{2}(1)\bigr),\ \bigl(2I(0,1),P(1,1)\bigr)\\
I(n,1):\ & \bigl(I(n-1,2),R_{1}^{(-n+3)\bmod 2+1}(1)\bigr),\ \bigl(I(n,3),R_{0}^{(-n+3)\bmod 2+1}(1)\bigr),\ \bigl(I(n,4),R_{\infty}^{(-n+3)\bmod 2+1}(1)\bigr)\\
 & \bigl((v+1)I,vP\bigr),\ n>2\\
\end{align*}
\end{fleqn}
\subsubsection*{Modules of the form $I(n,2)$}

Defect: $\partial I(n,2) = 1$, for $n\ge 0$.

\begin{fleqn}
\begin{align*}
I(0,2):\ & \bigl(I(0,3),R_{\infty}^{1}(1)\bigr),\ \bigl(I(0,4),R_{0}^{1}(1)\bigr),\ \bigl(I(0,5),P(0,2)\bigr)\\
I(1,2):\ & \bigl(I(0,1),R_{1}^{1}(1)\bigr),\ \bigl(I(1,3),R_{\infty}^{2}(1)\bigr),\ \bigl(I(1,4),R_{0}^{2}(1)\bigr)\\
I(2,2):\ & \bigl(I(1,1),R_{1}^{2}(1)\bigr),\ \bigl(I(2,3),R_{\infty}^{1}(1)\bigr),\ \bigl(I(2,4),R_{0}^{1}(1)\bigr),\ \bigl(2I(0,2),P(1,2)\bigr)\\
I(n,2):\ & \bigl(I(n-1,1),R_{1}^{(-n+3)\bmod 2+1}(1)\bigr),\ \bigl(I(n,3),R_{\infty}^{(-n+2)\bmod 2+1}(1)\bigr),\ \bigl(I(n,4),R_{0}^{(-n+2)\bmod 2+1}(1)\bigr)\\
 & \bigl((v+1)I,vP\bigr),\ n>2\\
\end{align*}
\end{fleqn}
\subsubsection*{Modules of the form $I(n,3)$}

Defect: $\partial I(n,3) = 1$, for $n\ge 0$.

\begin{fleqn}
\begin{align*}
I(0,3):\ & - \\
I(1,3):\ & \bigl(I(0,4),R_{1}^{1}(1)\bigr)\\
I(2,3):\ & \bigl(I(0,1),R_{0}^{1}(1)\bigr),\ \bigl(I(0,2),R_{\infty}^{2}(1)\bigr),\ \bigl(I(1,4),R_{1}^{2}(1)\bigr),\ \bigl(2I(0,3),P(1,3)\bigr)\\
I(n,3):\ & \bigl(I(n-2,1),R_{0}^{(-n+2)\bmod 2+1}(1)\bigr),\ \bigl(I(n-2,2),R_{\infty}^{(-n+3)\bmod 2+1}(1)\bigr),\ \bigl(I(n-1,4),R_{1}^{(-n+3)\bmod 2+1}(1)\bigr)\\
 & \bigl((v+1)I,vP\bigr),\ n>2\\
\end{align*}
\end{fleqn}
\subsubsection*{Modules of the form $I(n,4)$}

Defect: $\partial I(n,4) = 1$, for $n\ge 0$.

\begin{fleqn}
\begin{align*}
I(0,4):\ & - \\
I(1,4):\ & \bigl(I(0,3),R_{1}^{1}(1)\bigr)\\
I(2,4):\ & \bigl(I(0,1),R_{\infty}^{1}(1)\bigr),\ \bigl(I(0,2),R_{0}^{2}(1)\bigr),\ \bigl(I(1,3),R_{1}^{2}(1)\bigr),\ \bigl(2I(0,4),P(1,4)\bigr)\\
I(n,4):\ & \bigl(I(n-2,1),R_{\infty}^{(-n+2)\bmod 2+1}(1)\bigr),\ \bigl(I(n-2,2),R_{0}^{(-n+3)\bmod 2+1}(1)\bigr),\ \bigl(I(n-1,3),R_{1}^{(-n+3)\bmod 2+1}(1)\bigr)\\
 & \bigl((v+1)I,vP\bigr),\ n>2\\
\end{align*}
\end{fleqn}
\subsubsection*{Modules of the form $I(n,5)$}

Defect: $\partial I(n,5) = 2$, for $n\ge 0$.

\begin{fleqn}
\begin{align*}
I(0,5):\ & \bigl(I(0,3),I(1,3)\bigr),\ \bigl(I(0,4),I(1,4)\bigr)\\
I(1,5):\ & \bigl(I(0,1),I(1,1)\bigr),\ \bigl(I(0,2),I(1,2)\bigr),\ \bigl(I(1,3),I(2,3)\bigr),\ \bigl(I(1,4),I(2,4)\bigr)\\
I(n,5):\ & \bigl(I(n-1,1),I(n,1)\bigr),\ \bigl(I(n-1,2),I(n,2)\bigr),\ \bigl(I(n,3),I(n+1,3)\bigr)\\
 & \bigl(I(n,4),I(n+1,4)\bigr),\ n>1\\
\end{align*}
\end{fleqn}
\subsubsection{Schofield pairs associated to regular exceptional modules}

\subsubsection*{The non-homogeneous tube $\mathcal{T}_{1}^{\Delta(\widetilde{\mathbb{D}}_{4})}$}

\begin{figure}[ht]



\begin{center}
\captionof{figure}{\vspace*{-10pt}$\mathcal{T}_{1}^{\Delta(\widetilde{\mathbb{D}}_{4})}$}
\begin{scaletikzpicturetowidth}{0.40000000000000002\textwidth}

\end{scaletikzpicturetowidth}
\end{center}

\end{figure}
\begin{fleqn}
\begin{align*}
R_{1}^{1}(1):\ & - \\
R_{1}^{2}(1):\ & \bigl(I(0,2),P(0,1)\bigr),\ \bigl(I(0,1),P(0,2)\bigr),\ \bigl(I(0,4),P(0,3)\bigr),\ \bigl(I(0,3),P(0,4)\bigr)\\
\end{align*}
\end{fleqn}
\subsubsection*{The non-homogeneous tube $\mathcal{T}_{\infty}^{\Delta(\widetilde{\mathbb{D}}_{4})}$}

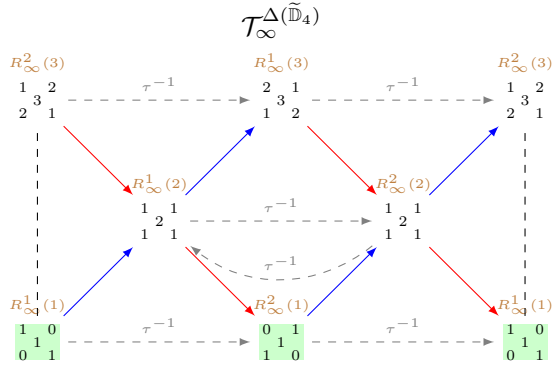
\begin{figure}[ht]



\begin{center}
\captionof{figure}{\vspace*{-10pt}$\mathcal{T}_{\infty}^{\Delta(\widetilde{\mathbb{D}}_{4})}$}
\begin{scaletikzpicturetowidth}{0.40000000000000002\textwidth}

\end{scaletikzpicturetowidth}
\end{center}

\end{figure}
\begin{fleqn}
\begin{align*}
R_{\infty}^{1}(1):\ & \bigl(I(1,3),P(0,2)\bigr),\ \bigl(I(0,4),P(1,1)\bigr)\\
R_{\infty}^{2}(1):\ & \bigl(I(1,4),P(0,1)\bigr),\ \bigl(I(0,3),P(1,2)\bigr)\\
\end{align*}
\end{fleqn}
\subsubsection*{The non-homogeneous tube $\mathcal{T}_{0}^{\Delta(\widetilde{\mathbb{D}}_{4})}$}

\begin{figure}[ht]



\begin{center}
\captionof{figure}{\vspace*{-10pt}$\mathcal{T}_{0}^{\Delta(\widetilde{\mathbb{D}}_{4})}$}
\begin{scaletikzpicturetowidth}{0.40000000000000002\textwidth}
 $\end{lrbox}
\subsection[Schofield pairs for the quiver $\Delta(\D_5)$]{Schofield pairs for the quiver $\Delta(\D_5)$ -- {\usebox{\boxDDDDDD}}}
\[
\vcenter{\hbox{\xymatrix{2 &   &   & 3\ar[dl]\\
  & 5\ar[dl]\ar[ul] & 6\ar[l] &  \\
1 &   &   & 4\ar[ul]}}}\qquad C_{\Delta(\widetilde{\mathbb{D}}_{5})} = %
\]
\subsubsection{Schofield pairs associated to preprojective exceptional modules}

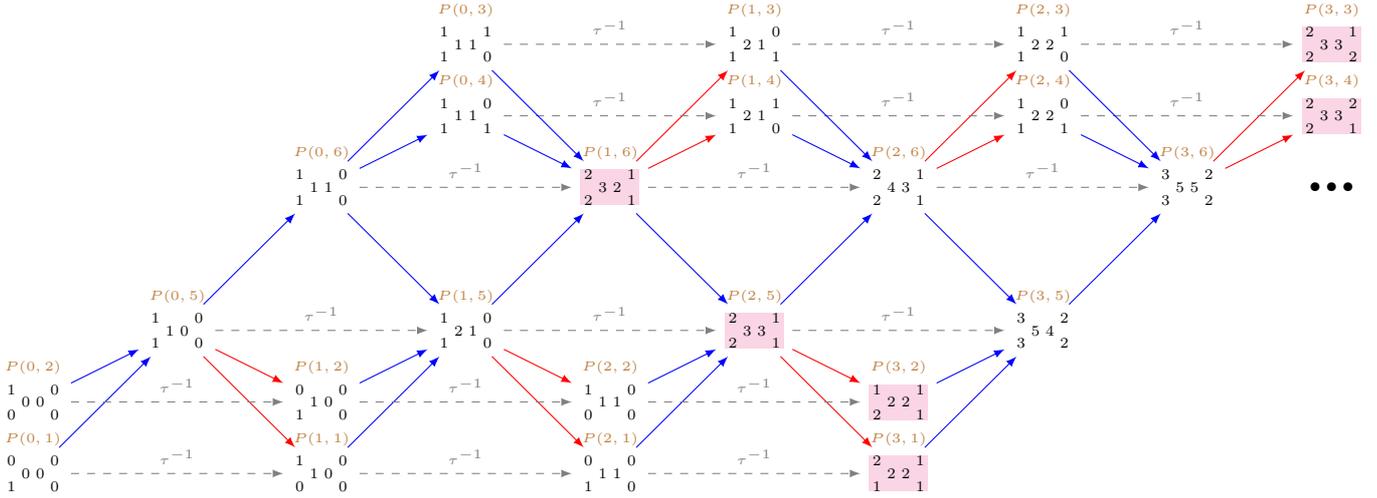
\begin{figure}[ht]


\begin{center}
\begin{scaletikzpicturetowidth}{\textwidth}

\end{scaletikzpicturetowidth}
\end{center}


\end{figure}
\medskip{}
\subsubsection*{Modules of the form $P(n,1)$}

Defect: $\partial P(n,1) = -1$, for $n\ge 0$.

\begin{fleqn}
\begin{align*}
P(0,1):\ & - \\
P(1,1):\ & \bigl(R_{1}^{3}(1),P(0,2)\bigr)\\
P(2,1):\ & \bigl(R_{1}^{3}(2),P(0,1)\bigr),\ \bigl(R_{1}^{1}(1),P(1,2)\bigr)\\
P(3,1):\ & \bigl(R_{\infty}^{1}(1),P(0,3)\bigr),\ \bigl(R_{0}^{1}(1),P(0,4)\bigr),\ \bigl(R_{1}^{1}(2),P(1,1)\bigr),\ \bigl(R_{1}^{2}(1),P(2,2)\bigr),\ \bigl(I(2,1),2P(0,2)\bigr)\\
P(n,1):\ & \bigl(R_{\infty}^{(n-3)\bmod 2+1}(1),P(n-3,3)\bigr),\ \bigl(R_{0}^{(n-3)\bmod 2+1}(1),P(n-3,4)\bigr),\ \bigl(R_{1}^{(n-3)\bmod 3+1}(2),P(n-2,1)\bigr)\\
 & \bigl(R_{1}^{(n-2)\bmod 3+1}(1),P(n-1,2)\bigr),\ \bigl(uI,(u+1)P\bigr),\ n>3\\
\end{align*}
\end{fleqn}
\subsubsection*{Modules of the form $P(n,2)$}

Defect: $\partial P(n,2) = -1$, for $n\ge 0$.

\begin{fleqn}
\begin{align*}
P(0,2):\ & - \\
P(1,2):\ & \bigl(R_{1}^{3}(1),P(0,1)\bigr)\\
P(2,2):\ & \bigl(R_{1}^{3}(2),P(0,2)\bigr),\ \bigl(R_{1}^{1}(1),P(1,1)\bigr)\\
P(3,2):\ & \bigl(R_{0}^{2}(1),P(0,3)\bigr),\ \bigl(R_{\infty}^{2}(1),P(0,4)\bigr),\ \bigl(R_{1}^{1}(2),P(1,2)\bigr),\ \bigl(R_{1}^{2}(1),P(2,1)\bigr),\ \bigl(I(2,2),2P(0,1)\bigr)\\
P(n,2):\ & \bigl(R_{0}^{(n-2)\bmod 2+1}(1),P(n-3,3)\bigr),\ \bigl(R_{\infty}^{(n-2)\bmod 2+1}(1),P(n-3,4)\bigr),\ \bigl(R_{1}^{(n-3)\bmod 3+1}(2),P(n-2,2)\bigr)\\
 & \bigl(R_{1}^{(n-2)\bmod 3+1}(1),P(n-1,1)\bigr),\ \bigl(uI,(u+1)P\bigr),\ n>3\\
\end{align*}
\end{fleqn}
\subsubsection*{Modules of the form $P(n,3)$}

Defect: $\partial P(n,3) = -1$, for $n\ge 0$.

\begin{fleqn}
\begin{align*}
P(0,3):\ & \bigl(R_{0}^{1}(1),P(0,1)\bigr),\ \bigl(R_{\infty}^{2}(1),P(0,2)\bigr),\ \bigl(I(1,4),P(0,5)\bigr),\ \bigl(I(0,3),P(0,6)\bigr)\\
P(1,3):\ & \bigl(R_{1}^{3}(1),P(0,4)\bigr),\ \bigl(R_{0}^{2}(1),P(1,1)\bigr),\ \bigl(R_{\infty}^{1}(1),P(1,2)\bigr),\ \bigl(I(0,4),P(1,5)\bigr)\\
P(2,3):\ & \bigl(R_{1}^{3}(2),P(0,3)\bigr),\ \bigl(R_{1}^{1}(1),P(1,4)\bigr),\ \bigl(R_{0}^{1}(1),P(2,1)\bigr),\ \bigl(R_{\infty}^{2}(1),P(2,2)\bigr)\\
P(3,3):\ & \bigl(R_{1}^{1}(2),P(1,3)\bigr),\ \bigl(R_{1}^{2}(1),P(2,4)\bigr),\ \bigl(R_{0}^{2}(1),P(3,1)\bigr),\ \bigl(R_{\infty}^{1}(1),P(3,2)\bigr),\ \bigl(I(2,3),2P(0,4)\bigr)\\
P(n,3):\ & \bigl(R_{1}^{(n-3)\bmod 3+1}(2),P(n-2,3)\bigr),\ \bigl(R_{1}^{(n-2)\bmod 3+1}(1),P(n-1,4)\bigr),\ \bigl(R_{0}^{(n-2)\bmod 2+1}(1),P(n,1)\bigr)\\
 & \bigl(R_{\infty}^{(n-3)\bmod 2+1}(1),P(n,2)\bigr),\ \bigl(uI,(u+1)P\bigr),\ n>3\\
\end{align*}
\end{fleqn}
\subsubsection*{Modules of the form $P(n,4)$}

Defect: $\partial P(n,4) = -1$, for $n\ge 0$.

\begin{fleqn}
\begin{align*}
P(0,4):\ & \bigl(R_{\infty}^{1}(1),P(0,1)\bigr),\ \bigl(R_{0}^{2}(1),P(0,2)\bigr),\ \bigl(I(1,3),P(0,5)\bigr),\ \bigl(I(0,4),P(0,6)\bigr)\\
P(1,4):\ & \bigl(R_{1}^{3}(1),P(0,3)\bigr),\ \bigl(R_{\infty}^{2}(1),P(1,1)\bigr),\ \bigl(R_{0}^{1}(1),P(1,2)\bigr),\ \bigl(I(0,3),P(1,5)\bigr)\\
P(2,4):\ & \bigl(R_{1}^{3}(2),P(0,4)\bigr),\ \bigl(R_{1}^{1}(1),P(1,3)\bigr),\ \bigl(R_{\infty}^{1}(1),P(2,1)\bigr),\ \bigl(R_{0}^{2}(1),P(2,2)\bigr)\\
P(3,4):\ & \bigl(R_{1}^{1}(2),P(1,4)\bigr),\ \bigl(R_{1}^{2}(1),P(2,3)\bigr),\ \bigl(R_{\infty}^{2}(1),P(3,1)\bigr),\ \bigl(R_{0}^{1}(1),P(3,2)\bigr),\ \bigl(I(2,4),2P(0,3)\bigr)\\
P(n,4):\ & \bigl(R_{1}^{(n-3)\bmod 3+1}(2),P(n-2,4)\bigr),\ \bigl(R_{1}^{(n-2)\bmod 3+1}(1),P(n-1,3)\bigr),\ \bigl(R_{\infty}^{(n-2)\bmod 2+1}(1),P(n,1)\bigr)\\
 & \bigl(R_{0}^{(n-3)\bmod 2+1}(1),P(n,2)\bigr),\ \bigl(uI,(u+1)P\bigr),\ n>3\\
\end{align*}
\end{fleqn}
\subsubsection*{Modules of the form $P(n,5)$}

Defect: $\partial P(n,5) = -2$, for $n\ge 0$.

\begin{fleqn}
\begin{align*}
P(0,5):\ & \bigl(P(1,1),P(0,1)\bigr),\ \bigl(P(1,2),P(0,2)\bigr)\\
P(1,5):\ & \bigl(R_{1}^{3}(1),P(0,6)\bigr),\ \bigl(P(2,1),P(1,1)\bigr),\ \bigl(P(2,2),P(1,2)\bigr)\\
P(2,5):\ & \bigl(P(2,4),P(0,3)\bigr),\ \bigl(P(2,3),P(0,4)\bigr),\ \bigl(R_{1}^{1}(1),P(1,6)\bigr),\ \bigl(P(3,1),P(2,1)\bigr),\ \bigl(P(3,2),P(2,2)\bigr)\\
P(n,5):\ & \bigl(P(n,4),P(n-2,3)\bigr),\ \bigl(P(n,3),P(n-2,4)\bigr),\ \bigl(R_{1}^{(n-2)\bmod 3+1}(1),P(n-1,6)\bigr)\\
 & \bigl(P(n+1,1),P(n,1)\bigr),\ \bigl(P(n+1,2),P(n,2)\bigr),\ n>2\\
\end{align*}
\end{fleqn}
\subsubsection*{Modules of the form $P(n,6)$}

Defect: $\partial P(n,6) = -2$, for $n\ge 0$.

\begin{fleqn}
\begin{align*}
P(0,6):\ & \bigl(P(2,2),P(0,1)\bigr),\ \bigl(P(2,1),P(0,2)\bigr),\ \bigl(R_{1}^{1}(1),P(0,5)\bigr)\\
P(1,6):\ & \bigl(P(1,3),P(0,3)\bigr),\ \bigl(P(1,4),P(0,4)\bigr),\ \bigl(P(3,2),P(1,1)\bigr),\ \bigl(P(3,1),P(1,2)\bigr),\ \bigl(R_{1}^{2}(1),P(1,5)\bigr)\\
P(n,6):\ & \bigl(P(n,3),P(n-1,3)\bigr),\ \bigl(P(n,4),P(n-1,4)\bigr),\ \bigl(P(n+2,2),P(n,1)\bigr)\\
 & \bigl(P(n+2,1),P(n,2)\bigr),\ \bigl(R_{1}^{n\bmod 3+1}(1),P(n,5)\bigr),\ n>1\\
\end{align*}
\end{fleqn}
\subsubsection{Schofield pairs associated to preinjective exceptional modules}

\begin{figure}[ht]


\begin{center}
\begin{scaletikzpicturetowidth}{\textwidth}

\end{scaletikzpicturetowidth}
\end{center}

\end{figure}
\medskip{}
\subsubsection*{Modules of the form $I(n,1)$}

Defect: $\partial I(n,1) = 1$, for $n\ge 0$.

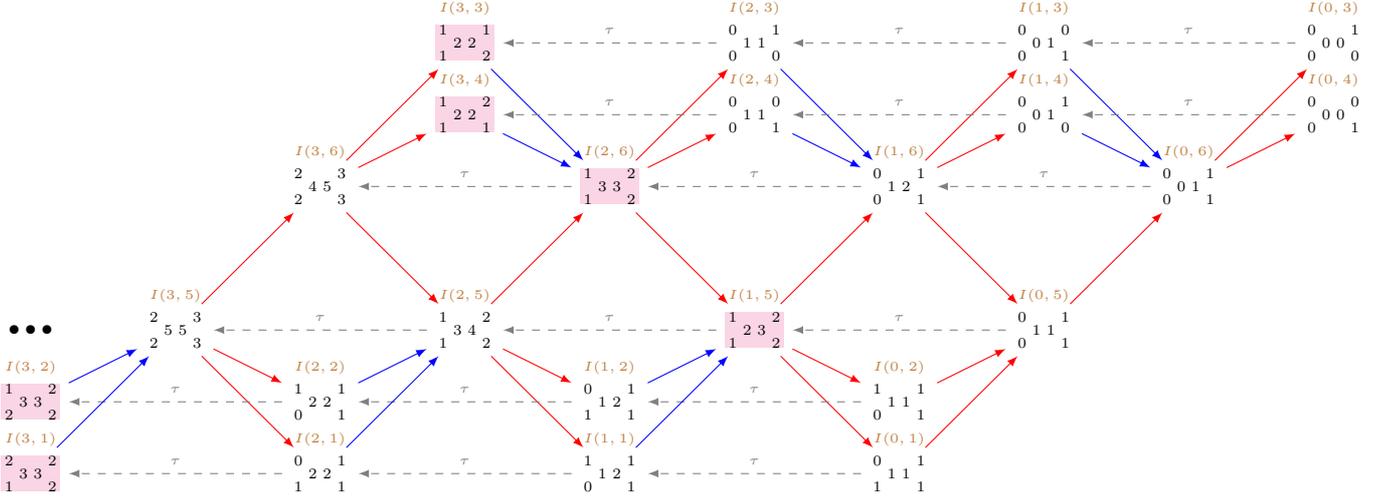
\begin{fleqn}
\begin{align*}
I(0,1):\ & \bigl(I(0,3),R_{0}^{2}(1)\bigr),\ \bigl(I(0,4),R_{\infty}^{2}(1)\bigr),\ \bigl(I(0,5),P(0,1)\bigr),\ \bigl(I(0,6),P(1,2)\bigr)\\
I(1,1):\ & \bigl(I(0,2),R_{1}^{1}(1)\bigr),\ \bigl(I(1,3),R_{0}^{1}(1)\bigr),\ \bigl(I(1,4),R_{\infty}^{1}(1)\bigr),\ \bigl(I(1,6),P(0,2)\bigr)\\
I(2,1):\ & \bigl(I(0,1),R_{1}^{3}(2)\bigr),\ \bigl(I(1,2),R_{1}^{3}(1)\bigr),\ \bigl(I(2,3),R_{0}^{2}(1)\bigr),\ \bigl(I(2,4),R_{\infty}^{2}(1)\bigr)\\
I(3,1):\ & \bigl(I(1,1),R_{1}^{2}(2)\bigr),\ \bigl(I(2,2),R_{1}^{2}(1)\bigr),\ \bigl(I(3,3),R_{0}^{1}(1)\bigr),\ \bigl(I(3,4),R_{\infty}^{1}(1)\bigr),\ \bigl(2I(0,2),P(2,1)\bigr)\\
I(n,1):\ & \bigl(I(n-2,1),R_{1}^{(-n+4)\bmod 3+1}(2)\bigr),\ \bigl(I(n-1,2),R_{1}^{(-n+4)\bmod 3+1}(1)\bigr),\ \bigl(I(n,3),R_{0}^{(-n+3)\bmod 2+1}(1)\bigr)\\
 & \bigl(I(n,4),R_{\infty}^{(-n+3)\bmod 2+1}(1)\bigr),\ \bigl((v+1)I,vP\bigr),\ n>3\\
\end{align*}
\end{fleqn}
\subsubsection*{Modules of the form $I(n,2)$}

Defect: $\partial I(n,2) = 1$, for $n\ge 0$.

\begin{fleqn}
\begin{align*}
I(0,2):\ & \bigl(I(0,3),R_{\infty}^{1}(1)\bigr),\ \bigl(I(0,4),R_{0}^{1}(1)\bigr),\ \bigl(I(0,5),P(0,2)\bigr),\ \bigl(I(0,6),P(1,1)\bigr)\\
I(1,2):\ & \bigl(I(0,1),R_{1}^{1}(1)\bigr),\ \bigl(I(1,3),R_{\infty}^{2}(1)\bigr),\ \bigl(I(1,4),R_{0}^{2}(1)\bigr),\ \bigl(I(1,6),P(0,1)\bigr)\\
I(2,2):\ & \bigl(I(0,2),R_{1}^{3}(2)\bigr),\ \bigl(I(1,1),R_{1}^{3}(1)\bigr),\ \bigl(I(2,3),R_{\infty}^{1}(1)\bigr),\ \bigl(I(2,4),R_{0}^{1}(1)\bigr)\\
I(3,2):\ & \bigl(I(1,2),R_{1}^{2}(2)\bigr),\ \bigl(I(2,1),R_{1}^{2}(1)\bigr),\ \bigl(I(3,3),R_{\infty}^{2}(1)\bigr),\ \bigl(I(3,4),R_{0}^{2}(1)\bigr),\ \bigl(2I(0,1),P(2,2)\bigr)\\
I(n,2):\ & \bigl(I(n-2,2),R_{1}^{(-n+4)\bmod 3+1}(2)\bigr),\ \bigl(I(n-1,1),R_{1}^{(-n+4)\bmod 3+1}(1)\bigr),\ \bigl(I(n,3),R_{\infty}^{(-n+4)\bmod 2+1}(1)\bigr)\\
 & \bigl(I(n,4),R_{0}^{(-n+4)\bmod 2+1}(1)\bigr),\ \bigl((v+1)I,vP\bigr),\ n>3\\
\end{align*}
\end{fleqn}
\subsubsection*{Modules of the form $I(n,3)$}

Defect: $\partial I(n,3) = 1$, for $n\ge 0$.

\begin{fleqn}
\begin{align*}
I(0,3):\ & - \\
I(1,3):\ & \bigl(I(0,4),R_{1}^{1}(1)\bigr)\\
I(2,3):\ & \bigl(I(0,3),R_{1}^{3}(2)\bigr),\ \bigl(I(1,4),R_{1}^{3}(1)\bigr)\\
I(3,3):\ & \bigl(I(0,1),R_{\infty}^{1}(1)\bigr),\ \bigl(I(0,2),R_{0}^{2}(1)\bigr),\ \bigl(I(1,3),R_{1}^{2}(2)\bigr),\ \bigl(I(2,4),R_{1}^{2}(1)\bigr),\ \bigl(2I(0,4),P(2,3)\bigr)\\
I(n,3):\ & \bigl(I(n-3,1),R_{\infty}^{(-n+3)\bmod 2+1}(1)\bigr),\ \bigl(I(n-3,2),R_{0}^{(-n+4)\bmod 2+1}(1)\bigr),\ \bigl(I(n-2,3),R_{1}^{(-n+4)\bmod 3+1}(2)\bigr)\\
 & \bigl(I(n-1,4),R_{1}^{(-n+4)\bmod 3+1}(1)\bigr),\ \bigl((v+1)I,vP\bigr),\ n>3\\
\end{align*}
\end{fleqn}
\subsubsection*{Modules of the form $I(n,4)$}

Defect: $\partial I(n,4) = 1$, for $n\ge 0$.

\begin{fleqn}
\begin{align*}
I(0,4):\ & - \\
I(1,4):\ & \bigl(I(0,3),R_{1}^{1}(1)\bigr)\\
I(2,4):\ & \bigl(I(0,4),R_{1}^{3}(2)\bigr),\ \bigl(I(1,3),R_{1}^{3}(1)\bigr)\\
I(3,4):\ & \bigl(I(0,1),R_{0}^{1}(1)\bigr),\ \bigl(I(0,2),R_{\infty}^{2}(1)\bigr),\ \bigl(I(1,4),R_{1}^{2}(2)\bigr),\ \bigl(I(2,3),R_{1}^{2}(1)\bigr),\ \bigl(2I(0,3),P(2,4)\bigr)\\
I(n,4):\ & \bigl(I(n-3,1),R_{0}^{(-n+3)\bmod 2+1}(1)\bigr),\ \bigl(I(n-3,2),R_{\infty}^{(-n+4)\bmod 2+1}(1)\bigr),\ \bigl(I(n-2,4),R_{1}^{(-n+4)\bmod 3+1}(2)\bigr)\\
 & \bigl(I(n-1,3),R_{1}^{(-n+4)\bmod 3+1}(1)\bigr),\ \bigl((v+1)I,vP\bigr),\ n>3\\
\end{align*}
\end{fleqn}
\subsubsection*{Modules of the form $I(n,5)$}

Defect: $\partial I(n,5) = 2$, for $n\ge 0$.

\begin{fleqn}
\begin{align*}
I(0,5):\ & \bigl(I(0,3),I(2,4)\bigr),\ \bigl(I(0,4),I(2,3)\bigr),\ \bigl(I(0,6),R_{1}^{3}(1)\bigr)\\
I(1,5):\ & \bigl(I(0,1),I(1,1)\bigr),\ \bigl(I(0,2),I(1,2)\bigr),\ \bigl(I(1,3),I(3,4)\bigr),\ \bigl(I(1,4),I(3,3)\bigr),\ \bigl(I(1,6),R_{1}^{2}(1)\bigr)\\
I(n,5):\ & \bigl(I(n-1,1),I(n,1)\bigr),\ \bigl(I(n-1,2),I(n,2)\bigr),\ \bigl(I(n,3),I(n+2,4)\bigr)\\
 & \bigl(I(n,4),I(n+2,3)\bigr),\ \bigl(I(n,6),R_{1}^{(-n+2)\bmod 3+1}(1)\bigr),\ n>1\\
\end{align*}
\end{fleqn}
\subsubsection*{Modules of the form $I(n,6)$}

Defect: $\partial I(n,6) = 2$, for $n\ge 0$.

\begin{fleqn}
\begin{align*}
I(0,6):\ & \bigl(I(0,3),I(1,3)\bigr),\ \bigl(I(0,4),I(1,4)\bigr)\\
I(1,6):\ & \bigl(I(0,5),R_{1}^{1}(1)\bigr),\ \bigl(I(1,3),I(2,3)\bigr),\ \bigl(I(1,4),I(2,4)\bigr)\\
I(2,6):\ & \bigl(I(0,1),I(2,2)\bigr),\ \bigl(I(0,2),I(2,1)\bigr),\ \bigl(I(1,5),R_{1}^{3}(1)\bigr),\ \bigl(I(2,3),I(3,3)\bigr),\ \bigl(I(2,4),I(3,4)\bigr)\\
I(n,6):\ & \bigl(I(n-2,1),I(n,2)\bigr),\ \bigl(I(n-2,2),I(n,1)\bigr),\ \bigl(I(n-1,5),R_{1}^{(-n+4)\bmod 3+1}(1)\bigr)\\
 & \bigl(I(n,3),I(n+1,3)\bigr),\ \bigl(I(n,4),I(n+1,4)\bigr),\ n>2\\
\end{align*}
\end{fleqn}
\subsubsection{Schofield pairs associated to regular exceptional modules}

\subsubsection*{The non-homogeneous tube $\mathcal{T}_{1}^{\Delta(\widetilde{\mathbb{D}}_{5})}$}

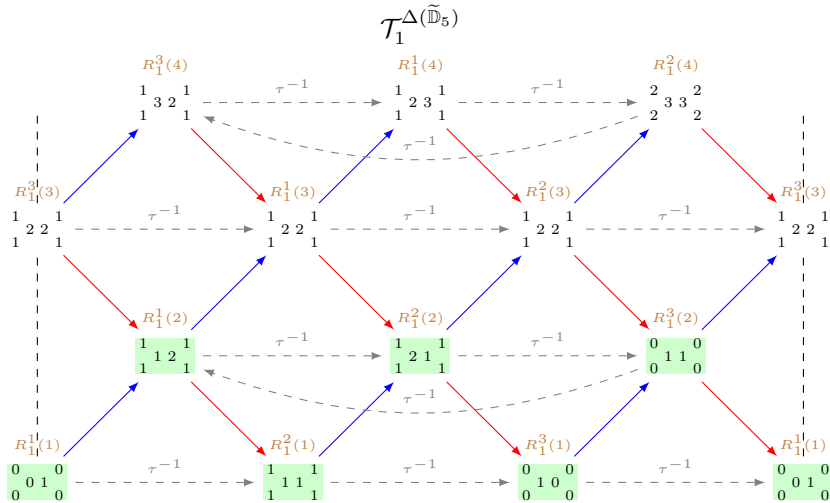
\begin{figure}[ht]



\begin{center}
\captionof{figure}{\vspace*{-10pt}$\mathcal{T}_{1}^{\Delta(\widetilde{\mathbb{D}}_{5})}$}
\begin{scaletikzpicturetowidth}{0.59999999999999998\textwidth}

\end{scaletikzpicturetowidth}
\end{center}

\end{figure}
\begin{fleqn}
\begin{align*}
R_{1}^{1}(1):\ & - \\
R_{1}^{1}(2):\ & \bigl(R_{1}^{2}(1),R_{1}^{1}(1)\bigr),\ \bigl(I(1,1),P(0,1)\bigr),\ \bigl(I(1,2),P(0,2)\bigr),\ \bigl(I(1,3),P(0,3)\bigr),\ \bigl(I(1,4),P(0,4)\bigr)\\
R_{1}^{2}(1):\ & \bigl(I(0,2),P(0,1)\bigr),\ \bigl(I(0,1),P(0,2)\bigr),\ \bigl(I(0,4),P(0,3)\bigr),\ \bigl(I(0,3),P(0,4)\bigr),\ \bigl(I(0,6),P(0,5)\bigr)\\
R_{1}^{2}(2):\ & \bigl(R_{1}^{3}(1),R_{1}^{2}(1)\bigr),\ \bigl(I(0,1),P(1,1)\bigr),\ \bigl(I(0,2),P(1,2)\bigr),\ \bigl(I(0,3),P(1,3)\bigr),\ \bigl(I(0,4),P(1,4)\bigr)\\
R_{1}^{3}(1):\ & - \\
R_{1}^{3}(2):\ & \bigl(R_{1}^{1}(1),R_{1}^{3}(1)\bigr)\\
\end{align*}
\end{fleqn}
\subsubsection*{The non-homogeneous tube $\mathcal{T}_{\infty}^{\Delta(\widetilde{\mathbb{D}}_{5})}$}

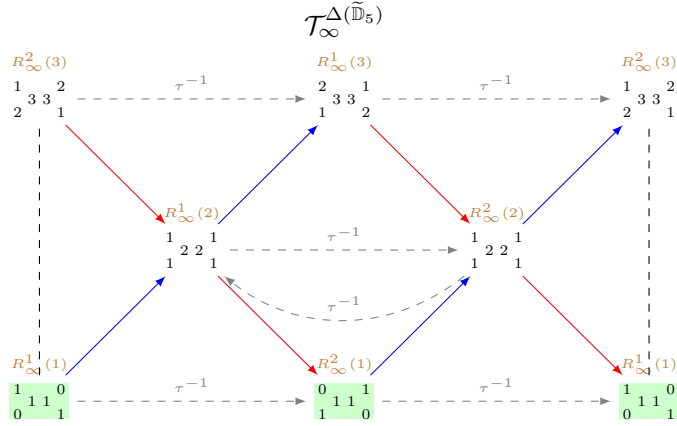
\begin{figure}[ht]



\begin{center}
\captionof{figure}{\vspace*{-10pt}$\mathcal{T}_{\infty}^{\Delta(\widetilde{\mathbb{D}}_{5})}$}
\begin{scaletikzpicturetowidth}{0.5\textwidth}

\end{scaletikzpicturetowidth}
\end{center}

\end{figure}
\begin{fleqn}
\begin{align*}
R_{\infty}^{1}(1):\ & \bigl(I(2,4),P(0,2)\bigr),\ \bigl(I(1,3),P(1,1)\bigr),\ \bigl(I(0,4),P(2,2)\bigr)\\
R_{\infty}^{2}(1):\ & \bigl(I(2,3),P(0,1)\bigr),\ \bigl(I(1,4),P(1,2)\bigr),\ \bigl(I(0,3),P(2,1)\bigr)\\
\end{align*}
\end{fleqn}
\subsubsection*{The non-homogeneous tube $\mathcal{T}_{0}^{\Delta(\widetilde{\mathbb{D}}_{5})}$}

\begin{figure}[ht]



\begin{center}
\captionof{figure}{\vspace*{-10pt}$\mathcal{T}_{0}^{\Delta(\widetilde{\mathbb{D}}_{5})}$}
\begin{scaletikzpicturetowidth}{0.5\textwidth}
 $\end{lrbox}
\subsection[Schofield pairs for the quiver $\Delta(\D_6)$]{Schofield pairs for the quiver $\Delta(\D_6)$ -- {\usebox{\boxDDDDDDD}}}
\[
\vcenter{\hbox{\xymatrix{2 &   &   &   & 3\ar[dl]\\
  & 5\ar[dl]\ar[ul] & 6\ar[l] & 7\ar[l] &  \\
1 &   &   &   & 4\ar[ul]}}}\qquad C_{\Delta(\widetilde{\mathbb{D}}_{6})} = %
\]
\subsubsection{Schofield pairs associated to preprojective exceptional modules}

\begin{figure}[ht]


\begin{center}
\begin{scaletikzpicturetowidth}{\textwidth}

\end{scaletikzpicturetowidth}
\end{center}


\end{figure}
\begin{figure}[ht]


\begin{center}
\begin{scaletikzpicturetowidth}{\textwidth}

\end{scaletikzpicturetowidth}
\end{center}


\end{figure}
\medskip{}
\subsubsection*{Modules of the form $P(n,1)$}

Defect: $\partial P(n,1) = -1$, for $n\ge 0$.

\begin{fleqn}
\begin{align*}
P(0,1):\ & - \\
P(1,1):\ & \bigl(R_{1}^{3}(1),P(0,2)\bigr)\\
P(2,1):\ & \bigl(R_{1}^{3}(2),P(0,1)\bigr),\ \bigl(R_{1}^{4}(1),P(1,2)\bigr)\\
P(3,1):\ & \bigl(R_{1}^{3}(3),P(0,2)\bigr),\ \bigl(R_{1}^{4}(2),P(1,1)\bigr),\ \bigl(R_{1}^{1}(1),P(2,2)\bigr)\\
P(4,1):\ & \bigl(R_{0}^{2}(1),P(0,3)\bigr),\ \bigl(R_{\infty}^{2}(1),P(0,4)\bigr),\ \bigl(R_{1}^{4}(3),P(1,2)\bigr),\ \bigl(R_{1}^{1}(2),P(2,1)\bigr),\ \bigl(R_{1}^{2}(1),P(3,2)\bigr)\\
 & \bigl(I(3,1),2P(0,1)\bigr)\\
P(n,1):\ & \bigl(R_{0}^{(n-3)\bmod 2+1}(1),P(n-4,3)\bigr),\ \bigl(R_{\infty}^{(n-3)\bmod 2+1}(1),P(n-4,4)\bigr),\ \bigl(R_{1}^{(n-1)\bmod 4+1}(3),P(n-3,2)\bigr)\\
 & \bigl(R_{1}^{(n-4)\bmod 4+1}(2),P(n-2,1)\bigr),\ \bigl(R_{1}^{(n-3)\bmod 4+1}(1),P(n-1,2)\bigr),\ \bigl(uI,(u+1)P\bigr),\ n>4\\
\end{align*}
\end{fleqn}
\subsubsection*{Modules of the form $P(n,2)$}

Defect: $\partial P(n,2) = -1$, for $n\ge 0$.

\begin{fleqn}
\begin{align*}
P(0,2):\ & - \\
P(1,2):\ & \bigl(R_{1}^{3}(1),P(0,1)\bigr)\\
P(2,2):\ & \bigl(R_{1}^{3}(2),P(0,2)\bigr),\ \bigl(R_{1}^{4}(1),P(1,1)\bigr)\\
P(3,2):\ & \bigl(R_{1}^{3}(3),P(0,1)\bigr),\ \bigl(R_{1}^{4}(2),P(1,2)\bigr),\ \bigl(R_{1}^{1}(1),P(2,1)\bigr)\\
P(4,2):\ & \bigl(R_{\infty}^{1}(1),P(0,3)\bigr),\ \bigl(R_{0}^{1}(1),P(0,4)\bigr),\ \bigl(R_{1}^{4}(3),P(1,1)\bigr),\ \bigl(R_{1}^{1}(2),P(2,2)\bigr),\ \bigl(R_{1}^{2}(1),P(3,1)\bigr)\\
 & \bigl(I(3,2),2P(0,2)\bigr)\\
P(n,2):\ & \bigl(R_{\infty}^{(n-4)\bmod 2+1}(1),P(n-4,3)\bigr),\ \bigl(R_{0}^{(n-4)\bmod 2+1}(1),P(n-4,4)\bigr),\ \bigl(R_{1}^{(n-1)\bmod 4+1}(3),P(n-3,1)\bigr)\\
 & \bigl(R_{1}^{(n-4)\bmod 4+1}(2),P(n-2,2)\bigr),\ \bigl(R_{1}^{(n-3)\bmod 4+1}(1),P(n-1,1)\bigr),\ \bigl(uI,(u+1)P\bigr),\ n>4\\
\end{align*}
\end{fleqn}
\subsubsection*{Modules of the form $P(n,3)$}

Defect: $\partial P(n,3) = -1$, for $n\ge 0$.

\begin{fleqn}
\begin{align*}
P(0,3):\ & \bigl(R_{0}^{1}(1),P(0,1)\bigr),\ \bigl(R_{\infty}^{2}(1),P(0,2)\bigr),\ \bigl(I(2,3),P(0,5)\bigr),\ \bigl(I(1,4),P(0,6)\bigr),\ \bigl(I(0,3),P(0,7)\bigr)\\
P(1,3):\ & \bigl(R_{1}^{3}(1),P(0,4)\bigr),\ \bigl(R_{0}^{2}(1),P(1,1)\bigr),\ \bigl(R_{\infty}^{1}(1),P(1,2)\bigr),\ \bigl(I(1,3),P(1,5)\bigr),\ \bigl(I(0,4),P(1,6)\bigr)\\
P(2,3):\ & \bigl(R_{1}^{3}(2),P(0,3)\bigr),\ \bigl(R_{1}^{4}(1),P(1,4)\bigr),\ \bigl(R_{0}^{1}(1),P(2,1)\bigr),\ \bigl(R_{\infty}^{2}(1),P(2,2)\bigr),\ \bigl(I(0,3),P(2,5)\bigr)\\
P(3,3):\ & \bigl(R_{1}^{3}(3),P(0,4)\bigr),\ \bigl(R_{1}^{4}(2),P(1,3)\bigr),\ \bigl(R_{1}^{1}(1),P(2,4)\bigr),\ \bigl(R_{0}^{2}(1),P(3,1)\bigr),\ \bigl(R_{\infty}^{1}(1),P(3,2)\bigr)\\
P(4,3):\ & \bigl(R_{1}^{4}(3),P(1,4)\bigr),\ \bigl(R_{1}^{1}(2),P(2,3)\bigr),\ \bigl(R_{1}^{2}(1),P(3,4)\bigr),\ \bigl(R_{0}^{1}(1),P(4,1)\bigr),\ \bigl(R_{\infty}^{2}(1),P(4,2)\bigr)\\
 & \bigl(I(3,3),2P(0,3)\bigr)\\
P(n,3):\ & \bigl(R_{1}^{(n-1)\bmod 4+1}(3),P(n-3,4)\bigr),\ \bigl(R_{1}^{(n-4)\bmod 4+1}(2),P(n-2,3)\bigr),\ \bigl(R_{1}^{(n-3)\bmod 4+1}(1),P(n-1,4)\bigr)\\
 & \bigl(R_{0}^{(n-4)\bmod 2+1}(1),P(n,1)\bigr),\ \bigl(R_{\infty}^{(n-3)\bmod 2+1}(1),P(n,2)\bigr),\ \bigl(uI,(u+1)P\bigr),\ n>4\\
\end{align*}
\end{fleqn}
\subsubsection*{Modules of the form $P(n,4)$}

Defect: $\partial P(n,4) = -1$, for $n\ge 0$.

\begin{fleqn}
\begin{align*}
P(0,4):\ & \bigl(R_{\infty}^{1}(1),P(0,1)\bigr),\ \bigl(R_{0}^{2}(1),P(0,2)\bigr),\ \bigl(I(2,4),P(0,5)\bigr),\ \bigl(I(1,3),P(0,6)\bigr),\ \bigl(I(0,4),P(0,7)\bigr)\\
P(1,4):\ & \bigl(R_{1}^{3}(1),P(0,3)\bigr),\ \bigl(R_{\infty}^{2}(1),P(1,1)\bigr),\ \bigl(R_{0}^{1}(1),P(1,2)\bigr),\ \bigl(I(1,4),P(1,5)\bigr),\ \bigl(I(0,3),P(1,6)\bigr)\\
P(2,4):\ & \bigl(R_{1}^{3}(2),P(0,4)\bigr),\ \bigl(R_{1}^{4}(1),P(1,3)\bigr),\ \bigl(R_{\infty}^{1}(1),P(2,1)\bigr),\ \bigl(R_{0}^{2}(1),P(2,2)\bigr),\ \bigl(I(0,4),P(2,5)\bigr)\\
P(3,4):\ & \bigl(R_{1}^{3}(3),P(0,3)\bigr),\ \bigl(R_{1}^{4}(2),P(1,4)\bigr),\ \bigl(R_{1}^{1}(1),P(2,3)\bigr),\ \bigl(R_{\infty}^{2}(1),P(3,1)\bigr),\ \bigl(R_{0}^{1}(1),P(3,2)\bigr)\\
P(4,4):\ & \bigl(R_{1}^{4}(3),P(1,3)\bigr),\ \bigl(R_{1}^{1}(2),P(2,4)\bigr),\ \bigl(R_{1}^{2}(1),P(3,3)\bigr),\ \bigl(R_{\infty}^{1}(1),P(4,1)\bigr),\ \bigl(R_{0}^{2}(1),P(4,2)\bigr)\\
 & \bigl(I(3,4),2P(0,4)\bigr)\\
P(n,4):\ & \bigl(R_{1}^{(n-1)\bmod 4+1}(3),P(n-3,3)\bigr),\ \bigl(R_{1}^{(n-4)\bmod 4+1}(2),P(n-2,4)\bigr),\ \bigl(R_{1}^{(n-3)\bmod 4+1}(1),P(n-1,3)\bigr)\\
 & \bigl(R_{\infty}^{(n-4)\bmod 2+1}(1),P(n,1)\bigr),\ \bigl(R_{0}^{(n-3)\bmod 2+1}(1),P(n,2)\bigr),\ \bigl(uI,(u+1)P\bigr),\ n>4\\
\end{align*}
\end{fleqn}
\subsubsection*{Modules of the form $P(n,5)$}

Defect: $\partial P(n,5) = -2$, for $n\ge 0$.

\begin{fleqn}
\begin{align*}
P(0,5):\ & \bigl(P(1,1),P(0,1)\bigr),\ \bigl(P(1,2),P(0,2)\bigr)\\
P(1,5):\ & \bigl(R_{1}^{3}(1),P(0,6)\bigr),\ \bigl(P(2,1),P(1,1)\bigr),\ \bigl(P(2,2),P(1,2)\bigr)\\
P(2,5):\ & \bigl(R_{1}^{3}(2),P(0,7)\bigr),\ \bigl(R_{1}^{4}(1),P(1,6)\bigr),\ \bigl(P(3,1),P(2,1)\bigr),\ \bigl(P(3,2),P(2,2)\bigr)\\
P(3,5):\ & \bigl(P(3,3),P(0,3)\bigr),\ \bigl(P(3,4),P(0,4)\bigr),\ \bigl(R_{1}^{4}(2),P(1,7)\bigr),\ \bigl(R_{1}^{1}(1),P(2,6)\bigr),\ \bigl(P(4,1),P(3,1)\bigr)\\
 & \bigl(P(4,2),P(3,2)\bigr)\\
P(n,5):\ & \bigl(P(n,3),P(n-3,3)\bigr),\ \bigl(P(n,4),P(n-3,4)\bigr),\ \bigl(R_{1}^{n\bmod 4+1}(2),P(n-2,7)\bigr)\\
 & \bigl(R_{1}^{(n-3)\bmod 4+1}(1),P(n-1,6)\bigr),\ \bigl(P(n+1,1),P(n,1)\bigr),\ \bigl(P(n+1,2),P(n,2)\bigr),\ n>3\\
\end{align*}
\end{fleqn}
\subsubsection*{Modules of the form $P(n,6)$}

Defect: $\partial P(n,6) = -2$, for $n\ge 0$.

\begin{fleqn}
\begin{align*}
P(0,6):\ & \bigl(P(2,2),P(0,1)\bigr),\ \bigl(P(2,1),P(0,2)\bigr),\ \bigl(R_{1}^{4}(1),P(0,5)\bigr)\\
P(1,6):\ & \bigl(R_{1}^{3}(1),P(0,7)\bigr),\ \bigl(P(3,2),P(1,1)\bigr),\ \bigl(P(3,1),P(1,2)\bigr),\ \bigl(R_{1}^{1}(1),P(1,5)\bigr)\\
P(2,6):\ & \bigl(P(2,4),P(0,3)\bigr),\ \bigl(P(2,3),P(0,4)\bigr),\ \bigl(R_{1}^{4}(1),P(1,7)\bigr),\ \bigl(P(4,2),P(2,1)\bigr),\ \bigl(P(4,1),P(2,2)\bigr)\\
 & \bigl(R_{1}^{2}(1),P(2,5)\bigr)\\
P(n,6):\ & \bigl(P(n,4),P(n-2,3)\bigr),\ \bigl(P(n,3),P(n-2,4)\bigr),\ \bigl(R_{1}^{(n+1)\bmod 4+1}(1),P(n-1,7)\bigr)\\
 & \bigl(P(n+2,2),P(n,1)\bigr),\ \bigl(P(n+2,1),P(n,2)\bigr),\ \bigl(R_{1}^{(n-1)\bmod 4+1}(1),P(n,5)\bigr),\ n>2\\
\end{align*}
\end{fleqn}
\subsubsection*{Modules of the form $P(n,7)$}

Defect: $\partial P(n,7) = -2$, for $n\ge 0$.

\begin{fleqn}
\begin{align*}
P(0,7):\ & \bigl(P(3,1),P(0,1)\bigr),\ \bigl(P(3,2),P(0,2)\bigr),\ \bigl(R_{1}^{4}(2),P(0,5)\bigr),\ \bigl(R_{1}^{1}(1),P(0,6)\bigr)\\
P(1,7):\ & \bigl(P(1,3),P(0,3)\bigr),\ \bigl(P(1,4),P(0,4)\bigr),\ \bigl(P(4,1),P(1,1)\bigr),\ \bigl(P(4,2),P(1,2)\bigr),\ \bigl(R_{1}^{1}(2),P(1,5)\bigr)\\
 & \bigl(R_{1}^{2}(1),P(1,6)\bigr)\\
P(n,7):\ & \bigl(P(n,3),P(n-1,3)\bigr),\ \bigl(P(n,4),P(n-1,4)\bigr),\ \bigl(P(n+3,1),P(n,1)\bigr)\\
 & \bigl(P(n+3,2),P(n,2)\bigr),\ \bigl(R_{1}^{(n-1)\bmod 4+1}(2),P(n,5)\bigr),\ \bigl(R_{1}^{n\bmod 4+1}(1),P(n,6)\bigr),\ n>1\\
\end{align*}
\end{fleqn}
\subsubsection{Schofield pairs associated to preinjective exceptional modules}

\begin{figure}[ht]


\begin{center}
\begin{scaletikzpicturetowidth}{\textwidth}

\end{scaletikzpicturetowidth}
\end{center}

\end{figure}
\begin{figure}[ht]


\begin{center}
\begin{scaletikzpicturetowidth}{\textwidth}

\end{scaletikzpicturetowidth}
\end{center}

\end{figure}
\medskip{}
\subsubsection*{Modules of the form $I(n,1)$}

Defect: $\partial I(n,1) = 1$, for $n\ge 0$.

\begin{fleqn}
\begin{align*}
I(0,1):\ & \bigl(I(0,3),R_{0}^{2}(1)\bigr),\ \bigl(I(0,4),R_{\infty}^{2}(1)\bigr),\ \bigl(I(0,5),P(0,1)\bigr),\ \bigl(I(0,6),P(1,2)\bigr),\ \bigl(I(0,7),P(2,1)\bigr)\\
I(1,1):\ & \bigl(I(0,2),R_{1}^{1}(1)\bigr),\ \bigl(I(1,3),R_{0}^{1}(1)\bigr),\ \bigl(I(1,4),R_{\infty}^{1}(1)\bigr),\ \bigl(I(1,6),P(0,2)\bigr),\ \bigl(I(1,7),P(1,1)\bigr)\\
I(2,1):\ & \bigl(I(0,1),R_{1}^{4}(2)\bigr),\ \bigl(I(1,2),R_{1}^{4}(1)\bigr),\ \bigl(I(2,3),R_{0}^{2}(1)\bigr),\ \bigl(I(2,4),R_{\infty}^{2}(1)\bigr),\ \bigl(I(2,7),P(0,1)\bigr)\\
I(3,1):\ & \bigl(I(0,2),R_{1}^{3}(3)\bigr),\ \bigl(I(1,1),R_{1}^{3}(2)\bigr),\ \bigl(I(2,2),R_{1}^{3}(1)\bigr),\ \bigl(I(3,3),R_{0}^{1}(1)\bigr),\ \bigl(I(3,4),R_{\infty}^{1}(1)\bigr)\\
I(4,1):\ & \bigl(I(1,2),R_{1}^{2}(3)\bigr),\ \bigl(I(2,1),R_{1}^{2}(2)\bigr),\ \bigl(I(3,2),R_{1}^{2}(1)\bigr),\ \bigl(I(4,3),R_{0}^{2}(1)\bigr),\ \bigl(I(4,4),R_{\infty}^{2}(1)\bigr)\\
 & \bigl(2I(0,1),P(3,1)\bigr)\\
I(n,1):\ & \bigl(I(n-3,2),R_{1}^{(-n+5)\bmod 4+1}(3)\bigr),\ \bigl(I(n-2,1),R_{1}^{(-n+5)\bmod 4+1}(2)\bigr),\ \bigl(I(n-1,2),R_{1}^{(-n+5)\bmod 4+1}(1)\bigr)\\
 & \bigl(I(n,3),R_{0}^{(-n+5)\bmod 2+1}(1)\bigr),\ \bigl(I(n,4),R_{\infty}^{(-n+5)\bmod 2+1}(1)\bigr),\ \bigl((v+1)I,vP\bigr),\ n>4\\
\end{align*}
\end{fleqn}
\subsubsection*{Modules of the form $I(n,2)$}

Defect: $\partial I(n,2) = 1$, for $n\ge 0$.

\begin{fleqn}
\begin{align*}
I(0,2):\ & \bigl(I(0,3),R_{\infty}^{1}(1)\bigr),\ \bigl(I(0,4),R_{0}^{1}(1)\bigr),\ \bigl(I(0,5),P(0,2)\bigr),\ \bigl(I(0,6),P(1,1)\bigr),\ \bigl(I(0,7),P(2,2)\bigr)\\
I(1,2):\ & \bigl(I(0,1),R_{1}^{1}(1)\bigr),\ \bigl(I(1,3),R_{\infty}^{2}(1)\bigr),\ \bigl(I(1,4),R_{0}^{2}(1)\bigr),\ \bigl(I(1,6),P(0,1)\bigr),\ \bigl(I(1,7),P(1,2)\bigr)\\
I(2,2):\ & \bigl(I(0,2),R_{1}^{4}(2)\bigr),\ \bigl(I(1,1),R_{1}^{4}(1)\bigr),\ \bigl(I(2,3),R_{\infty}^{1}(1)\bigr),\ \bigl(I(2,4),R_{0}^{1}(1)\bigr),\ \bigl(I(2,7),P(0,2)\bigr)\\
I(3,2):\ & \bigl(I(0,1),R_{1}^{3}(3)\bigr),\ \bigl(I(1,2),R_{1}^{3}(2)\bigr),\ \bigl(I(2,1),R_{1}^{3}(1)\bigr),\ \bigl(I(3,3),R_{\infty}^{2}(1)\bigr),\ \bigl(I(3,4),R_{0}^{2}(1)\bigr)\\
I(4,2):\ & \bigl(I(1,1),R_{1}^{2}(3)\bigr),\ \bigl(I(2,2),R_{1}^{2}(2)\bigr),\ \bigl(I(3,1),R_{1}^{2}(1)\bigr),\ \bigl(I(4,3),R_{\infty}^{1}(1)\bigr),\ \bigl(I(4,4),R_{0}^{1}(1)\bigr)\\
 & \bigl(2I(0,2),P(3,2)\bigr)\\
I(n,2):\ & \bigl(I(n-3,1),R_{1}^{(-n+5)\bmod 4+1}(3)\bigr),\ \bigl(I(n-2,2),R_{1}^{(-n+5)\bmod 4+1}(2)\bigr),\ \bigl(I(n-1,1),R_{1}^{(-n+5)\bmod 4+1}(1)\bigr)\\
 & \bigl(I(n,3),R_{\infty}^{(-n+4)\bmod 2+1}(1)\bigr),\ \bigl(I(n,4),R_{0}^{(-n+4)\bmod 2+1}(1)\bigr),\ \bigl((v+1)I,vP\bigr),\ n>4\\
\end{align*}
\end{fleqn}
\subsubsection*{Modules of the form $I(n,3)$}

Defect: $\partial I(n,3) = 1$, for $n\ge 0$.

\begin{fleqn}
\begin{align*}
I(0,3):\ & - \\
I(1,3):\ & \bigl(I(0,4),R_{1}^{1}(1)\bigr)\\
I(2,3):\ & \bigl(I(0,3),R_{1}^{4}(2)\bigr),\ \bigl(I(1,4),R_{1}^{4}(1)\bigr)\\
I(3,3):\ & \bigl(I(0,4),R_{1}^{3}(3)\bigr),\ \bigl(I(1,3),R_{1}^{3}(2)\bigr),\ \bigl(I(2,4),R_{1}^{3}(1)\bigr)\\
I(4,3):\ & \bigl(I(0,1),R_{0}^{1}(1)\bigr),\ \bigl(I(0,2),R_{\infty}^{2}(1)\bigr),\ \bigl(I(1,4),R_{1}^{2}(3)\bigr),\ \bigl(I(2,3),R_{1}^{2}(2)\bigr),\ \bigl(I(3,4),R_{1}^{2}(1)\bigr)\\
 & \bigl(2I(0,3),P(3,3)\bigr)\\
I(n,3):\ & \bigl(I(n-4,1),R_{0}^{(-n+4)\bmod 2+1}(1)\bigr),\ \bigl(I(n-4,2),R_{\infty}^{(-n+5)\bmod 2+1}(1)\bigr),\ \bigl(I(n-3,4),R_{1}^{(-n+5)\bmod 4+1}(3)\bigr)\\
 & \bigl(I(n-2,3),R_{1}^{(-n+5)\bmod 4+1}(2)\bigr),\ \bigl(I(n-1,4),R_{1}^{(-n+5)\bmod 4+1}(1)\bigr),\ \bigl((v+1)I,vP\bigr),\ n>4\\
\end{align*}
\end{fleqn}
\subsubsection*{Modules of the form $I(n,4)$}

Defect: $\partial I(n,4) = 1$, for $n\ge 0$.

\begin{fleqn}
\begin{align*}
I(0,4):\ & - \\
I(1,4):\ & \bigl(I(0,3),R_{1}^{1}(1)\bigr)\\
I(2,4):\ & \bigl(I(0,4),R_{1}^{4}(2)\bigr),\ \bigl(I(1,3),R_{1}^{4}(1)\bigr)\\
I(3,4):\ & \bigl(I(0,3),R_{1}^{3}(3)\bigr),\ \bigl(I(1,4),R_{1}^{3}(2)\bigr),\ \bigl(I(2,3),R_{1}^{3}(1)\bigr)\\
I(4,4):\ & \bigl(I(0,1),R_{\infty}^{1}(1)\bigr),\ \bigl(I(0,2),R_{0}^{2}(1)\bigr),\ \bigl(I(1,3),R_{1}^{2}(3)\bigr),\ \bigl(I(2,4),R_{1}^{2}(2)\bigr),\ \bigl(I(3,3),R_{1}^{2}(1)\bigr)\\
 & \bigl(2I(0,4),P(3,4)\bigr)\\
I(n,4):\ & \bigl(I(n-4,1),R_{\infty}^{(-n+4)\bmod 2+1}(1)\bigr),\ \bigl(I(n-4,2),R_{0}^{(-n+5)\bmod 2+1}(1)\bigr),\ \bigl(I(n-3,3),R_{1}^{(-n+5)\bmod 4+1}(3)\bigr)\\
 & \bigl(I(n-2,4),R_{1}^{(-n+5)\bmod 4+1}(2)\bigr),\ \bigl(I(n-1,3),R_{1}^{(-n+5)\bmod 4+1}(1)\bigr),\ \bigl((v+1)I,vP\bigr),\ n>4\\
\end{align*}
\end{fleqn}
\subsubsection*{Modules of the form $I(n,5)$}

Defect: $\partial I(n,5) = 2$, for $n\ge 0$.

\begin{fleqn}
\begin{align*}
I(0,5):\ & \bigl(I(0,3),I(3,3)\bigr),\ \bigl(I(0,4),I(3,4)\bigr),\ \bigl(I(0,6),R_{1}^{3}(1)\bigr),\ \bigl(I(0,7),R_{1}^{3}(2)\bigr)\\
I(1,5):\ & \bigl(I(0,1),I(1,1)\bigr),\ \bigl(I(0,2),I(1,2)\bigr),\ \bigl(I(1,3),I(4,3)\bigr),\ \bigl(I(1,4),I(4,4)\bigr),\ \bigl(I(1,6),R_{1}^{2}(1)\bigr)\\
 & \bigl(I(1,7),R_{1}^{2}(2)\bigr)\\
I(n,5):\ & \bigl(I(n-1,1),I(n,1)\bigr),\ \bigl(I(n-1,2),I(n,2)\bigr),\ \bigl(I(n,3),I(n+3,3)\bigr)\\
 & \bigl(I(n,4),I(n+3,4)\bigr),\ \bigl(I(n,6),R_{1}^{(-n+2)\bmod 4+1}(1)\bigr),\ \bigl(I(n,7),R_{1}^{(-n+2)\bmod 4+1}(2)\bigr),\ n>1\\
\end{align*}
\end{fleqn}
\subsubsection*{Modules of the form $I(n,6)$}

Defect: $\partial I(n,6) = 2$, for $n\ge 0$.

\begin{fleqn}
\begin{align*}
I(0,6):\ & \bigl(I(0,3),I(2,4)\bigr),\ \bigl(I(0,4),I(2,3)\bigr),\ \bigl(I(0,7),R_{1}^{4}(1)\bigr)\\
I(1,6):\ & \bigl(I(0,5),R_{1}^{1}(1)\bigr),\ \bigl(I(1,3),I(3,4)\bigr),\ \bigl(I(1,4),I(3,3)\bigr),\ \bigl(I(1,7),R_{1}^{3}(1)\bigr)\\
I(2,6):\ & \bigl(I(0,1),I(2,2)\bigr),\ \bigl(I(0,2),I(2,1)\bigr),\ \bigl(I(1,5),R_{1}^{4}(1)\bigr),\ \bigl(I(2,3),I(4,4)\bigr),\ \bigl(I(2,4),I(4,3)\bigr)\\
 & \bigl(I(2,7),R_{1}^{2}(1)\bigr)\\
I(n,6):\ & \bigl(I(n-2,1),I(n,2)\bigr),\ \bigl(I(n-2,2),I(n,1)\bigr),\ \bigl(I(n-1,5),R_{1}^{(-n+5)\bmod 4+1}(1)\bigr)\\
 & \bigl(I(n,3),I(n+2,4)\bigr),\ \bigl(I(n,4),I(n+2,3)\bigr),\ \bigl(I(n,7),R_{1}^{(-n+3)\bmod 4+1}(1)\bigr),\ n>2\\
\end{align*}
\end{fleqn}
\subsubsection*{Modules of the form $I(n,7)$}

Defect: $\partial I(n,7) = 2$, for $n\ge 0$.

\begin{fleqn}
\begin{align*}
I(0,7):\ & \bigl(I(0,3),I(1,3)\bigr),\ \bigl(I(0,4),I(1,4)\bigr)\\
I(1,7):\ & \bigl(I(0,6),R_{1}^{1}(1)\bigr),\ \bigl(I(1,3),I(2,3)\bigr),\ \bigl(I(1,4),I(2,4)\bigr)\\
I(2,7):\ & \bigl(I(0,5),R_{1}^{4}(2)\bigr),\ \bigl(I(1,6),R_{1}^{4}(1)\bigr),\ \bigl(I(2,3),I(3,3)\bigr),\ \bigl(I(2,4),I(3,4)\bigr)\\
I(3,7):\ & \bigl(I(0,1),I(3,1)\bigr),\ \bigl(I(0,2),I(3,2)\bigr),\ \bigl(I(1,5),R_{1}^{3}(2)\bigr),\ \bigl(I(2,6),R_{1}^{3}(1)\bigr),\ \bigl(I(3,3),I(4,3)\bigr)\\
 & \bigl(I(3,4),I(4,4)\bigr)\\
I(n,7):\ & \bigl(I(n-3,1),I(n,1)\bigr),\ \bigl(I(n-3,2),I(n,2)\bigr),\ \bigl(I(n-2,5),R_{1}^{(-n+5)\bmod 4+1}(2)\bigr)\\
 & \bigl(I(n-1,6),R_{1}^{(-n+5)\bmod 4+1}(1)\bigr),\ \bigl(I(n,3),I(n+1,3)\bigr),\ \bigl(I(n,4),I(n+1,4)\bigr),\ n>3\\
\end{align*}
\end{fleqn}
\subsubsection{Schofield pairs associated to regular exceptional modules}

\subsubsection*{The non-homogeneous tube $\mathcal{T}_{1}^{\Delta(\widetilde{\mathbb{D}}_{6})}$}

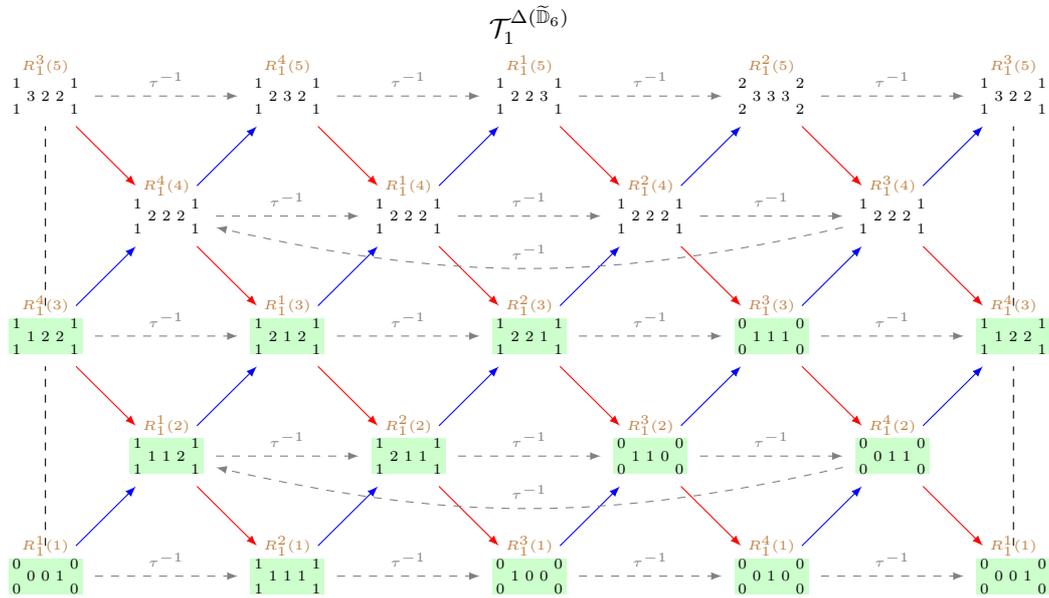
\begin{figure}[ht]



\begin{center}
\captionof{figure}{\vspace*{-10pt}$\mathcal{T}_{1}^{\Delta(\widetilde{\mathbb{D}}_{6})}$}
\begin{scaletikzpicturetowidth}{0.75\textwidth}

\end{scaletikzpicturetowidth}
\end{center}

\end{figure}
\begin{fleqn}
\begin{align*}
R_{1}^{1}(1):\ & - \\
R_{1}^{1}(2):\ & \bigl(R_{1}^{2}(1),R_{1}^{1}(1)\bigr),\ \bigl(I(1,1),P(0,1)\bigr),\ \bigl(I(1,2),P(0,2)\bigr),\ \bigl(I(1,3),P(0,3)\bigr),\ \bigl(I(1,4),P(0,4)\bigr)\\
 & \bigl(I(1,7),P(0,5)\bigr)\\
R_{1}^{2}(1):\ & \bigl(I(0,2),P(0,1)\bigr),\ \bigl(I(0,1),P(0,2)\bigr),\ \bigl(I(0,4),P(0,3)\bigr),\ \bigl(I(0,3),P(0,4)\bigr),\ \bigl(I(0,6),P(0,5)\bigr)\\
 & \bigl(I(0,7),P(0,6)\bigr)\\
R_{1}^{2}(2):\ & \bigl(R_{1}^{3}(1),R_{1}^{2}(1)\bigr),\ \bigl(I(0,1),P(1,1)\bigr),\ \bigl(I(0,2),P(1,2)\bigr),\ \bigl(I(0,3),P(1,3)\bigr),\ \bigl(I(0,4),P(1,4)\bigr)\\
 & \bigl(I(0,7),P(1,5)\bigr)\\
R_{1}^{3}(1):\ & - \\
R_{1}^{3}(2):\ & \bigl(R_{1}^{4}(1),R_{1}^{3}(1)\bigr)\\
R_{1}^{4}(1):\ & - \\
R_{1}^{4}(2):\ & \bigl(R_{1}^{1}(1),R_{1}^{4}(1)\bigr)\\
R_{1}^{4}(3):\ & \bigl(R_{1}^{1}(2),R_{1}^{4}(1)\bigr),\ \bigl(R_{1}^{2}(1),R_{1}^{4}(2)\bigr),\ \bigl(I(2,2),P(0,1)\bigr),\ \bigl(I(2,1),P(0,2)\bigr),\ \bigl(I(2,4),P(0,3)\bigr)\\
 & \bigl(I(2,3),P(0,4)\bigr)\\
R_{1}^{1}(3):\ & \bigl(R_{1}^{2}(2),R_{1}^{1}(1)\bigr),\ \bigl(R_{1}^{3}(1),R_{1}^{1}(2)\bigr),\ \bigl(I(1,2),P(1,1)\bigr),\ \bigl(I(1,1),P(1,2)\bigr),\ \bigl(I(1,4),P(1,3)\bigr)\\
 & \bigl(I(1,3),P(1,4)\bigr)\\
R_{1}^{2}(3):\ & \bigl(R_{1}^{3}(2),R_{1}^{2}(1)\bigr),\ \bigl(R_{1}^{4}(1),R_{1}^{2}(2)\bigr),\ \bigl(I(0,2),P(2,1)\bigr),\ \bigl(I(0,1),P(2,2)\bigr),\ \bigl(I(0,4),P(2,3)\bigr)\\
 & \bigl(I(0,3),P(2,4)\bigr)\\
R_{1}^{3}(3):\ & \bigl(R_{1}^{4}(2),R_{1}^{3}(1)\bigr),\ \bigl(R_{1}^{1}(1),R_{1}^{3}(2)\bigr)\\
\end{align*}
\end{fleqn}
\subsubsection*{The non-homogeneous tube $\mathcal{T}_{\infty}^{\Delta(\widetilde{\mathbb{D}}_{6})}$}

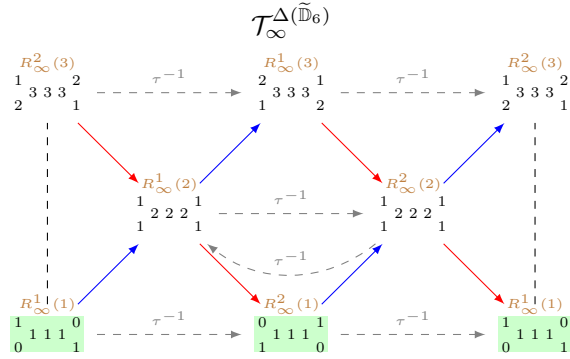
\begin{figure}[ht]



\begin{center}
\captionof{figure}{\vspace*{-10pt}$\mathcal{T}_{\infty}^{\Delta(\widetilde{\mathbb{D}}_{6})}$}
\begin{scaletikzpicturetowidth}{0.40000000000000002\textwidth}

\end{scaletikzpicturetowidth}
\end{center}

\end{figure}
\begin{fleqn}
\begin{align*}
R_{\infty}^{1}(1):\ & \bigl(I(3,3),P(0,2)\bigr),\ \bigl(I(2,4),P(1,1)\bigr),\ \bigl(I(1,3),P(2,2)\bigr),\ \bigl(I(0,4),P(3,1)\bigr)\\
R_{\infty}^{2}(1):\ & \bigl(I(3,4),P(0,1)\bigr),\ \bigl(I(2,3),P(1,2)\bigr),\ \bigl(I(1,4),P(2,1)\bigr),\ \bigl(I(0,3),P(3,2)\bigr)\\
\end{align*}
\end{fleqn}
\subsubsection*{The non-homogeneous tube $\mathcal{T}_{0}^{\Delta(\widetilde{\mathbb{D}}_{6})}$}

\begin{figure}[ht]



\begin{center}
\captionof{figure}{\vspace*{-10pt}$\mathcal{T}_{0}^{\Delta(\widetilde{\mathbb{D}}_{6})}$}
\begin{scaletikzpicturetowidth}{0.40000000000000002\textwidth}
 $\end{lrbox}
\subsection[Schofield pairs for the quiver $\Delta(\D_m)$]{Schofield pairs for the quiver $\Delta(\D_m)$ -- {\usebox{\boxDm}}, for $m\ge 4$}
\[
\xymatrix{2 &  &  &  &  & 3\ar[dl]\\
 & 5\ar[dl]\ar[ul] & 6\ar[l] & \ar[l]\cdots & \ar[l]m+1\\
1 &  &  &  &  & 4\ar[ul]
}
\]

\[
C_{\Delta(\widetilde{\mathbb{D}}_{m})}=%
\]

\subsubsection{Schofield pairs associated to preprojective exceptional modules}

$\dimz P(0,j)=j^{th}\text{ column of }C_{\Delta(\widetilde{\mathbb{D}}_{m})}$,\quad{}$\dimz I(0,i)=i^{th}\text{ row of }C_{\Delta(\widetilde{\mathbb{D}}_{m})}$\medskip{}

Throughout the lists below, in Schofield pairs of the form $(R,P)$
associated to the preprojective indecomposable $P'$, the preprojective
indecomposable $P$ is given by an explicit formula, while $R$ always
denotes the non-homogeneous regular with dimension vector $\dimz R=\dimz P'-\dimz P$
(for each such pair separately).

\subsubsection*{Modules of the form $P(n,1)$}

Defect: $\partial P(n,1)=-1$, for $n\ge0$. In the formulas below
$0\leq n_{1}<m-2$ and $n_{2}>m-2$:

\begin{fleqn}

\begin{align*}
P(n_{1},1):\  & \bigl(R,P(i,(n_{1}+i)\bmod2+1)\bigr),\ 0\leq i\leq n_{1}-1\\
P(m-2,1):\  & \bigl(R,P(i,(m+i)\bmod2+1)\bigr),\ 1\leq i\leq m-3\\
 & \bigl(R,P(0,3)\bigr),\ \bigl(R,P(0,4)\bigr),\ \bigl(I(m-3,1),2P(0,m\bmod2+1)\bigr)\\
P(n_{2},1):\  & \bigl(R,P(n_{2}-m+i+2,(m+i)\bmod2+1)\bigr),\ 1\leq i\leq m-3\\
 & \bigl(R,P(n_{2}-m+2,3)\bigr),\ \bigl(R,P(n_{2}-m+2,4)\bigr),\ \bigl(uI,(u+1)P\bigr)
\end{align*}

\end{fleqn}

\subsubsection*{Modules of the form $P(n,2)$}

Defect: $\partial P(n,2)=-1$, for $n\ge0$. In the formulas below
$0\leq n_{1}<m-2$ and $n_{2}>m-2$:

\begin{fleqn}

\begin{align*}
P(n_{1},2):\  & \bigl(R,P(i,(n_{1}+i+1)\bmod2+1)\bigr),\ 0\leq i\leq n_{1}-1\\
P(m-2,2):\  & \bigl(R,P(i,(m+i+1)\bmod2+1)\bigr),\ 1\leq i\leq m-3\\
 & \bigl(R,P(0,3)\bigr),\ \bigl(R,P(0,4)\bigr),\ \bigl(I(m-3,2),2P(0,(m+1)\bmod2+1)\bigr)\\
P(n_{2},2):\  & \bigl(R,P(n_{2}-m+i+2,(m+i+1)\bmod2+1)\bigr),\ 1\leq i\leq m-3\\
 & \bigl(R,P(n_{2}-m+2,3)\bigr),\ \bigl(R,P(n_{2}-m+2,4)\bigr),\ \bigl(uI,(u+1)P\bigr)
\end{align*}

\end{fleqn}

\subsubsection*{Modules of the form $P(n,3)$}

Defect: $\partial P(n,3)=-1$, for $n\ge0$. In the formulas below
$0\leq n_{1}<m-2$ and $n_{2}>m-2$:

\begin{fleqn}
\begin{align*}
P(n_{1},3):\  & \bigl(R,P(i,(n_{1}+i)\bmod2+3)\bigr),\ 0\leq i\leq n_{1}-1\\
 & \bigl(I(m-n_{1}-i+1,(m+i+1)\bmod2+3),P(n_{1},i)\bigr),\ 5\leq i\leq m-n_{1}+1\\
 & \bigl(R,P(n_{1},1)\bigr),\ \bigl(R,P(n_{1},2)\bigr)\\
P(m-2,3):\  & \bigl(R,P(i,(m+i)\bmod2+3)\bigr),\ 1\leq i\leq m-3\\
 & \bigl(R,P(m-2,1)\bigr),\ \bigl(R,P(m-2,2)\bigr),\ \bigl(I(m-3,3),2P(0,m\bmod2+3)\bigr)\\
P(n_{2},3):\  & \bigl(R,P(n_{2}-m+i+2,(m+i)\bmod2+3)\bigr),\ 1\leq i\leq m-3\\
 & \bigl(R,P(n_{2},1)\bigr),\ \bigl(R,P(n_{2},2)\bigr),\ \bigl(uI,(u+1)P\bigr)
\end{align*}

\end{fleqn}

\subsubsection*{Modules of the form $P(n,4)$}

Defect: $\partial P(n,4)=-1$, for $n\ge0$. In the formulas below
$0\leq n_{1}<m-2$ and $n_{2}>m-2$:

\begin{fleqn}
\begin{align*}
P(n_{1},4):\  & \bigl(R,P(i,(n_{1}+i+1)\bmod2+3)\bigr),\ 0\leq i\leq n_{1}-1\\
 & \bigl(I(m-n_{1}-i+1,(m+i)\bmod2+3),P(n_{1},i)\bigr),\ 5\leq i\leq m-n_{1}+1\\
 & \bigl(R,P(n_{1},1)\bigr),\ \bigl(R,P(n_{1},2)\bigr)\\
P(m-2,4):\  & \bigl(R,P(i,(m+i+1)\bmod2+3)\bigr),\ 1\leq i\leq m-3\\
 & \bigl(R,P(m-2,1)\bigr),\ \bigl(R,P(m-2,2)\bigr),\ \bigl(I(m-3,4),2P(0,(m+1)\bmod2+3)\bigr)\\
P(n_{2},4):\  & \bigl(R,P(n_{2}-m+i+2,(m+i+1)\bmod2+3)\bigr),\ 1\leq i\leq m-3\\
 & \bigl(R,P(n_{2},1)\bigr),\ \bigl(R,P(n_{2},2)\bigr),\ \bigl(uI,(u+1)P\bigr)
\end{align*}

\end{fleqn}

\subsubsection*{Modules of the form $P(n,j)$, for $5\protect\leq j\protect\leq m+1$}

Defect: $\partial P(n,j)=-2$, for $n\ge0$. In the formulas below
$0\leq n_{1}<m-j+1$ and $n_{2}>m-j+1$:

\begin{fleqn}

\begin{align*}
P(n_{1},j):\  & \bigl(R,P(i,j+n_{1}-i)\bigr),\ 0\leq i\leq n_{1}-1\\
 & \bigl(R,P(n_{1},i)\bigr),\ 5\leq i\leq j-1\\
 & \bigl(P(j-4+n_{1},(j+1)\bmod2+1),P(n_{1},1)\bigr),\ \bigl(P(j-4+n_{1},j\bmod2+1),P(n_{1},2)\bigr)\\
P(m-j+2,j):\  & \bigl(R,P(i+1,m-i+1)\bigr),\ 0\leq i\leq m-j\\
 & \bigl(R,P(m-j+2,i)\bigr),\ 5\leq i\leq j-1\\
 & \bigl(P(m-j+2,(m+j+1)\mod2+3),P(0,3)\bigr),\ \bigl(P(m-j+2,4-(m+j+1)\bmod2),P(0,4)\bigr),\\
 & \bigl(P(m-2,(j+1)\bmod2+1),P(m-j+2,1)\bigr),\ \bigl(P(m-2,j\bmod2+1),P(m-j+2,2)\bigr)\\
P(n_{2},j):\  & \bigl(R,P(n_{2}-m+j+i,m-i+1)\bigr),\ 0\leq i\leq m-j\\
 & \bigl(R,P(n_{2}+1,i)\bigr),\ 5\leq i\leq j-1\\
 & \bigl(P(n_{2}+1,(m+j+1)\mod2+3),P(n_{2}-m+j-1,3)\bigr),\\
 & \bigl(P(n_{2}+1,4-(m+j+1)\bmod2),P(n_{2}-m+j-1,4)\bigr),\\
 & \bigl(P(n_{2}+j-3,(j+1)\bmod2+1),P(n_{2}+1,1)\bigr),\\
 & \bigl(P(n_{2}+j-3,j\bmod2+1),P(n_{2}+1,2)\bigr)
\end{align*}

\end{fleqn}

\subsubsection{Schofield pairs associated to preinjective exceptional modules}

$\dimz P(0,j)=j^{th}\text{ column of }C_{\Delta(\widetilde{\mathbb{D}}_{m})}$,\quad{}$\dimz I(0,i)=i^{th}\text{ row of }C_{\Delta(\widetilde{\mathbb{D}}_{m})}$\medskip{}

Throughout the lists below, in Schofield pairs of the form $(I,R)$
associated to the preinjective indecomposable $I'$, the preinjective
indecomposable $I$ is given by an explicit formula, while $R$ always
denotes the non-homogeneous regular with dimension vector $\dimz R=\dimz I'-\dimz I$
(for each such pair separately).

\subsubsection*{Modules of the form $I(n,1)$}

Defect: $\partial I(n,1)=1$, for $n\ge0$. In the formulas below
$0\leq n_{1}<m-2$ and $n_{2}>m-2$:

\begin{fleqn}

\begin{align*}
I(n_{1},1):\  & \bigl(I(i,(n_{1}+i)\bmod2+1),R\bigr),\ 0\leq i\leq n_{1}-1\\
 & \bigl(I(n_{1},i),P(i-n_{1}-5,(i+1)\bmod2+1)\bigr),\ 5+n_{1}\leq i\leq m+1\\
 & \bigl(I(n_{1},3),R\bigr),\ \bigl(I(n_{1},4),R\bigr)\\
I(m-2,1):\  & \bigl(I(i,(m+i)\bmod2+1),R\bigr),\ 1\leq i\leq m-3\\
 & \bigl(I(m-2,3),R\bigr),\ \bigl(I(m-2,4),R\bigr),\ \bigl(2I(0,m\bmod2+1),P(m-3,1)\bigr)\\
I(n_{2},1):\  & \bigl(I(n_{2}-m+i+2,(m+i)\bmod2+1),R\bigr),\ 1\leq i\leq m-3\\
 & \bigl(I(n_{2},3),R\bigr),\ \bigl(I(n_{2},4),R\bigr),\ \bigl((v+1)I,vP\bigr)
\end{align*}

\end{fleqn}

\subsubsection*{Modules of the form $I(n,2)$}

Defect: $\partial I(n,2)=1$, for $n\ge0$. In the formulas below
$0\leq n_{1}<m-2$ and $n_{2}>m-2$:

\begin{fleqn}

\begin{align*}
I(n_{1},2):\  & \bigl(I(i,(n_{1}+i+1)\bmod2+1),R\bigr),\ 0\leq i\leq n_{1}-1\\
 & \bigl(I(n_{1},i),P(i-n_{1}-5,i\bmod2+1)\bigr),\ 5+n_{1}\leq i\leq m+1\\
 & \bigl(I(n_{1},3),R\bigr),\ \bigl(I(n_{1},4)\bigr)\\
I(m-2,2):\  & \bigl(I(i,(m+i+1)\bmod2+1),R\bigr),\ 1\leq i\leq m-3\\
 & \bigl(I(m-2,3),R\bigr),\ \bigl(I(m-2,4),R\bigr),\ \bigl(2I(0,(m+1)\bmod2+1),P(m-3,2)\bigr)\\
I(n_{2},2):\  & \bigl(I(n_{2}-m+i+2,(m+i+1)\bmod2+1),R\bigr),\ 1\leq i\leq m-3\\
 & \bigl(I(n_{2},3),R\bigr),\ \bigl(I(n_{2},4),R\bigr),\ \bigl((v+1)I,vP\bigr)
\end{align*}

\end{fleqn}

\subsubsection*{Modules of the form $I(n,3)$}

Defect: $\partial I(n,3)=1$, for $n\ge0$. In the formulas below
$0\leq n_{1}<m-2$ and $n_{2}>m-2$:

\begin{fleqn}

\begin{align*}
I(n_{1},3):\  & \bigl(I(i,(n_{1}+i)\bmod2+3),R\bigr),\ 0\leq i\leq n_{1}-1\\
I(m-2,3):\  & \bigl(I(i,(m+i)\bmod2+3),R\bigr),\ 1\leq i\leq m-3\\
 & \bigl(I(0,1),R\bigr),\ \bigl(I(0,2),R\bigr),\ \bigl(2I(0,m\bmod2+3),P(m-3,3)\bigr)\\
I(n_{2},3):\  & \bigl(I(n_{2}-m+i+2,(m+i)\bmod2+3),R\bigr),\ 1\leq i\leq m-3\\
 & \bigl(I(n_{2}-m+2,1),R\bigr),\ \bigl(I(n_{2}-m+2,2),R\bigr),\ \bigl((v+1)I,vP\bigr)
\end{align*}

\end{fleqn}

\subsubsection*{Modules of the form $I(n,4)$}

Defect: $\partial I(n,4)=1$, for $n\ge0$. In the formulas below
$0\leq n_{1}<m-2$ and $n_{2}>m-2$:

\begin{fleqn}

\begin{align*}
I(n_{1},4):\  & \bigl(I(i,(n_{1}+i+1)\bmod2+3),R\bigr),\ 0\leq i\leq n_{1}-1\\
I(m-2,4):\  & \bigl(I(i,(m+i+1)\bmod2+3),R\bigr),\ 1\leq i\leq m-3\\
 & \bigl(I(0,1),R\bigr),\ \bigl(I(0,2),R\bigr),\ \bigl(2I(0,(m+1)\bmod2+3),P(m-3,4)\bigr)\\
I(n_{2},4):\  & \bigl(I(n_{2}-m+i+2,(m+i+1)\bmod2+3),R\bigr),\ 1\leq i\leq m-3\\
 & \bigl(I(n_{2}-m+2,1),R\bigr),\ \bigl(I(n_{2}-m+2,2),R\bigr),\ \bigl((v+1)I,vP\bigr)
\end{align*}

\end{fleqn}

\subsubsection*{Modules of the form $I(n,j)$, for $5\protect\leq j\protect\leq m+1$}

Defect: $\partial I(n,j)=2$, for $n\ge0$. In the formulas below
$0\leq n_{1}<j-4$ and $n_{2}>j-4$:

\begin{fleqn}

\begin{align*}
I(n_{1},j):\  & \bigl(I(n_{1},i),R\bigr),\ j+1\leq i\leq m+1\\
 & \bigl(I(i-j+n_{1},i),R\bigr),\ j-n_{1}\leq i\leq j-1\\
 & \bigl(I(n_{1},3),I(m-j+n_{1}+2,(m+j+1)\bmod2+3)\bigr),\ \bigl(I(n_{1},4),I(m-j+n_{1}+2,(m+j)\bmod2+3)\bigr)\\
I(j-4,j):\  & \bigl(I(j-4,i),R\bigr),\ j+1\leq i\leq m+1\\
 & \bigl(I(i-4,i),R\bigr),\ 5\leq i\leq j-1\\
 & \bigl(I(j-4,3),I(m-2,(m+j+1)\bmod2+3)\bigr),\ \bigl(I(j-4,4),I(m-2,(m+j)\bmod2+3)\bigr),\\
 & \bigl(I(0,1),I(j-4,(j+1)\bmod2+1)\bigr),\ \bigl(I(0,2),I(j-4,j\bmod2+1)\bigr)\\
I(n_{2},j):\  & \bigl(I(n_{2},i),R\bigr),\ j+1\leq i\leq m+1\\
 & \bigl(I(n_{2}-j+i,i),R\bigr),\ 5\leq i\leq j-1\\
 & \bigl(I(n_{2},3),I(m-2,(m+j+1)\bmod2+3)\bigr),\ \bigl(I(n_{2},4),I(m-2,(m+j)\bmod2+3)\bigr),\\
 & \bigl(I(n_{2}-j+4,1),I(j-4,(j+1)\bmod2+1)\bigr),\ \bigl(I(n_{2}-j+4,2),I(j-4,j\bmod2+1)\bigr)
\end{align*}

\end{fleqn}

\subsubsection{Schofield pairs associated to regular exceptional modules}

\subsubsection*{The non-homogeneous tube $\mathcal{T}_{1}^{\Delta(\widetilde{\mathbb{D}}_{m})}$}

$\dimz R_{1}^{1}(1)=\begin{smallmatrix}0 &  &  &  &  &  & 0\\
 & 0 & 0 & \cdots & 0 & 1\\
0 &  &  &  &  &  & 0
\end{smallmatrix}$,$\quad\dimz R_{1}^{2}(1)=\begin{smallmatrix}1 &  &  &  &  &  & 1\\
 & 1 & 1 & 1 & \cdots & 1\\
1 &  &  &  &  &  & 1
\end{smallmatrix}$,$\quad\dimz R_{1}^{3}(1)=\begin{smallmatrix}0 &  &  &  &  &  & 0\\
 & 1 & 0 & 0 & \cdots & 0\\
0 &  &  &  &  &  & 0
\end{smallmatrix}$,\newline\newline$\dimz R_{1}^{4}(1)=\begin{smallmatrix}0 &  &  &  &  &  & 0\\
 & 0 & 1 & 0 & \cdots & 0\\
0 &  &  &  &  &  & 0
\end{smallmatrix}$,$\quad\dots$,$\quad\dimz R_{1}^{m-2}(1)=\begin{smallmatrix}0 &  &  &  &  &  & 0\\
 & 0 & \cdots & 0 & 1 & 0\\
0 &  &  &  &  &  & 0
\end{smallmatrix}$\newline

For all $1\leq t\leq m-3$ and $1\leq l\leq m-2$ we have the following
Schofield pairs for non-homogeneous regulars of the form $R_{1}^{l}(t)$,
where $t'=(l-3)\bmod(m-2)+1$:

\begin{fleqn}
\begin{align*}
R_{1}^{l}(t):\  & \bigl(R_{1}^{(l-1+i)\bmod(m-2)+1}(t-i),R_{1}^{l}(i)\bigr),\ 1\leq i\leq t-1\\
 & \bigl(I,P(t'-m+t+1,j)\bigr),\ 1\leq j\leq m-t+1,\ \text{if}\ t'-m+t+1\ge0,\\
 & \text{where }I\text{ is a preinjective exceptional with }\dimz I=\dimz R_{1}^{l}(t)-\dimz P(t'-m+t+1,j)
\end{align*}
\end{fleqn}

\subsubsection*{The non-homogeneous tube $\mathcal{T}_{\infty}^{\Delta(\widetilde{\mathbb{D}}_{m})}$}

$\dimz R_{\infty}^{1}(1)=\begin{smallmatrix}1 &  &  &  &  & 0\\
 & 1 & 1 & \cdots & 1\\
0 &  &  &  &  & 1
\end{smallmatrix}$,$\quad\dimz R_{\infty}^{2}(1)=\begin{smallmatrix}0 &  &  &  &  & 1\\
 & 1 & 1 & \cdots & 1\\
1 &  &  &  &  & 0
\end{smallmatrix}$

\begin{fleqn}
\begin{align*}
R_{\infty}^{1}(1):\  & \bigl(I(m-i-3,(m+i)\bmod2+3),P(i,2-(i\bmod2))\bigr),\ 0\leq i\leq m-3\\
R_{\infty}^{2}(1):\  & \bigl(I(m-i-3,(m+i+1)\bmod2+3),P(i,2-((i+1)\bmod2))\bigr),\ 0\leq i\leq m-3
\end{align*}
\end{fleqn}

\subsubsection*{The non-homogeneous tube $\mathcal{T}_{0}^{\Delta(\widetilde{\mathbb{D}}_{m})}$}

$\dimz R_{0}^{1}(1)=\begin{smallmatrix}1 &  &  &  &  & 1\\
 & 1 & 1 & \cdots & 1\\
0 &  &  &  &  & 0
\end{smallmatrix}$,$\quad\dimz R_{0}^{2}(0)=\begin{smallmatrix}0 &  &  &  &  & 0\\
 & 1 & 1 & \cdots & 1\\
1 &  &  &  &  & 1
\end{smallmatrix}$

\begin{fleqn}
\begin{align*}
R_{0}^{1}(1):\  & \bigl(I(m-i-3,(m+i+1)\bmod2+3),P(i,2-(i\bmod2))\bigr),\ 0\leq i\leq m-3\\
R_{0}^{2}(1):\  & \bigl(I(m-i-3,(m+i)\bmod2+3),P(i,2-((i+1)\bmod2))\bigr),\ 0\leq i\leq m-3
\end{align*}
\end{fleqn}

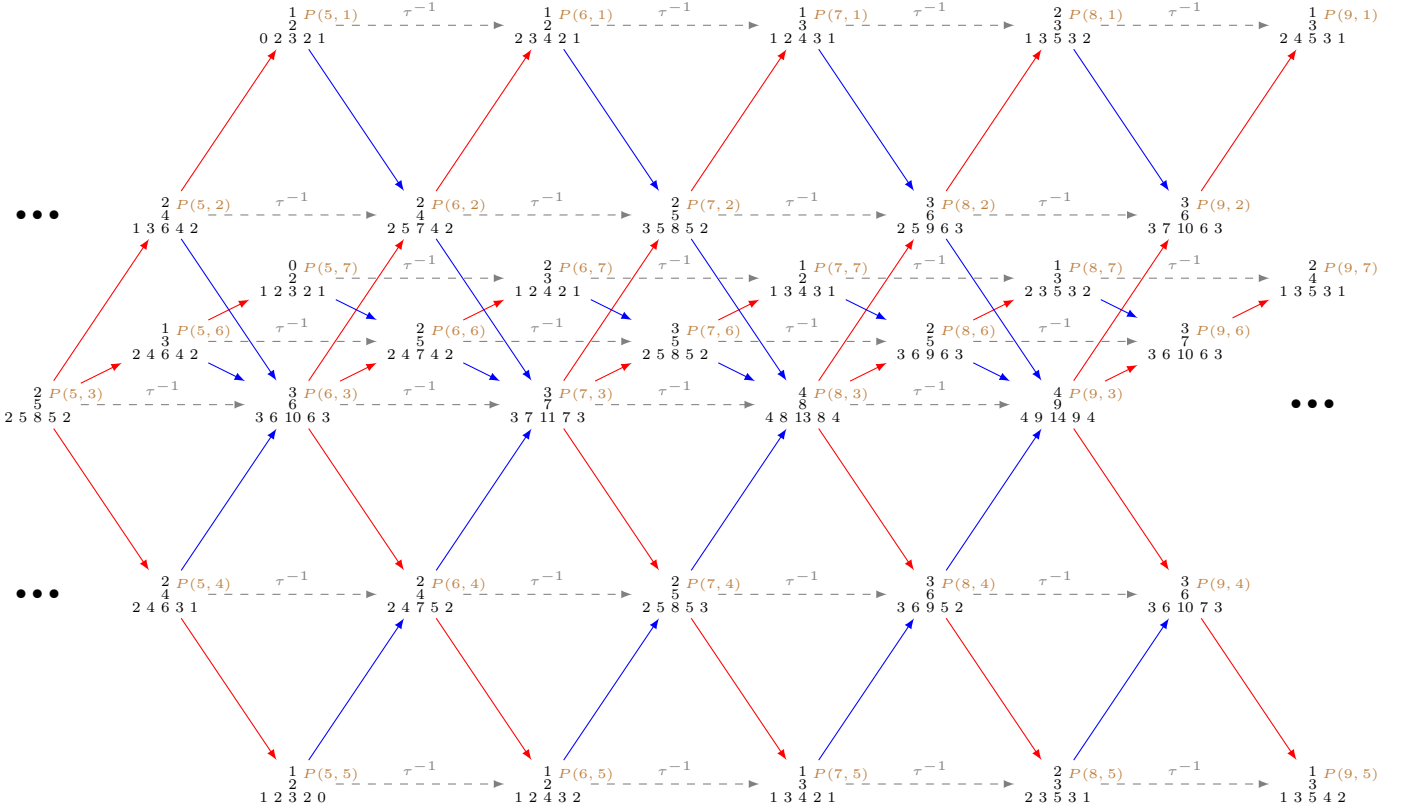
\begin{lrbox}{\boxEEEEEEE}$\delta =\begin{smallmatrix}&&1\\&&2\\1&2&3&2&1\end{smallmatrix} $\end{lrbox}
\subsection[Schofield pairs for the quiver $\Delta(\E_6)$]{Schofield pairs for the quiver $\Delta(\E_6)$ -- {\usebox{\boxEEEEEEE}}}
\[
\vcenter{\hbox{\xymatrix{  &   & 7\ar[d] &   &  \\
  &   & 6\ar[d] &   &  \\
1\ar[r] & 2\ar[r] & 3 & 4\ar[l] & 5\ar[l]}}}\qquad C_{\Delta(\widetilde{\mathbb{E}}_{6})} = \begin{bmatrix}1 & 0 & 0 & 0 & 0 & 0 & 0\\
1 & 1 & 0 & 0 & 0 & 0 & 0\\
1 & 1 & 1 & 1 & 1 & 1 & 1\\
0 & 0 & 0 & 1 & 1 & 0 & 0\\
0 & 0 & 0 & 0 & 1 & 0 & 0\\
0 & 0 & 0 & 0 & 0 & 1 & 1\\
0 & 0 & 0 & 0 & 0 & 0 & 1\end{bmatrix}\quad\Phi_{\Delta(\widetilde{\mathbb{E}}_{6})} = \begin{bmatrix}0 & 0 & -1 & 1 & 0 & 1 & 0\\
1 & 0 & -1 & 1 & 0 & 1 & 0\\
0 & 1 & -1 & 1 & 0 & 1 & 0\\
0 & 1 & -1 & 0 & 1 & 1 & 0\\
0 & 1 & -1 & 0 & 0 & 1 & 0\\
0 & 1 & -1 & 1 & 0 & 0 & 1\\
0 & 1 & -1 & 1 & 0 & 0 & 0\end{bmatrix}
\]
\subsubsection{Schofield pairs associated to preprojective exceptional modules}

\begin{figure}[ht]


\begin{center}
\begin{scaletikzpicturetowidth}{\textwidth}

\end{scaletikzpicturetowidth}
\end{center}


\end{figure}
\begin{figure}[ht]


\begin{center}
\begin{scaletikzpicturetowidth}{\textwidth}

\end{scaletikzpicturetowidth}
\end{center}


\end{figure}
\medskip{}
\subsubsection*{Modules of the form $P(n,1)$}

Defect: $\partial P(n,1) = -1$, for $n\ge 0$.

\begin{fleqn}
\begin{align*}
P(0,1):\ & \bigl(I(0,1),P(0,2)\bigr),\ \bigl(I(0,2),P(0,3)\bigr)\\
P(1,1):\ & \bigl(I(1,7),P(0,4)\bigr),\ \bigl(I(1,5),P(0,6)\bigr)\\
P(2,1):\ & \bigl(I(2,1),P(0,2)\bigr),\ \bigl(R_{0}^{2}(1),P(0,5)\bigr),\ \bigl(R_{1}^{2}(1),P(0,7)\bigr),\ \bigl(I(0,7),P(1,4)\bigr),\ \bigl(I(0,5),P(1,6)\bigr)\\
P(3,1):\ & \bigl(R_{\infty}^{1}(1),P(0,1)\bigr),\ \bigl(I(1,1),P(1,2)\bigr),\ \bigl(R_{0}^{3}(1),P(1,5)\bigr),\ \bigl(R_{1}^{3}(1),P(1,7)\bigr)\\
P(4,1):\ & \bigl(R_{1}^{3}(2),P(0,5)\bigr),\ \bigl(R_{0}^{3}(2),P(0,7)\bigr),\ \bigl(R_{\infty}^{2}(1),P(1,1)\bigr),\ \bigl(I(0,1),P(2,2)\bigr),\ \bigl(R_{0}^{1}(1),P(2,5)\bigr)\\
 & \bigl(R_{1}^{1}(1),P(2,7)\bigr)\\
P(n,1):\ & \bigl(R_{1}^{(n-2)\bmod 3+1}(2),P(n-4,5)\bigr),\ \bigl(R_{0}^{(n-2)\bmod 3+1}(2),P(n-4,7)\bigr),\ \bigl(R_{\infty}^{(n-3)\bmod 2+1}(1),P(n-3,1)\bigr)\\
 & \bigl(I(-n+4,1),P(n-2,2)\bigr),\ \bigl(R_{0}^{(n-4)\bmod 3+1}(1),P(n-2,5)\bigr),\ \bigl(R_{1}^{(n-4)\bmod 3+1}(1),P(n-2,7)\bigr),\ n>4\\
\end{align*}
\end{fleqn}
\subsubsection*{Modules of the form $P(n,2)$}

Defect: $\partial P(n,2) = -2$, for $n\ge 0$.

\begin{fleqn}
\begin{align*}
P(0,2):\ & \bigl(I(1,1),P(0,3)\bigr)\\
P(1,2):\ & \bigl(P(1,1),P(0,1)\bigr),\ \bigl(R_{1}^{3}(1),P(0,4)\bigr),\ \bigl(R_{0}^{3}(1),P(0,6)\bigr),\ \bigl(I(0,1),P(1,3)\bigr)\\
P(2,2):\ & \bigl(P(3,7),P(0,5)\bigr),\ \bigl(P(3,5),P(0,7)\bigr),\ \bigl(P(2,1),P(1,1)\bigr),\ \bigl(R_{1}^{1}(1),P(1,4)\bigr),\ \bigl(R_{0}^{1}(1),P(1,6)\bigr)\\
P(3,2):\ & \bigl(P(5,1),P(0,1)\bigr),\ \bigl(P(4,7),P(1,5)\bigr),\ \bigl(P(4,5),P(1,7)\bigr),\ \bigl(P(3,1),P(2,1)\bigr),\ \bigl(R_{1}^{2}(1),P(2,4)\bigr)\\
 & \bigl(R_{0}^{2}(1),P(2,6)\bigr)\\
P(n,2):\ & \bigl(P(n+2,1),P(n-3,1)\bigr),\ \bigl(P(n+1,7),P(n-2,5)\bigr),\ \bigl(P(n+1,5),P(n-2,7)\bigr)\\
 & \bigl(P(n,1),P(n-1,1)\bigr),\ \bigl(R_{1}^{(n-2)\bmod 3+1}(1),P(n-1,4)\bigr),\ \bigl(R_{0}^{(n-2)\bmod 3+1}(1),P(n-1,6)\bigr),\ n>3\\
\end{align*}
\end{fleqn}
\subsubsection*{Modules of the form $P(n,3)$}

Defect: $\partial P(n,3) = -3$, for $n\ge 0$.

\begin{fleqn}
\begin{align*}
P(0,3):\ & - \\
P(1,3):\ & \bigl(P(1,1),P(0,2)\bigr),\ \bigl(P(1,5),P(0,4)\bigr),\ \bigl(P(1,7),P(0,6)\bigr)\\
P(2,3):\ & \bigl(P(2,2),P(0,1)\bigr),\ \bigl(P(2,4),P(0,5)\bigr),\ \bigl(P(2,6),P(0,7)\bigr),\ \bigl(P(2,1),P(1,2)\bigr),\ \bigl(P(2,5),P(1,4)\bigr)\\
 & \bigl(P(2,7),P(1,6)\bigr)\\
P(n,3):\ & \bigl(P(n,2),P(n-2,1)\bigr),\ \bigl(P(n,4),P(n-2,5)\bigr),\ \bigl(P(n,6),P(n-2,7)\bigr)\\
 & \bigl(P(n,1),P(n-1,2)\bigr),\ \bigl(P(n,5),P(n-1,4)\bigr),\ \bigl(P(n,7),P(n-1,6)\bigr),\ n>2\\
\end{align*}
\end{fleqn}
\subsubsection*{Modules of the form $P(n,4)$}

Defect: $\partial P(n,4) = -2$, for $n\ge 0$.

\begin{fleqn}
\begin{align*}
P(0,4):\ & \bigl(I(1,5),P(0,3)\bigr)\\
P(1,4):\ & \bigl(R_{0}^{1}(1),P(0,2)\bigr),\ \bigl(P(1,5),P(0,5)\bigr),\ \bigl(R_{1}^{2}(1),P(0,6)\bigr),\ \bigl(I(0,5),P(1,3)\bigr)\\
P(2,4):\ & \bigl(P(3,7),P(0,1)\bigr),\ \bigl(P(3,1),P(0,7)\bigr),\ \bigl(R_{0}^{2}(1),P(1,2)\bigr),\ \bigl(P(2,5),P(1,5)\bigr),\ \bigl(R_{1}^{3}(1),P(1,6)\bigr)\\
P(3,4):\ & \bigl(P(5,5),P(0,5)\bigr),\ \bigl(P(4,7),P(1,1)\bigr),\ \bigl(P(4,1),P(1,7)\bigr),\ \bigl(R_{0}^{3}(1),P(2,2)\bigr),\ \bigl(P(3,5),P(2,5)\bigr)\\
 & \bigl(R_{1}^{1}(1),P(2,6)\bigr)\\
P(n,4):\ & \bigl(P(n+2,5),P(n-3,5)\bigr),\ \bigl(P(n+1,7),P(n-2,1)\bigr),\ \bigl(P(n+1,1),P(n-2,7)\bigr)\\
 & \bigl(R_{0}^{(n-1)\bmod 3+1}(1),P(n-1,2)\bigr),\ \bigl(P(n,5),P(n-1,5)\bigr),\ \bigl(R_{1}^{(n-3)\bmod 3+1}(1),P(n-1,6)\bigr),\ n>3\\
\end{align*}
\end{fleqn}
\subsubsection*{Modules of the form $P(n,5)$}

Defect: $\partial P(n,5) = -1$, for $n\ge 0$.

\begin{fleqn}
\begin{align*}
P(0,5):\ & \bigl(I(0,4),P(0,3)\bigr),\ \bigl(I(0,5),P(0,4)\bigr)\\
P(1,5):\ & \bigl(I(1,7),P(0,2)\bigr),\ \bigl(I(1,1),P(0,6)\bigr)\\
P(2,5):\ & \bigl(R_{1}^{1}(1),P(0,1)\bigr),\ \bigl(I(2,5),P(0,4)\bigr),\ \bigl(R_{0}^{3}(1),P(0,7)\bigr),\ \bigl(I(0,7),P(1,2)\bigr),\ \bigl(I(0,1),P(1,6)\bigr)\\
P(3,5):\ & \bigl(R_{\infty}^{1}(1),P(0,5)\bigr),\ \bigl(R_{1}^{2}(1),P(1,1)\bigr),\ \bigl(I(1,5),P(1,4)\bigr),\ \bigl(R_{0}^{1}(1),P(1,7)\bigr)\\
P(4,5):\ & \bigl(R_{0}^{1}(2),P(0,1)\bigr),\ \bigl(R_{1}^{2}(2),P(0,7)\bigr),\ \bigl(R_{\infty}^{2}(1),P(1,5)\bigr),\ \bigl(R_{1}^{3}(1),P(2,1)\bigr),\ \bigl(I(0,5),P(2,4)\bigr)\\
 & \bigl(R_{0}^{2}(1),P(2,7)\bigr)\\
P(n,5):\ & \bigl(R_{0}^{(n-4)\bmod 3+1}(2),P(n-4,1)\bigr),\ \bigl(R_{1}^{(n-3)\bmod 3+1}(2),P(n-4,7)\bigr),\ \bigl(R_{\infty}^{(n-3)\bmod 2+1}(1),P(n-3,5)\bigr)\\
 & \bigl(R_{1}^{(n-2)\bmod 3+1}(1),P(n-2,1)\bigr),\ \bigl(I(-n+4,5),P(n-2,4)\bigr),\ \bigl(R_{0}^{(n-3)\bmod 3+1}(1),P(n-2,7)\bigr),\ n>4\\
\end{align*}
\end{fleqn}
\subsubsection*{Modules of the form $P(n,6)$}

Defect: $\partial P(n,6) = -2$, for $n\ge 0$.

\begin{fleqn}
\begin{align*}
P(0,6):\ & \bigl(I(1,7),P(0,3)\bigr)\\
P(1,6):\ & \bigl(R_{1}^{1}(1),P(0,2)\bigr),\ \bigl(R_{0}^{2}(1),P(0,4)\bigr),\ \bigl(P(1,7),P(0,7)\bigr),\ \bigl(I(0,7),P(1,3)\bigr)\\
P(2,6):\ & \bigl(P(3,5),P(0,1)\bigr),\ \bigl(P(3,1),P(0,5)\bigr),\ \bigl(R_{1}^{2}(1),P(1,2)\bigr),\ \bigl(R_{0}^{3}(1),P(1,4)\bigr),\ \bigl(P(2,7),P(1,7)\bigr)\\
P(3,6):\ & \bigl(P(5,7),P(0,7)\bigr),\ \bigl(P(4,5),P(1,1)\bigr),\ \bigl(P(4,1),P(1,5)\bigr),\ \bigl(R_{1}^{3}(1),P(2,2)\bigr),\ \bigl(R_{0}^{1}(1),P(2,4)\bigr)\\
 & \bigl(P(3,7),P(2,7)\bigr)\\
P(n,6):\ & \bigl(P(n+2,7),P(n-3,7)\bigr),\ \bigl(P(n+1,5),P(n-2,1)\bigr),\ \bigl(P(n+1,1),P(n-2,5)\bigr)\\
 & \bigl(R_{1}^{(n-1)\bmod 3+1}(1),P(n-1,2)\bigr),\ \bigl(R_{0}^{(n-3)\bmod 3+1}(1),P(n-1,4)\bigr),\ \bigl(P(n,7),P(n-1,7)\bigr),\ n>3\\
\end{align*}
\end{fleqn}
\subsubsection*{Modules of the form $P(n,7)$}

Defect: $\partial P(n,7) = -1$, for $n\ge 0$.

\begin{fleqn}
\begin{align*}
P(0,7):\ & \bigl(I(0,6),P(0,3)\bigr),\ \bigl(I(0,7),P(0,6)\bigr)\\
P(1,7):\ & \bigl(I(1,5),P(0,2)\bigr),\ \bigl(I(1,1),P(0,4)\bigr)\\
P(2,7):\ & \bigl(R_{0}^{1}(1),P(0,1)\bigr),\ \bigl(R_{1}^{3}(1),P(0,5)\bigr),\ \bigl(I(2,7),P(0,6)\bigr),\ \bigl(I(0,5),P(1,2)\bigr),\ \bigl(I(0,1),P(1,4)\bigr)\\
P(3,7):\ & \bigl(R_{\infty}^{1}(1),P(0,7)\bigr),\ \bigl(R_{0}^{2}(1),P(1,1)\bigr),\ \bigl(R_{1}^{1}(1),P(1,5)\bigr),\ \bigl(I(1,7),P(1,6)\bigr)\\
P(4,7):\ & \bigl(R_{1}^{1}(2),P(0,1)\bigr),\ \bigl(R_{0}^{2}(2),P(0,5)\bigr),\ \bigl(R_{\infty}^{2}(1),P(1,7)\bigr),\ \bigl(R_{0}^{3}(1),P(2,1)\bigr),\ \bigl(R_{1}^{2}(1),P(2,5)\bigr)\\
 & \bigl(I(0,7),P(2,6)\bigr)\\
P(n,7):\ & \bigl(R_{1}^{(n-4)\bmod 3+1}(2),P(n-4,1)\bigr),\ \bigl(R_{0}^{(n-3)\bmod 3+1}(2),P(n-4,5)\bigr),\ \bigl(R_{\infty}^{(n-3)\bmod 2+1}(1),P(n-3,7)\bigr)\\
 & \bigl(R_{0}^{(n-2)\bmod 3+1}(1),P(n-2,1)\bigr),\ \bigl(R_{1}^{(n-3)\bmod 3+1}(1),P(n-2,5)\bigr),\ \bigl(I(-n+4,7),P(n-2,6)\bigr),\ n>4\\
\end{align*}
\end{fleqn}
\subsubsection{Schofield pairs associated to preinjective exceptional modules}

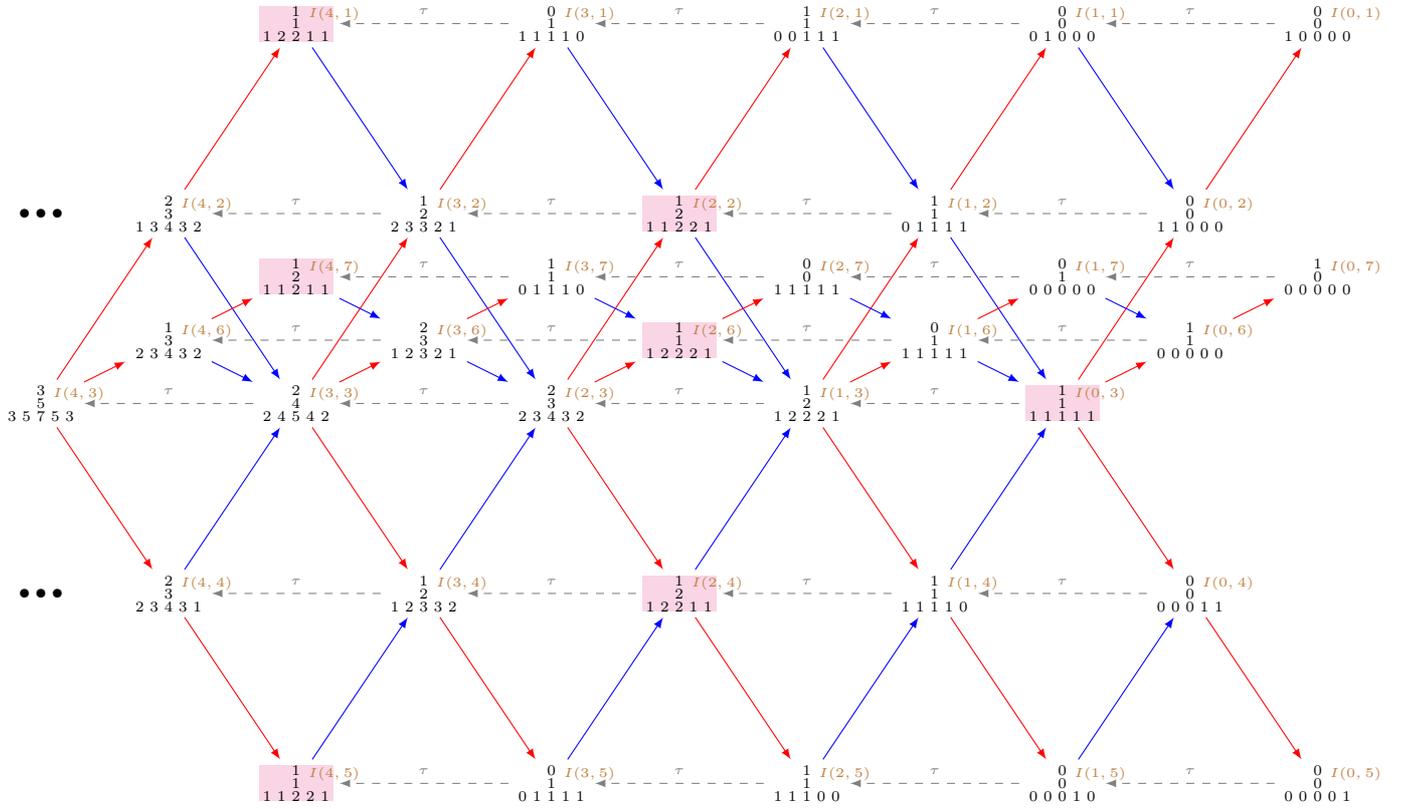
\begin{figure}[ht]


\begin{center}
\begin{scaletikzpicturetowidth}{\textwidth}

\end{scaletikzpicturetowidth}
\end{center}

\end{figure}
\begin{figure}[ht]


\begin{center}
\begin{scaletikzpicturetowidth}{\textwidth}

\end{scaletikzpicturetowidth}
\end{center}

\end{figure}
\medskip{}
\subsubsection*{Modules of the form $I(n,1)$}

Defect: $\partial I(n,1) = 1$, for $n\ge 0$.

\begin{fleqn}
\begin{align*}
I(0,1):\ & - \\
I(1,1):\ & - \\
I(2,1):\ & \bigl(I(0,4),P(0,7)\bigr),\ \bigl(I(0,5),R_{1}^{1}(1)\bigr),\ \bigl(I(0,6),P(0,5)\bigr),\ \bigl(I(0,7),R_{0}^{1}(1)\bigr)\\
I(3,1):\ & \bigl(I(0,1),R_{\infty}^{1}(1)\bigr),\ \bigl(I(0,2),P(1,1)\bigr),\ \bigl(I(1,5),R_{1}^{3}(1)\bigr),\ \bigl(I(1,7),R_{0}^{3}(1)\bigr)\\
I(4,1):\ & \bigl(I(0,5),R_{0}^{2}(2)\bigr),\ \bigl(I(0,7),R_{1}^{2}(2)\bigr),\ \bigl(I(1,1),R_{\infty}^{2}(1)\bigr),\ \bigl(I(1,2),P(0,1)\bigr),\ \bigl(I(2,5),R_{1}^{2}(1)\bigr)\\
 & \bigl(I(2,7),R_{0}^{2}(1)\bigr)\\
I(n,1):\ & \bigl(I(n-4,5),R_{0}^{(-n+5)\bmod 3+1}(2)\bigr),\ \bigl(I(n-4,7),R_{1}^{(-n+5)\bmod 3+1}(2)\bigr),\ \bigl(I(n-3,1),R_{\infty}^{(-n+5)\bmod 2+1}(1)\bigr)\\
 & \bigl(I(n-3,2),P(-n+4,1)\bigr),\ \bigl(I(n-2,5),R_{1}^{(-n+5)\bmod 3+1}(1)\bigr),\ \bigl(I(n-2,7),R_{0}^{(-n+5)\bmod 3+1}(1)\bigr),\ n>4\\
\end{align*}
\end{fleqn}
\subsubsection*{Modules of the form $I(n,2)$}

Defect: $\partial I(n,2) = 2$, for $n\ge 0$.

\begin{fleqn}
\begin{align*}
I(0,2):\ & \bigl(I(0,1),I(1,1)\bigr)\\
I(1,2):\ & \bigl(I(0,4),R_{0}^{2}(1)\bigr),\ \bigl(I(0,5),I(3,7)\bigr),\ \bigl(I(0,6),R_{1}^{2}(1)\bigr),\ \bigl(I(0,7),I(3,5)\bigr),\ \bigl(I(1,1),I(2,1)\bigr)\\
I(2,2):\ & \bigl(I(0,1),I(5,1)\bigr),\ \bigl(I(1,4),R_{0}^{1}(1)\bigr),\ \bigl(I(1,5),I(4,7)\bigr),\ \bigl(I(1,6),R_{1}^{1}(1)\bigr),\ \bigl(I(1,7),I(4,5)\bigr)\\
 & \bigl(I(2,1),I(3,1)\bigr)\\
I(n,2):\ & \bigl(I(n-2,1),I(n+3,1)\bigr),\ \bigl(I(n-1,4),R_{0}^{(-n+2)\bmod 3+1}(1)\bigr),\ \bigl(I(n-1,5),I(n+2,7)\bigr)\\
 & \bigl(I(n-1,6),R_{1}^{(-n+2)\bmod 3+1}(1)\bigr),\ \bigl(I(n-1,7),I(n+2,5)\bigr),\ \bigl(I(n,1),I(n+1,1)\bigr),\ n>2\\
\end{align*}
\end{fleqn}
\subsubsection*{Modules of the form $I(n,3)$}

Defect: $\partial I(n,3) = 3$, for $n\ge 0$.

\begin{fleqn}
\begin{align*}
I(0,3):\ & \bigl(I(0,1),I(1,2)\bigr),\ \bigl(I(0,2),I(2,1)\bigr),\ \bigl(I(0,4),I(2,5)\bigr),\ \bigl(I(0,5),I(1,4)\bigr),\ \bigl(I(0,6),I(2,7)\bigr)\\
 & \bigl(I(0,7),I(1,6)\bigr)\\
I(n,3):\ & \bigl(I(n,1),I(n+1,2)\bigr),\ \bigl(I(n,2),I(n+2,1)\bigr),\ \bigl(I(n,4),I(n+2,5)\bigr)\\
 & \bigl(I(n,5),I(n+1,4)\bigr),\ \bigl(I(n,6),I(n+2,7)\bigr),\ \bigl(I(n,7),I(n+1,6)\bigr),\ n>0\\
\end{align*}
\end{fleqn}
\subsubsection*{Modules of the form $I(n,4)$}

Defect: $\partial I(n,4) = 2$, for $n\ge 0$.

\begin{fleqn}
\begin{align*}
I(0,4):\ & \bigl(I(0,5),I(1,5)\bigr)\\
I(1,4):\ & \bigl(I(0,1),I(3,7)\bigr),\ \bigl(I(0,2),R_{1}^{1}(1)\bigr),\ \bigl(I(0,6),R_{0}^{3}(1)\bigr),\ \bigl(I(0,7),I(3,1)\bigr),\ \bigl(I(1,5),I(2,5)\bigr)\\
I(2,4):\ & \bigl(I(0,5),I(5,5)\bigr),\ \bigl(I(1,1),I(4,7)\bigr),\ \bigl(I(1,2),R_{1}^{3}(1)\bigr),\ \bigl(I(1,6),R_{0}^{2}(1)\bigr),\ \bigl(I(1,7),I(4,1)\bigr)\\
 & \bigl(I(2,5),I(3,5)\bigr)\\
I(n,4):\ & \bigl(I(n-2,5),I(n+3,5)\bigr),\ \bigl(I(n-1,1),I(n+2,7)\bigr),\ \bigl(I(n-1,2),R_{1}^{(-n+4)\bmod 3+1}(1)\bigr)\\
 & \bigl(I(n-1,6),R_{0}^{(-n+3)\bmod 3+1}(1)\bigr),\ \bigl(I(n-1,7),I(n+2,1)\bigr),\ \bigl(I(n,5),I(n+1,5)\bigr),\ n>2\\
\end{align*}
\end{fleqn}
\subsubsection*{Modules of the form $I(n,5)$}

Defect: $\partial I(n,5) = 1$, for $n\ge 0$.

\begin{fleqn}
\begin{align*}
I(0,5):\ & - \\
I(1,5):\ & - \\
I(2,5):\ & \bigl(I(0,1),R_{0}^{2}(1)\bigr),\ \bigl(I(0,2),P(0,7)\bigr),\ \bigl(I(0,6),P(0,1)\bigr),\ \bigl(I(0,7),R_{1}^{3}(1)\bigr)\\
I(3,5):\ & \bigl(I(0,4),P(1,5)\bigr),\ \bigl(I(0,5),R_{\infty}^{1}(1)\bigr),\ \bigl(I(1,1),R_{0}^{1}(1)\bigr),\ \bigl(I(1,7),R_{1}^{2}(1)\bigr)\\
I(4,5):\ & \bigl(I(0,1),R_{1}^{1}(2)\bigr),\ \bigl(I(0,7),R_{0}^{3}(2)\bigr),\ \bigl(I(1,4),P(0,5)\bigr),\ \bigl(I(1,5),R_{\infty}^{2}(1)\bigr),\ \bigl(I(2,1),R_{0}^{3}(1)\bigr)\\
 & \bigl(I(2,7),R_{1}^{1}(1)\bigr)\\
I(n,5):\ & \bigl(I(n-4,1),R_{1}^{(-n+4)\bmod 3+1}(2)\bigr),\ \bigl(I(n-4,7),R_{0}^{(-n+6)\bmod 3+1}(2)\bigr),\ \bigl(I(n-3,4),P(-n+4,5)\bigr)\\
 & \bigl(I(n-3,5),R_{\infty}^{(-n+5)\bmod 2+1}(1)\bigr),\ \bigl(I(n-2,1),R_{0}^{(-n+6)\bmod 3+1}(1)\bigr),\ \bigl(I(n-2,7),R_{1}^{(-n+4)\bmod 3+1}(1)\bigr),\ n>4\\
\end{align*}
\end{fleqn}
\subsubsection*{Modules of the form $I(n,6)$}

Defect: $\partial I(n,6) = 2$, for $n\ge 0$.

\begin{fleqn}
\begin{align*}
I(0,6):\ & \bigl(I(0,7),I(1,7)\bigr)\\
I(1,6):\ & \bigl(I(0,1),I(3,5)\bigr),\ \bigl(I(0,2),R_{0}^{1}(1)\bigr),\ \bigl(I(0,4),R_{1}^{3}(1)\bigr),\ \bigl(I(0,5),I(3,1)\bigr),\ \bigl(I(1,7),I(2,7)\bigr)\\
I(2,6):\ & \bigl(I(0,7),I(5,7)\bigr),\ \bigl(I(1,1),I(4,5)\bigr),\ \bigl(I(1,2),R_{0}^{3}(1)\bigr),\ \bigl(I(1,4),R_{1}^{2}(1)\bigr),\ \bigl(I(1,5),I(4,1)\bigr)\\
 & \bigl(I(2,7),I(3,7)\bigr)\\
I(n,6):\ & \bigl(I(n-2,7),I(n+3,7)\bigr),\ \bigl(I(n-1,1),I(n+2,5)\bigr),\ \bigl(I(n-1,2),R_{0}^{(-n+4)\bmod 3+1}(1)\bigr)\\
 & \bigl(I(n-1,4),R_{1}^{(-n+3)\bmod 3+1}(1)\bigr),\ \bigl(I(n-1,5),I(n+2,1)\bigr),\ \bigl(I(n,7),I(n+1,7)\bigr),\ n>2\\
\end{align*}
\end{fleqn}
\subsubsection*{Modules of the form $I(n,7)$}

Defect: $\partial I(n,7) = 1$, for $n\ge 0$.

\begin{fleqn}
\begin{align*}
I(0,7):\ & - \\
I(1,7):\ & - \\
I(2,7):\ & \bigl(I(0,1),R_{1}^{2}(1)\bigr),\ \bigl(I(0,2),P(0,5)\bigr),\ \bigl(I(0,4),P(0,1)\bigr),\ \bigl(I(0,5),R_{0}^{3}(1)\bigr)\\
I(3,7):\ & \bigl(I(0,6),P(1,7)\bigr),\ \bigl(I(0,7),R_{\infty}^{1}(1)\bigr),\ \bigl(I(1,1),R_{1}^{1}(1)\bigr),\ \bigl(I(1,5),R_{0}^{2}(1)\bigr)\\
I(4,7):\ & \bigl(I(0,1),R_{0}^{1}(2)\bigr),\ \bigl(I(0,5),R_{1}^{3}(2)\bigr),\ \bigl(I(1,6),P(0,7)\bigr),\ \bigl(I(1,7),R_{\infty}^{2}(1)\bigr),\ \bigl(I(2,1),R_{1}^{3}(1)\bigr)\\
 & \bigl(I(2,5),R_{0}^{1}(1)\bigr)\\
I(n,7):\ & \bigl(I(n-4,1),R_{0}^{(-n+4)\bmod 3+1}(2)\bigr),\ \bigl(I(n-4,5),R_{1}^{(-n+6)\bmod 3+1}(2)\bigr),\ \bigl(I(n-3,6),P(-n+4,7)\bigr)\\
 & \bigl(I(n-3,7),R_{\infty}^{(-n+5)\bmod 2+1}(1)\bigr),\ \bigl(I(n-2,1),R_{1}^{(-n+6)\bmod 3+1}(1)\bigr),\ \bigl(I(n-2,5),R_{0}^{(-n+4)\bmod 3+1}(1)\bigr),\ n>4\\
\end{align*}
\end{fleqn}
\subsubsection{Schofield pairs associated to regular exceptional modules}

\subsubsection*{The non-homogeneous tube $\mathcal{T}_{1}^{\Delta(\widetilde{\mathbb{E}}_{6})}$}

\begin{figure}[ht]



\begin{center}
\captionof{figure}{\vspace*{-10pt}$\mathcal{T}_{1}^{\Delta(\widetilde{\mathbb{E}}_{6})}$}
\begin{scaletikzpicturetowidth}{0.69999999999999996\textwidth}

\end{scaletikzpicturetowidth}
\end{center}

\end{figure}
\begin{fleqn}
\begin{align*}
R_{1}^{1}(1):\ & \bigl(I(0,6),P(0,4)\bigr),\ \bigl(I(1,5),P(0,7)\bigr),\ \bigl(I(0,7),P(1,1)\bigr)\\
R_{1}^{1}(2):\ & \bigl(R_{1}^{2}(1),R_{1}^{1}(1)\bigr),\ \bigl(I(3,7),P(0,5)\bigr),\ \bigl(I(2,1),P(1,7)\bigr),\ \bigl(I(1,5),P(2,1)\bigr),\ \bigl(I(0,7),P(3,5)\bigr)\\
R_{1}^{2}(1):\ & \bigl(I(0,4),P(0,2)\bigr),\ \bigl(I(1,1),P(0,5)\bigr),\ \bigl(I(0,5),P(1,7)\bigr)\\
R_{1}^{2}(2):\ & \bigl(R_{1}^{3}(1),R_{1}^{2}(1)\bigr),\ \bigl(I(3,5),P(0,1)\bigr),\ \bigl(I(2,7),P(1,5)\bigr),\ \bigl(I(1,1),P(2,7)\bigr),\ \bigl(I(0,5),P(3,1)\bigr)\\
R_{1}^{3}(1):\ & \bigl(I(1,7),P(0,1)\bigr),\ \bigl(I(0,2),P(0,6)\bigr),\ \bigl(I(0,1),P(1,5)\bigr)\\
R_{1}^{3}(2):\ & \bigl(R_{1}^{1}(1),R_{1}^{3}(1)\bigr),\ \bigl(I(3,1),P(0,7)\bigr),\ \bigl(I(2,5),P(1,1)\bigr),\ \bigl(I(1,7),P(2,5)\bigr),\ \bigl(I(0,1),P(3,7)\bigr)\\
\end{align*}
\end{fleqn}
\subsubsection*{The non-homogeneous tube $\mathcal{T}_{0}^{\Delta(\widetilde{\mathbb{E}}_{6})}$}

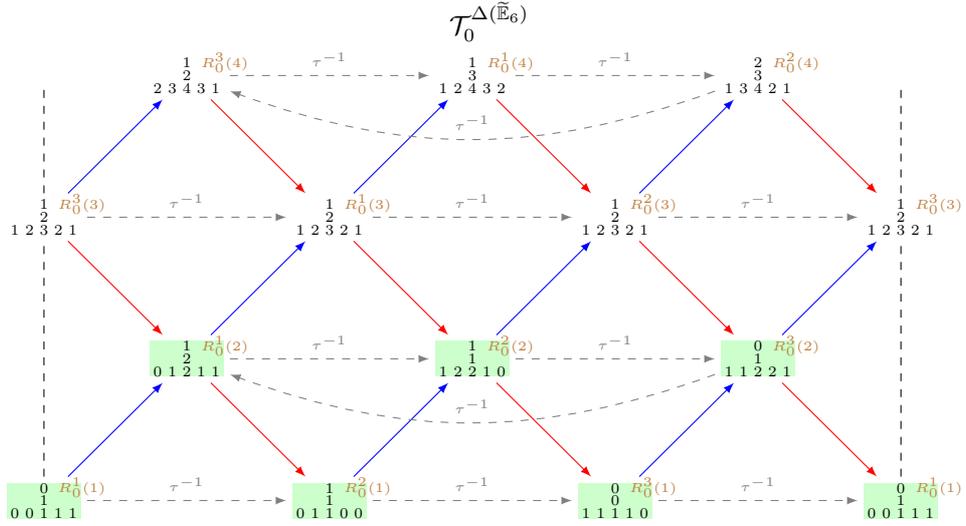
\begin{figure}[ht]



\begin{center}
\captionof{figure}{\vspace*{-10pt}$\mathcal{T}_{0}^{\Delta(\widetilde{\mathbb{E}}_{6})}$}
\begin{scaletikzpicturetowidth}{0.69999999999999996\textwidth}

\end{scaletikzpicturetowidth}
\end{center}

\end{figure}
\begin{fleqn}
\begin{align*}
R_{0}^{1}(1):\ & \bigl(I(1,7),P(0,5)\bigr),\ \bigl(I(0,4),P(0,6)\bigr),\ \bigl(I(0,5),P(1,1)\bigr)\\
R_{0}^{1}(2):\ & \bigl(R_{0}^{2}(1),R_{0}^{1}(1)\bigr),\ \bigl(I(3,5),P(0,7)\bigr),\ \bigl(I(2,1),P(1,5)\bigr),\ \bigl(I(1,7),P(2,1)\bigr),\ \bigl(I(0,5),P(3,7)\bigr)\\
R_{0}^{2}(1):\ & \bigl(I(0,6),P(0,2)\bigr),\ \bigl(I(1,1),P(0,7)\bigr),\ \bigl(I(0,7),P(1,5)\bigr)\\
R_{0}^{2}(2):\ & \bigl(R_{0}^{3}(1),R_{0}^{2}(1)\bigr),\ \bigl(I(3,7),P(0,1)\bigr),\ \bigl(I(2,5),P(1,7)\bigr),\ \bigl(I(1,1),P(2,5)\bigr),\ \bigl(I(0,7),P(3,1)\bigr)\\
R_{0}^{3}(1):\ & \bigl(I(1,5),P(0,1)\bigr),\ \bigl(I(0,2),P(0,4)\bigr),\ \bigl(I(0,1),P(1,7)\bigr)\\
R_{0}^{3}(2):\ & \bigl(R_{0}^{1}(1),R_{0}^{3}(1)\bigr),\ \bigl(I(3,1),P(0,5)\bigr),\ \bigl(I(2,7),P(1,1)\bigr),\ \bigl(I(1,5),P(2,7)\bigr),\ \bigl(I(0,1),P(3,5)\bigr)\\
\end{align*}
\end{fleqn}
\subsubsection*{The non-homogeneous tube $\mathcal{T}_{\infty}^{\Delta(\widetilde{\mathbb{E}}_{6})}$}

\begin{figure}[ht]



\begin{center}
\captionof{figure}{\vspace*{-10pt}$\mathcal{T}_{\infty}^{\Delta(\widetilde{\mathbb{E}}_{6})}$}
\begin{scaletikzpicturetowidth}{0.59999999999999998\textwidth}
 $\end{lrbox}
\subsection[Schofield pairs for the quiver $\Delta(\E_7)$]{Schofield pairs for the quiver $\Delta(\E_7)$ -- {\usebox{\boxEEEEEEEE}}}
\[
\xymatrix{  &   &   & 8\ar[d] &   &   &  \\
1\ar[r] & 2\ar[r] & 3\ar[r] & 4 & 5\ar[l] & 6\ar[l] & 7\ar[l]}
\]
\medskip{}
\[
 C_{\Delta(\widetilde{\mathbb{E}}_{7})} = %
\]

\medskip{}\subsubsection{Schofield pairs associated to preprojective exceptional modules}

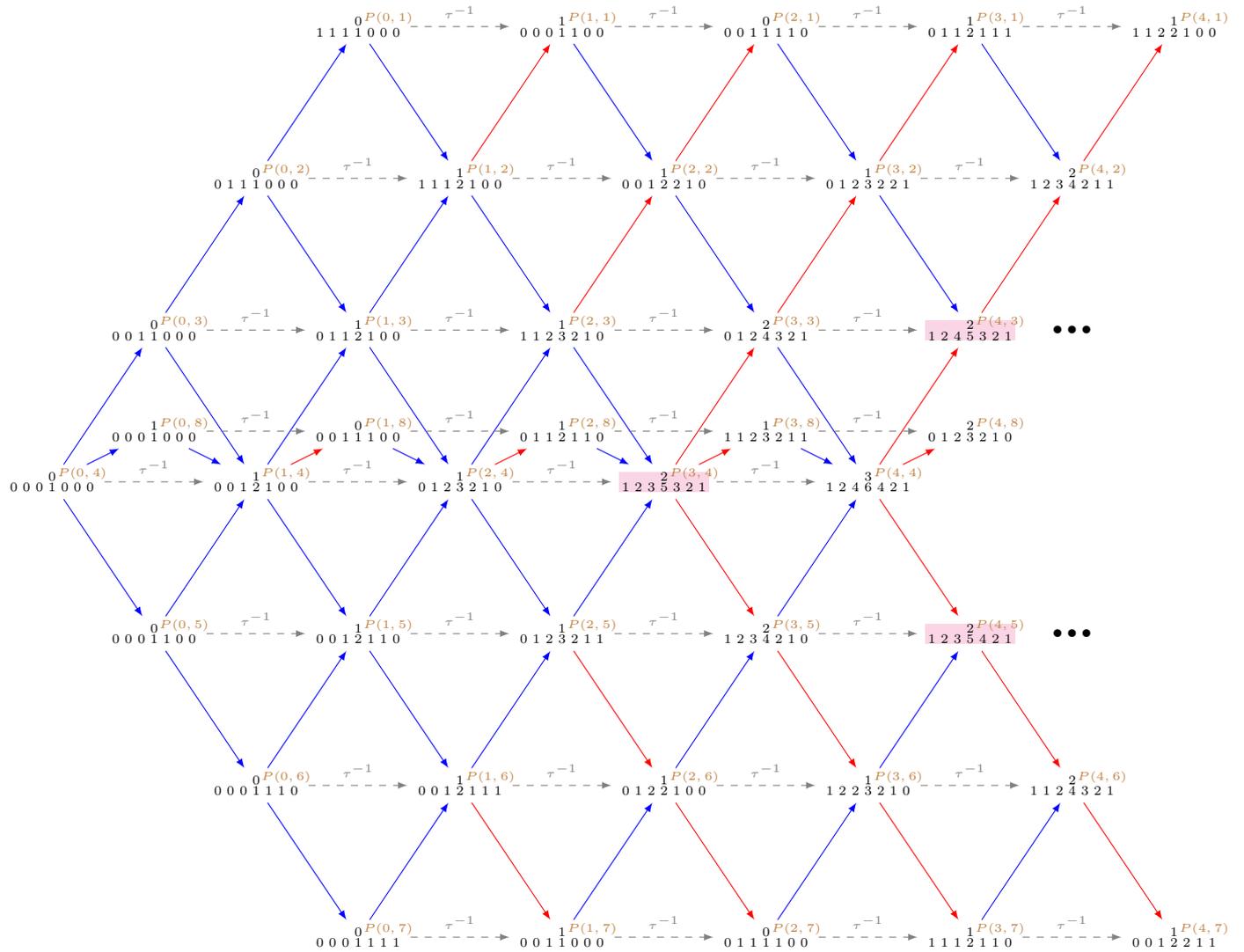
\begin{figure}[ht]


\begin{center}
\begin{scaletikzpicturetowidth}{\textwidth}

\end{scaletikzpicturetowidth}
\end{center}


\end{figure}
\begin{figure}[ht]


\begin{center}
\begin{scaletikzpicturetowidth}{\textwidth}

\end{scaletikzpicturetowidth}
\end{center}


\end{figure}
\medskip{}
\subsubsection*{Modules of the form $P(n,1)$}

Defect: $\partial P(n,1) = -1$, for $n\ge 0$.

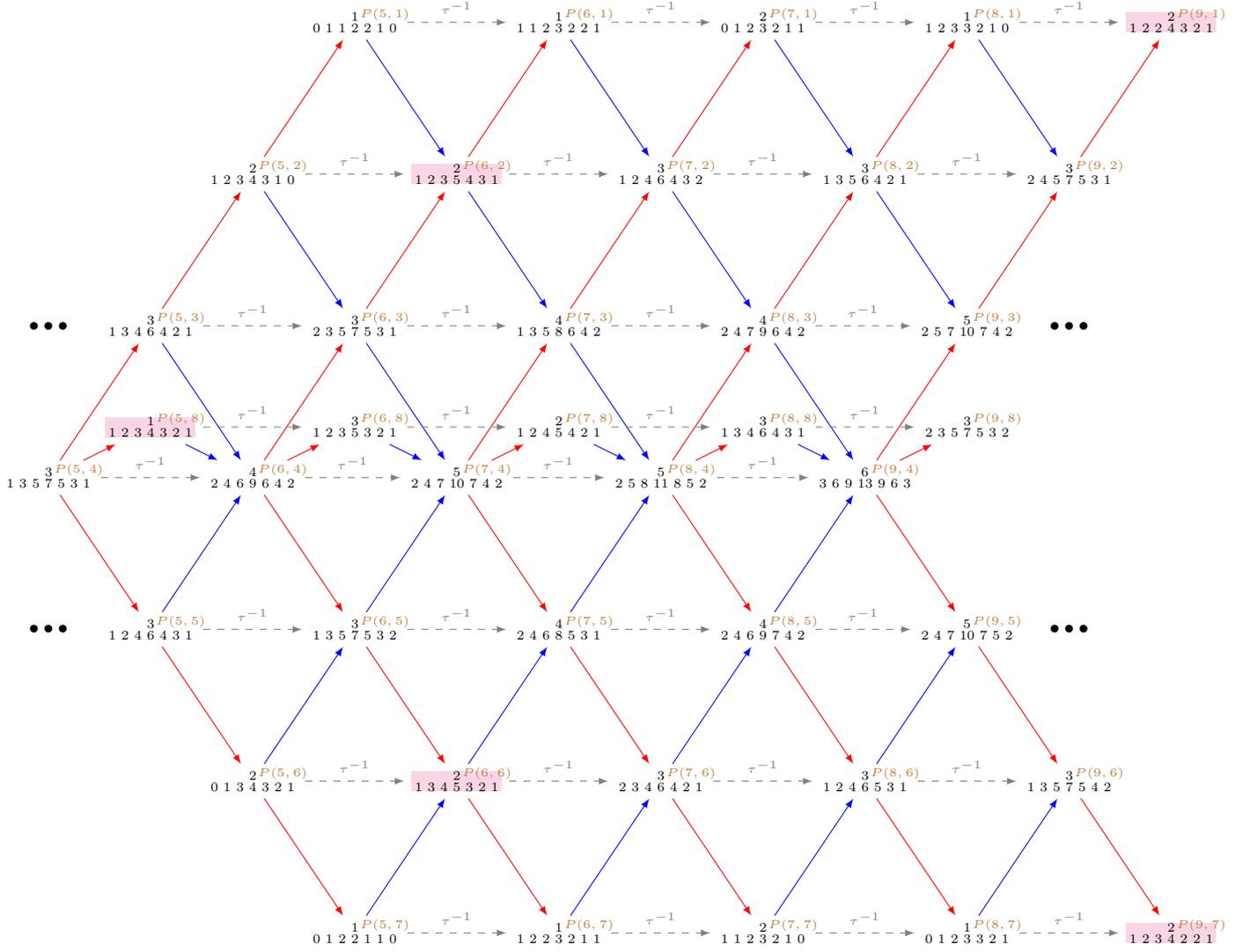
\begin{fleqn}
\begin{align*}
P(0,1):\ & \bigl(I(0,1),P(0,2)\bigr),\ \bigl(I(0,2),P(0,3)\bigr),\ \bigl(I(0,3),P(0,4)\bigr)\\
P(1,1):\ & \bigl(I(0,8),P(0,5)\bigr),\ \bigl(I(2,7),P(0,8)\bigr)\\
P(2,1):\ & \bigl(I(1,6),P(0,3)\bigr),\ \bigl(I(2,1),P(0,6)\bigr),\ \bigl(I(1,7),P(1,8)\bigr)\\
P(3,1):\ & \bigl(I(3,1),P(0,2)\bigr),\ \bigl(R_{1}^{3}(1),P(0,7)\bigr),\ \bigl(I(4,7),P(0,8)\bigr),\ \bigl(I(0,6),P(1,3)\bigr),\ \bigl(I(1,1),P(1,6)\bigr)\\
 & \bigl(I(0,7),P(2,8)\bigr)\\
P(4,1):\ & \bigl(R_{0}^{1}(1),P(0,1)\bigr),\ \bigl(I(2,1),P(1,2)\bigr),\ \bigl(R_{1}^{4}(1),P(1,7)\bigr),\ \bigl(I(3,7),P(1,8)\bigr),\ \bigl(I(0,1),P(2,6)\bigr)\\
P(5,1):\ & \bigl(I(5,1),P(0,6)\bigr),\ \bigl(R_{0}^{2}(1),P(1,1)\bigr),\ \bigl(I(1,1),P(2,2)\bigr),\ \bigl(R_{1}^{1}(1),P(2,7)\bigr),\ \bigl(I(2,7),P(2,8)\bigr)\\
P(6,1):\ & \bigl(R_{1}^{1}(2),P(0,1)\bigr),\ \bigl(R_{\infty}^{2}(1),P(0,7)\bigr),\ \bigl(I(4,1),P(1,6)\bigr),\ \bigl(R_{0}^{3}(1),P(2,1)\bigr),\ \bigl(I(0,1),P(3,2)\bigr)\\
 & \bigl(R_{1}^{2}(1),P(3,7)\bigr),\ \bigl(I(1,7),P(3,8)\bigr)\\
P(7,1):\ & \bigl(R_{1}^{2}(2),P(1,1)\bigr),\ \bigl(R_{\infty}^{1}(1),P(1,7)\bigr),\ \bigl(I(3,1),P(2,6)\bigr),\ \bigl(R_{0}^{1}(1),P(3,1)\bigr),\ \bigl(R_{1}^{3}(1),P(4,7)\bigr)\\
 & \bigl(I(0,7),P(4,8)\bigr)\\
P(8,1):\ & \bigl(R_{0}^{1}(2),P(0,1)\bigr),\ \bigl(R_{1}^{3}(2),P(2,1)\bigr),\ \bigl(R_{\infty}^{2}(1),P(2,7)\bigr),\ \bigl(I(2,1),P(3,6)\bigr),\ \bigl(R_{0}^{2}(1),P(4,1)\bigr)\\
 & \bigl(R_{1}^{4}(1),P(5,7)\bigr)\\
P(9,1):\ & \bigl(R_{1}^{3}(3),P(0,7)\bigr),\ \bigl(R_{0}^{2}(2),P(1,1)\bigr),\ \bigl(R_{1}^{4}(2),P(3,1)\bigr),\ \bigl(R_{\infty}^{1}(1),P(3,7)\bigr),\ \bigl(I(1,1),P(4,6)\bigr)\\
 & \bigl(R_{0}^{3}(1),P(5,1)\bigr),\ \bigl(R_{1}^{1}(1),P(6,7)\bigr)\\
P(n,1):\ & \bigl(R_{1}^{(n-7)\bmod 4+1}(3),P(n-9,7)\bigr),\ \bigl(R_{0}^{(n-8)\bmod 3+1}(2),P(n-8,1)\bigr),\ \bigl(R_{1}^{(n-6)\bmod 4+1}(2),P(n-6,1)\bigr)\\
 & \bigl(R_{\infty}^{(n-9)\bmod 2+1}(1),P(n-6,7)\bigr),\ \bigl(I(-n+10,1),P(n-5,6)\bigr),\ \bigl(R_{0}^{(n-7)\bmod 3+1}(1),P(n-4,1)\bigr)\\
 & \bigl(R_{1}^{(n-9)\bmod 4+1}(1),P(n-3,7)\bigr),\ n>9\\
\end{align*}
\end{fleqn}
\subsubsection*{Modules of the form $P(n,2)$}

Defect: $\partial P(n,2) = -2$, for $n\ge 0$.

\begin{fleqn}
\begin{align*}
P(0,2):\ & \bigl(I(1,1),P(0,3)\bigr),\ \bigl(I(1,2),P(0,4)\bigr)\\
P(1,2):\ & \bigl(P(1,1),P(0,1)\bigr),\ \bigl(I(3,7),P(0,5)\bigr),\ \bigl(R_{1}^{4}(1),P(0,8)\bigr),\ \bigl(I(0,1),P(1,3)\bigr),\ \bigl(I(0,2),P(1,4)\bigr)\\
P(2,2):\ & \bigl(R_{0}^{1}(1),P(0,6)\bigr),\ \bigl(P(2,1),P(1,1)\bigr),\ \bigl(I(2,7),P(1,5)\bigr),\ \bigl(R_{1}^{1}(1),P(1,8)\bigr)\\
P(3,2):\ & \bigl(R_{1}^{1}(2),P(0,2)\bigr),\ \bigl(P(5,7),P(0,7)\bigr),\ \bigl(R_{0}^{2}(1),P(1,6)\bigr),\ \bigl(P(3,1),P(2,1)\bigr),\ \bigl(I(1,7),P(2,5)\bigr)\\
 & \bigl(R_{1}^{2}(1),P(2,8)\bigr)\\
P(4,2):\ & \bigl(P(7,1),P(0,1)\bigr),\ \bigl(R_{1}^{2}(2),P(1,2)\bigr),\ \bigl(P(6,7),P(1,7)\bigr),\ \bigl(R_{0}^{3}(1),P(2,6)\bigr),\ \bigl(P(4,1),P(3,1)\bigr)\\
 & \bigl(I(0,7),P(3,5)\bigr),\ \bigl(R_{1}^{3}(1),P(3,8)\bigr)\\
P(5,2):\ & \bigl(P(8,1),P(1,1)\bigr),\ \bigl(R_{1}^{3}(2),P(2,2)\bigr),\ \bigl(P(7,7),P(2,7)\bigr),\ \bigl(R_{0}^{1}(1),P(3,6)\bigr),\ \bigl(P(5,1),P(4,1)\bigr)\\
 & \bigl(R_{1}^{4}(1),P(4,8)\bigr)\\
P(6,2):\ & \bigl(P(11,7),P(0,7)\bigr),\ \bigl(P(9,1),P(2,1)\bigr),\ \bigl(R_{1}^{4}(2),P(3,2)\bigr),\ \bigl(P(8,7),P(3,7)\bigr),\ \bigl(R_{0}^{2}(1),P(4,6)\bigr)\\
 & \bigl(P(6,1),P(5,1)\bigr),\ \bigl(R_{1}^{1}(1),P(5,8)\bigr)\\
P(n,2):\ & \bigl(P(n+5,7),P(n-6,7)\bigr),\ \bigl(P(n+3,1),P(n-4,1)\bigr),\ \bigl(R_{1}^{(n-3)\bmod 4+1}(2),P(n-3,2)\bigr)\\
 & \bigl(P(n+2,7),P(n-3,7)\bigr),\ \bigl(R_{0}^{(n-5)\bmod 3+1}(1),P(n-2,6)\bigr),\ \bigl(P(n,1),P(n-1,1)\bigr)\\
 & \bigl(R_{1}^{(n-6)\bmod 4+1}(1),P(n-1,8)\bigr),\ n>6\\
\end{align*}
\end{fleqn}
\subsubsection*{Modules of the form $P(n,3)$}

Defect: $\partial P(n,3) = -3$, for $n\ge 0$.

\begin{fleqn}
\begin{align*}
P(0,3):\ & \bigl(I(2,1),P(0,4)\bigr)\\
P(1,3):\ & \bigl(P(1,1),P(0,2)\bigr),\ \bigl(R_{1}^{3}(1),P(0,5)\bigr),\ \bigl(P(2,7),P(0,8)\bigr),\ \bigl(I(1,1),P(1,4)\bigr)\\
P(2,3):\ & \bigl(P(2,2),P(0,1)\bigr),\ \bigl(P(4,1),P(0,6)\bigr),\ \bigl(P(2,1),P(1,2)\bigr),\ \bigl(R_{1}^{4}(1),P(1,5)\bigr),\ \bigl(P(3,7),P(1,8)\bigr)\\
 & \bigl(I(0,1),P(2,4)\bigr)\\
P(3,3):\ & \bigl(P(4,8),P(0,7)\bigr),\ \bigl(P(3,2),P(1,1)\bigr),\ \bigl(P(5,1),P(1,6)\bigr),\ \bigl(P(3,1),P(2,2)\bigr),\ \bigl(R_{1}^{1}(1),P(2,5)\bigr)\\
 & \bigl(P(4,7),P(2,8)\bigr)\\
P(4,3):\ & \bigl(P(5,6),P(0,1)\bigr),\ \bigl(P(5,8),P(1,7)\bigr),\ \bigl(P(4,2),P(2,1)\bigr),\ \bigl(P(6,1),P(2,6)\bigr),\ \bigl(P(4,1),P(3,2)\bigr)\\
 & \bigl(R_{1}^{2}(1),P(3,5)\bigr),\ \bigl(P(5,7),P(3,8)\bigr)\\
P(n,3):\ & \bigl(P(n+1,6),P(n-4,1)\bigr),\ \bigl(P(n+1,8),P(n-3,7)\bigr),\ \bigl(P(n,2),P(n-2,1)\bigr)\\
 & \bigl(P(n+2,1),P(n-2,6)\bigr),\ \bigl(P(n,1),P(n-1,2)\bigr),\ \bigl(R_{1}^{(n-3)\bmod 4+1}(1),P(n-1,5)\bigr)\\
 & \bigl(P(n+1,7),P(n-1,8)\bigr),\ n>4\\
\end{align*}
\end{fleqn}
\subsubsection*{Modules of the form $P(n,4)$}

Defect: $\partial P(n,4) = -4$, for $n\ge 0$.

\begin{fleqn}
\begin{align*}
P(0,4):\ & - \\
P(1,4):\ & \bigl(P(1,1),P(0,3)\bigr),\ \bigl(P(1,7),P(0,5)\bigr),\ \bigl(P(1,8),P(0,8)\bigr)\\
P(2,4):\ & \bigl(P(2,2),P(0,2)\bigr),\ \bigl(P(2,6),P(0,6)\bigr),\ \bigl(P(2,1),P(1,3)\bigr),\ \bigl(P(2,7),P(1,5)\bigr),\ \bigl(P(2,8),P(1,8)\bigr)\\
P(3,4):\ & \bigl(P(3,3),P(0,1)\bigr),\ \bigl(P(3,5),P(0,7)\bigr),\ \bigl(P(3,2),P(1,2)\bigr),\ \bigl(P(3,6),P(1,6)\bigr),\ \bigl(P(3,1),P(2,3)\bigr)\\
 & \bigl(P(3,7),P(2,5)\bigr),\ \bigl(P(3,8),P(2,8)\bigr)\\
P(n,4):\ & \bigl(P(n,3),P(n-3,1)\bigr),\ \bigl(P(n,5),P(n-3,7)\bigr),\ \bigl(P(n,2),P(n-2,2)\bigr)\\
 & \bigl(P(n,6),P(n-2,6)\bigr),\ \bigl(P(n,1),P(n-1,3)\bigr),\ \bigl(P(n,7),P(n-1,5)\bigr)\\
 & \bigl(P(n,8),P(n-1,8)\bigr),\ n>3\\
\end{align*}
\end{fleqn}
\subsubsection*{Modules of the form $P(n,5)$}

Defect: $\partial P(n,5) = -3$, for $n\ge 0$.

\begin{fleqn}
\begin{align*}
P(0,5):\ & \bigl(I(2,7),P(0,4)\bigr)\\
P(1,5):\ & \bigl(R_{1}^{1}(1),P(0,3)\bigr),\ \bigl(P(1,7),P(0,6)\bigr),\ \bigl(P(2,1),P(0,8)\bigr),\ \bigl(I(1,7),P(1,4)\bigr)\\
P(2,5):\ & \bigl(P(4,7),P(0,2)\bigr),\ \bigl(P(2,6),P(0,7)\bigr),\ \bigl(R_{1}^{2}(1),P(1,3)\bigr),\ \bigl(P(2,7),P(1,6)\bigr),\ \bigl(P(3,1),P(1,8)\bigr)\\
 & \bigl(I(0,7),P(2,4)\bigr)\\
P(3,5):\ & \bigl(P(4,8),P(0,1)\bigr),\ \bigl(P(5,7),P(1,2)\bigr),\ \bigl(P(3,6),P(1,7)\bigr),\ \bigl(R_{1}^{3}(1),P(2,3)\bigr),\ \bigl(P(3,7),P(2,6)\bigr)\\
 & \bigl(P(4,1),P(2,8)\bigr)\\
P(4,5):\ & \bigl(P(5,2),P(0,7)\bigr),\ \bigl(P(5,8),P(1,1)\bigr),\ \bigl(P(6,7),P(2,2)\bigr),\ \bigl(P(4,6),P(2,7)\bigr),\ \bigl(R_{1}^{4}(1),P(3,3)\bigr)\\
 & \bigl(P(4,7),P(3,6)\bigr),\ \bigl(P(5,1),P(3,8)\bigr)\\
P(n,5):\ & \bigl(P(n+1,2),P(n-4,7)\bigr),\ \bigl(P(n+1,8),P(n-3,1)\bigr),\ \bigl(P(n+2,7),P(n-2,2)\bigr)\\
 & \bigl(P(n,6),P(n-2,7)\bigr),\ \bigl(R_{1}^{(n-1)\bmod 4+1}(1),P(n-1,3)\bigr),\ \bigl(P(n,7),P(n-1,6)\bigr)\\
 & \bigl(P(n+1,1),P(n-1,8)\bigr),\ n>4\\
\end{align*}
\end{fleqn}
\subsubsection*{Modules of the form $P(n,6)$}

Defect: $\partial P(n,6) = -2$, for $n\ge 0$.

\begin{fleqn}
\begin{align*}
P(0,6):\ & \bigl(I(1,6),P(0,4)\bigr),\ \bigl(I(1,7),P(0,5)\bigr)\\
P(1,6):\ & \bigl(I(3,1),P(0,3)\bigr),\ \bigl(P(1,7),P(0,7)\bigr),\ \bigl(R_{1}^{2}(1),P(0,8)\bigr),\ \bigl(I(0,6),P(1,4)\bigr),\ \bigl(I(0,7),P(1,5)\bigr)\\
P(2,6):\ & \bigl(R_{0}^{1}(1),P(0,2)\bigr),\ \bigl(I(2,1),P(1,3)\bigr),\ \bigl(P(2,7),P(1,7)\bigr),\ \bigl(R_{1}^{3}(1),P(1,8)\bigr)\\
P(3,6):\ & \bigl(P(5,1),P(0,1)\bigr),\ \bigl(R_{1}^{3}(2),P(0,6)\bigr),\ \bigl(R_{0}^{2}(1),P(1,2)\bigr),\ \bigl(I(1,1),P(2,3)\bigr),\ \bigl(P(3,7),P(2,7)\bigr)\\
 & \bigl(R_{1}^{4}(1),P(2,8)\bigr)\\
P(4,6):\ & \bigl(P(7,7),P(0,7)\bigr),\ \bigl(P(6,1),P(1,1)\bigr),\ \bigl(R_{1}^{4}(2),P(1,6)\bigr),\ \bigl(R_{0}^{3}(1),P(2,2)\bigr),\ \bigl(I(0,1),P(3,3)\bigr)\\
 & \bigl(P(4,7),P(3,7)\bigr),\ \bigl(R_{1}^{1}(1),P(3,8)\bigr)\\
P(5,6):\ & \bigl(P(8,7),P(1,7)\bigr),\ \bigl(P(7,1),P(2,1)\bigr),\ \bigl(R_{1}^{1}(2),P(2,6)\bigr),\ \bigl(R_{0}^{1}(1),P(3,2)\bigr),\ \bigl(P(5,7),P(4,7)\bigr)\\
 & \bigl(R_{1}^{2}(1),P(4,8)\bigr)\\
P(6,6):\ & \bigl(P(11,1),P(0,1)\bigr),\ \bigl(P(9,7),P(2,7)\bigr),\ \bigl(P(8,1),P(3,1)\bigr),\ \bigl(R_{1}^{2}(2),P(3,6)\bigr),\ \bigl(R_{0}^{2}(1),P(4,2)\bigr)\\
 & \bigl(P(6,7),P(5,7)\bigr),\ \bigl(R_{1}^{3}(1),P(5,8)\bigr)\\
P(n,6):\ & \bigl(P(n+5,1),P(n-6,1)\bigr),\ \bigl(P(n+3,7),P(n-4,7)\bigr),\ \bigl(P(n+2,1),P(n-3,1)\bigr)\\
 & \bigl(R_{1}^{(n-5)\bmod 4+1}(2),P(n-3,6)\bigr),\ \bigl(R_{0}^{(n-5)\bmod 3+1}(1),P(n-2,2)\bigr),\ \bigl(P(n,7),P(n-1,7)\bigr)\\
 & \bigl(R_{1}^{(n-4)\bmod 4+1}(1),P(n-1,8)\bigr),\ n>6\\
\end{align*}
\end{fleqn}
\subsubsection*{Modules of the form $P(n,7)$}

Defect: $\partial P(n,7) = -1$, for $n\ge 0$.

\begin{fleqn}
\begin{align*}
P(0,7):\ & \bigl(I(0,5),P(0,4)\bigr),\ \bigl(I(0,6),P(0,5)\bigr),\ \bigl(I(0,7),P(0,6)\bigr)\\
P(1,7):\ & \bigl(I(0,8),P(0,3)\bigr),\ \bigl(I(2,1),P(0,8)\bigr)\\
P(2,7):\ & \bigl(I(2,7),P(0,2)\bigr),\ \bigl(I(1,2),P(0,5)\bigr),\ \bigl(I(1,1),P(1,8)\bigr)\\
P(3,7):\ & \bigl(R_{1}^{1}(1),P(0,1)\bigr),\ \bigl(I(3,7),P(0,6)\bigr),\ \bigl(I(4,1),P(0,8)\bigr),\ \bigl(I(1,7),P(1,2)\bigr),\ \bigl(I(0,2),P(1,5)\bigr)\\
 & \bigl(I(0,1),P(2,8)\bigr)\\
P(4,7):\ & \bigl(R_{0}^{1}(1),P(0,7)\bigr),\ \bigl(R_{1}^{2}(1),P(1,1)\bigr),\ \bigl(I(2,7),P(1,6)\bigr),\ \bigl(I(3,1),P(1,8)\bigr),\ \bigl(I(0,7),P(2,2)\bigr)\\
P(5,7):\ & \bigl(I(5,7),P(0,2)\bigr),\ \bigl(R_{0}^{2}(1),P(1,7)\bigr),\ \bigl(R_{1}^{3}(1),P(2,1)\bigr),\ \bigl(I(1,7),P(2,6)\bigr),\ \bigl(I(2,1),P(2,8)\bigr)\\
P(6,7):\ & \bigl(R_{\infty}^{1}(1),P(0,1)\bigr),\ \bigl(R_{1}^{3}(2),P(0,7)\bigr),\ \bigl(I(4,7),P(1,2)\bigr),\ \bigl(R_{0}^{3}(1),P(2,7)\bigr),\ \bigl(R_{1}^{4}(1),P(3,1)\bigr)\\
 & \bigl(I(0,7),P(3,6)\bigr),\ \bigl(I(1,1),P(3,8)\bigr)\\
P(7,7):\ & \bigl(R_{\infty}^{2}(1),P(1,1)\bigr),\ \bigl(R_{1}^{4}(2),P(1,7)\bigr),\ \bigl(I(3,7),P(2,2)\bigr),\ \bigl(R_{0}^{1}(1),P(3,7)\bigr),\ \bigl(R_{1}^{1}(1),P(4,1)\bigr)\\
 & \bigl(I(0,1),P(4,8)\bigr)\\
P(8,7):\ & \bigl(R_{0}^{1}(2),P(0,7)\bigr),\ \bigl(R_{\infty}^{1}(1),P(2,1)\bigr),\ \bigl(R_{1}^{1}(2),P(2,7)\bigr),\ \bigl(I(2,7),P(3,2)\bigr),\ \bigl(R_{0}^{2}(1),P(4,7)\bigr)\\
 & \bigl(R_{1}^{2}(1),P(5,1)\bigr)\\
P(9,7):\ & \bigl(R_{1}^{1}(3),P(0,1)\bigr),\ \bigl(R_{0}^{2}(2),P(1,7)\bigr),\ \bigl(R_{\infty}^{2}(1),P(3,1)\bigr),\ \bigl(R_{1}^{2}(2),P(3,7)\bigr),\ \bigl(I(1,7),P(4,2)\bigr)\\
 & \bigl(R_{0}^{3}(1),P(5,7)\bigr),\ \bigl(R_{1}^{3}(1),P(6,1)\bigr)\\
P(n,7):\ & \bigl(R_{1}^{(n-9)\bmod 4+1}(3),P(n-9,1)\bigr),\ \bigl(R_{0}^{(n-8)\bmod 3+1}(2),P(n-8,7)\bigr),\ \bigl(R_{\infty}^{(n-8)\bmod 2+1}(1),P(n-6,1)\bigr)\\
 & \bigl(R_{1}^{(n-8)\bmod 4+1}(2),P(n-6,7)\bigr),\ \bigl(I(-n+10,7),P(n-5,2)\bigr),\ \bigl(R_{0}^{(n-7)\bmod 3+1}(1),P(n-4,7)\bigr)\\
 & \bigl(R_{1}^{(n-7)\bmod 4+1}(1),P(n-3,1)\bigr),\ n>9\\
\end{align*}
\end{fleqn}
\subsubsection*{Modules of the form $P(n,8)$}

Defect: $\partial P(n,8) = -2$, for $n\ge 0$.

\begin{fleqn}
\begin{align*}
P(0,8):\ & \bigl(I(0,8),P(0,4)\bigr)\\
P(1,8):\ & \bigl(I(2,7),P(0,3)\bigr),\ \bigl(I(2,1),P(0,5)\bigr)\\
P(2,8):\ & \bigl(R_{1}^{1}(1),P(0,2)\bigr),\ \bigl(R_{1}^{3}(1),P(0,6)\bigr),\ \bigl(R_{0}^{2}(1),P(0,8)\bigr),\ \bigl(I(1,7),P(1,3)\bigr),\ \bigl(I(1,1),P(1,5)\bigr)\\
P(3,8):\ & \bigl(P(4,7),P(0,1)\bigr),\ \bigl(P(4,1),P(0,7)\bigr),\ \bigl(R_{1}^{2}(1),P(1,2)\bigr),\ \bigl(R_{1}^{4}(1),P(1,6)\bigr),\ \bigl(R_{0}^{3}(1),P(1,8)\bigr)\\
 & \bigl(I(0,7),P(2,3)\bigr),\ \bigl(I(0,1),P(2,5)\bigr)\\
P(4,8):\ & \bigl(P(5,7),P(1,1)\bigr),\ \bigl(P(5,1),P(1,7)\bigr),\ \bigl(R_{1}^{3}(1),P(2,2)\bigr),\ \bigl(R_{1}^{1}(1),P(2,6)\bigr),\ \bigl(R_{0}^{1}(1),P(2,8)\bigr)\\
P(5,8):\ & \bigl(P(8,7),P(0,1)\bigr),\ \bigl(P(8,1),P(0,7)\bigr),\ \bigl(P(6,7),P(2,1)\bigr),\ \bigl(P(6,1),P(2,7)\bigr),\ \bigl(R_{1}^{4}(1),P(3,2)\bigr)\\
 & \bigl(R_{1}^{2}(1),P(3,6)\bigr),\ \bigl(R_{0}^{2}(1),P(3,8)\bigr)\\
P(n,8):\ & \bigl(P(n+3,7),P(n-5,1)\bigr),\ \bigl(P(n+3,1),P(n-5,7)\bigr),\ \bigl(P(n+1,7),P(n-3,1)\bigr)\\
 & \bigl(P(n+1,1),P(n-3,7)\bigr),\ \bigl(R_{1}^{(n-2)\bmod 4+1}(1),P(n-2,2)\bigr),\ \bigl(R_{1}^{(n-4)\bmod 4+1}(1),P(n-2,6)\bigr)\\
 & \bigl(R_{0}^{(n-4)\bmod 3+1}(1),P(n-2,8)\bigr),\ n>5\\
\end{align*}
\end{fleqn}
\subsubsection{Schofield pairs associated to preinjective exceptional modules}

\begin{figure}[ht]


\begin{center}
\begin{scaletikzpicturetowidth}{\textwidth}

\end{scaletikzpicturetowidth}
\end{center}

\end{figure}
\begin{figure}[ht]


\begin{center}
\begin{scaletikzpicturetowidth}{\textwidth}

\end{scaletikzpicturetowidth}
\end{center}

\end{figure}
\medskip{}
\subsubsection*{Modules of the form $I(n,1)$}

Defect: $\partial I(n,1) = 1$, for $n\ge 0$.

\begin{fleqn}
\begin{align*}
I(0,1):\ & - \\
I(1,1):\ & - \\
I(2,1):\ & - \\
I(3,1):\ & \bigl(I(0,5),P(0,8)\bigr),\ \bigl(I(0,6),P(1,1)\bigr),\ \bigl(I(0,7),R_{1}^{1}(1)\bigr),\ \bigl(I(0,8),P(0,7)\bigr)\\
I(4,1):\ & \bigl(I(0,1),R_{0}^{2}(1)\bigr),\ \bigl(I(0,2),P(2,1)\bigr),\ \bigl(I(0,3),P(0,6)\bigr),\ \bigl(I(1,6),P(0,1)\bigr),\ \bigl(I(1,7),R_{1}^{4}(1)\bigr)\\
I(5,1):\ & \bigl(I(0,8),P(2,7)\bigr),\ \bigl(I(1,1),R_{0}^{1}(1)\bigr),\ \bigl(I(1,2),P(1,1)\bigr),\ \bigl(I(2,7),R_{1}^{3}(1)\bigr)\\
I(6,1):\ & \bigl(I(0,1),R_{1}^{2}(2)\bigr),\ \bigl(I(0,6),P(4,1)\bigr),\ \bigl(I(0,7),R_{\infty}^{2}(1)\bigr),\ \bigl(I(1,8),P(1,7)\bigr),\ \bigl(I(2,1),R_{0}^{3}(1)\bigr)\\
 & \bigl(I(2,2),P(0,1)\bigr),\ \bigl(I(3,7),R_{1}^{2}(1)\bigr)\\
I(7,1):\ & \bigl(I(1,1),R_{1}^{1}(2)\bigr),\ \bigl(I(1,6),P(3,1)\bigr),\ \bigl(I(1,7),R_{\infty}^{1}(1)\bigr),\ \bigl(I(2,8),P(0,7)\bigr),\ \bigl(I(3,1),R_{0}^{2}(1)\bigr)\\
 & \bigl(I(4,7),R_{1}^{1}(1)\bigr)\\
I(8,1):\ & \bigl(I(0,1),R_{0}^{1}(2)\bigr),\ \bigl(I(2,1),R_{1}^{4}(2)\bigr),\ \bigl(I(2,6),P(2,1)\bigr),\ \bigl(I(2,7),R_{\infty}^{2}(1)\bigr),\ \bigl(I(4,1),R_{0}^{1}(1)\bigr)\\
 & \bigl(I(5,7),R_{1}^{4}(1)\bigr)\\
I(9,1):\ & \bigl(I(0,7),R_{1}^{3}(3)\bigr),\ \bigl(I(1,1),R_{0}^{3}(2)\bigr),\ \bigl(I(3,1),R_{1}^{3}(2)\bigr),\ \bigl(I(3,6),P(1,1)\bigr),\ \bigl(I(3,7),R_{\infty}^{1}(1)\bigr)\\
 & \bigl(I(5,1),R_{0}^{3}(1)\bigr),\ \bigl(I(6,7),R_{1}^{3}(1)\bigr)\\
I(n,1):\ & \bigl(I(n-9,7),R_{1}^{(-n+11)\bmod 4+1}(3)\bigr),\ \bigl(I(n-8,1),R_{0}^{(-n+11)\bmod 3+1}(2)\bigr),\ \bigl(I(n-6,1),R_{1}^{(-n+11)\bmod 4+1}(2)\bigr)\\
 & \bigl(I(n-6,6),P(-n+10,1)\bigr),\ \bigl(I(n-6,7),R_{\infty}^{(-n+9)\bmod 2+1}(1)\bigr),\ \bigl(I(n-4,1),R_{0}^{(-n+11)\bmod 3+1}(1)\bigr)\\
 & \bigl(I(n-3,7),R_{1}^{(-n+11)\bmod 4+1}(1)\bigr),\ n>9\\
\end{align*}
\end{fleqn}
\subsubsection*{Modules of the form $I(n,2)$}

Defect: $\partial I(n,2) = 2$, for $n\ge 0$.

\begin{fleqn}
\begin{align*}
I(0,2):\ & \bigl(I(0,1),I(1,1)\bigr)\\
I(1,2):\ & \bigl(I(1,1),I(2,1)\bigr)\\
I(2,2):\ & \bigl(I(0,5),P(1,7)\bigr),\ \bigl(I(0,6),R_{0}^{1}(1)\bigr),\ \bigl(I(0,7),I(5,7)\bigr),\ \bigl(I(0,8),R_{1}^{2}(1)\bigr),\ \bigl(I(2,1),I(3,1)\bigr)\\
I(3,2):\ & \bigl(I(0,1),I(7,1)\bigr),\ \bigl(I(0,2),R_{1}^{1}(2)\bigr),\ \bigl(I(1,5),P(0,7)\bigr),\ \bigl(I(1,6),R_{0}^{3}(1)\bigr),\ \bigl(I(1,7),I(6,7)\bigr)\\
 & \bigl(I(1,8),R_{1}^{1}(1)\bigr),\ \bigl(I(3,1),I(4,1)\bigr)\\
I(4,2):\ & \bigl(I(1,1),I(8,1)\bigr),\ \bigl(I(1,2),R_{1}^{4}(2)\bigr),\ \bigl(I(2,6),R_{0}^{2}(1)\bigr),\ \bigl(I(2,7),I(7,7)\bigr),\ \bigl(I(2,8),R_{1}^{4}(1)\bigr)\\
 & \bigl(I(4,1),I(5,1)\bigr)\\
I(5,2):\ & \bigl(I(0,7),I(11,7)\bigr),\ \bigl(I(2,1),I(9,1)\bigr),\ \bigl(I(2,2),R_{1}^{3}(2)\bigr),\ \bigl(I(3,6),R_{0}^{1}(1)\bigr),\ \bigl(I(3,7),I(8,7)\bigr)\\
 & \bigl(I(3,8),R_{1}^{3}(1)\bigr),\ \bigl(I(5,1),I(6,1)\bigr)\\
I(n,2):\ & \bigl(I(n-5,7),I(n+6,7)\bigr),\ \bigl(I(n-3,1),I(n+4,1)\bigr),\ \bigl(I(n-3,2),R_{1}^{(-n+7)\bmod 4+1}(2)\bigr)\\
 & \bigl(I(n-2,6),R_{0}^{(-n+5)\bmod 3+1}(1)\bigr),\ \bigl(I(n-2,7),I(n+3,7)\bigr),\ \bigl(I(n-2,8),R_{1}^{(-n+7)\bmod 4+1}(1)\bigr)\\
 & \bigl(I(n,1),I(n+1,1)\bigr),\ n>5\\
\end{align*}
\end{fleqn}
\subsubsection*{Modules of the form $I(n,3)$}

Defect: $\partial I(n,3) = 3$, for $n\ge 0$.

\begin{fleqn}
\begin{align*}
I(0,3):\ & \bigl(I(0,1),I(1,2)\bigr),\ \bigl(I(0,2),I(2,1)\bigr)\\
I(1,3):\ & \bigl(I(0,5),R_{1}^{3}(1)\bigr),\ \bigl(I(0,6),I(5,1)\bigr),\ \bigl(I(0,7),I(2,8)\bigr),\ \bigl(I(0,8),I(4,7)\bigr),\ \bigl(I(1,1),I(2,2)\bigr)\\
 & \bigl(I(1,2),I(3,1)\bigr)\\
I(2,3):\ & \bigl(I(0,1),I(4,6)\bigr),\ \bigl(I(1,5),R_{1}^{2}(1)\bigr),\ \bigl(I(1,6),I(6,1)\bigr),\ \bigl(I(1,7),I(3,8)\bigr),\ \bigl(I(1,8),I(5,7)\bigr)\\
 & \bigl(I(2,1),I(3,2)\bigr),\ \bigl(I(2,2),I(4,1)\bigr)\\
I(n,3):\ & \bigl(I(n-2,1),I(n+2,6)\bigr),\ \bigl(I(n-1,5),R_{1}^{(-n+3)\bmod 4+1}(1)\bigr),\ \bigl(I(n-1,6),I(n+4,1)\bigr)\\
 & \bigl(I(n-1,7),I(n+1,8)\bigr),\ \bigl(I(n-1,8),I(n+3,7)\bigr),\ \bigl(I(n,1),I(n+1,2)\bigr)\\
 & \bigl(I(n,2),I(n+2,1)\bigr),\ n>2\\
\end{align*}
\end{fleqn}
\subsubsection*{Modules of the form $I(n,4)$}

Defect: $\partial I(n,4) = 4$, for $n\ge 0$.

\begin{fleqn}
\begin{align*}
I(0,4):\ & \bigl(I(0,1),I(1,3)\bigr),\ \bigl(I(0,2),I(2,2)\bigr),\ \bigl(I(0,3),I(3,1)\bigr),\ \bigl(I(0,5),I(3,7)\bigr),\ \bigl(I(0,6),I(2,6)\bigr)\\
 & \bigl(I(0,7),I(1,5)\bigr),\ \bigl(I(0,8),I(1,8)\bigr)\\
I(n,4):\ & \bigl(I(n,1),I(n+1,3)\bigr),\ \bigl(I(n,2),I(n+2,2)\bigr),\ \bigl(I(n,3),I(n+3,1)\bigr)\\
 & \bigl(I(n,5),I(n+3,7)\bigr),\ \bigl(I(n,6),I(n+2,6)\bigr),\ \bigl(I(n,7),I(n+1,5)\bigr)\\
 & \bigl(I(n,8),I(n+1,8)\bigr),\ n>0\\
\end{align*}
\end{fleqn}
\subsubsection*{Modules of the form $I(n,5)$}

Defect: $\partial I(n,5) = 3$, for $n\ge 0$.

\begin{fleqn}
\begin{align*}
I(0,5):\ & \bigl(I(0,6),I(2,7)\bigr),\ \bigl(I(0,7),I(1,6)\bigr)\\
I(1,5):\ & \bigl(I(0,1),I(2,8)\bigr),\ \bigl(I(0,2),I(5,7)\bigr),\ \bigl(I(0,3),R_{1}^{1}(1)\bigr),\ \bigl(I(0,8),I(4,1)\bigr),\ \bigl(I(1,6),I(3,7)\bigr)\\
 & \bigl(I(1,7),I(2,6)\bigr)\\
I(2,5):\ & \bigl(I(0,7),I(4,2)\bigr),\ \bigl(I(1,1),I(3,8)\bigr),\ \bigl(I(1,2),I(6,7)\bigr),\ \bigl(I(1,3),R_{1}^{4}(1)\bigr),\ \bigl(I(1,8),I(5,1)\bigr)\\
 & \bigl(I(2,6),I(4,7)\bigr),\ \bigl(I(2,7),I(3,6)\bigr)\\
I(n,5):\ & \bigl(I(n-2,7),I(n+2,2)\bigr),\ \bigl(I(n-1,1),I(n+1,8)\bigr),\ \bigl(I(n-1,2),I(n+4,7)\bigr)\\
 & \bigl(I(n-1,3),R_{1}^{(-n+5)\bmod 4+1}(1)\bigr),\ \bigl(I(n-1,8),I(n+3,1)\bigr),\ \bigl(I(n,6),I(n+2,7)\bigr)\\
 & \bigl(I(n,7),I(n+1,6)\bigr),\ n>2\\
\end{align*}
\end{fleqn}
\subsubsection*{Modules of the form $I(n,6)$}

Defect: $\partial I(n,6) = 2$, for $n\ge 0$.

\begin{fleqn}
\begin{align*}
I(0,6):\ & \bigl(I(0,7),I(1,7)\bigr)\\
I(1,6):\ & \bigl(I(1,7),I(2,7)\bigr)\\
I(2,6):\ & \bigl(I(0,1),I(5,1)\bigr),\ \bigl(I(0,2),R_{0}^{1}(1)\bigr),\ \bigl(I(0,3),P(1,1)\bigr),\ \bigl(I(0,8),R_{1}^{4}(1)\bigr),\ \bigl(I(2,7),I(3,7)\bigr)\\
I(3,6):\ & \bigl(I(0,6),R_{1}^{3}(2)\bigr),\ \bigl(I(0,7),I(7,7)\bigr),\ \bigl(I(1,1),I(6,1)\bigr),\ \bigl(I(1,2),R_{0}^{3}(1)\bigr),\ \bigl(I(1,3),P(0,1)\bigr)\\
 & \bigl(I(1,8),R_{1}^{3}(1)\bigr),\ \bigl(I(3,7),I(4,7)\bigr)\\
I(4,6):\ & \bigl(I(1,6),R_{1}^{2}(2)\bigr),\ \bigl(I(1,7),I(8,7)\bigr),\ \bigl(I(2,1),I(7,1)\bigr),\ \bigl(I(2,2),R_{0}^{2}(1)\bigr),\ \bigl(I(2,8),R_{1}^{2}(1)\bigr)\\
 & \bigl(I(4,7),I(5,7)\bigr)\\
I(5,6):\ & \bigl(I(0,1),I(11,1)\bigr),\ \bigl(I(2,6),R_{1}^{1}(2)\bigr),\ \bigl(I(2,7),I(9,7)\bigr),\ \bigl(I(3,1),I(8,1)\bigr),\ \bigl(I(3,2),R_{0}^{1}(1)\bigr)\\
 & \bigl(I(3,8),R_{1}^{1}(1)\bigr),\ \bigl(I(5,7),I(6,7)\bigr)\\
I(n,6):\ & \bigl(I(n-5,1),I(n+6,1)\bigr),\ \bigl(I(n-3,6),R_{1}^{(-n+5)\bmod 4+1}(2)\bigr),\ \bigl(I(n-3,7),I(n+4,7)\bigr)\\
 & \bigl(I(n-2,1),I(n+3,1)\bigr),\ \bigl(I(n-2,2),R_{0}^{(-n+5)\bmod 3+1}(1)\bigr),\ \bigl(I(n-2,8),R_{1}^{(-n+5)\bmod 4+1}(1)\bigr)\\
 & \bigl(I(n,7),I(n+1,7)\bigr),\ n>5\\
\end{align*}
\end{fleqn}
\subsubsection*{Modules of the form $I(n,7)$}

Defect: $\partial I(n,7) = 1$, for $n\ge 0$.

\begin{fleqn}
\begin{align*}
I(0,7):\ & - \\
I(1,7):\ & - \\
I(2,7):\ & - \\
I(3,7):\ & \bigl(I(0,1),R_{1}^{3}(1)\bigr),\ \bigl(I(0,2),P(1,7)\bigr),\ \bigl(I(0,3),P(0,8)\bigr),\ \bigl(I(0,8),P(0,1)\bigr)\\
I(4,7):\ & \bigl(I(0,5),P(0,2)\bigr),\ \bigl(I(0,6),P(2,7)\bigr),\ \bigl(I(0,7),R_{0}^{2}(1)\bigr),\ \bigl(I(1,1),R_{1}^{2}(1)\bigr),\ \bigl(I(1,2),P(0,7)\bigr)\\
I(5,7):\ & \bigl(I(0,8),P(2,1)\bigr),\ \bigl(I(1,6),P(1,7)\bigr),\ \bigl(I(1,7),R_{0}^{1}(1)\bigr),\ \bigl(I(2,1),R_{1}^{1}(1)\bigr)\\
I(6,7):\ & \bigl(I(0,1),R_{\infty}^{1}(1)\bigr),\ \bigl(I(0,2),P(4,7)\bigr),\ \bigl(I(0,7),R_{1}^{4}(2)\bigr),\ \bigl(I(1,8),P(1,1)\bigr),\ \bigl(I(2,6),P(0,7)\bigr)\\
 & \bigl(I(2,7),R_{0}^{3}(1)\bigr),\ \bigl(I(3,1),R_{1}^{4}(1)\bigr)\\
I(7,7):\ & \bigl(I(1,1),R_{\infty}^{2}(1)\bigr),\ \bigl(I(1,2),P(3,7)\bigr),\ \bigl(I(1,7),R_{1}^{3}(2)\bigr),\ \bigl(I(2,8),P(0,1)\bigr),\ \bigl(I(3,7),R_{0}^{2}(1)\bigr)\\
 & \bigl(I(4,1),R_{1}^{3}(1)\bigr)\\
I(8,7):\ & \bigl(I(0,7),R_{0}^{1}(2)\bigr),\ \bigl(I(2,1),R_{\infty}^{1}(1)\bigr),\ \bigl(I(2,2),P(2,7)\bigr),\ \bigl(I(2,7),R_{1}^{2}(2)\bigr),\ \bigl(I(4,7),R_{0}^{1}(1)\bigr)\\
 & \bigl(I(5,1),R_{1}^{2}(1)\bigr)\\
I(9,7):\ & \bigl(I(0,1),R_{1}^{1}(3)\bigr),\ \bigl(I(1,7),R_{0}^{3}(2)\bigr),\ \bigl(I(3,1),R_{\infty}^{2}(1)\bigr),\ \bigl(I(3,2),P(1,7)\bigr),\ \bigl(I(3,7),R_{1}^{1}(2)\bigr)\\
 & \bigl(I(5,7),R_{0}^{3}(1)\bigr),\ \bigl(I(6,1),R_{1}^{1}(1)\bigr)\\
I(n,7):\ & \bigl(I(n-9,1),R_{1}^{(-n+9)\bmod 4+1}(3)\bigr),\ \bigl(I(n-8,7),R_{0}^{(-n+11)\bmod 3+1}(2)\bigr),\ \bigl(I(n-6,1),R_{\infty}^{(-n+10)\bmod 2+1}(1)\bigr)\\
 & \bigl(I(n-6,2),P(-n+10,7)\bigr),\ \bigl(I(n-6,7),R_{1}^{(-n+9)\bmod 4+1}(2)\bigr),\ \bigl(I(n-4,7),R_{0}^{(-n+11)\bmod 3+1}(1)\bigr)\\
 & \bigl(I(n-3,1),R_{1}^{(-n+9)\bmod 4+1}(1)\bigr),\ n>9\\
\end{align*}
\end{fleqn}
\subsubsection*{Modules of the form $I(n,8)$}

Defect: $\partial I(n,8) = 2$, for $n\ge 0$.

\begin{fleqn}
\begin{align*}
I(0,8):\ & - \\
I(1,8):\ & \bigl(I(0,1),I(4,7)\bigr),\ \bigl(I(0,2),R_{1}^{2}(1)\bigr),\ \bigl(I(0,3),P(0,7)\bigr),\ \bigl(I(0,5),P(0,1)\bigr),\ \bigl(I(0,6),R_{1}^{4}(1)\bigr)\\
 & \bigl(I(0,7),I(4,1)\bigr)\\
I(2,8):\ & \bigl(I(0,8),R_{0}^{2}(1)\bigr),\ \bigl(I(1,1),I(5,7)\bigr),\ \bigl(I(1,2),R_{1}^{1}(1)\bigr),\ \bigl(I(1,6),R_{1}^{3}(1)\bigr),\ \bigl(I(1,7),I(5,1)\bigr)\\
I(3,8):\ & \bigl(I(0,1),I(8,7)\bigr),\ \bigl(I(0,7),I(8,1)\bigr),\ \bigl(I(1,8),R_{0}^{1}(1)\bigr),\ \bigl(I(2,1),I(6,7)\bigr),\ \bigl(I(2,2),R_{1}^{4}(1)\bigr)\\
 & \bigl(I(2,6),R_{1}^{2}(1)\bigr),\ \bigl(I(2,7),I(6,1)\bigr)\\
I(n,8):\ & \bigl(I(n-3,1),I(n+5,7)\bigr),\ \bigl(I(n-3,7),I(n+5,1)\bigr),\ \bigl(I(n-2,8),R_{0}^{(-n+3)\bmod 3+1}(1)\bigr)\\
 & \bigl(I(n-1,1),I(n+3,7)\bigr),\ \bigl(I(n-1,2),R_{1}^{(-n+6)\bmod 4+1}(1)\bigr),\ \bigl(I(n-1,6),R_{1}^{(-n+4)\bmod 4+1}(1)\bigr)\\
 & \bigl(I(n-1,7),I(n+3,1)\bigr),\ n>3\\
\end{align*}
\end{fleqn}
\subsubsection{Schofield pairs associated to regular exceptional modules}

\subsubsection*{The non-homogeneous tube $\mathcal{T}_{1}^{\Delta(\widetilde{\mathbb{E}}_{7})}$}

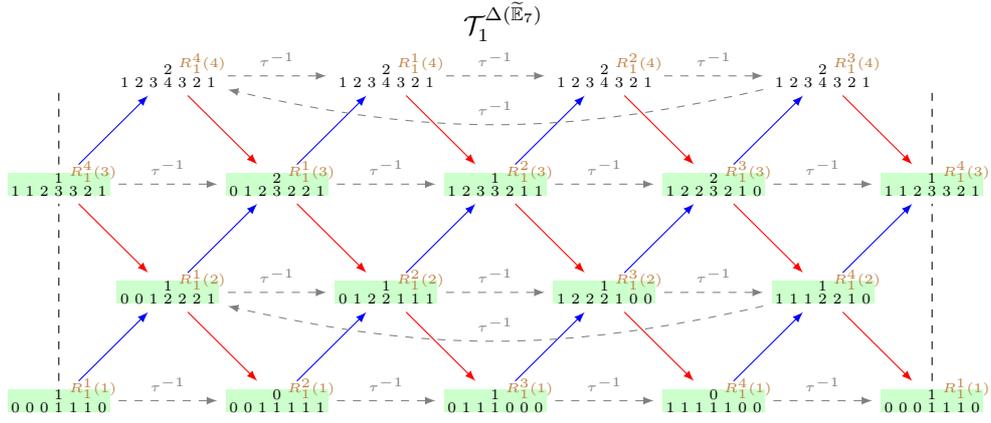
\begin{figure}[ht]



\begin{center}
\captionof{figure}{\vspace*{-10pt}$\mathcal{T}_{1}^{\Delta(\widetilde{\mathbb{E}}_{7})}$}
\begin{scaletikzpicturetowidth}{0.69999999999999996\textwidth}

\end{scaletikzpicturetowidth}
\end{center}

\end{figure}
\begin{fleqn}
\begin{align*}
R_{1}^{1}(1):\ & \bigl(I(0,8),P(0,6)\bigr),\ \bigl(I(1,6),P(0,8)\bigr),\ \bigl(I(1,7),P(1,1)\bigr)\\
R_{1}^{1}(2):\ & \bigl(R_{1}^{2}(1),R_{1}^{1}(1)\bigr),\ \bigl(I(5,7),P(0,7)\bigr),\ \bigl(I(1,6),P(1,6)\bigr),\ \bigl(I(3,1),P(2,1)\bigr),\ \bigl(I(1,7),P(4,7)\bigr)\\
R_{1}^{2}(1):\ & \bigl(I(0,5),P(0,3)\bigr),\ \bigl(I(2,1),P(0,7)\bigr),\ \bigl(I(0,6),P(1,8)\bigr),\ \bigl(I(0,7),P(2,1)\bigr)\\
R_{1}^{2}(2):\ & \bigl(R_{1}^{3}(1),R_{1}^{2}(1)\bigr),\ \bigl(I(2,2),P(0,2)\bigr),\ \bigl(I(4,7),P(1,7)\bigr),\ \bigl(I(0,6),P(2,6)\bigr),\ \bigl(I(2,1),P(3,1)\bigr)\\
 & \bigl(I(0,7),P(5,7)\bigr)\\
R_{1}^{3}(1):\ & \bigl(I(0,8),P(0,2)\bigr),\ \bigl(I(1,2),P(0,8)\bigr),\ \bigl(I(1,1),P(1,7)\bigr)\\
R_{1}^{3}(2):\ & \bigl(R_{1}^{4}(1),R_{1}^{3}(1)\bigr),\ \bigl(I(5,1),P(0,1)\bigr),\ \bigl(I(1,2),P(1,2)\bigr),\ \bigl(I(3,7),P(2,7)\bigr),\ \bigl(I(1,1),P(4,1)\bigr)\\
R_{1}^{4}(1):\ & \bigl(I(2,7),P(0,1)\bigr),\ \bigl(I(0,3),P(0,5)\bigr),\ \bigl(I(0,2),P(1,8)\bigr),\ \bigl(I(0,1),P(2,7)\bigr)\\
R_{1}^{4}(2):\ & \bigl(R_{1}^{1}(1),R_{1}^{4}(1)\bigr),\ \bigl(I(2,6),P(0,6)\bigr),\ \bigl(I(4,1),P(1,1)\bigr),\ \bigl(I(0,2),P(2,2)\bigr),\ \bigl(I(2,7),P(3,7)\bigr)\\
 & \bigl(I(0,1),P(5,1)\bigr)\\
R_{1}^{4}(3):\ & \bigl(R_{1}^{1}(2),R_{1}^{4}(1)\bigr),\ \bigl(R_{1}^{2}(1),R_{1}^{4}(2)\bigr),\ \bigl(I(8,1),P(0,7)\bigr),\ \bigl(I(6,7),P(2,1)\bigr),\ \bigl(I(4,1),P(4,7)\bigr)\\
 & \bigl(I(2,7),P(6,1)\bigr),\ \bigl(I(0,1),P(8,7)\bigr)\\
R_{1}^{1}(3):\ & \bigl(R_{1}^{2}(2),R_{1}^{1}(1)\bigr),\ \bigl(R_{1}^{3}(1),R_{1}^{1}(2)\bigr),\ \bigl(I(7,1),P(1,7)\bigr),\ \bigl(I(5,7),P(3,1)\bigr),\ \bigl(I(3,1),P(5,7)\bigr)\\
 & \bigl(I(1,7),P(7,1)\bigr)\\
R_{1}^{2}(3):\ & \bigl(R_{1}^{3}(2),R_{1}^{2}(1)\bigr),\ \bigl(R_{1}^{4}(1),R_{1}^{2}(2)\bigr),\ \bigl(I(8,7),P(0,1)\bigr),\ \bigl(I(6,1),P(2,7)\bigr),\ \bigl(I(4,7),P(4,1)\bigr)\\
 & \bigl(I(2,1),P(6,7)\bigr),\ \bigl(I(0,7),P(8,1)\bigr)\\
R_{1}^{3}(3):\ & \bigl(R_{1}^{4}(2),R_{1}^{3}(1)\bigr),\ \bigl(R_{1}^{1}(1),R_{1}^{3}(2)\bigr),\ \bigl(I(7,7),P(1,1)\bigr),\ \bigl(I(5,1),P(3,7)\bigr),\ \bigl(I(3,7),P(5,1)\bigr)\\
 & \bigl(I(1,1),P(7,7)\bigr)\\
\end{align*}
\end{fleqn}
\subsubsection*{The non-homogeneous tube $\mathcal{T}_{0}^{\Delta(\widetilde{\mathbb{E}}_{7})}$}

\begin{figure}[ht]



\begin{center}
\captionof{figure}{\vspace*{-10pt}$\mathcal{T}_{0}^{\Delta(\widetilde{\mathbb{E}}_{7})}$}
\begin{scaletikzpicturetowidth}{0.69999999999999996\textwidth}

\end{scaletikzpicturetowidth}
\end{center}

\end{figure}
\begin{fleqn}
\begin{align*}
R_{0}^{1}(1):\ & \bigl(I(2,1),P(1,1)\bigr),\ \bigl(I(2,7),P(1,7)\bigr),\ \bigl(I(0,8),P(1,8)\bigr)\\
R_{0}^{1}(2):\ & \bigl(R_{0}^{2}(1),R_{0}^{1}(1)\bigr),\ \bigl(I(5,1),P(2,1)\bigr),\ \bigl(I(5,7),P(2,7)\bigr),\ \bigl(I(2,1),P(5,1)\bigr),\ \bigl(I(2,7),P(5,7)\bigr)\\
R_{0}^{2}(1):\ & \bigl(I(1,6),P(0,2)\bigr),\ \bigl(I(1,2),P(0,6)\bigr),\ \bigl(I(1,1),P(2,1)\bigr),\ \bigl(I(1,7),P(2,7)\bigr)\\
R_{0}^{2}(2):\ & \bigl(R_{0}^{3}(1),R_{0}^{2}(1)\bigr),\ \bigl(I(7,1),P(0,1)\bigr),\ \bigl(I(7,7),P(0,7)\bigr),\ \bigl(I(4,1),P(3,1)\bigr),\ \bigl(I(4,7),P(3,7)\bigr)\\
 & \bigl(I(1,1),P(6,1)\bigr),\ \bigl(I(1,7),P(6,7)\bigr)\\
R_{0}^{3}(1):\ & \bigl(I(3,1),P(0,1)\bigr),\ \bigl(I(3,7),P(0,7)\bigr),\ \bigl(I(1,8),P(0,8)\bigr),\ \bigl(I(0,6),P(1,2)\bigr),\ \bigl(I(0,2),P(1,6)\bigr)\\
 & \bigl(I(0,1),P(3,1)\bigr),\ \bigl(I(0,7),P(3,7)\bigr)\\
R_{0}^{3}(2):\ & \bigl(R_{0}^{1}(1),R_{0}^{3}(1)\bigr),\ \bigl(I(6,1),P(1,1)\bigr),\ \bigl(I(6,7),P(1,7)\bigr),\ \bigl(I(3,1),P(4,1)\bigr),\ \bigl(I(3,7),P(4,7)\bigr)\\
 & \bigl(I(0,1),P(7,1)\bigr),\ \bigl(I(0,7),P(7,7)\bigr)\\
\end{align*}
\end{fleqn}
\subsubsection*{The non-homogeneous tube $\mathcal{T}_{\infty}^{\Delta(\widetilde{\mathbb{E}}_{7})}$}

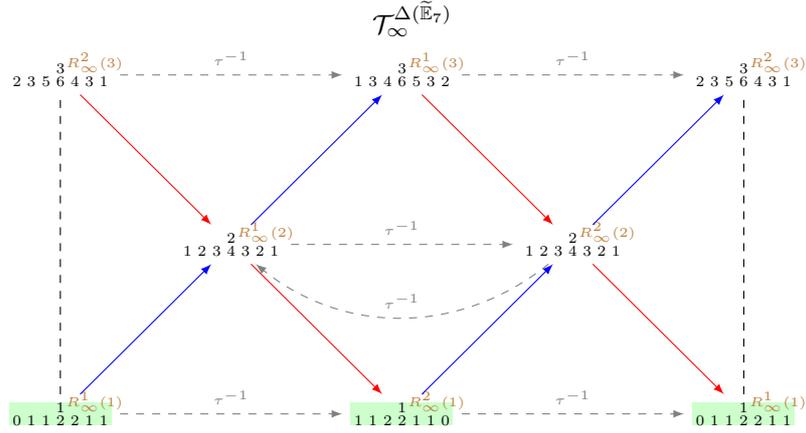
\begin{figure}[ht]



\begin{center}
\captionof{figure}{\vspace*{-10pt}$\mathcal{T}_{\infty}^{\Delta(\widetilde{\mathbb{E}}_{7})}$}
\begin{scaletikzpicturetowidth}{0.59999999999999998\textwidth}

\end{scaletikzpicturetowidth}
\end{center}

\end{figure}
\begin{fleqn}
\begin{align*}
R_{\infty}^{1}(1):\ & \bigl(I(5,1),P(0,7)\bigr),\ \bigl(I(4,7),P(1,1)\bigr),\ \bigl(I(3,1),P(2,7)\bigr),\ \bigl(I(2,7),P(3,1)\bigr),\ \bigl(I(1,1),P(4,7)\bigr)\\
 & \bigl(I(0,7),P(5,1)\bigr)\\
R_{\infty}^{2}(1):\ & \bigl(I(5,7),P(0,1)\bigr),\ \bigl(I(4,1),P(1,7)\bigr),\ \bigl(I(3,7),P(2,1)\bigr),\ \bigl(I(2,1),P(3,7)\bigr),\ \bigl(I(1,7),P(4,1)\bigr)\\
 & \bigl(I(0,1),P(5,7)\bigr)\\
\end{align*}
\end{fleqn}
\clearpage

\include{sp_E8}

\end{document}